\input amstex
\documentstyle{amsppt}
\magnification\magstep1
\NoRunningHeads
\hsize6.9truein
\pageheight{23 truecm}
\baselineskip=15pt
\input btxmac.tex               
\bibliographystyle{plain}
\input epsf


\def\ec#1{e_{\text{top}}(#1)}
\def\t#1{\tilde{#1}}
\def\ol#1{\overline{#1}}
\def\ul#1{\underline{#1}}

\def\sb{\ol{\Sigma}}

\def\ideal#1.{I_{#1}}
\def\ring#1.{\Cal O_{#1}}
\def\proj#1.{\Bbb P(#1)}
\def\pr #1.{\Bbb P^{#1}}
\def\af #1.{\Bbb A^{#1}}
\def\Hz #1.{\Bbb F_{#1}}
\def\Hbz #1.{\overline{\Bbb F}_{#1}}
\def\fb#1.{\underset #1 \to \times}
\def\res#1.{\underset {\ \ring #1.} \to \otimes}
\def\au#1.{\operatorname {Aut}\,(#1)}
\def\deg#1.{\operatorname {deg } (#1)}
\def\pic#1.{\operatorname {Pic}\,(#1)}
\def\pico#1.{\operatorname{Pic}^0(#1)}
\def\ner#1.{NS (#1)}
\define\rdown#1.{\llcorner#1\lrcorner}\
\define\rup#1.{\ulcorner#1\urcorner}

\def\kts{K_{\t{S}_0}^2}
\def\list#1.#2.{{#1}_1,{#1}_2,\dots,{#1}_{#2}}
\def\loc#1.#2.{\Cal O_{#1,#2}}
\def\fderiv#1.#2.{\frac {\partial #1}{\partial #2}}
\def\map#1.#2.{#1 \longrightarrow #2}
\def\rmap#1.#2.{#1 \dasharrow #2}
\def\emb#1.#2.{#1 \hookrightarrow #2}
\def\non#1.#2.{\text {Spec }#1[\epsilon]/(\epsilon)^{#2}}
\def\H#1.#2.{\text {Hilb}^{#1}(#2)}
\define\sym#1.#2.{\operatorname {Sym}^{#1}(#2)}
\def\Hb#1.#2.{\text {Hilb}_{#1}(#2)}
\def\Hm#1.#2.{\Hom_{#1}(#2)}
\def\prd#1.#2.{{#1}_1\cdot {#1}_2\cdots {#1}_{#2}}
\def\OO{\Cal O}


\def\alist#1.#2.#3.{#1_1\,#2\,#1_2#2\,\dots\,#1_{#3}}
\def\lmap#1.#2.#3.{#1 \overset #2\to \longrightarrow #3}
\def\ses#1.#2.#3.{0\longrightarrow #1 \longrightarrow #2 \longrightarrow #3 
\longrightarrow 0}
\def\les#1.#2.#3.{0\longrightarrow #1 \longrightarrow #2 \longrightarrow #3}

\def\Hi#1.#2.#3.{\text {Hilb}^{#1}_{#2}(#3)}


\define\Diff{\operatorname{Diff}}

\define\Hom{\operatorname{Hom}}

\def\Supp{\operatorname{Supp}}

\def\di{\operatorname{dim}}

\def\deg{\operatorname{deg}}

\def\sg{\operatorname{Sing}}

\def\Pic{\operatorname{Pic}}


\def\dd{\Cal D}
\def\e{\Cal E}
\def\ZZ{\Bbb Z}
\def\gH{\goth H}
\def\alg{\pi_1^{\text {alg}}}
\def\aos{{\Bbb A}^1_*}


\def\mapdown#1{\big\downarrow\rlap{$\vcenter
{\hbox{$\scriptstyle#1$}}$}}

\def\mapse#1{
{\vcenter{\hbox{$\mathop{\smash{\raise1pt\hbox{$\diagdown$}\!\lower7pt
\hbox{$\searrow$}}\vphantom{p}}\limits_{#1}\vphantom{\mapdown{}}$}}}}

\nologo
\topmatter
\title Rational curves on quasi-projective surfaces \endtitle
\author Se\'an Keel and James M\raise 1.6pt \hbox{\text {\smc c}}Kernan
 \endauthor
\address Department of Mathematics, 
University of Texas, at Austin, Austin TX 78712 
Department of Mathematics,
University of California at Santa Barbara
Santa Barbara, CA 93101
\endaddress
\endtopmatter


\heading Contents \endheading

\item{\S 1} Introduction and Statement of Results			\hfill page 2\ \
\item{\S 2} Glossary of notation and conventions			\hfill page 15\
\item{\S 3} Gorenstein del Pezzos					\hfill page 17\
\item{\S 4} Koll\'ar's Bug-Eyed cover					\hfill page 22\
\item{\S 5} Deformation theory of log pairs				\hfill page 28\
\item{\S 6} Criteria for log uniruledness				\hfill page 35\
\item{\S 7} The algebraic fundamental group				\hfill page 45\
\item{\S 8} Flushness and Preparation for the Hunt			\hfill page 47\
\item{\S 9} Bogomolov Bound						\hfill page 70\
\item{\S 10} Surfaces with small coefficient				\hfill page 72\
\item{\S 11} Partial classification of $K_T$-contractions		\hfill page 76\
\item{\S 12} The linear system $|K_S + A|$				\hfill page 92\
\item{\S 13} Classification of Bananas and Fences			\hfill page 96\
\item{\S 14} $T_1$ a net						\hfill page 101\
\item{\S 15} $g(A_1) > 1$						\hfill page 104
\item{\S 16} $A_1$ has a simple cusp					\hfill page 111
\item{\S 17} $A_1$ has a simple node					\hfill page 114
\item{\S 18} $A_1$ smooth						\hfill page 122
\item{\S 19} The smooth banana						\hfill page 128
\item{\S 20} Proof of (1.1) and corollaries				\hfill page 133
\item{\S 21} A surface with $\alg(S^0) = \{1\}$ but no tiger		\hfill page 135
\item{\S 22} Complements and toric pairs				\hfill page 139
\item{\S 23} Classification of almost all rank one log del Pezzos	\hfill page 143
\item{\S L}  Appendix: Log terminal singularities and adjunction	\hfill page 152
\item{\S N}  Appendix: Normalisation of an algebraic space		\hfill page 155
\item{} Index to Obscure or unconventional notation 			\hfill page 156
\item{} References							\hfill page 157

\heading  \S 1 Introduction and Statement of Results \endheading

In classifying smooth projective varieties, one looks for an intrinsic
map to projective space. The natural intrinsic line bundles are the 
various tensor powers of $\omega_X$, and so one is led to consider 
the pluricanonical linear series $|mK_X|$. The obvious question arises, 
what if these series are all empty? The answer is striking, one of 
the principal accomplishments of Mori's program: 
\proclaim{1.0 Theorem (Miyaoka, Mori, Kawamata, Shokurov, and others)} Let 
$X$ be a smooth projective variety of dimension at most three. 

Then $|mK_X|$ is empty for all $m > 0$ iff $X$ is covered by images of $\pr 1.$.
\endproclaim

(We say $X$ is covered by images of $\pr 1.$ if for every point $p$ of 
$X$, we may find a morphism $f:\map {\pr 1.}.X.$ such that the point 
$p$ belongs to $f(\pr 1.)$.)

In this paper we consider the analogous question for quasi-projective (or
more generally log) varieties. If $U$ is a smooth quasi-projective variety, 
then, following Iitaka, one picks a smooth compactification $U \subset X$ 
such that the complement $D=X\setminus U$ is a divisor with normal 
crossings. The linear series $|m(K_X+D)|$ turn out to depend only on $U$, 
not on $X$ or $D$, and is thus the natural analogue of $|mK_X|$. 
$|m(K_X + D)|$ is called the log pluricanonical series, and the problem is to 
characterise those smooth quasi-projective varieties for which the log pluricanonical
series are all empty. 
  
When $U$ is a curve the solution is elementary, $|m(K_X+D)|$ is 
empty for all $m > 0$ iff $U = \af 1.$ or $\pr 1.$.

 In dimension two the problem is already surprisingly subtle, and has received
considerable attention.   An import special case was  settled
by Miyanishi and Tsunoda, \cite{MT84a},\cite{MT84b}, and further results 
have been obtained by Zhang, \cite{Zhang88}. Here our main goal is a 
complete, and
self-contained solution (throughout the paper, everything takes place over
${\Bbb C}$):

\proclaim{1.1 Theorem} Let $U$ be a smooth quasi-projective variety of 
dimension at most two. Then $|m(K_X+D)|$ is empty for all $m>0$ 
iff $U$ is dominated by images of $\af 1.$. 
\endproclaim

We say $U$ is {\bf dominated by images of a curve $C$}, if there is a dense open 
subset $V$ of $U$ and for every point $p$ of $V$, we can find a non-constant 
morphism $f:\map C.U.$ such that $p\in f(C)$.  We say $U$ is dominated by
rational curves, if it is dominated by images of $\pr 1.$. 

We note that the reverse direction of the implication in (1.1) is fairly straightforward,
and holds in all dimensions, see (5.11). Thus (1.1) can be viewed as saying, either
$U$ has a log pluricanonical section, or there is a clear geometric reason why
it cannot. 

The other main results of this paper are very strong partial classifications
of log del Pezzo surfaces, that is projective surfaces with quotient
singularities and $-K_S$ ample. This includes a classification of
all but a bounded family of log del Pezzo surfaces of Picard number one. We will
explain this classification at the end of this introduction.
We believe that with sufficient effort the methods of the paper would yield
a complete classification.

Log del Pezzo surfaces are of interest for several reasons, 
independently of (1.1). They are obviously important for the study of open 
surfaces. Beyond this, the log category is important even for the study of 
projective varieties with, for example, log surfaces (surfaces with
boundary) playing an intermediary r\^ole between surfaces and
threefolds. Good examples of this are Kawamata's proof of the
Abundance conjecture, and Shokurov's program for Flips. Log del Pezzos
also occur as the centres of 4-fold log flips.

In addition log del Pezzos play an essential r\^ole in the compactification
of the moduli space of surfaces of general type. They occur at
the boundary of the moduli space of surfaces of general type (just as 
rational curves occur at the boundary of the moduli space of curves), and are 
the main object of study in the proofs of Alexeev's boundedness theorems 
\cite{Alexeev94}, used to show the moduli space is projective.

 We will use the following terminology. Given a variety $X$, let $X^0$ denote 
the smooth locus of $X$.

\definition{1.2 Definition} Let $X$ be a variety, $C$ a curve
in $X$ and $D$ a subset of $X$. We say that $C$ meets $D$ 
$k$ times if the inverse image of $D$ on the normalisation of
$C$ is a set of $k$ points.
\enddefinition

For example a smooth point or a unibranch singularity of $C$ counts once,
and more generally a singularity with $r$ branches counts $r$ times, 
regardless of the order of contact.

We note also that throughout the paper, by a rational curve, we
mean a complete rational curve, that is an image of $\pr 1.$. By uniruled,
we mean dominated by rational curves.

\remark{1.2.1 Remarks} If a complete variety is dominated by rational
curves, then it is in fact covered by rational curves (that is through every
point there is a rational curve), since the degeneration of a rational
curve is again rational. 

In the definition of dominating, one can equivalently require that the maps 
$f:\map C.U.$ form a flat family (see (IV.1.3.5) of \cite{Kollar96}).
\endremark

\subhead Connection between (1.1) and log del Pezzo surfaces \endsubhead
In proving (1.1) one tries to simplify the situation by a birational
contraction $\pi:\map X.Y.$ of an irreducible curve $E \subset X$. 
In order to preserve the dimension of $|m(K_X + D)|$, it is sufficient that 
$(K_X+D)\cdot E<0$. Under the contraction, the Picard number goes down by one, and 
from this point of view, $Y$ is simpler than $X$. The complication is that
$Y$ can have quotient singularities (if $E$ has self-intersection at most $-2$). 
We replace $(X,D)$ by $(Y,D_Y= \pi_*(D))$, and continue, by looking for
a $(K_Y+D_Y)$-negative contraction. Such a series of contactions is 
called the $(K_X+D)$-minimal model program (MMP for short).

If at the start $|m(K_X+D)|=\emptyset$ for all $m>0$, then the Log Abundance 
Theorem implies that either $X$ has a fibration $\map X.C.$ whose
general fibre is $\pr 1.$ and meets $D$ at most once, or, we have a birational morphism
$\pi:\map X.S.$, a composition of contractions as above, with
$S$ of Picard number one, and $K_S + D_S$ anti-ample and log terminal (discussed
below). In the first case
of course we take the general fibre. In the second case, $X$ is dominated by
rational curves meeting $D$ at most once, iff $S$ is dominated by
rational curves meeting $\pi(D)$ at most once. The one dimensional part of
$\pi(D)$ is $D_S$, the zero dimensional part, $V$, is a union of components of $D$
contracted by $\pi$. In this way, (1.1) is reduced (and in fact equivalent) to:

\proclaim{1.3 Theorem} Let $S$ be a normal projective surface of Picard number 
one, with quotient singularities. Suppose $D \subset S$ is a reduced curve, such 
that $K_S+D$ is log terminal, and $-(K_S+D)$ is ample. Let $V \subset S$ be
any finite set of points. Then $S$ is dominated by rational curves, meeting
$D \cup V$ at most once.
\endproclaim

An elementary discussion of the notion of log terminal is given in Appendix L. 
We note here (so that the reader may have some idea of the term's meaning) that
if $D$ is irreducible then $K_S+D$ is log terminal iff the pair $(S,D)$ has 
quotient singularities, that is locally analytically the quotient of a smooth 
pair $(S',D')$ by a finite group. For the $(K_S+D)$-MMP, log terminal is the 
correct generalisation of normal crossings, in the sense that it is preserved by 
the birational contractions that occur in the $(K_S+D)$-minimal model program. 

In fact in (1.3) it is enough to consider the case $V = \sg(S)$. The main point behind 
this reduction is Mori's observation that on a smooth space, the general member of a 
dominating family of images of $\af 1.$ deforms freely (as an image of $\af 1.$) and 
so in particular misses any fixed codimension two subset, see (5.5) and (5.8). Thus 
(1.3) can be equivalently (and somewhat more aesthetically) stated as:

\proclaim{1.3.1 Theorem} Let $S$ be a 
normal projective surface of Picard number 
one, with quotient singularities. Suppose $D$ is a reduced curve, such 
that $K_S+D$ is log terminal, and $-(K_S+D)$ is ample. 
Then $S^0$ is dominated by rational curves which meet $D$ at most once.
\endproclaim

As we will explain shortly, the main work is proving (1.3) in the case when $D$ is 
empty, and is thus to show that the smooth locus of a log del Pezzo surface 
is uniruled (that is dominated by complete rational curves). 

Our proof of (1.3) proceeds roughly as follows: For one class of rank one log del Pezzo surfaces, 
those with a tiger (defined below), we give a short and elegant deformation theoretic proof of 
(1.3), see (6.1). Boundedness results of Alexeev and Koll\'ar imply that rank one log del Pezzos 
without a tiger are bounded, see (23.1). We complete the proof by constructing an explicit 
finite list of families of rank one log del Pezzos which includes any with no tiger and whose 
smooth locus has trivial algebraic fundamental group (it is easy to reduce (1.1) to this case, 
see \S 7). For each of the surfaces in this list, we directly construct a dominating family of
rational curves. 

We will explain all this in much greater detail below, but first we present some additional 
results which are of independent interest. Most are corollaries of (1.1). Proofs are given 
in \S 20. We will use the following notation. 

$\aos = \af 1. \setminus \{0\}$. We say that a ${\Bbb Q}$-Weil divisor is {\bf effective} 
if some positive multiple is an integral Weil divisor linearly equivalent to an effective Weil divisor.
Thus a ${\Bbb Q}$-Weil divisor is effective iff it has non-negative Kodaira dimension.

In view of (5.11), it is natural to think of $\af 1.$ and $\aos$ as the open analogues of
$\pr 1.$ and elliptic curves, since the existence of 
dominating rational, or elliptic families has analogous implications
on ordinary Kodaira dimension. 

\proclaim{1.4 Proposition} Let $B \subset S$ be a reduced curve on a normal projective surface and
set $U = S \setminus (B \cup \sg(S))$. Consider the following conditions:
\roster
\item $K_S + B$ is not effective.
\item $K_S + B$ is numerically trivial, but not log canonical.
\item $K_S + B$ is numerically trivial, and $B \neq \emptyset$.
\item $K_S$ is numerically trivial, $B = \emptyset$, and 
$S$ has a singularity which is {\bf not} a quotient singularity.
\endroster

If any of the above hold, then $U$ is dominated by images of $\aos$, and by images of
$\af 1.$ if (1) or (2) holds.
\endproclaim

We note one interesting implication of (1.4.2) and (5.11):
The complement of a integral plane cubic, $B$, is dominated by images of $\af 1.$
iff $B$ has a cusp, or $B$ is the union of a smooth conic and  a  tangent line, or $B$ is
the union of three lines meeting at a point (indeed this is just a list of the non log 
canonical possibilities).

Next we have a version of (1.3) for any boundary:

\proclaim{1.5 Corollary} Suppose the pair $(S, \Delta)$ consists of
a projective surface $S$ and boundary $\Delta$, such that $K_S+\Delta$ is 
log canonical. Then either
\roster
\item $|m(K_S+\Delta)|$ is non-empty for some $m>0$, or
\item There is a covering family of rational curves $C$, such that 
$(K_S+\Delta)\cdot C<0$. 
\endroster
\endproclaim 

(1.5) thus says that either some multiple of $K_S+\Delta$ has a section, or
there is a good geometric reason why it does not. 

It may be tempting to
believe that (1.5) follows automatically from the MMP: One can assume
the $(K_S+\Delta)$-MMP gives a composition of contractions
$f:\map S.S'.$ such that $-(K_S' + \Delta')$ is ample. Every
curve meets $(K_{S'}+\Delta')$-negatively.  However, if the curve
meets the exceptional locus of $f$, its strict transform may be
$(K_S+\Delta)$-positive. Thus in order to prove (1.5) along these lines 
some result such as (1.3) is required. Such an attempt
(in dimension three !) was the original impetus for this paper.

\proclaim{1.6 Corollary} If $(S,\Delta)$ is a log Fano surface 
then
$S^0$ is rationally connected, and $U = S \setminus \rdown \Delta .$ is connected
by images of $\aos$. In particular, $\pi_1(S^0)$ is finite, and 
$\pi_1(U)$ is almost Abelian.
\endproclaim

(Almost abelian is defined in \S 7. $\rdown \Delta.$ is the support of $\Delta$). 
The main issue in (1.6) is the uniruledness 
of the smooth locus of a rank one log del Pezzo. This was
conjectured by Miyanishi and Tsunoda, \cite{MT82}. 
The finiteness of the fundamental group $\pi_1(S^0)$ has been established previously in 
\cite{Zhang94b}-\cite {Zhang94a} and separately in \cite{FKL93}. The smooth
locus of a log Fano threefold (or surface) has  finite algebraic fundamental group
by a boundedness result of \cite{Borisov96}.

\proclaim{1.7 Corollary} Suppose $K_S+D+\Delta$ is log terminal for some 
effective ${\Bbb Q}$-divisor $\Delta$. Then every co-extremal ray 
(cf. \cite{Batryev89}) is spanned by classes $[C]$, for $C$ the general member 
of a covering family of rational curves, contained in the smooth locus, and 
meeting $D$ in at most one point. In particular if $-(K_S+D+\Delta)$ is ample,
then the nef cone of $S$ is polyhedral and spanned by such classes. 
\endproclaim

The fact that the nef cone in (1.7) is polyhedral was proved in 
\cite{Batryev89}, our contribution is the description of the generators. 
(1.7) should be compared to the cone theorem, which describes generators for 
the cone of effective curves.

(1.7) in turn gives:

\proclaim{1.8 Corollary} Suppose $(S,D+\Delta)$ is log Fano. 
Then $T_S(-\log D)$ is generically semi-positive, in other words for a general member 
$C$ of a sufficiently ample linear series, the vector bundle $T_S(-\log D)|_C$ 
has no quotient line bundles of negative degree.
\endproclaim

The fact that (1.8) is implied by (1.7) was pointed out to us by Mori.

We note that for terminal Fano threefolds, the analogue of (1.6) implies the 
analogue of (1.8). Thus if one could prove terminal Fano threefolds have 
uniruled smooth locus, it would follow as in  \cite{Kawamata92} that 
the set of terminal Fano threefolds is bounded.

(1.8) itself has some interesting corollaries:

\proclaim{1.8.1 Corollary} Let $(S,B)$ be a log Fano pair (with $S$ a projective
surface, and $B\subset S$ a reduced curve). 
$$
\sum_{p \in \text{Sing}(S \setminus B)} \frac{r_p -1}{r_p} \leq 2 + \rho(S) - \ec{B}
$$
where $r_p$ is the order of the local fundamental group, $\ec{B}$ is the 
topological Euler characteristic of $B$, and $\rho$ is the Picard number.
\endproclaim

Note when $B$ is empty then (1.8.1) implies that log del Pezzo surface
can have at most $4 + 2 \cdot \rho$ singularities.
As far as we know, no such bound has been previously observed. 

For $\rho =1$, we prove a stronger version of (1.8.1) (it applies to any log
terminal surface $S$), following an argument of \cite{Kawamata92}, see
(9.2). The proof is independent of (1.1) (and indeed, anything else in
the paper). We call this intermediate result the {\bf Bogomolov Bound} 
and use it repeatedly in the proof of (1.3). The Bogomolov Bound has
an interesting corollary, see (9.3). In \cite{Alexeev94}, using different methods, 
Alexeev obtains results related to (but much deeper than) (9.3). 

In the light of (1.0) and (1.1), we propose the following:

\proclaim{1.9 Conjecture} Let $U$ be a smooth quasi-projective variety. 
Then, either
\roster
\item $\bar \kappa (U) \geq 0$, or
\item $U$ is dominated by images of $\af 1.$.
\endroster
\endproclaim

Note that (1.9) is equivalent to the following: 

\proclaim{1.10 Conjecture} Let the pair $(X,D)$ consist of a smooth variety
$X$, and reduced normal crossings divisor $D$. Then, either
\roster
\item $|m(K_X+D)|$ is non-empty for some $m>0$, (that is $K_X + D$ is effective), or
\item there is a dominating family of rational curves meeting $D$
in at most one point.
\endroster
\endproclaim

By (5.11), (1) and (2) of (1.10) are mutually exclusive, and one can view the conjecture as 
saying either $X$ has a log pluricanonical section, or there is a clear geometric reason why 
none can exist.

There are natural analogues of (1.4-7) in higher dimensions, and similar 
implications will hold between them, if one has the MMP. We wish however to 
note one difference; in dimension two (1.1) implies (1.6). This follows from
deformation theory, and the 
fact that log terminal surface singularities are quotient singularities,
see Appendix L and \S 5. Three-fold log terminal singularities need not be quotient
singularities. 

In higher dimensions there is not much evidence for (1.9-10). We can prove 
(1.9) in the case $-(K_X+D)$ is ample, see (5.4), and a version of (1.3) for 
a threefold $X$ in the case $D$ has two components, see (6.7).

Under the assumptions of (1.9), if $X$ has dimension at most three, then the
MMP and log Abundance \cite {KMM94} imply that $X$ is covered by (not 
necessarily rational) curves $C$ with $(K_X + D) \cdot C < 0$. Thus one is 
moved to consider:

\proclaim{1.11 Conjecture: Log Bend and Break} (notation as in 1.10). If $(K_X +D) \cdot C < 0$,
and $C \not \subset D$ then through a general point of $C$ there is a rational
curve meeting $D$ at most once.
\endproclaim

Mori's famous bend and break argument proves (1.11) in the case
$D$ empty. We don't see any way to extend the proof to the 
case when $D$ is not empty. In Mori's argument, $K_X  \cdot C <0$
is used to show that (mod $p$) some multiple of $C$ moves with
a fixed point. If a curve moves with a point fixed, then rigidity
implies there is a rational curve through the fixed point. However
it is not clear (at least to us) how to control how the generated rational curve intersects
$D$. In (6.9) we give examples indicating some of the difficulties. 

In any case, the rational curves we use to prove (1.1) are obtained by an 
entirely different procedure, which we discuss below. However, once we have 
(1.1) we obtain a version of (1.11) as a:
\proclaim{1.12 Corollary} Assumptions as in (1.10). If $X$ has dimension 
at most two, or $X$ has dimension three and $K_X + D$ {\bf is} effective, 
then there is a rational curve through a general point of $C$, meeting $D$ at most once.
\endproclaim

As already mentioned, our second main result is a classification of {\it most}
rank one log del Pezzo surfaces. 

\subhead Outline of the proof of (1.3) \endsubhead

Now we turn to an outline of the proof of (1.3), both to indicate its logical structure, and to 
high light many of issues which arise that are of independent interest. The logical order of the 
proof is also indicated by the flowchart (1.18) below.

The proof divides into two cases. We will use the notation of (1.3). Thus $S$ is a normal 
projective surface with Picard number one with quotient singularities. $D$ is a curve,
such that $K_S+D$ is log terminal and $-(K_S+D)$ is ample. 

Case I: $D \neq \emptyset$. The case of non-empty $D$ is proved in (6.2), using
deformation theory and Koll\'ar's Bug-Eyed cover. 
Here is a sketch: For simplicity, suppose $D$ is irreducible.
By adjunction $D$ itself is a smooth rational curve, and we 
obtain a covering family of rational curves by deforming a high multiple of 
$D$. The difficulty of course is the presence of singularities. For this 
we use the Bug-Eyed cover:  Given a normal surface $S$ with 
quotient singularities, there exists a unique smooth (but non-separated) 
algebraic space $\flat:\map S^{\flat}.S.$ such that $\flat$ is a universal
homeomorphism, and an isomorphism over $S^0$. See \S 4.
For (1.3) the important point is that in terms of
$\Hom(\Sigma,S^{\flat})$, for a proper curve $\Sigma$, $S^{\flat}$ behaves 
exactly like a smooth variety. The idea then is to lift the problem to $S^{\flat}$. 
The condition that $K_S+D$ is log terminal is equivalent to the condition that 
$\flat^{-1}(D) \subset S^{\flat}$ is smooth, and the condition that $K_S+D$ is 
anti-ample, is equivalent to the existence of a endomorphism 
$f:\map {\pr 1.}.D.= \pr 1.$ and a commutative diagram
$$
\CD
\pr 1. @>{g}>> S^{\flat} \\
@V{f}VV         @V{\flat}VV \\
D @>{}>>    S ,
\endCD
$$
see (4.14). Now deformation theory implies we can deform $g$ away from the singularities 
of $S$, maintaining a single point of (necessarily high order) contact with 
$D$. 

The case of $D$ non-empty has been previously proved by Miyanishi and 
Tsunoda \cite{MT84a} and \cite{MT84b}. Their proof is based on a 
classification of pairs $(S,D)$. Deformation 
theory together with the Bug-Eyed cover 
gives a simple, transparent proof, as well as some generalisations to higher 
dimensions.

Case two: 
When $D$ is empty, the idea is to replace $S$ by $S_1$, birational to $S$,
which contains a non-empty $D$, in such a way that  (1.3) for
$(S_1,D)$ implies $S^0$ is uniruled. $(S_1,D)$ will come from extracting
a divisor $E$ via a blow up $f:\map T.S.$, and then blowing
down in a different direction $\pi:\map T.S_1.$, 
$D = \pi(E)$. To find $(S_1,D)$ essentially reduces to finding
$E$. For this we introduce the following, the most important technical definition of
the paper:

\definition{1.13 Definition} Let $X$ be proper, normal variety. Let $\Delta$ be an effective
${\Bbb Q}$-Weil divisor on $X$. A {\bf special tiger} for 
$K_X+\Delta$ is an effective $\Bbb Q$-divisor $\alpha$ such that 
$K_X+\Delta+\alpha$ is numerically  trivial, but not Kawamata log terminal 
(klt). 
\enddefinition
For some $m$, $m\alpha$ is a very singular element of $|-m(K_X+\Delta)|$ and 
so a sort of antithesis to Reid's general elephant. Hence the terminology. 
By a special tiger for $X$, or a special tiger (without further reference), 
we mean a special tiger for $K_X$.

 Note that if we have a special tiger there is at least one divisor 
$E$ of discrepancy (with respect to $K_S+\alpha$) at most $-1$. We will call any such 
$E$ a {\bf tiger}. In fact we are much more interested in $E$ (and the resulting
new log Mori fibre space structure) than in $\alpha$. Sometimes we will 
be a little sloppy in our notation and use the word tiger to mean either 
$E$ or $\alpha$. This should not cause any harm, because normally we are only 
interested in when $X$ does or does not have a tiger, $E$, which is 
equivalent to the existence of a special tiger, $\alpha$.

We will define klt, and explain the motivation behind (1.13) in (1.15) below.

The next two results indicate the usefulness of tigers for (1.3):

\proclaim{1.14 Proposition} If $S$ is a projective surface with
quotient singularities and $S$ has a tiger, then $S^0$ is uniruled.
\endproclaim

\proclaim{1.15 Proposition} The collection of rank one log del Pezzo surfaces
which do not have a special tiger is bounded.
\endproclaim

Note (1.14) is implied by (6.1).

(1.15) follows from (9.3) and quite general boundedness principals. We prove a stronger 
statement in \S 23. Of course, together Case I and (1.14-15) imply (1.3) in all but a 
bounded collection of cases. We will turn to the question of classifying these cases  
(that is surfaces without tiger) in a moment. But first some general remarks on tigers.

For the rest of the introduction, unless otherwise noted, $S$ will
indicate a log del Pezzo surface of rank one.

Shokurov has independently considered tigers, in relation to complements. 
He has observed that if $S$ has a tiger, then $K_S$ is $1$, $2$, $3$, $4$ 
or $6$-complemented ($K_S$ is $n$-complemented if there is a member 
$M \in |-nK_S|$ such that $K + 1/n M$ is log canonical). We include 
Shokurov's proof of this in \S 22. Note Shokurov's result, and (1.15) 
imply there is a uniform $N > 0$ such that 
$|-NK_S|$ is non-empty for all $S$. This result is 
interesting in view of the fact that there is no bound on the index of $K_S$.

Complements in dimension $n$ give information on extremal neighbourhoods in 
dimension $n+2$ (see \cite{Kollaretal}). They play an important  r\^ole in
Shokurov's program for log flips. These observations were pointed 
out to us by Alessio Corti. 

Tigers have the following relation to the notion of affine-ruled (as
defined by Miyanishi). The proof is given in (21.4).
\proclaim{Lemma} If $S$ is affine-ruled (that is contains
a product neighbourhood $U \times \af 1.$)  then $S$ has a tiger.
\endproclaim

There are $S$ with simply connected smooth locus (cf. \S 21), but no tiger. 
These give counter-examples to  Miyanishi's  conjecture (see 
\cite{Zhang94a} ) that  the smooth locus of any rank one log del Pezzo
has a finite \'etale cover, that is affine-ruled. 

It is a fairly simple matter to reduce the proof of (1.3) to the case when 
$\pi_1^{alg}(S^0)$ is trivial (cf. \S 7). We will assume this for the rest of 
the introduction. We will sometimes abuse notation and say that $S^0$ (or even $S$) is 
simply connected. Of course, a posteriori, by (1.6), the two are even equivalent. 

As remarked above, not every $S$ has a tiger.
In \S 15-19 we generate a finite set of families of
surfaces $\goth F$
which includes all (simply connected) $S$ without tigers. For each $S \in {\goth F}$ we 
have an explicit description of the minimal desingularisation, $\t{S}$,
in terms of blow ups of $\pr 2.$, and we check in each case that
$S^0$ is uniruled, by explicitly exhibiting a dominating family of 
rational curves.  We note here that it is quite possible that
${\goth F}$ is too big, in the sense that some of the  surfaces in
${\goth F}$ actually have tigers, or non-simply connected smooth locus.
Given an explicit 
surface it is usually easier for us to check that the smooth locus is
uniruled, than to check either that the smooth locus is simply connected, 
or that the surface has no tiger, as we have rather robust tools
for the first, and only rather ad-hoc techniques for the second and
third. Of course for (1.3) this approach is sufficient. 

It is perhaps surprising, that even with an explicit description
of $\t{S}$, it  can still be very difficult to show $S^0$ is dominated by rational 
curves. Typically the strict transform of these covering curves have 
rather high degree in $\pr 2.$, often over a hundred. 

To prove uniruledness in a specific case, we try to find rational
curves, $Z \subset S$, with an endomorphism lifting to $S^b$. In practical
situations this is only possible when $Z$ meets the singular locus
at most twice, see (4.9.4). Then we try to deform the lifted endomorphism. We
have the following sufficient condition (6.5):

Suppose $Z$ meets the singular locus twice, and the local analytic index 
of its two branches (meeting the singular locus) are $u$ and $v$. Then 
there is a surjection $f:\map {\pr 1.}.Z.$ of degree $lcm(u,v)$ lifting 
to $S^b$, and if 
$$
-K_S \cdot Z \geq 1/u + 1/v
$$ 
then $f$ deforms to cover $S^0$. We usually use the criterion when $Z$ 
has two branches meeting at one point, and $u=v$. In this case
$uZ$ is Cartier, and so the criterion only fails when
$-K_S \cdot Z = 1/u$. Example (6.8) indicates that the criterion
is essentially sharp.

In order to apply the criterion to a particular case, we have to find rational curves meeting the 
singular locus twice. The $Z$ we use occasionally come from interesting geometric configurations
(see for example 17.5.1-3). For concrete example applications of the criterion, see (6.10) and (8.1). 

There are around sixty surfaces in $\goth F$. The generation of ${\goth F}$ 
is one main focus of the paper, and accounts for most of its volume.

{\bf The hunt:} Our construction of $\goth F$ is based on a simple idea. 
To explain it (and the notion of tiger) we use the following:

\definition{1.16 Definition} Let $X$ be a normal surface, and $\Delta$ an 
effective $\Bbb Q$-Weil divisor (that is $\Delta = \sum a_i D_i$
with $a_i \in {\Bbb Q}$, $a_i > 0$). Let $\pi :\map Y.X.$ be any birational 
morphism. Let $\t{\Delta} = \sum a_i \t{D}_i$ be the
strict transform. There is a unique ${\Bbb Q}$-Weil divisor $F$ 
supported on the exceptional locus such that 
$$
K_Y + \t{\Delta} + F = \pi ^*(K_X+\Delta).
$$
$\Gamma = \t{\Delta} + F$ is called the {\bf log pullback} of $\Delta$.
The {\bf coefficient} $e(E,K_X + \Delta)$ 
of any irreducible divisor $E$ on $Y$ is just the coefficient as it appears in 
$\Gamma$, it depends only on $(X,\Delta)$ and the discrete valuation associated to $E$.
The {\bf coefficient} $e(X,\Delta)$  of the pair $(X, \Delta)$ is the largest coefficient
of any divisor (or discrete valuation). We will write $e(X)$ for the coefficient of $(X, \emptyset)$.
\enddefinition

\remark{Remark-Definition} $e(E,K_X + \Delta)$ is just the negative of its 
discrepancy. In particular, $K_X+\Delta$ is log canonical iff the coefficient 
of $(X,\Delta)$ is at most one, and log terminal iff in addition the
coefficient of every exceptional divisor is less than one. It is
Kawamata log terminal iff its coefficient is less than one.
\endremark

We note one trivial, but useful, fact: {\it coefficients are invariant under log
pullback}.

Now let us explain the motivation behind definition (1.13), specifically
the definition of a special tiger for $S$. In view
of case I, when $D$ in (1.3) is empty, one
naturally wonders if there is some curve $C$ on $S$ with $K_S+C$
anti-ample and log terminal. It turns out this fails for infinitely many
families of $S$, so one looks for weaker conditions. 
Note that if $K_S+C$ is anti-ample, then
we can add on some effective $\Bbb Q$-Weil divisor $\beta$ so that
$K_S+C+\beta$ is trivial. A general philosophy of the
log category is to treat exceptional divisors and divisors on
$S$ uniformly. Applying the philosophy to $\alpha  = C + \beta$,
leads to the notion of a special tiger (for $K_S$)  which can be defined as
an effective $\alpha$ with $K_S + \alpha$ numerically trivial, 
with coefficient at least one. We note also the following equivalent formulation
of the definition of tiger (which follows readily from the definitions):

\proclaim{Lemma} Let $E$ be an exceptional divisor over $S$. Let $\map T.S.$ be the
extraction (of relative Picard number one) of $E$. $E$ is a tiger iff
$-(K_T + E)$ is effective.
\endproclaim

Since a  tiger is a divisor of coefficient at least one,
in hunting for a tiger, it is natural to extract the exceptional divisor from the
minimal desingularisation which has 
maximal coefficient. This is the idea behind the hunt. 
We will use the tiger/hunt metaphor in
various notations throughout the paper, occasionally to a tiresome degree.

Beginning  with $(S_0,\Delta_0) = (S, \emptyset)$,
we inductively construct a sequence of pairs $(S_i,\Delta_i)$ of a rank one
log del Pezzo with a boundary, such that:
\roster
\item $-(K_{S_i}+\Delta_i)$ is ample.
\item If $(S_i,\Delta_i)$ has a tiger, then so does $S$.
\endroster

The construction is by a sequence of a $K$-positive extraction (that is a blow up)
$f_i:\map T_{i+1}.S_i.$ followed by a $K$-negative contraction, $\pi_{i+1}$, each
of relative Picard number one.
$\pi_{i+1}$ is either a $\pr 1.$-fibration, or a blow down
$\pi_{i+1}:\map T_{i+1}.S_{i+1}.$. In the first case we say $T_{i+1}$ is a 
{\bf net} and the process stops. We give the details in \S 8.
Such sequences are frequently studied in the MMP, see for example \cite{Shokurov93}.
The only choice in the sequence is which divisor is extracted by $f_i$. The hunt is a sequence 
given by extracting an exceptional divisor $E_{i+1}$, of the minimal desingularisation
of $S_{i}$, for which the coefficient $e(E_{i+1},K_{S_i} + \Delta_i)$ is maximal.

To generate the collection ${\goth F}$ ( and complete the proof of (1.3) ) 
we classify all possibilities for the hunt for which we are unable to find 
a tiger. Let us give a few remarks to explain why such a classification
is possible:

If the coefficient $e(S)$ is sufficiently close to one 
(cf. (21.1)), then by (5.4) of \cite{Kollar92b}, $E_1$ is a tiger. 
By (9.3) the collection of $S$ with $e(S) < 1 - \epsilon$ is bounded. Thus in all but a 
bounded number of cases, the hunt finds a tiger at the very first step, 
and what is needed is an efficient means of dealing with the exceptions. Our
choice for the hunt, that is always extracting a divisor of maximal coefficient (which
is a natural choice, from the point of view of tigers) 
turns out to have remarkably strong geometric consequences. We will explain this
in considerable detail in the introduction to \S 8. It is these consequences which make
feasible an explicit classification of the exceptional cases. In (8.4.7) we give a 
detailed breakdown describing possibilities for the hunt. We then complete the proof by 
analysing each of the possibilities. 

\subhead 1.17 Classification of all but a bounded collection of $S$ \endsubhead

As discussed above, a detailed analysis of the hunt yields a collection ${\goth F}$ 
containing all $S$ (with $S^0$ simply connected) without tiger. We introduced the hunt for
exactly this purpose. Somewhat surprisingly, the hunt is also a useful tool for
classification at the other extreme:

In \S 23 we classify
rank one log del Pezzos $S=S_0$ (no assumption on fundamental group) such that $E_1$
(of the first hunt step) is a tiger. As we remarked above, this includes all $S$
with sufficiently large coefficient $e(S)$, which in turn includes all but a bounded
collection of $S$. 

Our classification, which is independent of the hunt analysis \S 14--19,
is of the following sort:

First we classify abstract pairs 
$(S_1,A_1)$ of a rank one log
del Pezzo surface containing an integral rational curve $A_1$ such that
$K_{S_1} + A_1$ is anti-nef, and log terminal at singular points of $S_1$.
They fall into a short number of series. We note that a similar classification of 
pairs  is obtained in \cite{MT84a} and \cite{MT84b}. Our argument is based on quite 
different ideas, and is considerably simpler.

We apply this classification to the first hunt step.
Let $A_1 \subset S_1$ be the image of $E_1$ (assuming $\pi_1$ is birational).
It is easy to show that if $E_1$ is a tiger, then either $T_1$ is a $\pr 1.$-fibration, 
or $(S_1,A_1)$ is as
in the preceding paragraph. The first case is easy to classify. 
To classify $S$, it remains to classify 
possibilities for the transformation $S_1 \dasharrow S$.
This amounts to classifying possibilities for $\pi_1:\map T_1.S_1.$
such that $E_1$ (the strict transform of $A_1$) is contractible, and contracts
to a log terminal singularity. We indicate how this can be done in \S 23. It
is easy and elementary, but we do not actually list the possibilities, as it would
be notationally too involved. 

Observe that the existence of such a classification is at least to some degree 
counter-intuitive. One might have expected a simple classification of $S$ with mild 
singularities, with a progressively less tractable list of possibilities as more 
complicated singularities are allowed. Indeed, $\pr 2.$ is the only smooth $S$, 
Gorenstein $S$ are classified in \cite{Furushima86}, and log del Pezzos with 
index (of $K$) at most two are classified in \cite{AlekseevNikulin89}.  
However, comparison of our two main classification results -the collection $\goth F$, 
and our classification of $S$ with large coefficient-gives an indication in the 
opposite direction: Surfaces of small coefficient, though bounded, appear 
(at least to us) rather sporadic, while surfaces with sufficiently large coefficient 
exhibit relatively uniform behaviour.

Given the effectiveness of our techniques at either singularity extreme, we believe 
that repeated application of the same methods would eventually yield a complete 
classification of rank one log del Pezzo surfaces. 


Several people have asked the following question; is it true that $\{K_S^2\}$ (over all log del 
Pezzos) satisfies ACC for bounded subsequences ? The question is suggested by analogy with a 
result of \cite{Alexeev94} that the corresponding set over minimal log terminal surfaces of 
general type satisfies DCC. The answer to this question is a resounding no, see (22.5).

Thanks: We would like to thank A. Baragar, A. Corti, D. Morrison, S. Mori,
P. Deligne, D. Huybrechts, T. Hsu, M. Lustig, R. Morelli, Y. Petridis and V. Shokurov for 
advice and suggestions. We would like to particularly thank J\'anos Koll\'ar, 
who suggested (1.5) as a natural companion to Log Abundance, and who has 
assisted us considerably throughout our research.

\heading 1.18 Logical Structure of Proof of (1.1) \endheading
\heading Assume $\kappa(K_X + B) = -\infty$. Goal: $(X,B)$ log
uniruled \endheading
$$
\aligned
&\boxed {\aligned &\text {Run $(K_X + B)$-negative MMP} \\            
                  &\text {$\map (X,B).(S,D).$, $S$ rank one ldp } \endaligned} \\
&\qquad\qquad \Downarrow \\
&\left (\aligned &\text {Goal: $(S,D)$ log } \\
                 &\text {uniruled; two cases: }  \endaligned \right )  \Rightarrow 
                 \boxed{\aligned &\text {If $D \neq \emptyset$: $(S,D)$ log }  \\
                                &\text {uniruled (6.2) \qed}\endaligned }\Leftarrow 
                  \left\{ \aligned &\text {Bug-Eyed cover, \S 4,  and} \\
		       &\text {Def. Theory \S 5 } \endaligned \right.  \\
&\qquad\qquad \Downarrow \\
&\left (\aligned &\text {If $D = \emptyset$: Goal $S^0$  } \\ 
  		 &\text {uniruled; Two cases:} \endaligned \right ) \Rightarrow 
	          \boxed{\aligned &\text {If $K_S$ has a tiger, $E$:}  \\
                                 &\text{transform $(S,\emptyset) \dasharrow (S_1,E)$} \\
                                 &\text{$S_1^0 \setminus E$ dominated by ${\Bbb A}^1_*s$ } \\
                               &\text {$\Rightarrow S^0$ uniruled. (6.1) \qed} \endaligned } 
                  \Leftarrow \left \{\aligned &\text {Bug-Eyed Cover, } \\
		  &\text {Def. Theory, and }  \\
                  &\text {Gorenstein ldps, \S 3} \endaligned \right.    \\
&\qquad\qquad \Downarrow \\
&\boxed {\aligned &\text{ If $K_S$ has no tiger:} \\
                  &\text{ \hskip.25truecmRun the Hunt } \endaligned}  \\
&\qquad\qquad \Downarrow \\
&\boxed {\aligned &\text{ $S \in \goth F$ } \\
                 &\text{ $\exists$ special rat $Z \subset S$} \endaligned} \Rightarrow 
       \boxed {\text{$S^0$ dominated by rats $\simeq mZ$ \qed}} \Leftarrow 
                           \left \{ \aligned &\text{ Def. Theory, } \\
                                            &\text{Bug-Eyed cover, and} \\
      			   &\text{criteria (6.5-6). } \endaligned \right. \\
&\qquad\qquad \Uparrow \\
&\boxed {\aligned &\text{ Analysis of Hunt, \S 8-19} \\
                  &\text{ See flow chart 8.0.16} \endaligned}
\endaligned
$$

\heading \S 2 Glossary of notation and conventions \endheading

If $S$ is a normal surface we indicate by $\t{S}$ its minimal 
desingularisation. If $C \subset S$ is an effective divisor, then 
$\t{C} \subset \t{S}$ will indicate its strict transform. 

$\Hz n.={\Bbb P}(\ring. \oplus \ring.(n))$ denotes the unique minimal rational ruled 
surface, with a curve $\sigma_{\infty} \subset \Hz n.$ of self-intersection
$-n$, and $\Hbz n.$ denotes the log del Pezzo surface of rank one, obtained
by contracting this curve. Let $\sigma_n \subset \Hz n.$ indicate a section disjoint from 
$\sigma_{\infty}$ (all such sections are linearly equivalent, and have self-intersection $n$). 

Where no confusion arises, we will use the same notation to
indicate a divisor, and its strict transform under a 
birational transformation. For a singular point $p \in S$, let 
$\rho(p)$ count the number of exceptional divisor of $\map \t{S}. S.$ lying over $p$.
On the other hand let $\rho$ be the relative Picard number of $\map \t{S}.S.$. Clearly
$$
\rho=\sum _{p\in S}\rho (p).
$$

By a curve {\bf over} a point $p \in S$, we mean a curve in
the exceptional locus of $\map \t{S}.S.$ over $p$.

By a {\bf $k$-curve} on a normal surface $T$ we mean a complete
curve whose strict transform on $\t{T}$ has self-intersection
$k$. Note that if the curve is $K_T$-negative, and $k \leq 0$, then by
adjunction, its strict transform is a smooth rational curve.

A ${\Bbb Q}$-Weil divisor  $\sum a_i D_i$ is called
a {\bf subboundary} if $a_i \leq 1$, and a {\bf boundary}
if $0 \leq a_i \leq 1$. It is called {\bf pure} if
$a_i < 1$. $\sum a_i D_i$ is called reduced if all the non-zero
$a_i$ are equal to one.

When we express a ${\Bbb Q}$-Weil divisor as $\sum a_i D_i$ we 
assume, that the $D_i$ are distinct and irreducible.

By a component of the ${\Bbb Q}$-Weil divisor $\Delta = \sum a_i D_i$, with
$a_i \geq 0$, we mean a $D_i$ with $a_i > 0$. The support of $\Delta$,
$\Supp(\Delta)$, is the union of its components (with reduced structure). 

We will say that a $\Bbb Q$-Weil divisor  $\alpha$ is {\bf effective}, 
if there is an integer $m > 0$ such that $m \alpha$ is integral and 
$|m \alpha|$ is non-empty.

By $>$ or $<$ we mean strict inequality.

We write $\sum a_i D_i \geq \sum b_i D_i$ if $a_i \geq b_i$ for all
$i$. We write $\sum a_i D_i > \sum b_i D_i$ if 
in addition $a_i > b_i$ for some $i$.

A log resolution of a pair $(X,D)$ of a variety $X$ and a reduced
divisor $D$ is a birational map $f:\map Y. X.$ with divisorial exceptional
locus $E$, such that $Y$, and all components of $E + \t{D}$ are smooth,
and $E + \t{D}$ has normal crossings. A log resolution of
$(X,\Delta)$ is a log resolution of $(X,\Supp(\Delta))$. 
$(X,D)$ is called log smooth, if $X$ and every component of $D$ is smooth,
and $(X,D)$ has normal crossings.

Throughout the paper, unless otherwise noted, by a divisor, we
mean a $\Bbb Q$-Weil divisor, and by a component we mean an irreducible
component.

We will make frequent, and occasionally unremarked use
of the classification of two dimensional quotient
singularities, see Chapter 3 of \cite{Kollaretal}. For the readers
convenience we recall the most important facts in Appendix L. We also make 
frequent use of the formulas of Chapter 3 of \cite{Kollaretal} for the index
and discrepancies. We will often refer to a cyclic singularity as a 
{\bf chain} singularity, as the resolution graph is in this case a chain. 
For a non-cyclic singularity, there is exactly one vertex in the resolution 
graph meeting three edges. We call this the corresponding exceptional
divisor the {\bf central divisor}. Three chains meet this vertex. We
refer to these as the {\bf branches}. In a branch, we call the curve
meeting the central divisor, the {\bf adjacent} curve, and the other end
of the chain, the {\bf opposite} curve. 

We will say that a chain singularity is {\bf almost Du Val } if
its resolution graph has only one vertex of weight other 
than $2$, that of weight $3$ and occurring at the end of
the chain, that is the chain has form $(2,2,2,\dots,3)$ (where we
also allow $(3)$).

By the index of a quotient singularity, we mean the order of the
local fundamental group.

We will make frequent use of the 
main results of the (log) Minimal Model Program:
the cone theorem, the contraction theorem,
and the log abundance theorem. The first two of these are known in
all dimensions, and the last up to dimension three. However we only
use them in dimension two. For proofs in this case see
\cite{MT82}. For overviews of the general theory (as well
as further references) see \cite{CKM88}, \cite{Kollaretal}, \cite{KMM87}. We will
also make use of the standard definitions and notations of the program, as
in \cite{Kollaretal}.

We will say a pair $(X,D)$ of a normal variety and a reduced divisor
is {\bf log uniruled} if $X^0$ is covered by rational curves meeting
$D$ at most once (that is through a general point of $x$ there is
a complete rational curve, contained in the smooth locus, and meeting
$D$ at most once).  We will say $X$ is log uniruled if $(X,\emptyset)$
is log uniruled, or equivalently, if $X^0$ is uniruled, 
that is through a general point of $X$, there is a 
complete rational curve, contained in the smooth locus.



\heading \S 3 Gorenstein del Pezzo Surfaces \endheading

In this section we collect a number of results on
Gorenstein log del Pezzo surfaces that we will use at various points
in the paper. Especially useful will be (3.6-8) which together give a 
complete and simple picture of rank one Gorenstein log del Pezzo surfaces
whose smooth locus is algebraically simply connected.

{\bf Notation: Throughout the section $S$ denotes a rank one Gorenstein log del Pezzo.}

In \cite{Furushima86} possible singularities of $S$
are classified. We will make frequent use
of this classification, which we refer to as 
{\bf the list}. We  write for example $S(A_1 + A_ 3 + A_5)$ for a rank one Gorenstein
del Pezzo with those singularities, and for example
$\t{S}(A_1 + A_ 3 + A_5)$ for its desingularisation. 

For the reader's convenience, here is a copy of
the list:
$$
\align
&A_1,\quad A_1+A_2,\quad A_4,\quad 2A_1+A_3,\quad D_5,\quad A_1+A_5,\quad 
3A_2,\quad E_6,\quad \\
&3A_1+D_4,\quad A_7,\quad A_1+D_6,\quad E_7,\quad A_1+2A_3,\quad 
A_2+A_5,\quad \\
&D_8,\quad 2A_1+D_6,\quad E_8,\quad A_1+E_7,\quad A_1+A_7
,\quad 2A_4,\quad A_8,\quad \\
&A_1+A_2+A_5,\quad A_2+E_6,\quad A_3+D_5,\quad 4A_2,\quad 2A_1+2A_3,\quad 2D_4. \\
\tag{3.1}
\endalign
$$
By \cite{MZ88}, the following are the subset of possibilities with
$S^0$ algebraically simply connected, and any rank one Gorenstein del Pezzo
with these singularities has algebraically simply connected smooth locus:
$$
A_1,\quad A_1 + A_2,\quad A_4,\quad D_5,\quad E_6,\quad E_7,\quad E_8.
\tag{3.2}
$$
We will abuse notation and call this the
{\bf simply connected list}. 

We note that in general $S$ is not
determined by its singularities. However, for the simply connected list,
this is almost the case. See (3.10) below.

The $K$-negative MMP for Gorenstein surfaces is almost as simple
as that for smooth surfaces. 

\proclaim{3.3 Lemma} Let $g:\map T.W.$ be a birational $K_T$-contraction, of
relative Picard number one, with (irreducible)
exceptional divisor $\Sigma$, such that $T$ is Gorenstein along $\Sigma$. Along $\Sigma$, 
$T$ is either smooth, or has a unique singularity. In the latter case the 
singularity is an $A_r$ singularity, with $K_T+\Sigma$ log terminal, and the singularity 
is removed by the contraction. Furthermore the induced map 
$\map \pi_1^{alg}(T^0).\pi_1^{alg}(W^0).$ is an isomorphism.
\endproclaim
\demo{Proof} On the minimal desingularisation, $\Sigma$ is a $-1$-curve. 
Contractibility considerations give the description of the possible 
singularities. The last claim follows from (7.3). \qed \enddemo

Gorenstein $\pr 1.$-fibrations (of relative Picard number one) are
also quite simple. Here we state the possibilities. We will prove
a more precise result in \S 11, see (11.5.4). 

\proclaim{3.4 Lemma} Suppose $\pi: \map T. {C}.$ is a $\pr 1.$-fibration, of relative 
Picard number one, and $G \subset T$ is a fibre, contained in the Du Val locus of $T$. 
One of the following holds:
\roster
\item $T$ is smooth along $G$.
\item There are exactly two singularities, $A_1$ points, along $G$. $K_T + G$
is lt. 
\item There is a unique singularity, a $D_n$ point along $G$. 
\endroster
\endproclaim
\demo{Proof} Let $\map \t{T}.T.$ be the minimal desingularisation, and 
let $h:\map \t{T}.W.$ be a relative minimal model of $\map \t{T}.T.$,
thus $\map W.C.$ is smooth, and $h$ is a composition of blow ups at smooth
points. Because $\pi$ has relative Picard number one, each blow up in $h$ is at
a point along the $-1$-curve of the previous blow up. Now one checks easily
that the proposition gives all possibilities with Du Val singularities. For more
details see (11.5). \qed \enddemo

\proclaim{3.5 Lemma} If $T$ is a Gorenstein del Pezzo surface,
then $|\pi_1^{alg}(T^0)|\leq 8$. 
\endproclaim
\demo{Proof} Note if $\map T'.T.$ is a degree $n$ cover, \'etale in 
codimension one, then $T'$ is again a 
Gorenstein del Pezzo, and $K_{T'}^2 = n K_{T}^2 \leq 8$. 
Since $K_{T}^2$ is a positive integer, the result follows. \qed \enddemo

Recall that throughout the paper, a proper $K_T$-negative curve
on a normal surface $T$ is called 
a $-1$-curve, if $\t{C} \subset \t{T}$ is a $-1$-curve in the
usual sense. Note when $T$ is Gorenstein, such a $C$ is a $-1$-curve iff $C$ is rational,
$K_T \cdot C = -1$ and $\t{C}$ is smooth.

\proclaim{3.6 Lemma} Let $S = S(E_8)$. $S$ contains a unique $-1$-curve $D$. 
$\t{D} \subset \t{S}$ meets a unique exceptional curve, the opposite end of
the $A_4$ chain, with normal contact. 
$|-K_S|$ is one dimensional, and has a unique
basepoint, a smooth point of $S$. There are two possibilities for the collection
of rational members of $|-K_S|$. Either
\roster
\item $|-K_S|$ has exactly three rational members, $D$ and
two integral nodal rational curves $N_1$, $N_2 \subset S^0$, or
\item $|-K_S|$ has exactly two rational members, $D$, and 
a unique integral cuspidal rational curve, $C \subset S^0$.
\endroster
\endproclaim
\demo{Proof} Let $S=S(E_8)$ be any such surface. Obviously $S$ contains some $-1$-curve, $D$. 
$K_S^2 =1$, so $|-K_S|$ is one dimensional
by Riemann-Roch, and every member of $|-K_S|$ is reduced and irreducible.
$|-K_S|$ and has a smooth elliptic member $E \subset S^0$ by \cite{HW81}.
$\ring E.(E) = \ring E.(q)$ for a unique $q \in E$. 
Since $H^1(\ring S.) =0$, $q \in S^0$ is the unique basepoint of $|-K_S|$.
Let $\map T.S.$ blow up $q$. $|-K_S|$ yields an elliptic 
fibration $g:\map T.{\pr 1.}.$. Let $\t{g}:\map \t{T}.{\pr 1.}.$ be the
induced map. Since $S^0$ is simply connected, any $-1$-curve on $S$
is a member of $|-K_S|$. Since any $-1$-curve must pass through the singular point,
and $g$ has irreducible fibres, $D$ is the unique $-1$-curve. The fibre of $\t{g}$
containing $D$ has $9$ irreducible components, thus by 
Kodaira's classification of singular fibres, see page 150 of \cite{BPV84}, 
the fibre containing $D$ is
$\t{E}_8$, thus $D$ meets the opposite end of the $A_4$ chain as required. 
Note $e(\t{T})$ (the topological Euler characteristic) is $12$. By formula 
(11.4) on page 97 of \cite{BPV84}, the additional singular fibres (which we
know are reduced and irreducible) are either 
exactly two nodal rational curves, or exactly one cuspidal rational curve. 
\qed\enddemo

\proclaim{3.7 Lemma} There are exactly two isomorphism classes of $S(E_8)$ corresponding
to the two possibilities in (3.6).
\endproclaim
\demo{Proof} Let $S = S(E_8)$. Let $D$ be the unique $-1$-curve of (3.6). Let $B \subset S^0$
be any member of $|-K_S|$. Let $q \in B^0$ be the unique basepoint of $|-K_S|$, see (3.6).
Let $L$ be the $-2$-curve of the $A_1$ branch of the 
$E_8$ singularity. Let $\map T.S.$ extract $L$. $T$ has an $A_7$ singularity,
and $D \subset T$ contracts $\pi:\map T.{\pr 2.}.$, the image of $L$ is a flex line
to the cubic curve $B \subset \pr 2.$ at (the image of the) $q$.
The induced map $\map \t{S}.{\pr 2.}.$ is obtained by
blowing up $8$ times along $B$ over $q$. $B \subset \pr 2.$ is embedded by the full
linear system $|3q|$, thus 
$S$ depends only on $(B,q)$. The automorphism group of $B$ acts transitively on $B^0$, so
$S$ depends only on $B$. By (3.6) we 
can take for $B$ either a cuspidal rational curve $C$, or a nodal rational curve $N$. 
\qed \enddemo

\proclaim{3.8 Lemma} Assume $S^0$ has trivial
algebraic fundamental group. If $S$ is not $\pr 2.$ or $\Hbz 2.$ then
$S$ contains a unique $-1$-curve, $D$. $\t{D}$ has normal
crossings with the exceptional locus of $\map \t{S}.S.$, and meets
exactly one exceptional curve (of the minimal desingularisation)
over each singular point of $S$. If $S \neq \pr 2.$
then there is a $-2$-curve, $E$, of $\t{S}$, such that extracting
$E$ gives a $\pr 1.$-fibration, with $E$ a section, and $D$ (in case $S \neq \Hbz 2.$)
the only multiple fibre. $E$, and the singularities of $D$ are as
follows:
\roster
\item If $S= S(A_1) = \Hbz 2.$ then $E$ is the unique $-2$-curve and 
the fibration is smooth.
\item If $S=S(A_1 + A_2)$, then $K_S + D$ is log terminal, $E$ is the
$-2$-curve over the $A_2$ point which is disjoint from $\t{D}$.
\item If $S = S(A_4)$ with singularity $(2,2',2,\ul{2})$, then $D$
meets the primed curve and $E$ is the underlined curve.
\item If $S = S(D_5)$ then $D$ meets one of the $(2)$ branches,
and $E$ is the opposite (from the central curve) end of
the $A_2$ chain.
\item If $S= S(E_k)$ ($8 \geq k \geq 6$) then $D$ meets the opposite
end of an $A_{k-4}$ chain (this chain is unique except when $k=6$),
and $E$ is the opposite end of an $A_2$ chain, different, in the
case $k=6$, from the $A_2$ chain which $D$ meets.
\endroster

Furthermore, assume $1 < K_S^2 < 8$.  
Let $\map T.\t{S}.$ blow up a point of $D$ not
on any $-2$-curve. Then $T$ is the minimal desingularisation of a rank one log
del Pezzo, $S'$ with algebraically simply connected smooth locus, and
$K_{S'}^2 = K_S^2 -1$. If we repeat this process $K_S^2 -1$ times, we obtain
$S(E_8)$. The induced map $\map \t{S}(E_8).\t{S}.$ is canonical, contracting
at each stage the unique $-1$-curve. \endproclaim

\demo{Proof} We will prove that there exists a $-1$-curve, $D$, meeting
the singularities as prescribed. The final remarks are  then immediate from the
singularity description, and imply the uniqueness of the $-1$-curve by (3.6).
One also checks easily 
that extracting the indicated $-2$-curve $E$ gives the required $\pr 1.$-fibration.
Hence it is enough to prove the existence of $D$.

We can assume by (3.6) that $K_S^2 \geq 2$. 

Note if $C \subset S$ is any $-1$-curve, then $(K_S + C) \cdot C < 0$
and thus $C$ is smooth, and so meets at most one exceptional divisor over
each singular point, and the contact is normal (see for example (6.11)).

Let $f:\map T.S.$ extract a $-2$-curve which according to the statement of
the lemma is to have contact with $D$. In the case of $S(A_1 + A_2)$ let $V$ be 
the $-2$-curve of the $A_1$ point, otherwise choose any of the possible curves,
that is on $S(D_5)$ either of the $A_1$ branches, on $S(E_6)$ the opposite end
of either $A_2$ chain, and on $S(A_4)$ either of the interior $-2$-curves (in
the other cases $V$ is unique). Suppose first that $T$ has a $\pr 1.$-fibration. 
By (3.4) this is only possible if $S= S(E_6)$. In this case one checks
that the fibre through the $D_5$ point of $T$ is a $-1$-curve meeting the singularities
as prescribed (for details see (11.5.4)). So we can assume $T$ has a birational
contraction $\pi:\map T.S_1.$ of a $-1$-curve, $D$. Using (3.3) and the
list, one checks in each case that $D$ meets the singularities as prescribed:

If $S=S(D_5)$, $T$ has an $A_4$ singularity.
By (3.3) and the simply connected list $S_1$ is either $\pr 2.$,
$D \subset T$ contains the $A_4$ point, and $K_T + D$ is log terminal,
or $S' = S(A_4)$ and $D \subset T^0$. In the second case $D \subset S$
has the prescribed singularities. In the first either the 
singularities are as required, or $D$ meets the opposite end of
the $A_3$ chain. But the latter is impossible, for then $D^2 =0$ on
$S$. 

If $S=S(A_4)$. $T$ has singularities $A_1 + A_2$.
Either $D \subset T^0$, and $S_1 = S(A_1 + A_2)$, or $D$ meets one end of
the $A_2$ chain, and $S_1 = S(A_1)$. In the first case, either the singularities
are as prescribed, or $D \subset S$ meets an end of the $A_4$ chain. But the
latter is impossible, for in that case $D^2 < 0$ on $S$. 

If $S = S(A_4)$, $T$ has a single singularity, an $A_2$ point. $D \not \subset T^0$,
from the list, so $D \subset S$ has the prescribed singularities.

In the cases, $S(E_6)$, $S(E_7)$, 
$T$ has a single singularity, a non-cyclic singularity,
so by (3.3) $D \subset T^0$, and the singularities are as prescribed.
\qed \enddemo

\proclaim{3.9 Corollary} Let $W$ be a Gorenstein log del Pezzo surface of 
rank one, which contains two distinct rational curves $D_1$ and $D_2$ 
with $K_W \cdot D_i = -1$.  The following implications hold:
\roster
\item If $W$ is simply connected then $W = S(E_8)$, one of the 
$D_i$ is contained in $W^0$, and is a rational curve of arithmetic genus one.  
$D_1 \cap D_2$ is a single smooth point of $W$, 
the unique basepoint of the linear series $|-K_W|$. 
\item If $K^2_W\geq 4$, then $K^2_W=4$, $W=S(2A_1+A_3)$ 
and $D_1$ and $D_2$ each pass through one of the $A_1$ points and 
opposite ends of the $A_3$ point.
\endroster
\endproclaim
\demo{Proof} (1) is immediate from (3.6), (3.7) and (3.8).

Now suppose $K_W^2 \geq 4$. By (1) $W$ is not simply connected, and
so from the (full) list, $S=S(2A_1 + A_3)$. $K_W + D_1 + D_2$ is
anti-ample. In particular $D_i$ is smooth. By adjunction,
they meet at only one point, and neither can contain all three singular points.
Thus each  $D_i$ must contain exactly two singular points, and they cannot both be
$A_1$, for
otherwise $D_i^2 \leq 0$.  Thus each $D_i$ contains one
of the $A_1$ points, and they both contain the $A_3$ point. The described
configuration is now the only possibility, as one checks that otherwise after
extracting an appropriate curve over the $A_3$ point, both have
non-positive self-intersection, and at least one is contractible,
violating the Picard number. 
\qed\enddemo

\proclaim{3.10 Corollary} Let $S$ be a rank one Gorenstein log del Pezzo surface,
with algebraically simply connected smooth locus. $S$ is uniquely determined by
its singularities (and thus by $K_S^2$)  except in the case of $S(E_8)$
(that is $K_S^2 =1$) when there are two possibilities, as described in (3.7).
\endproclaim
\demo{Proof} We will write $S(C)$ (resp $S(N)$) for the $S(E_8)$ of (3.6-7) 
with a cuspidal member (resp. nodal member). 

$S(A_1)$ is unique by (3.8). Assume $1 < K_S^2 < 8$.

By the final remarks of (3.8) there are at most two possibilities for $\t{S}$: Start with either
$S(C)$ or $S(N)$ and contract repeatedly the unique $-1$-curve, $K_S^2 -1$ times.
We will show that the two possibilities give the same surface. Start with $S(C)$. One obtains 
$S$ with $C \subset S^0$, $C \in |-K_S|$. Let $f:\map {\pr 1.}.C.$ be the normalisation.
$\Hom(\pr 1., S)$ has dimension at least $4$ at $f$ by (II.1.2) of \cite{Kollar96}, 
and thus deforms to give a dominating family of rational curves. By (II.3.14) of \cite{Kollar96}  
the general member has only nodes, and so gives a nodal rational curve $N \in |-K_S|$, 
$N \subset S^0$. Choosing the point $x$, of the final remarks in (3.8), to be (repeatedly) 
$N \cap D$, we see we can also obtain $S$ starting with $S(N)$. \qed \enddemo


\def\aos{{\Bbb A}^1_*}
\heading \S 4 Bug Eyed Covers. \endheading
\subhead Brief Introduction \endsubhead

As indicated in the introduction, a key tool in the proof of (1.1) is 
Koll\'ar's bug-eyed cover. For a rank one log del Pezzo $S$, this is a
universal homeomorphism $\flat:\map S^{\flat}.S.$, with
$S^\flat$ a 
non-singular (but if $S$ is singular, non-separated) algebraic space. 
Informally, $S^{\flat}$ allows us
to do deformation theory on $S$, {\it as if $S$ were smooth}. More precisely,
we have a good deformation theory for a map $f:\map C.S.$ (from a proper smooth curve)
provided $f$ factors through ${\flat}$. Here our main goal is to present
criteria for the existence of such a lifting. We will consider the question of
how to deform $f$ in \S 5 and \S 6.

We note at the outset that all of the important ideas in this section are
due to Koll\'ar, and most of the statements below can be found implicitly
or explicitly in \cite{Kollar91}. We present them here rather than quote Koll\'ar
directly, as it will be convenient to have them in a slightly different form. 
Our only original contribution to the general theory is the uniqueness of the 
bug-eyed cover of a variety with isolated quotient singularities, (4.4). This is
stated in \cite{Kollar91}, but without proof.

The reader who is familiar with Koll\'ar's construction, and willing to take uniqueness
on faith, can skip directly to (4.5).

One point of notation: or us a curve is assumed to be a separated variety.
Recall that everything in this paper takes place over ${\Bbb C}$. 

We begin by recalling Koll\'ar's definition:

\definition{4.0 Definition} Let $f:\map Y.X.$ be a map between algebraic 
spaces. Let $U \subset Y$ be the locus where $f$ is \'etale. Assume $U$ is 
separated. The {\bf bug-eyed cover of $X$ associated to $f$} is the algebraic 
space $\map Y/R_f.X.$, where $R_f$ is the equivalence relation 
$R_f=(U \fb X. U \setminus \Delta_U)\cup \Delta_Y$. Let $B=U^c$.
\enddefinition

\remark{4.1 Remark} Informally, $\map Y/R_f.X.$ is obtained from 
$\map Y.X.$ by gluing together inverse image points where $f$ is \'etale. In 
particular, if $f$ is \'etale at every inverse image of $x$, then 
$\map Y/R_f.X.$ is an isomorphism over a neighbourhood of $x$. It is a bit
surprising to us  that such a simple idea turns out to be so useful.
\endremark

\proclaim{4.2 Lemma} Let $g:\map Y'.Y.$ be a separated \'etale surjection, 
which is one to one on closed points over $B$. Then
\roster 
\item $R_{f \circ g}=R_f\fb Y\times Y.(Y' \times Y')$, 
\item the induced map $\map Y'/R_{f \circ g}.Y/R.$ is an isomorphism, and
\item if further $f$ is an isomorphism on $U$ and $f$ is bijective on closed
points, then $\map Y.X.$ is the bug-eyed cover associated to $f\circ g$,
$Y'/R_{f \circ g}$.
\endroster
\endproclaim
\demo{Proof} Let $\amalg$ indicate a disjoint union into 
open and closed subsets. Note $U' = g^{-1}(U)$ ($U',B' \subset Y'$ are
defined with respect to $f \circ g$ in the obvious way). 
Since $g$ is separated and \'etale, there is a decomposition 
$Y'\fb Y. Y' = A \amalg \Delta_{Y'}$. By assumption
$B' \fb Y. B' \subset \Delta_{Y'}$, and so
$U' \fb Y. U' = A \amalg \Delta_{U'}$.  Thus we have
$$
\align
R_f \fb {Y \times Y}. (Y' \times Y') &= (U' \fb X. U' \setminus
U' \fb Y. U') \amalg Y' \fb Y.Y' \\
&= (U' \fb X. U' \setminus \Delta_{U'}) \amalg \Delta_{Y'} = R_{g \circ f}.
\endalign
$$
Hence (1). Since $g$ is an \'etale surjection, (2) follows. 

 Note that when $f$ is an isomorphism on $U$ and $f$ is bijective on closed 
points, $R_f=\Delta_Y$. This gives (3). \qed\enddemo

Now we turn to uniqueness. 
 We first consider what happens locally. Note that when $X$ is two dimensional, 
the case that most interests us, the local analytic possibilities for
$f$ are classified in \cite{Brieskorn68}.

\proclaim{4.3 Lemma} Let $(X,x)$ be a pointed normal variety of dimension 
$n\geq 2$. Let $f:\map (W,w).(X,x).$ be a quasi-finite surjection from a smooth 
algebraic space, with $w=f^{-1}(x)$ (as sets), \'etale away from $w$. Let $G$ be 
the local fundamental group of $X$ at $p$. 
\roster 
\item If $W$ is a separated variety, then locally analytically around $x$ and 
$w$, 
$$ 
f: \map W \setminus \{w\}.X \setminus \{x\}.
$$
is the universal cover. 
\item The extension of Hensilisations $\ring X,p.^h \subset \ring W,w.^h$ 
induces a Galois extension of quotient fields, with group $G$, and 
$\ring X,p.^h = (\ring W,w.^h)^G$. 
\endroster
\endproclaim
\demo{Proof} Replacing $W$ by an \'etale cover with a single point
lying over $w$, (which doesn't effect $\ring W,w.^h$) we may assume $W$ is an 
affine variety. Locally analytically 
$f:\map W \setminus \{w\}.X \setminus \{x\}.$ is a cover, and since 
$W \setminus \{w\}$ is simply connected (locally analytically), it is the 
universal cover. Hence (1). The analogous statement to (2), for analytic local 
rings follows. (2) now follows by (4.11) of \cite{Artin69}.
\qed\enddemo

\proclaim{4.4 Proposition} Let $X$ be a normal variety with isolated quotient 
singularities. Then there is a unique smooth algebraic space $f:\map W.X.$ 
such that $f$ is an isomorphism over $X^0$, and one to one on closed points. 
The map of Hensilisations is described by (4.3.2). $f$ is a universal
homeomorphism. In particular, when $X$ is complete, $W$ is universally closed.
\endproclaim
\demo{Proof} By definition there is a surjection $\map Y.X.$, 
\'etale outside the singular locus, where $Y$ is a smooth variety. By dropping 
points of $Y$ we can assume the map is one to one on points over the singular 
locus. For the existence of $W$ we take the associated bug-eyed cover. 

Uniqueness can be checked Zariski locally on $X$. Thus for the rest of the proof
we work in the category of pointed spaces. Suppose $f_i:\map (W_i,w_i).(X,x).$, 
$i=1,2$, are two candidates. By (4.2) there is a smooth affine variety $(W,w)$ 
and a pair of \'etale surjections $p_i:\map (W,w).(W_i,w_i).$. Dropping closed 
points, we can assume $w=p_i^{-1}(w_i)$. Since $X$ is separated, 
$f_1\circ p_1=f_2\circ p_2$. By (4.2.3), $W_i = W/R_{f_i \circ p_i}$.

Since $f$ is a set bijection, to show it is a universal homeomorphism, it is enough to show it 
is universally open. This is a local \'etale question on $W$ and $X$, thus we reduce to the case 
of a quotient by a finite group, where the result is known. Alternatively, one can apply 
Chevalley's criterion, see (15.4) of \cite{EGAIV}. \qed\enddemo

\definition{4.5 Definition} We will denote the unique $W$ of (4.4)
by $\flat: \map X^{\flat}.X.$, and refer to it as {\bf the bug-eyed cover of }
$X$. \enddefinition

We will be principally interested in the bug-eyed cover of a log del
Pezzo surface, $S$.

Next we turn to lifting criteria.  

\proclaim{4.6 Lemma} Let $c:\map C'.C.$ be a surjective generic isomorphism from
a non-singular one dimensional algebraic space to a smooth curve.
For each $Q \in C'$ such that $c$ is not an isomorphism over $p_Q = c(Q)$, 
let $r_Q$ be the ramification index. Then in the notation of \cite{Kollar91}, 
$C'= C\{r_Q,p_Q\}$. \endproclaim
\demo{Proof} Let $D= C\{r_Q,p_Q\}$. There is an obvious set bijection between 
$D$ and $C'$ and one can check locally that this is an isomorphism. Thus
one can drop points from fibres, and assume $c$ is a set bijection. In
this case the result follows from (4.2.3) above, and (4.3) of \cite{Kollar91}.
\qed\enddemo

We now present Koll\'ar's basic lifting criterion in a convenient intrinsic form:

\proclaim{4.7 Lemma} Let $c:\map C'.C.$ be as in (4.6). Let
$f:\map D.C.$ be a finite map, with $D$ a smooth curve. Then $f$ factors through
$c$ iff the following holds: for each point $y \in D$, there is an inverse
image point $x \in c^{-1}(f(y))$ such that the ramification index
$r_c(x)$ divides the ramification index $r_f(y)$. 
\endproclaim
\demo{Proof} Immediate from (4.6), and (4.2.2) of \cite{Kollar91}. \qed \enddemo

\proclaim{4.8 Lemma} Let $C'$ be a smooth one dimensional (not
necessarily separated) algebraic space. Let $C$ be a proper smooth curve,
and let $c:\map C'.{C}.$ be a generically isomorphic surjection.

If $K_{C'}$ has negative degree then $C = \pr 1.$ and there exists an endomorphism $f$ of $C$ 
which factors through $c$. 

If further $c$ is a set bijection, then the converse holds. If $C = \pr 1.$ and
there is an endomorphism $f$ factoring through $c$, then $K_{C'}$
has negative degree.  
\endproclaim
\demo{Proof} By (4.6), $C'= C\{r_i,p_i\}$. 
$K_{C'}$ has 
degree $-2 + \sum \frac{r_i -1}{r_i}$ by the Hurwitz formula.
Now apply (4.4.2) and (6.5) of \cite{Kollar91} (note the remark 
at the end of the proof of (6.5)). \qed\enddemo

\remark{4.9 Remarks}
\roster 
\item Notation as in (4.8). We note that by (4.6) above, (4.4.2) and  
(6.5) of \cite{Kollar91} there is always a quasi-finite surjection $d:\map D.C'.$ with $D$ 
a proper smooth curve, and the degree of a line bundle $L$ on $C'$ can be defined (as
a rational number) as $\deg (d^*(L)) / \deg (d)$.
\item Of course in (4.7) the condition on $f(y)$ holds if there is
a point over $f(y)$ where $c$ is \'etale. Thus one only needs to consider
points in $C$ over which $c$ is {\it no where \'etale}.
\item (4.7) implies in particular that if $C = \pr 1.$ and there are at
most two points in $C$ over which $c$ is no where \'etale, then there is
an endomorphism $f$ of $C$ factoring through $c$. 
\item  In the final paragraph of (4.8), when $C = \pr 1.$, one can express the necessary and 
sufficient condition for the existence of a lifting endomorphism as: There 
are at most three points in $C'$ where the map is not \'etale and if
$r_1$, $r_2$, $r_3$ are the ramification indices (which are allowed to be one), then
there exists a lifting endomorphism iff $2 > \sum \frac{r_i -1}{r_i}$. For three
ramification points, the 
possible $(r_1,r_2,r_3)$ are the so called {\it Platonic triples}, ubiquitous
in the study of lt surface singularities. 
\endroster
\endremark

\proclaim{4.10 Lemma} Let $f:\map C.S.$ be a proper map, and a generic 
embedding, where $C$ is a smooth curve and $S$ is a rank one log del Pezzo 
surface. Let $C'$ be the normalisation of $\flat^{-1}(f(C))$ 
(for the construction of $C'$ see \S 22).
\item{(1)} There is a commutative diagram
$$
\CD 
C'@>>> S^{\flat} \\
@VcVV     @VVV \\
C  @>{f} >> S.
\endCD
$$
\item{(2)} If $\map D.C.$ is a map from a smooth curve such that $\map D.S.$ 
factors through $h:\map D.S^{\flat}.$, then $h$ factors uniquely through
$C'$. 

 Let $p\in S$ and $g:\map (U,y).(V,p).$ express a local analytic neighbourhood
$V$ of $p$ as a quotient of smooth $U$ by the finite group $G$ (so $G$ acts on 
$U$ freely outside of $y$, and $G$ is the local fundamental group).
Fix $q \in f^{-1}(p)$. Let $Z$ be the corresponding local analytic
branch of $f(C)$ through $p$. Let $n$ be the number of analytic branches
of $g^{-1}(Z)$. Let $m$ be the index of $S$ at $p$ (that is the order of $G$)
and let $d$ be the local analytic index of $Z$. Let $t \in C'$ be 
a point of $c^{-1}(q)$. 

\item{(3)} $n$ is the number of points in $c^{-1}(q)$. $n| m$ and $r= m/n$
is the ramification index of $c$ at $t$. $d | r$. If $m$ is prime then 
$d=r=m$ unless $Z$ is Cartier at $p$. 
\endproclaim 
\demo{Proof} (1) and (2) are immediate from (N.3), the universal property
of normalisation. 

 (3) can be considered at the level of strict Hensilisations and so by (4.3)
we can replace $\flat:\map S^{\flat}.S.$ by $g:\map (U,y).(V,p).$ and 
$C'$ by the normalisation of $g^{-1}(Z)$. 

By (5.E.vi) of \cite{Matsumura80} $G$ acts transitively on the branches of 
$g^{-1}(Z)$, and so $n | m$. Clearly each branch of $C'$ has 
degree $m/n$ over $C$, thus $r = m/n$, since $q$ has a unique inverse image 
on each branch.

Let $h$ be a local equation for one branch of $g^{-1}(Z)$. Since $G$ acts 
transitively on the branches of $D$, $N_G(h)$ (the product of the 
translates), is a defining equation for $r \cdot g^{-1}(Z)$. Since 
$N_G(h)$ is $G$ invariant, it's a defining equation for $r Z \subset V$. 
Thus $r Z$ is Cartier, and so $d | r$ by the definition of the index.

Finally, if $m$ is prime, and $d \neq m$, then $d =1$, and $Z$ is 
Cartier. \qed\enddemo

It will be convenient for applications to have a criterion for when $Z$ of (4.10) fails to be 
Cartier. 

\proclaim{4.11 Lemma} Let $(A,m)$ be a one dimensional local Noetherian Domain, a 
$k$-algebra, with $k$ isomorphic to $A/m$. Assume the normalisation $\t{A}$ 
is finite over $A$ and local. Then
$$
\dim_k \frac{m}{m^2} \leq  \min_{h \in m \setminus 0}
\dim_k \frac{A}{h} =\dim_k \frac{\t{A}}{m \cdot \t{A}} .
$$
\endproclaim
\demo{Proof} By the projection formula \cite{Fulton84} 
$$
dim_k \frac{A}{h} = \dim_k \frac{\t{A}}{h \cdot \t{A}}.
$$
Since $\t{A}$ is a local DVR, the second equality follows. 
$$
dim_k \frac{m}{m^2} \leq dim_k \frac{m}{m^2 + h} + 1 
\leq dim_k \frac{m}{h} + 1 = \dim_k \frac{A}{h}
$$
which gives the first inequality.\qed\enddemo

\proclaim{4.12 Lemma} Let $(S,p)$ be a germ of a normal surface,
with a rational singularity, and $g:\map \t{S}.S.$ the minimal 
desingularisation, with reduced exceptional divisor $E$. 
Let $p \in C \subset S$ be an analytically irreducible curve, such that 
the strict transform $\t{C} \subset \t{S}$ is smooth. Assume $E$ is 
$g$-anti-nef.
\roster
\item $\dim T_p(C) \leq \t{C} \cdot E$.
\item $C$ is smooth iff $\t{C} \cdot E =1$. In particular if $\t{C}\cdot E=1$ 
then $C$ is not Cartier.
\item If $\t{C} \cdot E =2$ then $C$ has a hypersurface singularity.
If furthermore $C \subset S$ is Cartier, then $S$ is Du Val.
\endroster
\endproclaim
\demo{Proof} Since $E$ is anti-nef, there is a sum of analytic discs
$W$, meeting $E$ normally, and disjoint from $\t{C}$ such
that $E + W$ is cut out by $g^{-1}(h)$ for some
$h \in m_{S,p}$. (1) follows from (4.11).

Now suppose $C$ is smooth. Then if we take $h$ a generator
of $m_{C,p}$, then 
$$
1 = Z(h) \cdot \t{C} \geq E \cdot \t{C}
$$
(where $Z(h) \subset \t{S}$ is the scheme cut out by $g^*(h)$)
which gives (2). 

For (3) suppose $\t{C} \cdot E =2$. Then by (1-2), $C$ has two
dimensional tangent space and in particular $C$ is Gorenstein. Now if $C$ is 
Cartier, then $S$ is Gorenstein. Hence (3).\qed\enddemo

\remark{4.13 Remarks} Note the conditions of (4.12) hold for any cyclic 
quotient singularity, and for any non-cyclic quotient singularity for which 
the central curve has self-intersection at most $-3$. The
implication $C$ smooth $\implies$ $\t{C} \cdot E =1$ (and its proof)
holds for any curve on a normal surface germ, that is without assuming $E$ is 
$g$-anti-nef. The reverse implication, however, does not. One counter
example is the unique $-1$-curve on $S(E_8)$, see (3.8).
\endremark

\remark{4.14 Remark} The conditions for a pair $(S,D)$ of an irreducible (non-empty)
curve on a rank one log terminal surface to be log Fano have a nice expression
in terms of the bug-eyed cover ${\flat}:\map S^{\flat}.S.$:
\roster 
\item An irreducible curve $D$ on a log terminal surface $S$ is log terminal
iff the pullback  $\flat^{*}(D) \subset S$ is non-singular.
\item If $D$ is log terminal and $S$ has rank one, then $K_S + D$ is anti-ample
iff $D = \pr 1.$ and there exists an endomorphism of $D$ which factors through
${\flat}$.
\endroster
\endremark
\demo{Proof} For (1) see Appendix L. (2) follows from (4.8) (or (4.9.3)) and
the adjunction formula, see  (4.15) below. \qed \enddemo

\remark{4.15 Remark} For log terminal singularities, 
Shokurov's different (see Appendix L)
also has a nice interpretation in terms of the 
bug-eyed cover. Suppose we have an lt pair $(S,C)$, with $C$ reduced and
irreducible. 
Suppose the singular points of $S$ along 
$C$ are $\list p.n.$ and $\list r.n.$ are the 
corresponding orders of the local fundamental groups. $C'= {\flat}^{-1}(C)$ is 
non-singular as $C$ is lt, and $c:\map C'.C.$ is a set
bijection, totally ramified at each $p_i$ with ramification index $r_i$. In
the notation of \cite{Kollar91}, $C' = C\{r_i,p_i\}$.
Let $q_i$ be the unique point of $D$ over $p_i$. 
Since ${\flat}:\map S^{\flat}.S.$ is \'etale in codimension one
$$
\align
{\flat}^*(K_S + C)|_{C'} = (K_{S^{\flat}}+C')|_{C'}
&= K_{C'} \\
&= c^*(K_C)  + \sum (r_i -1)q_i \\
&= c^*(K_C + \sum \frac{r_i -1}{r_i} p_i).
\endalign
$$
Thus there is a canonical identification of $\Bbb Q$-Weil divisors
$$
(K_S + C)|_C = K_C + \sum \frac{r_i -1}{r_i} p_i.
$$
\endremark

\heading \S 5  Log Deformation theory \endheading

As indicated in the introduction, the proof of (1.3) involves deforming a
rational curve, while maintaining some prescribed contact with a divisor. 
Here we begin by developing the necessary deformation theory, (5.1-6), a straightforward
extension of the results of \cite{KMM92}. Then we will draw a number of
corollaries, all of which we believe are of some independent interest. 

In the applications, we will often be working on a bug-eyed cover, a
non-singular, but (usually) non-separated algebraic space. So we work
in this level of generality (which in fact does not in any way complicate
the treatment).

Below, by a space, we mean a (not necessarily separated) algebraic space.
In all of the  $\Hom$s we consider, the domain is always assumed to
be a projective scheme. See \cite{Kollar91} for a discussion of
$\Hom$ in this context. 

First a few remarks on pullbacks and notation (all of which are standard).  
Given a map $f: \map X. Y.$, and a subspace $K \subset Y$, $f^{-1}(K) \subset X$ 
indicates the scheme-theoretic inverse image (defined by pulling back the defining equations).
If $K$ is an effective Cartier divisor, then
$f^{-1}(K) \subset X$ is a locally principal subspace (that is its ideal sheaf is
locally principal), and if $X$ is integral and $f(X)$ is not a subset
of the support of $K$, then $f^{-1}(K) \subset X$ is again a Cartier divisor, which we will
often denote by the more conventional notation $f^*K$.  
In any case we can always pullback the line bundle $\ring Y.(K)$.
Given two subspaces $H$, $K$ of a space $Y$, 
$H\subset K$  means that the space $H$ is a subspace of  $K$ (as 
opposed to the weaker condition of set-theoretic inclusion).

\proclaim{5.1 Definition-Lemma} Let $D \subset X$ and $E \subset Y$ be 
effective Cartier divisors. Assume $X$ is purely one dimensional. 

There is a space $\Hom(X,Y,D \subset E)$ with the following universal
property. A $T$-point $g: \map T \times X.Y.$ of $\Hom(X,Y)$ is a point of 
$\Hom(X, Y, D \subset E)$ iff $T \times D \subset g^{-1}(E)$.
\endproclaim
\demo{Proof} Let 
$$
k:\map \Hom(X,Y) \times X.Y.
$$
be the tautological map. Let 
$$
G = \pi_*k^*\ring Y.(E).
$$
$G$ is obviously a vector bundle of rank $\deg {D}$. 
The composition
$$
\map {k^*\ring Y.(-E)}.{\ring {\Hom(X, Y) \times X}.}.
\twoheadrightarrow \ring \Hom(X, Y) \times D.
$$
induces a canonical section of $G$. Let 
$$
\Hom(X, Y, D \subset E) \subset \Hom(X, Y)
$$ 
be the zero locus. The universal property of the zero locus
of this section is obvious from the construction. 
\qed\enddemo

Let $\dd$ be a collection of effective Cartier divisors $\list D.n.$ on $X$ and $\e$ a 
collection of effective Cartier divisors $\list E.n.$ on $Y$. We let
$E = E_1 + \dots E_n$. We let 
$\gH=\Hom(X,Y, \dd \subset \e)$ denote the scheme-theoretic intersection of the
$\Hom(X,Y,D_i \subset E_i)$ inside $\Hom(X,Y)$. $\gH$ then inherits an obvious 
universal property from (5.1).

\proclaim{5.2 Lemma} Suppose $Y$ is smooth and the divisors $\list E.n.$ 
have normal crossings. The sequences
$$
\ses {T_Y(-\log E)}.T_Y.\oplus {\ring E_i.(E_i)}. \tag{1}
$$
and
$$
\ses {\ring E_1.}.T_Y(-\log E)|_{E_1}.
T_{E_1}(-\log [E_1 \cap (E-E_1)]). \tag{2}
$$
are exact.
\endproclaim
\demo{Proof} For (1) see Lecture 2  of \cite{Green94}. (2) follows from (1).\qed\enddemo

Next we extend (1.2) and (1.3) of \cite{KMM92}. 

\proclaim{5.3 Proposition} Let $X=C$ be a smooth projective curve. 
$\gH \subset \Hom(C,Y)$ is cut out by $\deg(D)$ equations. The 
dimension of any irreducible component of $\gH$ at a point $g$
is at least $\chi(g^*T_Y) - \deg(D)$. 

If furthermore 
$(Y,E)$ is log smooth, and 
$D_i=g^{-1}(E_i)$ for all $i$, then the Zariski tangent space at 
$[g]$ is $H^0(g^*T_Y(-\log E))$.
\endproclaim
\demo{Proof} We note first that (1.2) and (1.3) of \cite{KMM92} and their 
proofs hold when the target is not separated. $\gH\subset \Hom(C,Y)$ is cut out 
by $\deg (D)$ equations by (5.1). The dimension estimate
now follows from (1.2) of \cite{KMM92}. 

Next consider the commutative diagram
$$
\CD
0 @>>>  H^0(g^*T_Y(-\log E_i)) @>>> H^0(g^*T_Y) @>>> H^0(g^*\OO_{E_i}(E_i))  \\
@.      @VVV            @V{r}VV    @|   \\
0 @>>>  H^0((g|_{D_i})^*T_{E_i}) @>>> H^0((g|_{D_i})^*T_Y) @>z>> 
H^0((g|_{D_i})^*\OO_{E_i}(E_i)). 
\endCD
$$
(The first row is pulled back from (5.2.1) and remains exact as $E_i$ pulls 
back to a Cartier divisor). By the universal property, the tangent space to 
$\Hom(C,Y,D_i\subset E_i)$ is the kernel of the composition $ z \circ r$ and 
thus by the diagram equal to  $H^0(g^*T_Y(-\log E_i))$. The tangent space to 
$\gH$ is the intersection of the tangent spaces of $\Hom(C,Y,D_i \subset E_i)$
and so by (5.2.1) is equal to $H^0(g^*T_Y(-\log E))$.\qed \enddemo

Juggling techniques of Mori we obtain the following: 

\proclaim{5.4 Corollary} If $X$ is a smooth variety and $-(K_X+D)$ is ample, 
then $(X,D)$ is log uniruled. If $-K_X-D$ is pseudo-effective, and $X$ is 
uniruled, then through a general point of $X$ there is a rational curve meeting 
$D$ at most twice.
\endproclaim
\demo{Proof} Let $w$ be a zero dimensional scheme on a curve $C$. We say $w$ has 
{\bf size} $s$, if $s$ is the maximum length of a subscheme $z \subset w$ 
supported at a single point. Fix a polarisation $H$ and a general point $x$ of 
$X$. Among rational curves $f:\map C.X.$ through $x$ with $H\cdot C$ minimum, 
choose one such that $f^*D$ has greatest size $s$. 

Suppose the result fails. Then in the pseudo-effective case,
$s \leq C \cdot D -2$, and $s \leq C \cdot D -1$ in the ample case.
Let $z$ be a subscheme of $f^*D$ of length $s$, let $p$ be the 
support of $z$, and let $t$ be a point mapping to $x$. Note if $-K_X-D$ is 
pseudo-effective, $C \cdot D \leq -K_X \cdot C$, 
with strict inequality if $-K_X-D$ is
ample. Thus in any case $\Hom(C,X,z \subset D)$ has dimension at least
$-K_X\cdot C+n-s\geq n+2$ where $n$ is the dimension of $X$. Thus $f$ moves,
fixing $x$ and $z \subset f^*(D)$, in a two dimensional family, and so its 
image must move. Thus we can find a projective surface $S$, a smooth 
projective curve $B$, and a finite map $\map S.B \times X.$ such that 
the map $h:\map S.X.$ is generically finite, $\map S.B.$ has general 
fibre $\pr 1.$, and for some open set $U \subset B$, $h$ gives a $U$-point
of $\Hom(C,X,z \subset D, t \subset x)$, and $h_b: \map {\pr 1.}.X.$ is the 
map $f$. Since $C \cdot H$ is minimal, $\pi_1$ has reduced irreducible fibres,
and thus is a smooth $\pr 1.$-bundle. Let $Z \subset S$ be the closure of 
$p \times U$. Then $Z$ is a subscheme of $h^*D$ and a section of $\pi_1$. 
Let $T$ be the closure of $t \times U$. By assumption there is another
component $E$ of $h^*D$, necessarily disjoint from $Z \cup T$. $T$ has
negative self-intersection (as $h(T) = x$). Thus (since the
cone of $S$ is two dimensional), $Z$ and $E$ have positive self
intersection, and $Z \cdot E =0$, a contradiction. \qed\enddemo

(5.4) gives some evidence for (1.10), but it is not very useful for proving (1.1)
-assuming both the smoothness for $X$ and ampleness of $-(K_X + D)$ is 
too strong. The next two corollaries, on the other hand, will prove to be key:

\proclaim{5.5 Corollary} Let $g:\map {X=C=\pr 1.}.Y.$ be a morphism, where $Y$ 
is smooth and the divisors $\list E.r.$ have normal crossings. Assume 
$D_i=g^{-1}(E_i)$. The following hold:
\roster
\item If $D$ is supported in at most two points, and $g$ is in a component 
$Z$ of $\gH$ such that $\map {Z \times \pr 1.}.Y.$ is dominating, then 
$f^*T_Y(-\log E)$ is semipositive for a general point $f \in Z$.
\endroster
\smallbreak
If $g^*T_Y(-\log E)$ is semipositive then (near $g$)
\roster
\item[2] $\map \gH\times \{p\}.Y.$ is smooth for any $p \not \in D$ and 
\item $\map \gH \times \{p\}.E_j$ is smooth for 
$p \in D_j \setminus \cup_{i \neq j} D_i$. 
\endroster
\endproclaim
\demo{Proof} (1) follows from (5.3) as in the proof of (1.3) of \cite{KMM92}. 
(2) and (3) are clear from (5.3).\qed\enddemo

The following corollary holds with any number of components, but we 
present it for at most two components, as this simplifies notation and is 
sufficient for our needs.

\proclaim{5.6 Corollary} Let $g:\map {X=C=\pr 1.}.Y.$ be a morphism, where $Y$ 
is smooth and the divisors $E_1,E_2$ have normal crossings.
Suppose either 
\roster 
\item"(i)" $g(\pr 1.)\subset E_1$, $g^{-1}(E_2)=D_2$, 
$\deg g^*\ring Y.(E_1) =\deg (D_1)$, 
and 
$g^*T_{E_1}(-\log (E_1 \cap E_2))$ is semipositive. Set 
$$
J=\Hom(C,E_1,D_2 \subset E_2 \cap E_1),
$$
or
\item"(ii)" $g(\pr 1.)\subset E_1\cap E_2$, $\deg g^*\ring Y.(E_i)=\deg (D_i)$, 
and $g^*T_{E_1 \cap E_2}$ is semipositive. Set 
$$
J =\Hom(C,E_1,D_2 \subset E_2 \cap E_1)\cup \Hom(C,E_2,D_1\subset E_1\cap E_2).
$$ 
\endroster
Then $\gH$ is l.c.i and reduced at $g$ and $J$ is a nowhere dense subspace. 
Furthermore (near $g$)
\roster
\item $\map \gH \setminus J \times \{p\}.Y.$ is smooth for any $p \not \in D$.
\item $\map \gH\setminus J \times \{p\}. E_1.$ is smooth for any 
$p \in D_1 \setminus D_2 \cap D_1$.
\endroster
\endproclaim
\demo{Proof} We will only deal with case (i), case (ii) is very similar
(and not used in the proof of (1.1)).

 (5.2.1) and (5.2.2) show that $g^*T_Y(-\log E)$, and $g^*T_Y$
are both semipositive. Let $H=\Hom(\pr 1., Y)$. Observe that
$J = \gH \cap \Hom(\pr 1.,E_1) \subset H$. (1-2) are immediate from
(5.5.2-3). By (5.3), $H$ and $J$  are smooth at 
$g$ of dimension $\chi(g^*T_Y)$, and $\chi(g^*T_{E_1}(-\log (E_2 \cap E_1)))$ 
respectively. By (5.3) every component of $\gH$ has dimension at least 
$\chi(g^*T_Y) -\deg (g^*E)=\dim (J)+1$. Thus $J$ is nowhere dense. 
So by (1-2), $\gH \subset H$ has pure codimension
$\text{deg}(D)$. Since it is cut out by $\text{deg}(D)$ equations,
it is a complete intersection. In particular it is Cohen-Macaulay, and
so since it is generically reduced, it is reduced. \qed\enddemo

\example{5.7 Examples} 
As a simple example of (5.6.i), let $g$ embed $\pr 1.$ as a smooth conic $E= E_1$
in $\pr 2.$. Let $D = 4[p]$, for a fixed $p \in \pr 1.$. (5.6.1) 
says we can deform $g$ in a four dimensional family, while 
maintaining fourth order contact at $p$. Thus we 
have the elementary observation that given any conic in $\pr 2.$ there is a two dimensional 
family of conics which meets the given conic once. 
 
As another simple example of (5.6.i), 
suppose $B \subset S^0$ is a log canonical,
reduced, but not irreducible, member of $|-K_S|$ (here,
as always, $S$ is a rank one log del Pezzo). Then adjunction shows
that every irreducible component $E_1$ of $B$, is a smooth rational curve, 
and $E_1$ meets $E_2 = B - E_1$ normally at two points, say $p$ and $q$. Let
$g:\map {\pr 1.}.E_1.$ be the identity. $g^*T_{E_1}(-\log E_1 \cap E_2)$ 
is trivial by adjunction, so we may apply (5.6.i), with
$D_1 = (E_1^2)p$, $D_2 =p+q$, to conclude that $S^0$ is dominated
by rational curves (in fact deformations of $E_1$) meeting $B$ twice. 
In this case $p$ will be a basepoint of the dominating family. One can 
obtain the same result by counting dimensions in $|E_1|$, since any 
deformation of $E_1$ is again smooth rational. 
\endexample

We will frequently use the following consequence of (5.5) (sometimes
without reference): 

\proclaim{5.8 Corollary} Let $B\subset X$ be a reduced divisor on
a normal algebraic space. Assume $X$ is a separated variety in a neighbourhood
of $\sg(B) \cup \sg(X)$. 
Let $V \subset Z$ be a closed subset of codimension at least two, disjoint from
$\sg(B) \cup \sg(X)$.
 
 If $X$ is dominated by a family of rational curves meeting $B \cup \sg(X)$ at 
most $d$ times, with $d \leq 2$, then it is dominated by rational curves
disjoint from $V$, meeting $B \cup \sg (X)$ at most $d$ times, and algebraically
equivalent to members of the original family.
\endproclaim
\demo{Proof} After passing to a log resolution of $(X,B)$, replacing $B$ by its
total transform and  the dominating family by the strict
transform,  we may assume $(X,B)$ is log smooth. By (5.5.1-3), the general member
of the family deforms, maintaining the number of points of contact with $B$, with no 
basepoints outside of $\sg(B)$. \qed \enddemo

Next a method of proving $X^0$ is uniruled. We use the following:

\definition{Definition} We say a pair $(X,A)$ of a reduced divisor on a normal
variety has {\bf quotient singularities} if there is a quasi-finite surjection
$q:\map U.X.$, \'etale in codimension one, with $U$ and $q^{-1}(A)$ smooth. 
\enddefinition

\proclaim{5.9 Corollary} Let the pair $(X,A)$ consist of a normal variety and 
a divisor $A$. Suppose there is a dominating family of rational curves meeting 
$\sg (X)\cup A$ at most $d$ times, for $d \leq 2$. 
Suppose, either
\roster 
\item $(X,A)$ has quotient singularities, or
\item $X$ has terminal, LCIQ singularities (see \cite {Kollar91}), and 
$A$ is empty.  
\endroster
 Then $X^0$ is dominated by rational curves meeting $A$ at most $d$ times. 
In particular $X^0$ is uniruled.
\endproclaim
\demo{Proof} By (4.0), as in the first paragraph of the proof of (4.4), 
there is a bug-eyed cover 
$\flat:\map X^{\flat}.X.$ such that in case (2), $X^{\flat}$ is l.c.i, and
in case (1), $X^{\flat}$ and $A'= \flat^{-1}(A)$ are non-singular. By (4.9.3), $X^{\flat}$ 
is dominated by rational curves meeting $A'$ at most twice. 

For (1), apply (5.8) to $(X^\flat,A')$ with $V = \flat^{-1}(\sg(X))$.

In case (2), let $n$ be the dimension of $X$. Let $C = \pr 1.$ and $\map C.X^{\flat}.$
the general map in a dominating family, with $f:\map C.X.$ the induced map.
By (2.10) of
\cite{Kollar91}, $\Hom (C,X)$ has dimension at least 
$-f(C) \cdot K_X  + n$ at $[f]$. Let $[f] \in H \subset \Hom (C,X)$ 
be an irreducible component of that dimension. Let
$\map Y.X.$ be a resolution.  We may lift $f$ to $f':\map C.Y.$. Let $H'$ be the
irreducible component of $\Hom(C,Y)$ containing $f'$. 
Since maps in $H'$ dominate $Y$, $H'$ has dimension $- f'(C) \cdot K_Y + n$,
by (5.3) and (5.5).
Since $\map H'.H.$ is dominant, 
$f'(C) \cdot K_Y \leq f(C) \cdot K_X$. Now by the definition of terminal, 
it follows that $K_Y \cdot f'(C) = K_X \cdot f(C)$ and $f(C) \subset X^0$. \qed \enddemo

\remark{Remark} For cases where (5.9) applies,  
we note that a log terminal pair $(S,B)$ of an irreducible curve on
a surface has quotient singularities, see Appendix L, and that terminal three-fold
singularities are LCIQ, see \cite{Kollar91}. \endremark

For surfaces, we have the following combination of (5.8-9):

\proclaim{5.9.3 Corollary} Let $D \subset S$ be a reduced curve on a log
terminal surface. Let $V$ be any finite set of points at which $(S,D)$
is purely log terminal (or equivalently, has quotient singularities). If $S$ 
is dominated by rational curves meeting $D \cup \sg(S)$ at most $d$ times,
with $d \leq 2$, then $S$ is dominated by rational curves disjoint from $V$,
and meeting $D \cup \sg(S)$ at most $d$ times.
\endproclaim
\demo{Proof} Let $C$ be a general member of a dominating family of 
rational curves meeting $D \cup \sg(S)$ at most $d$ times. 	
After blowing up points away from $V$, taking the strict transform of $C$,
and then restricting to a neighbourhood of $C$, we can assume $(S,D)$ has
quotient singularities. Now apply (5.9). \qed \enddemo

\proclaim{5.10 Corollary} The smooth locus of a $\Bbb Q$-factorial projective 
toric variety is uniruled.
\endproclaim

\demo{Proof} $X$ has quotient singularities and $X^0$ is
obviously dominated by copies of $\aos$. Now apply (5.9.1). \qed\enddemo

The existence of dominating families of rational curves has implications for
the Kodaira dimension:

\proclaim{5.11 Lemma} Let $Y$ be a variety. Suppose that $K_Y+A$ is 
semi log canonical and that $C$ is the general member of a covering family of rational 
curves. Let $d$ be the number of times that $C$ meets the union of $A$ with
the locus of points where $Y$ is not canonical. 
\roster
\item If $d\leq 2$ then $(K_Y+A)\cdot C\leq 0$. In particular 
$K_Y+A$ is not big.  
\item If $d\leq 1$, then $(K_Y+A)\cdot C < 0$. In particular 
$K_Y+A$ is not effective.  
\endroster
\endproclaim
\demo{Proof} Passing to the normalisation of $Y$, we may assume that $Y$ is log 
canonical. Let $\pi:\map Z.Y.$ be a log resolution of 
the pair $(Y,A)$. By the definitions of canonical and log canonical, we may write 
$$
K_Y+B=\pi ^*(K_Z+A)+E
$$
where $E$ is effective and exceptional, $B$ is reduced, and $\pi(B)$ is contained
in the union of $A$ and the locus where $Y$ is not canonical. Thus the
strict transform of $C$ 
meets $B$ at most $d$ times. Replacing the pair $(Y,A)$ by the pair
$(Z,B)$ we may assume that $(Y,A)$ is log smooth.
Let $f:\map {\pr 1.}.Y.$ be the normalisation of $C$. By (5.5.1), after possibly
replacing $f$ by a deformation (with the same contact with $A$),
$f^*T_Y(-\log A)$ is semipositive. This proves (1). 

Now suppose $d=1$ (the case $d=0$ is easier and well known). Let $a= \deg f^{-1}(A)$. 
By assumption, there is a component $Z$ of $\Hom(\pr 1.,X,ap\subset A)$ which contains $f$ 
whose universal family dominates $Y$. The universal family has one
dimensional fibres over $Z$, but it has at least two dimensional fibres over $Y$, since 
$\af 1.$ has a two dimensional family of automorphisms. Thus the Zariski tangent
space at $[f]$ has dimension at least $\di(Y)+1$. (2) follows easily by Riemann Roch 
and (5.3). \qed\enddemo

\remark{5.12 Remark} 
Note by (5.9.1), if the smooth locus of a log terminal
surface $S$ is dominated by images of $\aos$ then in fact the smooth
locus is uniruled, that is $S$ is dominated by complete rational curves that
miss the singularities entirely. 

On the other hand, it is easy to construct rational surfaces with quotient singularities 
such that $K_S$ is ample (for example the surface in (19.3.2.2.1) with $k \geq 10$). Obviously 
these surfaces are uniruled, but by (5.11) any dominating family of rational curves meets
the singular locus at least three times. Similarly, the bug-eyed covers of these surfaces 
cannot be dominated by images of $\pr 1.$, for example by (5.11) and the proof of (5.9). \endremark

\heading \S 6 Log Uniruledness \endheading
\subhead Brief Introduction \endsubhead

In this section we develop our main techniques for proving log uniruledness.
First we will prove (1.3) in the case when the boundary $D$ is non-empty. Thus
for (1.3) we are reduced to showing that a rank one log del Pezzo has uniruled
smooth locus. In this section we provide two sorts of criteria.

First we will prove the general implication
$$
\text{existence of a tiger} \Rightarrow \text{$S^0$ uniruled}.
$$
As indicated in the introduction, this will complete the first half of the proof of (1.3), 
and proves (1.3) in all but a bounded collection of cases.

Secondly (6.5) and (6.6) give criteria for when a multiple of a rational
curve $Z \subset S$ deforms (as a rational curve) to dominate $S^0$.
We will finish the proof of (1.3) in \S 8-18 by generating a finite set
$\goth F$ of (families of) $S$, containing any simply connected $S$ without tiger, and then
for each $S \in {\goth F}$, finding a rational curve $Z \subset S$ to which 
we can apply one of (6.5) or (6.6).

At the end of this section we give a number of examples. (6.8) shows that our main criterion,
(6.5), is in some sense optimal. (6.10) contains a detailed example application of (6.5). In 
(6.11) we find the minimal (degree) dominating family of rational curves in $S^0$ in an 
interesting specific case. In (6.9) we give some examples to convince the reader that one 
cannot prove (1.1) by (at least a direct naive application of) the {\it classical} method  
of constructing rational curves, Mori's bend and break technique (we worry the reader might 
duck the rather involved analysis ahead if they believed there was a promising alternative).

\subhead Existence of a tiger implies uniruledness \endsubhead

The main result of this subsection is the following:

\proclaim{6.1 Proposition} Let $S$ be a normal surface and $\alpha$
an effective $\Bbb Q$-Weil divisor on $S$, such that $K_S+\alpha$ is 
numerically trivial, and not klt. Let $B$ be the support of $\rdown \alpha.$. 
Then $S$ is dominated by rational curves meeting $B\cup \sg(S)$ at most twice. 
\endproclaim

Note that when $S$ in (6.1) is log terminal, then, by (5.12), 
the conclusion implies in particular that $S^0$ is uniruled. Thus (6.1)
implies (1.14) of the introduction.

We will obtain a stronger form of (6.1) in \S 20 as a corollary of (1.1), see
(1.4).

The main point for (6.1) is the following, the non-empty boundary
case of (1.3).

\proclaim{6.2 Lemma} Let $S$ be a rank one log del Pezzo surface. 
If $B$ is a non-empty reduced curve, such that $K_S+B$ is 
log terminal  and $-(K_S+B)$ is ample, then $(S,B)$ is log-uniruled.
\endproclaim
\demo{Proof} Note by adjunction and the classification of lt singularities,
$B$ has at most two components, and if $B$ is reducible, the two components
meet exactly once, normally at a smooth point of $S$.

Let $\flat :\map S^{\flat}.S.$ be the bug-eyed cover of $S$, (4.5). 
By (4.14.1) the pair $(S^{\flat},B'=\flat ^{-1}(B))$ is log smooth. 

Let $C$ be a component of $B$ and let $C'$ be the corresponding component of 
$B'$. By (4.8) and (4.14), $C'$ is non-singular, $K_{C'}$ has negative degree,
$C = \pr 1.$ and there is a surjection $g:\map {\pr 1.}.{C'}.$. Note
$C^2 > 0$ since $S$ has Picard number one. There are two cases. 

 If $B=C$ has one component we can apply (i) of (5.6) with 
$E=B'$ and $D$ equal to $\deg (g^*B')p$ for 
$p$ a general  point of $\pr 1.$. By (5.6.1) $g$ deforms to give 
a dominating family of 
maps $h: \map {\pr 1.}.{S^\flat}.$ such that 
$$
h^{-1}(B')=(\flat \circ h)^{-1}(B)  \supset (\deg (\flat \circ h)^{-1}(B))p.
$$
Since the two Cartier divisors above have the same degree, we have equality, and thus
$h(\pr 1.)$ meets $B'$ exactly once, at $h(p)$. By (5.6.1) the general member misses
any codimension two subset, in particular $\flat^{-1}(\sg(S))$, and thus we obtain
a family of rational curves dominating $S^0$ and meeting $B$ exactly once.
 
 Otherwise $B$ has exactly two components. Let $B=B_1+B_2$, where $C=B_1$. 
$B_2$ meets $C$ normally at a single point, $x \in S^0$, and there is (by adjunction)
at most one singular point, $y$, of $S$ along $C$. By (4.9.3), we can
assume $\flat \circ g$ is an endomorphism of $\pr 1.$ totally ramified at $x$ 
and $y$. Let $p$ be the unique point in $g^{-1}(x)$.
Now we apply (5.6.i) with 
$$
E_1=B_1', \quad E_2 = B_2', \quad D_1= \deg (g^*B_1')p, \quad
D_2 = \deg (g^*B_2')p = g^{-1}(B_2').
$$ 
The existence of a dominating family follows as before. 
\qed\enddemo

 The next Lemma reduces (6.1) to a question about rank one Gorenstein log del
Pezzo surfaces, using (6.2) and some standard 
applications of the MMP. First a little notation. In the notation of (6.1), let 
$\beta$ be the largest boundary, contained inside $\alpha$. Note that $\alpha$ 
and $\beta$ have the same support, and agree on the complement of $B$.

\proclaim{6.3 Lemma} To prove (6.1) we may assume, 
\roster
\item There is no morphism $\pi :\map S.\Sigma.$, with general fibre $\pr 1.$,
\item $K_S+\beta$ is log canonical,
\item $B$ is non-empty, 
\item $S$ is a rank one log del Pezzo surface,
\item $K_S + B$ is log terminal at singular points,
\item $B = \alpha$ (or equivalently, $K_S + B$ is numerically trivial),
\item $B\in |-K_S|$, $B\subset S^0$ and $S$ is Gorenstein, 
\item $S^0$ is simply connected, and
\item $B$ is irreducible. 
\endroster
\endproclaim
\demo{Proof} Suppose there is a morphism $\pi:\map S.\Sigma.$, 
with general fibre, $F$, isomorphic
to $\pr 1.$. Then $F$ lies in the smooth locus and since $(K_S+B)$ is non-positive
on $F$, $F$ meets $B$ at most twice. Hence 
we have found our dominating family of rational curves. Thus 
we may assume (1). 

 Let $f:\map T.S.$ be a log terminal model of the pair $(S,\beta)$, see (17.10) 
of \cite{Kollaretal}. By definition of a log terminal model
we may write the log pullback (see (1.16) ) as 
$$
K_T+\tilde \beta + E =f^*(K_T+\beta),
$$
where $\rdown E.$ is the full reduced exceptional locus, 
$\tilde \beta$ is the strict transform 
of $\beta$, and $K_T + \tilde{\beta} + \rdown E.$ is log terminal.
As discrepancies
are not effected by log pullback,
we are free to replace $(S,\alpha)$, by $(T,E + \tilde\beta+ f^*(\alpha-\beta))$. 
Now $K_S + \beta$ is log terminal, in particular (2) holds. 

Suppose $B$ is empty. Then $\alpha$ and $\beta$ are equal, whilst $K_S+\beta$ is 
klt, a contradiction. Hence (3). 

Suppose $f:\map S.S'.$ is any contraction. Let $B'$ be the one dimensional 
part of $f(B)$. By (5.8), applied with
$V = f(B) \setminus B'$, it is enough to find 
a dominating family of rational curves $C$ which meets $B'\cup \sg(S')$ at 
most twice.

$K_S + (\beta - B)$ cannot be effective, so by (1), the  MMP with respect
to this divisor gives a birational contraction $f: \map S. S'.$ to a
rank one log del Pezzo with $K_{S'} + f_*(\beta - B)$ log terminal and 
anti-ample. Not every component of $B$ can be contracted by $f$, for otherwise
(since $\alpha - \beta$ is supported on $B$) $K_S + \alpha$ is the log
pullback of $K_{S'} + f_*(\beta - B)$, contradiction. Replacing $(S,\alpha)$ by 
$(S', f_*(\alpha))$ we have (1-4).

Now suppose $B$ is not log terminal at singular points. Then we 
may extract a divisor of the minimal desingularisation of coefficient one
(for  $K_S + B$) $f:\map T.S.$ and contract the
$K_T$-negative ray $\pi:\map T.S'.$  (replacing
the boundaries by pullback and pushforward as above). This process improves the
singularities so eventually terminates, with (1-5).

Now suppose $B \neq \alpha$. Then $K_S + B$ is anti-ample, and so by adjunction
and (5), $K_S + B$ is log terminal. Thus we may apply (6.2). Hence (1-6).

There is a connected cover, 
$h:\map T.S.$ \'etale in codimension one, such that $K_T + B' =h^*(K_S + B)$
is trivial (in particular Cartier). By 
(20.4) of \cite{Kollaretal}, and (5), $K_T + B'$ is log canonical, and 
log terminal at singular points. It follows, by the classification of log
terminal pairs,
that $B' \subset T^0$. Thus $T$ is a Gorenstein log del Pezzo. Taking a further
cover we may assume, by (3.5), that $T^0$ has trivial algebraic fundamental group.
Running the $K_T$-MMP to return to Picard number one, we obtain (1-8) by
(3.3). 

 Finally, by (5.7) we may assume that $B$ is irreducible. Hence 
(9). \qed\enddemo

We will use the following rather ad hoc result in the proof of (6.1):

\proclaim{6.4 Lemma} Let $S$ be a surface with quotient singularities,
and let $f:\map {\pr 1.}.S.$ be a generic embedding with image $Z$. 
Let $B \subset S^0$ be a reduced curve. Assume $Z$ meets $B$ exactly once and
$\sg(S)$ at most once, and $f^{-1}(B)$ has degree at least 
two.
If $(K_S+B) \cdot Z \leq 0$ then $S^0$ is dominated by rational curves, 
algebraically equivalent to a multiple of $Z$, meeting $B$ at most twice. 
\endproclaim
\demo{Proof} Let $C = \pr 1.$. Let $z = Z \cap B$ and if $Z$ meets
$\sg (S)$ then let 
$p =Z\cap \sg(S)$ and otherwise choose any smooth point $p$ of $Z$ other than $z$. 
We follow the notation of (4.10). Note that $f:\map C.Z.$, $\map C'.Z.$, and
$\map S^{\flat}.S.$ are all set bijections in a neighbourhood of $z$, and we will
also use $z$ to denote the its unique inverse image under any of these maps. 
Similarly, we'll use the same notation for $B \subset S^0$ and its inverse
image on $S^{\flat}$.

Let $\ol{h}:\map {\pr 1.}.C.$ be of degree $r$ ($r$ as in (4.10.3), which is $1$ if
$Z \subset S^0$) totally ramified over $z$ and $p$. By (4.7), $\ol{h}$ factors through
a map $h:\map {\pr 1.}.{C'}.$.  
Let $x=h^{-1}(z)$ and $y=h^{-1}(p)$. Note $h$ is \'etale at $y$ and 
totally ramified at $x$. We also indicate by $h$ the induced map 
$\map {\pr 1.}.S^{\flat}.$. Note $f^{-1}(B)=s[z]$, $h^{-1}(B)=rs[x]$, for some $s \geq 2$.

Let 
$$
H=\Hom(\pr 1.,S^{\flat}, (rs-1)[x] \subset B).
$$

By (5.8) it is enough to show that  $H$ gives a dominating family of rational 
curves (necessarily, by degree considerations, meeting $B$ at most
twice). Suppose this fails, or equivalently, that 
maps in $H$ have fixed image $\flat^{-1}(Z)$. Then by (4.10.2) for 
$$
J = \Hom(\pr 1.,C', (rs-1)[x] \subset s[z]),
$$
the natural quasi-finite map $\map J.H.$ is surjective near
$[h]$. 

Obviously, for any non-constant map $\gamma:\map {\pr 1.}.C'.$, the
length of $\gamma^{-1}(s[z])$ at any point in its support is divisible by
$s$. Thus, since $s \geq 2$, at least set theoretically
$J = \Hom(\pr 1.,C',r[x] \subset [z])$. By (5.3) the latter
has Zariski tangent space
$$
H^0(h^* T_{C'}(-\log z)) = H^0(T_{\pr 1.}(-\log x))= 
H^0(\ring {\pr 1.}.(1))
$$
at $[h]$ and thus has dimension at most two.
But by 
(5.3) the dimension of $H$ is 
at least $-r(K_S + B) \cdot Z + 3 \geq 3$ at $[h]$, a contradiction, \qed \enddemo

Note that (6.4) cannot be extended to the case $Z \cdot B =1$. For example 
take for $Z$ a $-1$-curve on a smooth surface, and for $B$ some curve meeting $Z$ 
normally at one point.

\demo{Proof of (6.1)} By (6.3), we may assume we have properties (6.3.1-9). 

Note that $B$ is either a smooth elliptic curve, or a rational curve with a node. If 
$S = \pr 2.$ take tangent lines to $B$. Otherwise
let $\pi: \map T. {\pr 1.}.$ be the $\pr 1.$-fibration of (3.8), obtained by extracting the 
$-2$-curve $E$. $B \subset T^0$ is a double cover. By the
Hurwitz formula, since $\pi$ has at most one multiple fibre, there
is a smooth fibre $F \subset T^0$ 
meeting $B$ exactly once, with $F$ and $B$ simply
tangent (at a point of $B^0$). Let $C \subset S$ be the image of $F$. Since $E$ is a section,
$C$ meets $\sg(S)$ once. Now apply (6.4). \qed \enddemo  

\subhead Deforming special rational curves $Z \subset S$ \endsubhead

Here we give some criteria for when a multiple of a rational curve $Z \subset S$ 
deforms, as a rational curve, to dominate $S^0$. 
The idea is to lift (an endomorphism of) $Z$ to the bug, and then apply 
deformation theory.  

Here is the main criteria, discussed in the introduction. 

\proclaim{6.5 Lemma} Let $S$ be a normal surface with quotient singularities. 
Let $f:\map C.S.$ be a generic embedding, with 
$C=\pr 1.$, such that $f^{-1}(\sg(S))$ is two points, $p$, $q$. Let $x$, 
$y$ be the local indices of $Z=f(C)$ at $f(p)$, $f(q)$ (when $f(p)=f(q)$ this
means the indices of the two branches). If $-K_S \cdot Z \geq 1/x + 1/y$ 
then $S^0$ is uniruled, dominated by rational curves algebraically
equivalent to a multiple of $Z$.
\endproclaim
\demo{Proof} Let $u$ and $v$ be the ramification indices of 
$c:\map C'.C.$ over $p$ and $q$, where $C'$ is the normalisation of 
$\flat ^{-1}(Z)$ . By (4.10.3) $-K_S \cdot C \geq 1/u+1/v$. 
Let $h :\map {\pr 1.}.C.$ be a map of degree $uv$, totally ramified above 
$p$, $q$. By (4.8), there is a lifting $h':\map {\pr 1.}.C'.$. 
By (5.3), $\Hom(\pr 1.,C')$ is smooth at $[h']$ 
of dimension $u+v+1$. By (5.3), $\Hom(\pr 1.,S^{\flat})$ has dimension at least 
$-uvK_S \cdot C +2$. The result now follows from (4.10.2) and 
(5.5).\qed\enddemo

We note that in the proof of (6.5) we actually obtain a bit more. It's enough
to get $K_S \cdot Z \geq 1/r(p) + 1/r(q)$ where $r$ is as in (4.10.3) (note
by (4.10.3), $x,y$ divide $r(p)$, $r(q)$, and they are not in general equal). 
In applications we will know how
$\t{Z} \subset \t{S}$ meets the exceptional locus. From this (and the 
description of the exceptional locus) it is easy to compute the local
analytic index of (some branch of) $Z$. The numbers $r(p)$ and $r(q)$ can
also be computed from this description (which after all determines the pair
$(S,Z)$ locally analytically), but the analysis is more involved. We use the
weaker criterion, and its partner below, as they turn out to suffice.

In most applications of (6.5), $x$ and $y$ will be equal, and each branch
of $Z$ will have the same index, the index of the singular point (for example
this holds when each branch is log terminal, or by (4.10.3) if neither
branch is Cartier and the singularity has prime index). If the inequality
in (6.5) fails, we will look for a second curve of the same sort, in order
to apply the next result:

\proclaim{6.6 Lemma} Let $S$ be a normal surface with quotient singularities. 
Let $f:\map {C=\pr 1.}.S.$, $g:\map {D=\pr 1.}.S.$ be 
two generic embeddings, such that $f(C)\cap \sg(S)=\{p\}$, 
$g(D)\cap \sg(S)=\{q\}$, $Z_1= f(C) \neq g(D) = Z_2$ and such that $f^{-1}(p)$ and 
$g^{-1}(q)$ are each a pair of points. Suppose further that each of the local 
branches of $Z_1$ (resp. $Z_2$) have the same index $n$ (resp. $m$) at $p$ 
(resp. $q$) where $n$ (resp. $m$) is the index of the surface at $p$  
(resp. $q$). Assume there are points $x$ and $y$ of $\pr 1.$ such that 
$f(x)=g(y)\in S^0$. 
 
 Then $S^0$ is uniruled.
\endproclaim
\demo{Proof} By (6.5) we may assume $nf(C)\cdot K_S=mg(C)\cdot K_S=-1$. As in 
the proof of (6.5), there are \'etale maps of degree $n$ and $m$, 
$h' :\map {\pr 1.}.C'.$ and $j':\map {\pr 1.}.D'.$, where $C'$ and 
$D'$ are the normalisations of $\flat ^{-1}(Z_1)$ and $\flat ^{-1}(Z_2)$. 
Let $N^{\flat}$ be the union (with normal crossing) of $C'$ and 
$D'$ at $h'(x)$ and $j'(y)$. Let $N$ be the union of two copies of 
$\pr 1.$ joined along $x$ and $y$. There is an induced morphism 
$i:\map N.N^{\flat}.$. Let $\map \Cal N.B.$ be a one dimensional deformation of 
$N$ to a smooth rational curve. By \cite{KMM92} 
$\Hm B.\Cal N, S^{\flat}\times B.$ has dimension at least $5$, while 
$\Hom(N,N^{\flat})$ has dimension four, since a deformation of $i$ is given 
by \'etale maps $\map {\pr 1.}.C^{\flat}.$ and $\map {\pr 1.}.D^{\flat}.$, with 
the image of $x$ and $y$ fixed. Thus $S^0$ is uniruled.\qed\enddemo

An example application of (6.5) is worked out in (6.10).

We present a three dimensional version of (6.2) (dlt is defined in \cite {Kollaretal}):

\proclaim{6.7 Lemma} Let $X$ be a projective threefold with a dlt divisor
$K_X+E_1+E_2$ where $E_1$ and $E_2$ are irreducible and normal.
Let $C=E_1 \cap E_2$. If $C \cdot E_i$ and $-C \cdot (K_X + E_1 + E_2)$
are both positive then $(X, E_1 + E_2)$ is log uniruled. 
\endproclaim
\demo{Proof} By Adjunction $C=E_1\cap E_2$ is a smooth rational curve. A local 
analytic description along $C$ is given in (16.15.2) of \cite{Kollaretal}. 
In particular $X$ has quotient singularities along $C$. Let  
$\flat :\map X^{\flat}.X.$ be a smooth bug-eyed cover. Then if we set
$E_i'=\flat ^{-1}(E_i)$ then $E_1'$ and $E_2'$ are smooth and cross normally 
along $C'=\flat ^{-1}(C)$. Adjunction shows that $K_{C'}$ is negative and so 
there is a surjection $g :\map {\pr 1.}.C'.$. Now apply (5.6), checking 
semipositivity as in the proof of (6.2). \qed\enddemo

\subhead Optimality of (6.5) \endsubhead

It is possible to have $f:\map {\pr 1.}.S^{\flat}.$, such that no multiple of $f$ deforms to 
cover $S^{\flat}$. Here by a multiple we mean a composition  $g=f \circ h$, for $h$ an 
endomorphism of $\pr 1.$. Thus it is not sufficient for proving uniruledness to find a rational 
curve lifting to $S^{\flat}$, we need some additional local information. For example:

\proclaim {6.8 Example} Suppose $C$ is a rational curve, smooth away from $\sg S$, meeting 
$\sg S$ twice, with each branch of $C$ log terminal, and such that each of the singular
points along $C$ has the same index, $m$. Let $f: \map {\pr 1.}. S^{\flat}.$ be the $m$ to one 
lift to the bug-eyed cover (as in the proof of (6.5)). 

If (6.5) fails, that is if $-K_S\cdot C = 1/m$ then no multiple of $f$ deforms to cover 
$S^{\flat}$. \endproclaim
\demo{Proof} ${\flat}^{-1}(C)$ has smooth branches, thus $j:\map C^{\flat}.S^{\flat}.$ 
is unramified, and $T_{C^{\flat}}$ is a subbundle of $j^*(T_{S^{\flat}})$. $f$ factors through an
\'etale map $\map {\pr 1.}. C^{\flat}.$ thus $T_{\pr 1.}$ is a subbundle of $f^*(T_{S^{\flat}})$. 
Since $f^*(T_{S^{\flat}})$ has degree $1$, the quotient the quotient is $\ring .(-1)$ and thus
$$
f^*(T_{S^{\flat}}) = \ring .(2) \oplus \ring. (-1).
$$
Let $h$ be an endomorphism of $\pr 1.$ of degree $d$, and let $g = f \circ h$. The Zariski tangent
space to $\Hom(\pr 1.,S^{\flat})$ at $g$ is
$$
H^0(g^*(T_{S^{\flat}}))=H^0(\ring .(2d) \oplus \ring.(-d))
$$
which has dimension $2d +1$. $\Hom(\pr 1., {\flat}^{-1}(C))$ has dimension at $g$ at least 
the dimension of $\Hom(\pr 1., \pr 1.)$ at $h$, which is also $2d + 1$. Thus 
$\Hom(\pr 1.,S^{\flat})$ is smooth at $g$, and equal to $\Hom(\pr 1., {\flat}^{-1}(C))$. 
\qed \enddemo
 
One instance  of (6.8) is the curve $L_{aa'}$ in (15.2).

\example{6.9 Examples: Trouble with Bend and Break}
\roster
\item If $C$ is a curve in the smooth locus of $\Hbz e.$ then 
$-K_{\Hbz e.}\cdot C \geq e+2$. Note in contrast that if $X$ is a smooth 
projective $n$-dimensional Fano variety, then the bend and break techniques 
give
that $X$ is covered by rational curves $C$ with $C \cdot -K_X \leq n+1$.
\item Let $S$ be the cone over a smooth elliptic curve. Then $S$ is a log 
canonical surface of Picard number one, with $-K_S$ ample, but $S^0$ contains
no rational curves. Note that one can choose $C \subset S^0$ any curve, and 
then apply Mori's Bend and Break argument to produce a rational curve, but no 
matter how generally one chooses $C$, the rational curve will meet the 
singular point.
\item Pick any threefold $Y$, and any smooth curve $C$ in $Y$. Blow up $Y$ 
along $C$. The exceptional divisor is a ruled surface $\pi:\map S.C.$. Now 
blow up again, along a section of $\pi$ with very high self-intersection, 
to obtain a smooth threefold $Z$. Let $E$ be the strict transform of $S$, 
and $F$ be the general fibre of $E$ over $C$. The cone of $Z$ over $Y$ is 
two dimensional and $F$ generates one edge. Since $(K_Z + E) \cdot F = -2$,
we may, by the contraction theorem, contract $E$ down to $C$ to obtain
$X$, a $\Bbb Q$-factorial threefold over $Y$. Since $K_Z \cdot F =0$,
$X$ is Gorenstein and canonical.
Now the exceptional divisor of $X$ 
over $Y$ is dominated by $K_X$-negative curves contained in $X^0$, 
but if $Y$ itself contains no rational curves 
then every rational curve in $X$ meets $\sg (X)$. 
\item In \cite{Zhang95b} Zhang gives examples of log terminal
$\pr 1.$-fibrations $f:\map T.{\pr 1.}.$ such that 
$\kappa(-K_T) =2$ but $T^0$ is {\it not} simply connected. In particular 
$T^0$ is not rationally connected. By 
(IV.3.10.3) of \cite{Kollar96} the fibres of $f$ are the unique
family of rational curves dominating $T^0$, while any two points
of $T^0$ are connected by a smooth complete $K_T$ 
negative curve contained in $T^0$ (take a high multiple of an ample divisor).
\endroster
\endexample 

\heading 6.10 Example Computation \endheading

Here we work through one example application of (6.5) in considerable detail. 
We hope that after digesting this example, the reader will have little trouble
following the analogous, but of necessity rather terse, analysis
of \S 17 and \S 19, where we apply (6.5) to each of the roughly
sixty surfaces $S \in \goth F$.

Choose a configuration in $\pr 2.$ of a conic $B$, a secant line $A$, and a 
tangent line $C$. Let $\{a,b\}=A\cap B$ and $c=B\cap C$ (see Figure 3, (8.1)).

Let $\map \t{S}.{\pr 2.}.$ be obtained by blowing up $3$ times at $c$, along 
$C$, (that is further blow ups are at the point of the strict transform of 
$C$ lying over $c$) $5$ times at $b$ along $B$ and $7$ times at $a$ along 
$A$. Let $\map \t{S}.S.$ contract the curves on which $K_{\t {S}}$ is 
nonnegative. The resulting surface has two singular points: 
$x=(2^C,7^A,2,2,2,2)$, index $57$, and $y=(2,2,4^B,2,2,2,2,2,2)$, index $52$, 
where for example $2^C$ indicates that this $-2$-curve is the strict transform of 
$C$ (the index can be calculated using (3.1.8) of \cite{Kollaretal}).

We use Lemma (6.5) to argue that $S^0$ is uniruled. We look for rational 
curves $Z \subset S$, such that $Z \cap \sg(S)$ is a single point, $p$, $Z$ 
has two analytic branches at this point, (we will abuse notation and say $Z$ 
has a node at $p$, even though either of these branches may be singular, and 
so the singularity may not actually be a double point) and each branch has the 
same index, $m$. If $Z$ is such a curve, then (6.5) can only fail if
$-K_S \cdot Z = 1/m$ (note $m$ is in this case the global 
index of $Z$, that is $m$ is the smallest positive integer such that $mZ$ is 
Cartier).

There are several natural choices for $Z$. We follow the notation of (4.10). Thus
$\map C.Z.$ is the normalisation. One could try the strict transform
of the tangent line to $B$ at $b$. This curve has a node at $x$, one branch 
is lt (so index 57), but the other branch is $(2,7,2,2',2,2)$, and one checks 
that the index of this branch is 19.  One computes that $-K_S \cdot Z = 2/57$,
and so we can't apply (6.5). 

Another reasonable choice would be to take $Z$ 
the strict transform of the secant line from $c$ to $a$. Then $Z$ has a node 
at $y$, with one lt branch (and so the index is $52$) and the other branch
$$
(2,2,4,2',2,2,2,2,2)= (A_2,4,2',A_5).
$$
One checks that the index of this non lt branch is $26$. 
One computes $-K_S \cdot Z =2/52$ so again we can't apply (6.5).

Now take instead for $Z$ the smooth conic which has third order contact 
with $B$ at $c$ and is tangent to $A$ at $a$. $Z$ is a zero curve on $S$, with 
two branches at $y$: one is $(A_2,4,2,2',A_4)$ and the other is singular,
$(2,2',4',A_6)$. By this last notation we mean, $Z$ meets two exceptional 
divisors of $\map \t{S}.S.$, the two marked curves, at their point of 
intersection. One checks that each of these branches has index 52. Also one 
computes that
$$
-K_S \cdot Z = 2 - 42/52 -28/52 - 30/52 = 4/52
$$
and thus $S$ is log uniruled by (6.5).

\subhead 6.11 Rational curves in $S^0(E_8)$ \endsubhead

 Here we work through an interesting and non-trivial example. In fact our proof of (6.1) does not generate 
the dominating family of minimal $-K_S$-degree. Here  we determine such a family for $S(E_8)$. 
This will not be used in the proof of (1.1). 

We follow the notation of (3.6-7). We indicate the two surfaces of (3.7) by
$S_N$ and $S_C$, corresponding to a nodal and cuspidal cubic, respectively.
Let $N,M \in |-K_{S_N}|$ be the rational curves $N_1$, $N_2$ of (3.6).

 The following result should be compared with (6.6). 

\proclaim{6.11.1 Proposition} There is a dominating family of rational curves 
$G \subset S_N^0$ with $G  \cdot K_{S_N} =-2$. 
\endproclaim
\demo{Proof} Let $S = S_N$. Let $\map \Cal C.D.$ be a one parameter 
deformation of a union of two $\pr 1.$s meeting at a point, to a smooth 
$\pr 1.$, with singular fibre $\Cal C_0$. Let $f:\map \Cal C_0.S.$ be the map 
given by the normalisations of $N$ and $M$. The dimension of 
$\Hom(\Cal C_0,S)$ at $[f]$ is four (neither branch can move away from 
$N \cup M$, thus the node must get sent to $q= N\cap M$, and then we can take any 
automorphism of either branch fixing the preimage of $q$). But according to
(1.2) of \cite{KMM92}, $Hom_D(\Cal C,D \times S)$ has dimension at least 
$\chi(f^*(T_S))+\dim (D)= 5$. Thus $f$ deforms to a map of a smooth $\pr 1.$ 
into $S^0$. \qed\enddemo

Fix a flex point $q \in B^0$ to an integral plane cubic $B$. Let $L$ be the flex
line. Blowing up $8$ times at $q$ along $B$ gives $\map \tilde{S}(E_8).{\pr 2.}.$,
see the proof of (3.7).
Let $S(B)$ be the corresponding copy of $S(E_8)$. By the construction
$B \subset S(B)^0$. We note that if $B_t$ is a flat family of such elliptic curves,
the $S(B_t)$ form a flat family. Choose such a family with $B_t$ nodal
for $t \neq 0$, and $B_0=C$ cuspidal. Then there is a family $D_t$
of members of $|-2K_{S_t}|$ with $D_t=M_t+N_t$ for $t \neq 0$,
and $D_0=2C$. 

\proclaim{6.11.2  Proposition} Let $S= S(E_8)$. There is a dominating family of 
rational curves $J \subset S^0$ with $J \cdot K_S=-2$. 
\endproclaim
\demo{Proof} By (6.11.1) we can assume $S=S_C$. We consider the family $S_t$ 
defined above. By (6.11.1) we have for each $t \neq 0$ a one dimensional family 
$\Cal G_t \subset |-2K_{S_t}|$ of rational curves contained in the smooth 
locus, which includes $N_t+M_t$. After a possible base change, these lie in a 
family for all $t$. Then the limit, $\Cal G_0 \subset |-2K_{S_C}|$ gives
a one dimensional family of rational members, containing $2C$. Since
$C \subset S_C^0$, these dominate $S_C^0$. 
\qed\enddemo

Note there are only finitely many rational members of $|-K_{S(E_8)}|$ so
the family in (6.11.2) is obviously of minimal degree.

\def\aos{{\Bbb A}^1_*}
\heading \S 7 Reduction to $\alg(S^0) = \{1\}$. \endheading

Our main goal is to reduce the proof of the empty boundary case in 
(1.3) to the case when $S^0$ has 
trivial algebraic fundamental group, (7.2). We also prove, (7.10), the Picard number
one, non-empty boundary, case of (1.6).  

Additionally, we include some easy lemmas describing how $\alg$ is effected
by birational transformations, and some implications
for $\pi_1$ of the existence of certain rational curves. These will be useful
in the hunt analysis. The results 
are well known, we include proofs for the reader's convenience. 

We start with an easy Lemma. 

\proclaim{7.1 Lemma} Let $X$ be a $\Bbb Q$-Gorenstein Fano variety of
dimension $n$.
 
 If $K_X^n> n^n$ then $X$ has a special tiger.
\endproclaim
\demo{Proof} By Riemann-Roch and Serre vanishing, the dimension of the linear 
system $|-mK_X|$ grows like $m^nK^n/n !$. On the other hand, if $p$ is smooth 
point of $X$, then it is at most $\binom {mr + n -1}n$ conditions on a linear
system to have a member of multiplicity at least $mr$ at $p$. So if the
given inequality holds, there will be a member
$D \in |-mK_X|$ with multiplicity $nm$ at $p$. Then if $E$ is the exceptional
divisor of the blow up at $p$, $e(E,K_X + 1/m D) \geq 1$.
\qed\enddemo

\proclaim{7.2 Proposition} If every rank one log del Pezzo with 
$\alg (S^0)=\{1\}$ is log uniruled then every rank one log del Pezzo is log 
uniruled.
\endproclaim
\demo{Proof} Suppose $S$ is any rank one log del Pezzo. If 
$\alg (S^0)\neq \{1\}$, then there is a cover $\pi :\map S'.S.$, where 
$K^2_{S'}\geq 2K^2_S$. Running the MMP, we may assume $S'$ is a rank one log 
del Pezzo, as $K^2_{S'}$ will only increase. Moreover $S$ is log uniruled 
iff $S'$ is log uniruled. 

 If this process does not terminate, then eventually $K^2_S>4=2^2$, and we may 
apply (7.1) \qed\enddemo

We also obtain some easy lemmas describing how $\alg$ is affected by the hunt.

\proclaim{7.3 Lemma} Let $g:\map T.S.$ be a birational map between normal 
surfaces, with exceptional locus $E$. There is a natural surjection
$\map \alg (S^0).\alg (T^0)$, which is an isomorphism if $g(E) \subset S^0$.
\endproclaim
\demo{Proof} We note that if $U \subset X$ is an open subset of a normal 
variety $X$,
then $\map \alg (U).\alg (X).$ is surjective (since a connected cover cannot 
become disconnected when restricted to $U$), and if $X$ is smooth and $U^c$ 
has codimension at least two, $\alg (X)=\alg (U)$ (by purity of the branch 
locus). \qed\enddemo

\proclaim{7.3.1 Corollary} In the hunt, if $\alg(S_i)$ is trivial, and 
$\pi(\Sigma_{i+1}) \in S_{i+1}^0$ then 
$$
\alg(T_{i+1})=\alg(S_{i+1})=\{1\}.
$$
\endproclaim

\proclaim{7.4 Lemma} Let $Y_1$, $Y_2$ be two complete $\Bbb Q$-Cartier 
divisors in an integral quasi-projective variety $X$ of dimension 
at least two. If $Y_i|{Y_i}$ is nef and big, $i=1$, $2$ then 
$Y_1 \cap Y_2 \neq \emptyset$. 
\endproclaim
\demo{Proof} Cutting by hyperplanes, compactifying and desingularising, we may 
assume $X$ is a smooth projective surface. The result now follows from the 
Hodge index Theorem. \qed\enddemo

\proclaim{7.4.1 Corollary} Let $Y$ be an effective complete integral Weil,
$\Bbb Q$-Cartier divisor in a connected integral quasi-projective variety $X$ 
of dimension at least two. If the normal bundle of $Y$ is big and nef, then 
$\map \alg(Y).\alg(X).$ is surjective.
\endproclaim
\demo{Proof} Its enough to show a connected finite \'etale cover of $X$ 
restricts to a connected cover of $Y$. Let $f:\map X'.X.$ be such a cover. 
Then $f^*(Y)$ has big nef normal bundle, and so is connected by 
(7.4). \qed\enddemo

\proclaim{7.5 Lemma} Let $S$ be a normal projective surface of Picard number 
one. If $S^0$ is uniruled, then it is rationally connected, and $\pi_1(S^0)$ 
is finite.
\endproclaim
\demo{Proof} Just copy the proof of (3.4) of \cite{KMM92b}, to deduce 
that $S^0$ is rationally connected. Now apply (7.8) below. \qed\enddemo

\definition{7.6 Definition} Let $C$, $X$ be varieties. We say that $X$ is
$C$-connected, if given any two general points $x,y \in X$, there is a 
generically finite map $\map C.X.$ whose image contains $x$ and $y$.
\enddefinition

\definition{7.7 Definition} A homomorphism of groups is called {\bf almost surjective}
if the image is a subgroup of finite index. A group is called {\bf almost Abelian} if
it has an Abelian subgroup of finite index.
\enddefinition

\proclaim{7.8 Lemma} If $X$ and $C$ are quasi-projective, and $X$ is $C$-connected, then there 
is a map $f:\map C.X.$ such that the induced map $f:\map \pi_1(C).\pi_1(X).$ is almost 
surjective.
\endproclaim
\demo{Proof} Let $x \in X$ be a general point. By quasi-projectivity of (the connected 
components of) $\Hom(C,X)$, there is a variety $U$, a dominant map $g: \map U \times C.X.$, 
and a multi-section $\Sigma \subset U \times C$ (for the first projection) such that 
$g(\Sigma) = x$. Replacing $U$ by $\Sigma$ (and pulling back), we can assume $\Sigma$ is a 
section. By \cite{Kollar93} a dominant map of varieties induces an almost surjection on 
fundamental groups. $\pi_1(U \times C)$ is generated by the images of $\pi_1(\Sigma)$ and 
$\pi_1(C)$ (where the latter is included via any fibre). Clearly the image of 
$\map \pi_1(\Sigma).\pi_1(X).$ is trivial. \qed \enddemo

\proclaim{7.9 Corollary} If $X$ is $\af 1.$-connected (resp. $\aos$-connected) then
$\pi_1(X)$ is finite (resp. almost Abelian).
\endproclaim

\proclaim{7.10 Proposition} Let $(S,B)$ be a log Fano pair, of Picard
number one, with $B$ a non-empty curve.
Let $U = S \setminus (B \cup \sg(S)).$ Then $U$ is
$\aos$-connected, in particular, by (7.9), $\pi_1(U)$ is almost Abelian.

Furthermore, if $U$ is not $\af 1.$-connected, then 
$g^*(T_S(-\log B)) = \ring.(1) \oplus \ring.$ for $g:\map {\pr 1.}.C.$ the
normalisation of the general member of any family of rational curves dominating
$S^0$ and meeting $B$ once.
\endproclaim
\demo{Proof} By (6.2), $S^0$ is dominated by a family of rational curves meeting
$B$ exactly once. 
Let $g: \map {\pr 1.}.C.$ be the normalisation
of a general member of any such family. Suppose $U$ is not $\af 1.$-connected.  We have
$g^{-1}(D) = d [q]$ for some point $q \in \pr 1.$ and some $d > 0$. 
Consider the dimension, $n$, of $H=\Hom(\pr 1.,S,d [q] \subset B)$ at $[g]$. By
(5.5) and (5.3), $n =\chi(g^*T_S(-\log B))$. 
Let $s \in \pr 1.$ be a general
point, and let $x = g(s)$. The subspace
$$
H' = \Hom(\pr 1.,S,d [q] \subset B, s \rightarrow x) \subset H
$$
(where the additional notation means that $s$ is required to map to $x$) is
cut out by $2$ equations. By assumption the image of a map in $H'$ is fixed (equal to $C$). Thus
$H'$ has dimension one (there is a one dimensional family of maps with fixed image, since
in $H'$ the images of $q$ and $s$ are fixed), and so $n \leq 3$.
By (5.5), $g^*(T_S(-\log B))$ is semi-positive. By assumption it has positive degree. 
Thus the bundle is
$\ring.(1) \oplus \ring.$. 

Assume $d \geq 2$. Then $\Hom(\pr 1., S, (d-1)[q] \subset B)$ has dimension
at least $4$ at $[g]$, and the subspace, $H''$, with $s \rightarrow x$, dimension at
least $2$. Thus (since in $H''$ we fix the images of two points, $q,s$) 
the images of maps in $H''$ is not fixed. By construction they meet
$B$ at most twice. Thus $U$ is $\aos$-connected.

The final possibility is that $d=1$. Let $C_1$, $C_2 \subset S^0$ be two general members
of the family. $C_1 \cap C_2 \neq \emptyset$, since the Picard number is one. 
As in the proof of (6.11.1), we can deform $C_1 + C_2$ to get
a dominating family of integral rational curves $D \subset S^0$ with
$B \cdot D = 2$, $-K \cdot D = 4$. The reasoning above shows that
$D$ deforms with a general point fixed (no condition on the contact with $B$)
--if it did not, then the pullback of $T_S(-\log D)$ to the normalisation of $D$ would be
$\ring. \oplus \ring.(2)$.
Since $D \cdot B =2$, $D \cap U$ gives a connecting family of images of $\aos$. \qed\enddemo

Note $U$ in (7.10) is not in general $\af 1.$-connected. For example take
$(S,B) = (\pr 2., L_1 + L_2)$ for two lines $L_1$, $L_2$. $U = \af 1. \times \aos$, which
has infinite fundamental group, and thus cannot be $\af 1.$-connected by (7.9)


\heading \S 8 Preparations for the Hunt \endheading

\subhead 8.0 Introduction \endsubhead

We now take up the second half of the proof of (1.3). 
As outlined in the
introduction, the proof amounts to explicitly determining a finite set ${\goth F}$ of
$S$, which includes all $S$ with $\alg(S^0) = \{1\}$ but no tiger, 
and for each $S \in {\goth F}$, explicitly exhibiting a family
of rational curves dominating the smooth locus. Once $S$ is in hand, 
the dominating family will be obtained by finding a special rational curve
$Z \subset S$ which satisfies the criteria of (6.5) or (6.6). 
We find $Z$ in each case by a straightforward
but rather ad-hoc analysis. Here we consider the problem of finding ${\goth F}$.
As noted in the introduction, this is done by analysing possibilities for a series of 
simplifying birational transformations,
which we call the hunt, and which we now describe. In the hunt definition that follows 
there are a few choices to be made, and a few points that require proof. We will give the
precise details in (8.2). 

The hunt generates a series of pairs $(S_i,\Delta_i)$ of rank one log del Pezzo
surfaces with boundary with the properties that
\roster
\item $-(K_{S_i}+\Delta_i)$ is ample.
\item If $K_{S_{i+1}} + \Delta_{i+1}$ has a tiger, then so does $K_{S_i} + \Delta_i$.
\endroster

The process is as follows:

Given $(S_i,\Delta_i)$, we extract an exceptional divisor $E_{i+1}$ of 
$\tilde{S_i}$ via $f_{i}:\map T_{i+1}.S_{i}.$. We will discuss the choice of $E_{i+1}$ 
below. Let $\Gamma_{i+1}$ be the log pullback of $\Delta_i$ (see \S 2).
$T_{i+1}$ has Picard number two. It admits a (unique) $K$-negative contraction
$\pi_{i+1}$. Define $\Gamma' > \Gamma$ so that $K_{T_{i+1}} + \Gamma'$ is 
$\pi_{i+1}$-trivial (there is some choice in defining $\Gamma'$). 
Either $\pi_{i+1}$ is a $\pr 1.$-fibration, or it is a birational
contraction $\pi_{i+1}:\map T_{i+1}.S_{i+1}.$ to a rank one log del Pezzo.
In the first case we call $T_{i+1}$ a {\bf net} and the process stops. In
the second we define $\Delta_{i+1} = \pi_{i+1}(\Gamma')$. 

Given $S$ we inductively carry out the above, starting with 
$(S_0,\Delta_0) = (S,\emptyset)$. To prove (1.3), we assume $S_0$ has no tiger, and 
$\alg(S^0)$ is trivial, and our goal is to give 
an explicit description of $\tilde{S}$ as a blow up of $\pr 2.$. Our method is 
to at some stage determine $S_i$ (or when there is a net, $T_i$) and then recover
$S$ by classifying possibilities for the inverse transformation $S_i \dasharrow S_0$. 
The remainder of (8.0) is devoted to a detailed discussion of the techniques 
which allow us to control the possibilities. Our aim is to give the reader a feeling for 
how and why these techniques work. These techniques will prove very efficient for 
classification; we will never need to consider $S_n$ beyond $S_3$,  and in fact the vast 
majority of the analysis will be on $S_1$ and $S_2$.

In the discussion, we will use the following notation, fixed for the remainder of the 
paper. Let $A= \pi(E_1)$, $S_1 \ni q = \pi(\Sigma_1)$.  Let $a$ be the
coefficient of $A$ in $\Delta_1$, that is $\Delta_1 = a A$. We let $B \subset S_2$ be
the image of $E_2$, and let $b$ be the coefficient of $B$ in $\Delta_2$.

\subhead Choice of $E_i$ \endsubhead

The hunt is a series of $K$ non-negative blow ups followed by $K$-negative
blow downs, with the blow up and the blow down each of relative Picard number one.
Such sequences of transformations are frequently studied in the MMP, see
for example \cite{Shokurov93}. The only choice in defining such a sequence of transformations
is which divisor $E_{i+1}$ to extract. The hunt is the series of transformations 
defined by choosing $E_{i+1}$ of maximal coefficient (in $K_{S_i} + \Delta_i$).
This choice, natural from the point of view of tigers, turns out to have remarkably 
strong geometric consequences, and it is these consequences, more
than anything, which makes it possible to classify outcomes of the hunt.
In order to explain these geometric consequences, we introduce some definitions:

\subhead 8.0.1 The Flush Condition \endsubhead
{\it Let $X$ be a normal quasi-projective variety and  $\Delta= \sum a_i D_i$ 
a boundary. Let  $m=m(\Delta)$ be the 
minimum of the non-zero $a_i$. If $\Delta$ is empty, we let $m(\Delta)=1$. }
\definition{8.0.2 Definition} Let $E$ be any exceptional divisor
over $X$. We say that $K_X+\Delta$, and also the pair $(X,\Delta)$ is:
\roster
\item {\bf flush} (resp. {\bf level)} at $E$ if 
$e(E,K_X+\Delta) < m$ (resp $\leq$).
\item {\bf flush} (resp. {\bf level)} if $e(E,K_X+\Delta)< m$ (resp. $\leq$) 
for all exceptional divisors $E$.
\endroster
\enddefinition

We also have local versions: 
{\it Let $Z$ be a subset of $X$. We say $K_X+\Delta$ is flush (resp. level) at, or
around, $Z$,
if there is some neighbourhood $U$ of $Z$, so that $K_U+\Delta|_U$ is flush 
(resp. level). We say $K_X+\Delta$ is flush away from, or outside of, $Z$ 
if $K_X+\Delta$ is flush at $X\setminus Z$.}

\remark{8.0.3 Warning-Remark} Neither (1) nor (2) in (8.0.2) 
can in general be checked locally, since $m$ can increase when you shrink
$X$ (if we throw away all the components of coefficient $m$). Of course
(1) can be checked around the center, $P$, of $E$ on $X$,
if there is a component of $\Delta$ of coefficient $m$ which meets $P$. 
\endremark

We will show that in the hunt, $(S_1,\Delta_1)$ is flush (an almost immediate 
consequence of the choice of $E_1$) and that something slightly
weaker, but with the same geometric consequences, holds for all hunt stages, see (8.4.5-6).

The flush condition controls both the singularities of $X$ and of 
$\rup \Delta.$, in a manner we will now discuss. Proofs are given in (8.3).

\subhead Some geometric consequences \endsubhead

Consider first the local situation, 
a flush pair $(S,\Delta)$, where $p \in S$ is a surface germ. 
Let $A = \rup \Delta.$. 

\proclaim{8.0.4 Lemma} If $p$ is singular, then $K_S+ A$ is log terminal at $p$.
\endproclaim

Thus from the classification of log terminal pairs, recalled in Appendix $L$, 
$p$ is a cyclic singularity,
$A$ is a smooth curve germ
and $\t{A}$ meets only one exceptional divisor of the minimal desingularisation over $p$, 
one of the ends of the chain. 

We introduce some convenient notation for describing this situation. The notation
is  fixed throughout the paper:

\subhead 8.0.5 Notation for flush pairs at a point \endsubhead

We indicate the singularity type of the triple $p \in A \subset S$ by a vector
of integers $\alpha = -(E_1^2,E_2^2,\dots,E_n^2)$, where the $E_i$ are the exceptional
divisors of $\t{S}$ (over $p$) and $\t{A}$ meets $E_1$.  We will say $(S,\Delta)$,
has a singularity of {\bf type} $\alpha$ at $p$. Let $a$ be the coefficient of $A$, that is
$\Delta = a A$. Of course $(S,\Delta)$ is uniquely determined
(locally analytically) by $\alpha$ and $a$.
We note that this vector
notation is the standard way of describing a cyclic singularity, the only new
feature here is our specification of the {\it adjacent curve}, that is, the curve
meeting $\t{A}$.

The flush condition implies much more than log terminality of $A$, 
the smaller the coefficient $a$ is, the
milder the singularity. The precise statement involves a definition:
\definition{8.0.6 Definition (Spectral Value)} We will say a vector of positive 
integers $\alpha$ (and also a flush pair $(S,\Delta)$) has {\bf spectral value} $k$ if 
(in the notation of (8.0.5)) the coefficient of $E_1$ (for $K_S$) has the form $k/r$, 
where $r$ is the index of (the cyclic singularity defined by) $\alpha$.
\enddefinition

\proclaim{8.0.7 Lemma} Suppose the pair $(S,\Delta)$ is flush. Then
\roster
\item There is a unique exceptional divisor of maximal coefficient (with respect to
$K_S + \Delta$). It is the exceptional divisor adjacent to $A$, that is $E_1$ of (8.0.5).
\item If the pair $(S,\Delta)$ has spectral value $k$ at $p$, 
then the coefficient $e$ of the pair $(S,\Delta)$ at $p$ is at least
$k/(k+1)$. In particular $a > k/(k+1)$.
\endroster 
\endproclaim

Singularities with small spectral value are very restricted. For the next lemma, we will 
call the vector $\beta$ the {\bf suspension} of $\alpha$ if $\beta$
is obtained from $\alpha$ by successively adding a two on the left. 

\proclaim{8.0.8 Lemma} If $\beta = (j,\alpha)$ then the difference of the
spectral value of $\beta$ and $\alpha$ is $(j-2)r$, where $r$ is the index of 
$\alpha$. In particular the spectral value is invariant under suspension, and 
if $\alpha$ has spectral value $k$, then $\alpha$ is the suspension of  
\roster
\item"(a)" the empty string, that is $\alpha = (2,\dots,2)$, when $k=0$.
\item"(b)" $(3)$, that is $\alpha = (2,\dots, 3)$ (or $(3)$), when $k=1$. 
\item"(c)"  $(4)$ or $(3,2)$, when $k=2$.
\endroster
\endproclaim

Note that without the flush assumption, that is if all we know is that $K_S + a A$ is
log terminal, then 
completely the opposite holds: the smaller $a$ is, the more possibilities
there are for $p$ and $A$. 

If the curve $A$ is singular at $p$, then by (8.0.7), $S$ is smooth
at $p$. The size of the coefficients of $\Delta$ 
will control the singularities of $A$. 
For simplicity of exposition, suppose $\Delta = a A$ (that is each component of $\Delta$ has 
the same coefficient). 
\proclaim{8.0.9 Lemma} Assume $(S,aA)$ is flush at $p$.
\roster
\item If $A$ has multiplicity $m$ then $(m-1)a < 1$.
\item If $A$ has two smooth branches meeting to order $g$, or a unibranch singularity of 
genus $g$ then $a < \frac{g}{2g-1}$ (for the definition of these singularities see (11.1) 
and (11.2)).
\item If $a \geq 4/5$ then $A$ has normal crossings at $p$.
\endroster
\endproclaim

Now consider the global situation, a pair $(S,A)$, of a curve on a surface, such
that $(S,aA)$ is flush (as will be the case in the hunt for $(S_1,A)$). The two
cases above can be played against each other. 
If $A$ has any singularities, then these singularities put bounds on  $a$ by (8.0.9). 
These bounds control the singularities of $S$ along $A$, by (8.0.7-8). 
More concretely, suppose $A$ has a singularity of multiplicity $m$ at some point, and let
$k$ be the spectral value of $(S,\Delta)$ at a singular point $q$. Then
by (8.0.7) and (8.0.9)
$$
k/(k+1) < 1/(m-1).
$$
Thus if $m \geq 3$, $k =0$, so, by (8.0.8), $S$ is Du Val along $A$. Similarly, if $A$ 
has a double point of arithmetic genus $g$, then
$$
k/(k+1)< g/(2g-1)
$$
so for example if $g \geq 3$, $k \leq 1$. Thus by (8.0.8) the singularities along
$A$ are at worst of type $(2,\dots,3)$.

Singularities of $A$ also control to a lesser
extend the singularities of $S$ off $A$, since $a$  bounds the coefficient
of $S$ at any singular point, and singularities with small coefficient are fairly 
restricted. For example, we list the possibilities with $e(S) < 3/5$ in (10.1). 

Of course we have the same control whenever we have an upper bound on $a$. 

\subhead Implications for the hunt \endsubhead
In the hunt, $-(K_{S_i} + \Delta_i)$ is ample. Thus the (near) flushness of 
$(\Delta_i,S_i)$ sets up a useful dichotomy: 

{\it The smaller the coefficients of $\Delta_i$, the stronger our control on
the pair $(S_i,\Delta_i)$; the larger the coefficients of $\Delta_i$, the closer we are to
a tiger. }

As noted above, one way to obtain bounds on the coefficients is from singularities of 
$\rup \Delta_i.$. Thus we have the general philosophy: When $A$ or $A+ B$ is quite 
singular, we expect strong control. As an easy example:

\proclaim{8.0.10 Sample Lemma} If, in the hunt, 
$A$ has a singularity of multiplicity at least $4$, then $S_0$ is
Du Val. Any Du Val $S$ has a tiger. 
\endproclaim
\demo{Proof} Note by the definition of the hunt $a > e_0 = e(S_0)$. Thus by (8.0.9.1), 
$e_0 < 1/3$. An easy coefficient calculation shows $S_0$ must be Du Val, see (10.1). The 
last remark is an easy consequence of Riemann Roch, see (10.4). \qed\enddemo

When $\Delta_2$ has two components, the possibilities for $(S_2,A_2 + B_2)$  will
divide into a small number of fairly simple geometric configurations, which 
we now introduce. We define the configurations in general (that is outside the
context of the hunt).

{\it Let $A$ and $B$ be two rational curves on $S$, such that $K_S+A+B$ is lt at 
any singular point of $S$. We emphasize that either $A$ or $B$ may 
have singularities (at smooth points of $S$). }

The following figure will hopefully clarify some of the notation:
\smallskip
\hskip 35pt
\epsfysize=300pt
\epsffile{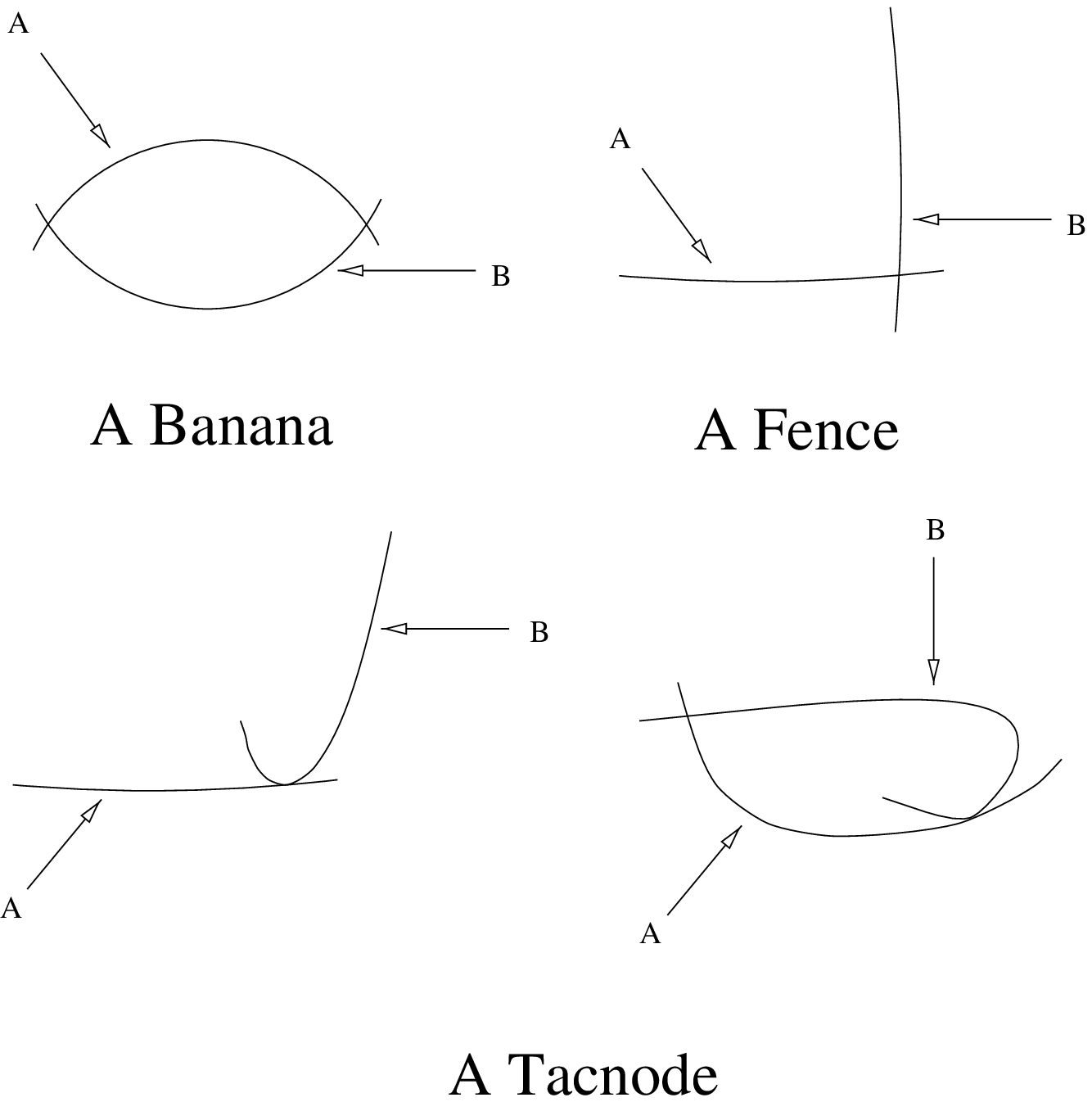}
\smallskip
\centerline{Figure 1}
\bigskip

\definition{8.0.11 Definition} We say that $(S,A+B)$ is a {\bf banana}, 
if $A$ and $B$ meet in exactly two points, and there normally.
\enddefinition

\definition{8.0.12 Definition} We say that $(S,A+B)$ is a {\bf fence} if 
$A$ and $B$ meet at exactly one point, and there normally.
\enddefinition

\definition{8.0.13  Definition} We say that $(S,A+B)$ is a {\bf tacnode} if 
$A \cap B$ is at most two points, there is one point $q \in A \cap B$ such 
that $A+B$ has a node of genus $g \geq 2$ at $q$ (for the definition, see (11.1)), and if 
there is a second point in $A \cap B$, then $A$ and $B$ meet there normally.
\enddefinition

Another tool which we use repeatedly is a partial classification of the contractions
$\pi_i: \map T_i.S_i.$ (locally analytically around $\Sigma_i$). The 
classification is easy and elementary. We carry it out in \S 11. With this classification,
once we have identified $(S_2,A+B)$ or $(S_1,A)$, we can reverse the process to
recover $S_0$.  

\subhead 8.0.14 Rough Sketch of the proof \endsubhead

The general idea of the second half of the proof of (1.3) is as follows (see
also the flow chart below). Assuming $S = S_0$ has no tiger (and simply connected smooth
locus), we run the hunt. If $A$ is sufficiently singular, we expect, by the
above philosophy, very strong control. In fact when $g(A) \geq 2$ (the arithmetic genus)
we are able to rule out all but one case, which we can explicitly construct,
(15.2). When $g(A)=1$ or $A$ is smooth we consider the next hunt step. Again we expect 
strong control if $A_2 + B_2$ has a singularity of arithmetic genus at least two,
and in fact we can rule out all possibilities except $(S_2,A_2 + B_2)$ either a smooth 
Banana, or $(S_2,A_2 + B_2)$ a fence with $B_2$ smooth and $A_2$ of arithmetic genus 1, with
a simple node. For such pairs the Bogomolov Bound (9.2) gives very strong global control, 
and we are able to give a  classification, an explicit construction of each possible case,
see (13.2) and (13.5). We also have a classification for the inverse transformation 
$S_2 \dasharrow S_0$ (see \S 11) and so can recover an explicit expression for all possible 
$S_0$. 

We note, in case the reader is puzzled by the occurrence of a node, but not a cusp, in 
the preceding paragraph, that the triviality of $\alg(S^0)$ allows us at various points 
to simplify the analysis, and in particular to rule out a fence where $A_2$ has a genus 1 
cusp. We should also note that, for simplicity of exposition, in the above sketch we have 
left out a few possibilities. Namely $T_1$ or $T_2$ could be a net, or $A$ might be 
contracted by $\pi_2$ (so $\Delta_2$ has only one component). Of course in the actual proof 
we analyse these cases, and as implicitly indicated in the sketch, rule them out.

\subhead 8.0.15 Content overview for \S 8-20 \endsubhead

The remainder of \S 8  is divided as follows. 

In (8.1) we give a concrete example of the hunt. In (8.2) we make precise the definition 
of the hunt;  there are a few choices to be made, and a few assertions that require proof.
In (8.3) we study the general properties of the flush condition. 
In (8.4) we determine to what degree flushness is preserved by the hunt. Then
we give a proposition (8.4.7) which gives a detailed description of
the possibilities for the first two hunt steps. The proof of (1.3) then proceeds
by an explicit analysis of these possibilities. This is carried out in
sections \S 14-19. 

 Next we give an overview of the sections, and then a flowchart which shows the logical 
order of the proof. \S 10-13 treat a collection of essentially independent technical issues:

In \S 10 we give a partial classification of $S$ with no tiger, and small coefficient.
The results, and notation, are used repeatedly in the hunt analysis. This will
be a main tool in the cases  $g(A) \geq 2$.

In \S 11 we give partial classifications for the contractions 
$\pi_i$, and partial classifications for $\pr 1.$-bundles. These will be
sufficient for classifying the inverse transformations $\map S_i.S_0.$
in the circumstances we will encounter in the hunt analysis.

In \S 12 we study the linear system $|K_{S_1} + A|$ in the case when $A$
is not smooth. In the cases when $g= 1,2$ the system turns
out to contain a surprising amount of geometric information. For $g =2$ it defines a 
$\pr 1.$-fibration, and for $g =1$ a birational transformation to a Gorenstein log del
Pezzo.

In \S 13 we give a partial classification of Bananas and Fences. These will
be sufficient for classifying $(S_2,A_2 + B_2)$, in the cases we will need to
consider. 

In \S 14-19 we analyse the possibilities for the hunt, according to the
breakdown of (8.4.7). These sections are essentially independent, and can be
read in any order. 

Finally in \S 20 we formally present the proof of (1.1), which
will be nothing but a series of references to the preceding sections. This section also
includes proofs of the corollaries stated in the introduction.

\newpage
\heading 8.0.16 Finding ${\goth F}$.  Assume $S$ has no tiger and $\alg(S^0) = \{1\}$. \\
Goal: Find $S$ and special rational curve  $Z \subset S$ \endheading
$$
{\aligned
&\boxed {\aligned &\text {Run Hunt with $(S_0,\Delta_0) = (S,\emptyset)$} \\
              &\text {General breakdown into cases in (8.4.7)} \endaligned} 
\Leftarrow \left \{ \text {Geometric Consequences of Flush, \S 8}\right. \\
&\qquad\qquad \Downarrow\\
&\left (\aligned &\text {First Hunt Step:} \\ 
                  &\text {$T_1$ is not a net, \S 14} \\
                  &\text {Divide further analysis according} \\
                   &\text  {to $g(A_1)$}    \endaligned \right ) 
{\aligned &\Leftarrow \left \{\text {Classification of $\pr 1.$-fibrations, (11.5)}\right. \\
                               &\Rightarrow \boxed{\aligned &\text {If $g(A) \geq 2$, \S 15} \\
                               &\text {$S$ given by (15.2)} \\
 			       &\text {$\exists Z \subset S$ \qed} \endaligned}
		        \Leftarrow \left \{\aligned &\text {Class. $S$ with small $e$ \S 10} \\
                        &\text {$|K_S + A|$ \quad \S 12  } \endaligned \right.  \endaligned} \\
& \qquad\qquad\Downarrow \qquad \qquad \qquad \qquad \qquad \searrow \\
&{\aligned &\boxed{\text{ If $A_1$ smooth, \S 18}} \\ 
           &\qquad\qquad \Downarrow \\
            &\boxed { \aligned &\text{ Second Hunt Step:} \\
                              &\text{ Rule out all but one case: } \\
                              &\text{ Smooth Banana. Classify pairs} \\
                              &\text{ $(S_2,A_2 + B_2)$, (13.2)} \endaligned } \\
           &\qquad\qquad \Downarrow \\
            &\boxed { \aligned &\text{ Classify $S$, \S 19} \\
                               &\text{ $\exists Z \subset S$ \qed }\endaligned} \\
           &\qquad\qquad \Uparrow \\
           &\boxed { \text{ Class. $\pi_1,\pi_2$, (11.1-2)}} \endaligned} \qquad \qquad 
          {\aligned &\boxed{\aligned &\text{If $g(A_1) =1$, \S 16-17} \\
         			     &\text {If $A_1$ has simple node, \S 17} \endaligned} 
                                   \Leftarrow \left \{ \aligned &\text{No Cusp, \S 16} \\
                  		     &\text{$|K_S + A|$ \quad \S 12}  \\
	      		             &\text{Du Val ldps, \S 3} \endaligned \right. \\
           &\qquad\qquad \Downarrow \\
            &\boxed{ \aligned &\text{ Second Hunt Step:} \\
                              &\text{ Rule out all but one case:} \\
                              &\text{ Fence. Classify pairs } \\
                              &\text{ $(S_2,A_2 + B_2)$, (13.5)} \endaligned} \\
           &\qquad\qquad \Downarrow \\
            &\boxed { \aligned &\text{ Classify $S$, (17.5-14)} \\
                               &\text{ $\exists Z \subset S$ \qed} \endaligned } \\
           &\qquad\qquad \Uparrow \\
           &\boxed { \text{ Class. $\pi_1,\pi_2$, (11.1-2)}} \endaligned}
\endaligned}
$$

\subhead 8.1 An example of the hunt \endsubhead

Fix a configuration in $\pr 2.$ of a tangent line, and a secant line, 
$D$ and $A$, to a conic $B$. Let $d = D \cap B$, and $\{a,b\} = A \cap B$. 
Let $h:\map \t{S}.{\pr 2.}.$ be given by blowing up $3$ times at $d$, always 
along the strict transform of $D$, $5$ times at $b$, always along the strict 
transform of $B$, and $5$ times at $a$, always along the strict transform of 
$A$. Let $\Sigma_1$, $\Sigma_2$, and $\Sigma_3$  be the $-1$-curves of 
$\t{S}$ over $a$, $b$, and $d$ respectively. 

Here, as throughout the paper, we will often use the same notation to denote 
both a curve and its strict transform under some birational transformation. 

Let $\map \t{S}.S.$ contract the $K_{\t{S}}$ non-negative curves. Then $S$ is a 
rank one log del Pezzo surface, and $\t{S}$ is its minimal desingularisation. 
$S$ has two singular points,
$$
x=x_0 = (2^D,5^A,2,2,2,2) = (2^D,5^A,A_4)
$$
and
$$
y = (A_2,4^B,A_4)
$$
where for example $5^A$ indicates the $-5$-curve of $\t{S}$ which is the strict 
transform of the secant line $A$. $x$ has index $37$ and $y$ has index $38$.
Here is a picture of the dual graph of $\t{S}$ over $x$

\bigskip
\hskip 15pt
\epsfysize=50pt
\epsffile{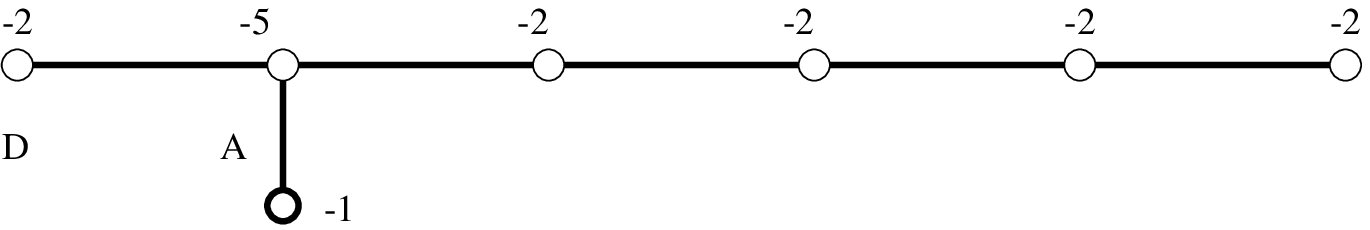}
\smallskip
\centerline{Figure 2}
\bigskip

We have chosen the notation to correspond with the outcome of
the hunt. 

$f_0:\map T_1.S_0=S.$ extracts $E_1 = A$, the $-5$-curve over $x_0$. 
$\pi_1:\map T_1.S_1.$ contracts the $-1$-curve $\Sigma_1$. $A_1 \subset S_1$ 
is a smooth $\pr 1.$, $K_{S_1} + A$ is purely log terminal, and there are 
three singularities of $S_1$ along $A$, 
$$
(2), (2,2,2,2) \text{ and } (\ul{3},2,2)
$$
where the underline indicates the curve which $A$ meets on the
minimal desingularisation. A picture of $S$, $\t{S}$ and the relevant loci
in $\pr 2.$ is given by:

\smallskip
\epsfysize=400pt
\epsffile{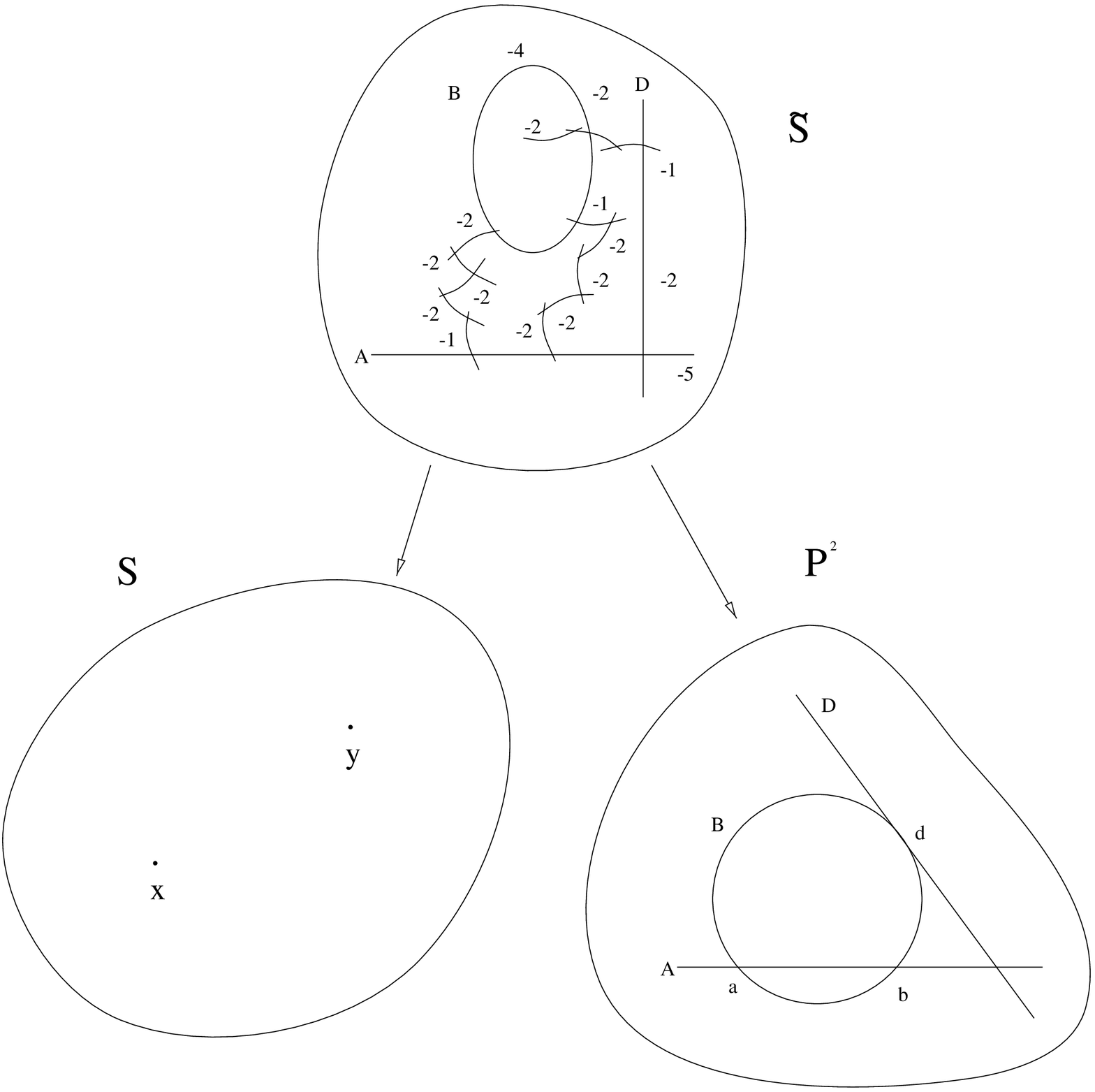}
\smallskip
\centerline{Figure 3}
\bigskip

$x_1 \in A$ is the $(3,2,2)$ point, and $f_1:\map T_2.S_1.$ extracts the 
$-3$-curve, $E_2$, which is the strict transform of $B$. $\pi_2$ contracts 
$\Sigma_2 \subset T_2$. $\Sigma_2$ meets $E_2$ at a smooth point, and meets 
$A$ at the $A_4$ point (and $A$ and $\Sigma_2$ meet ``opposite'' ends of 
the $A_4$ chain).

The configuration $(S_2,A + B)$ is a smooth banana (that is a banana in which both
components are smooth). Such configurations are classified in (13.2).
$S_2$ is a Du Val surface, with two singular points $A_1 \in A$ and $A_2 \in B$. 

$x_2$ is the $A_2$ point along $B$, and $E_3$ is the $-2$-curve meeting $B$. 
$\pi_3:\map T_3.S_3,$ contracts $\Sigma_3$. Let $C$ be the image of $E_3$. 
$S_3 =\Hbz 2.$ 
and $A + B+ C \subset S_3$ is a configuration of two sections and
a fibre, meeting normally. 

\remark{Remark} If instead on $S_2$ we chose $f_2:\map T_3.S_2.$ to blow up 
the $A_1$ point along $A$, then we would be extracting $E_3 = D$. Again 
$\Sigma_3$ is contractible, and $\pi_3$ would contract to $\pr 2.$, and
the composition of $\pi_1$, $\pi_2$, and $\pi_3$ would exactly reverse our 
construction of $S$. However this is not how the hunt proceeds, as one can 
check by computing coefficients, see (13.2).
\endremark

We now use the criterion of \S 6 to prove $S^0$ is uniruled.
We need to find rational curves on $S$ with two analytic branches
through one of the singular points, and otherwise contained in
the smooth locus. For this let $M_b \subset \pr 2.$ be the
tangent line to $B$ at $b$. On $S$, $M_b$ meets $\sg(S)$ twice, with
two analytic branches at $x$:
$$
(2',5,A_4) \text{ and } (2,5,2,2',A_2)
$$
where for example $(2,5,2,2',A_2)$ indicates the strict transform
of the branch meets only the marked $-2$-curve, and meets
it normally. Note the branch with singularity $(2',5,A_4)$ is
thus purely log terminal. Each branch has
index $37$: since the index of $x$ is prime, this can only
fail if the branch is Cartier, and this cannot happen as the
two local branches are in fact smooth, see (4.10.3) and (4.12.2). 
One computes
$$
-K_S \cdot M_b= 4/37
$$
and thus by the criterion (6.5) $S^0$ is uniruled (and dominated
by rational deformations of $37 M_b$).


\subhead 8.2 General Framework for Hunts \endsubhead

Here we make precise the definition of the hunt. A given step has three salient features:
\roster
\item The {\it underlying transformations} $f:\map T.S.$, $\pi:\map T.S_1.$.
\item The {\it scaling}, that is the definition of $\Gamma'$ (for which there
is some choice).
\item The choice of divisor $E$ extracted by $f$.
\endroster

It will be convenient to have a framework that separates the first two features from the 
third. 

\subhead 8.2.4 Scaling \endsubhead 
{\it
Begin with a pair $(S,\Delta)$ of a rank one log del Pezzo surface, and an effective
${\Bbb Q}$-divisor $\Delta$ (we do not assume $\Delta$ is a boundary).

Let $f:\map T.S.$ be an extraction of relative Picard number one,
of an irreducible divisor $E$ of the minimal desingularisation. Thus
$T$ has Picard number two, and so its cone of curves
$\ol{NE}_1(T)$ has two edges. One edge is generated by $E$. Let $R$
generate the other edge. Let $x = f(E)$.

Let $\Gamma$ be the log pullback of $\Delta$ (for the definition, see (1.16)).

We let $\Gamma_{\epsilon} = \Gamma + \epsilon E$. When 
$K_S + \Delta$ is numerically trivial, we take $\epsilon = 0$. When
$-(K_S+\Delta)$ is ample we assume $0 < \epsilon < < 1$, but we will never specify the
exact value. }

\proclaim{8.2.5 Definition-Lemma} Assume $-(K_S+\Delta)$ is ample (resp.
numerically trivial). Then
\roster
\item $R$ is $K_T$-negative, and contractible. Let $\pi$ be the associated map.
\item $K_T + \Gamma_{\epsilon}$ is anti-ample (resp. numerically trivial).
\item $\Gamma_{\epsilon}$ is $E$ negative (resp. non-positive).
\item There is a unique rational number $\lambda$ such that with
$\Gamma' = \lambda \Gamma_{\epsilon}$, $K_T + \Gamma'$ is $R$
trivial. $\lambda > 1$ (resp. $\lambda=1$).
\item $K_T + \Gamma'$ is $E$ negative (resp. trivial).
\item $\pi$ is either birational, or a $\pr 1.$-fibration.
\item If $\pi:\map T.S_1.$ is birational, 
and $\Delta_1 = \pi(\Gamma')$, then $K_{S_1} + \Delta_1$ is anti-ample
(resp. numerically trivial) and $S_1$ is a rank one log del Pezzo.
\item If $K_S + \Delta$ does not have a special tiger, then $K_S + \Delta$, $K_T + \Gamma'$, 
and $K_{S_1} + \Delta_1$ are all klt (in particular we have pure boundaries) and none of
the three has a special tiger.
\endroster
\endproclaim
\demo{Proof} We prove the anti-ample case. The numerically trivial 
case is similar, but
easier. If $K_T \cdot R \geq 0$, then $K_T$ and hence $K_S$ are
nef, a contradiction, hence $K_T \cdot R < 0$. $R$ is contractible
by the contraction theorem, see \cite{KMM87}. Hence (1). 
(2) is clear by intersecting with
$E$ and $R$, for example by Kleiman's criterion for ampleness. 
(3)
follows from (2), since $K_T \cdot E \geq 0$ by assumption. 
By (3), 
$\Gamma_{\epsilon}$ is $R$ positive (since it's effective), Now (4)
follows from (2). (6) is obvious, from (1). (5) follows from (2) and (3), and implies (7), since
$S_1$ has Picard number one. (8) is clear from the construction.
\qed \enddemo

\remark{8.2.6 Remarks} We introduce $\Gamma_{\epsilon}$ for technical reasons.
In the hunt, if $K_S + \Delta$ is flush, we want $K_T + \Gamma'$ to be
flush, but if we scale directly, without introducing $\epsilon$, we can
only in general guarantee level. See the proof of (8.4.1). The difference between the
geometric implications of flush and level are, we believe, sufficient to justify this
additional technicality, see (8.3.6.3).

The process of blowing up, scaling, and blowing down can be carried out
in much greater generality, any time one has a log minimal model program. See (2.1) of 
\cite{KMM94}.
\endremark

\subhead 8.2.7 Hunting for a tiger \endsubhead

\definition{8.2.8 Definition-Remark} We will call the transformations $f,\pi$ 
(with the associated divisors) a {\bf next hunt step} for $(S,\Delta)$ if 
$e(E,K_{S} + \Delta)$ is maximal among exceptional divisors of the minimal 
desingularisation.

Note, there may be several choices of $E$. If $x$ is a cyclic singularity, we allow any choice of 
maximal $E$, which is not a $-2$-curve (this is always possible by (10.11)) but 
if $x$ is a non-chain singularity, we require that $E$ be the central curve, which has 
maximal coefficient by (8.3.9).
\enddefinition

The main hunt we will consider is the series of transformations described in the 
introduction, beginning with $(S_0,\emptyset)$, where we assume $S_0$ does not have a 
tiger, and our goal is to explicitly determine $S_0$.
We will call this {\bf the hunt with $\Delta_0 = \emptyset$}. 
However occasionally we will run a hunt just to generate a series of (simplifying)
birational transformations, see for example the proof of (23.5). The first case will
satisfy some additional useful assumptions:

\definition{8.2.9 Definition} We will say we are {\bf hunting for a tiger} if (in addition
to the assumptions of (8.2.4)) the following hold: 
\roster
\item $K_S + \Delta$ is klt and anti-ample.
\item $(S,\Delta)$ is flush at singular points of $S$ lying in the support of $\Delta$.
\endroster
\enddefinition

The hunt with $\Delta_0 = \emptyset$ will indeed satisfy (8.2.9), see (8.4.5-7).

We will use still further notation for the hunt with $\Delta_0 = \emptyset$. Since we
have now introduced notation in several places, we will repeat it all here, so
the reader will have a convenient reference. It is fixed for the remainder of the paper.

\subheading{8.2.10 Glossary of notation for the hunt with $\Delta_0 = \emptyset$}

$f_i, \pi_{i+1}$ define the next hunt step for $(S_i,\Delta_i)$. 	
$$
S_{i} \ni x_i = f_i(E_{i+1}),\quad S_{i+1} \ni q_{i+1} = \pi_{i+1}(\Sigma_{i+1}).
$$
$\Gamma_{i+1}$ is defined by
$$
K_{T_{i+1}} + \Gamma_{i+1} = f^*(K_{S_i} + \Delta_i).
$$
$\Gamma_{\epsilon} = \Gamma_{i+1} + \epsilon E_{i+1}$.
$\Gamma' > \Gamma_{\epsilon}$ is defined in (8.2.5). It is $\pi_{i+1}$-trivial.
$\Delta_{i+1} = \pi_{i+1}(\Gamma')$ and satisfies
$$
K_{T_{i+1}} + \Gamma' = \pi_{i+1}^*(K_{S_{i+1}} + \Delta_{i+1}).
$$
$$
S_1 \supset A = \pi_1(E_1), \quad S_2 \supset B = \pi_2(E_2) \quad S_3 \supset C = \pi_3(E_3).
$$

If there seems a possibility of confusion, we will use a subscripts, $A_i,B_i, C_i$ to
indicate the strict transform of $A,B,C$ on $S_i$. 

Let $a,b,c$ be the coefficients of $A,B,C$ in $\Delta_1,\Delta_2,\Delta_3$. These
are also the coefficients of $E_1$, $E_2$, $E_3$ in $\Gamma'_1$,$\Gamma'_2$,$\Gamma'_3$.
Let $a_i,b_i,c_i$ be the coefficient of $A_i,B_i,C_i$ in $\Delta_i$.

Note that because of the scaling, the coefficients $a_i,b_i,c_i$ each 
strictly increase with $i$. Flushness and scaling imply $a_i > b_i > c_i$. .

Let $e_i$ be the coefficient of $E_{i+1}$ in $(S_i,\Delta_i)$. This is also
the coefficient of the pair $(S_i,\Delta_i)$, see (8.3.2.2) and (8.3.5.4).

We let $e_i$ be the coefficient of $K_{S_i} + \Delta_i$.

We let $E_1$ be a $-k$-curve and $E_2$ a $-j$-curve. 

If $\Sigma_i$ passes through a chain singularity $w$ of $T_i$, the notation 
$(\list a.n.)$ denotes the negative of the self-intersection of the curves 
in the minimal resolution of $w$, a prime above (resp. underline) $a_i$ means 
that $\Sigma_i$ (resp. the strict transform of $\Delta_i$) meets the corresponding curve on 
$\t{T}_i$. We let $\sb_i$ be the image of $\Sigma_i$ on $S_{i-1}$. We will write
$\sb$ for $\sb_1$.

\subhead 8.3 Properties of the flush condition \endsubhead

In this sub-section we have two main goals. We want to obtain the geometric consequences
of the flush condition discussed in (8.0), and to determine how the flush 
condition is effected by the operations of the hunt, namely extracting or contracting
a divisor, and scaling a boundary.

\noindent
{\bf 8.3.0 Notation:} {\it Throughout (8.3) we consider a boundary $\Delta$, with
support $D$.}

We begin with some simple properties that hold in any dimension.

\proclaim{8.3.1 Lemma} Let $\Delta' \geq \Delta \geq 0$, with $\Delta'$ a boundary, on a 
normal $\Bbb Q$-factorial variety $X$. Suppose $f:\map Y.X.$ is a birational morphism with 
irreducible divisorial exceptional locus, $E$. Assume $e=e(E,K_X + \Delta) \geq 0$. Let 
$\Gamma$ be the log pullback of $\Delta$,  that is $\Gamma = eE + \t{\Delta}$. Assume $\Gamma$ 
is a boundary.
\roster 
\item If $K_X+\Delta$ is level and $E$ has maximal coefficient (for $(X,\Delta)$)  then 
$K_Y+\Gamma$ is level.
\item If $K_Y+\Gamma$ is flush, then $K_X+\Delta$ is not flush for $F$ 
iff $F=E$, and the coefficient of $E$ in $\Gamma$  is at least as large
as the coefficient of some non exceptional component of $\Gamma$.
\item If $X$ is terminal and $K_X+\Delta'$ is flush then $K_X+\Delta$ is flush.
\endroster
\endproclaim
\demo{Proof} (1) and (2) are obvious.

We prove (3) by induction on the number of irreducible components in the
support of $\Delta$. If this is zero, the result holds by the definition
of terminal. If there is more than one component, then 
dropping a component clearly preserves flushness, since
$m$ (of (8.0.1)) will only go up, and the coefficient of any 
divisor will only go down. Thus we may assume $\Delta'$ and $\Delta$
have the same support, but are not equal.
Let $\Theta$ be the boundary $\Delta'-\Delta$, and $\Delta(t)$ the divisor
$\Delta+t\Theta$. Let $F$ be any exceptional divisor. Then 
$f(t)=e(F,K_S+\Delta(t))$ and $g(t) = m(\Delta(t))$, the smallest coefficient of any
component of $\Delta(t)$, are both 
affine functions of $t$. Consider the smallest value of $t$, $t_0$, such that 
$\Delta(t)$ is a boundary. $t_0 < 0$ since $0< \Theta$.
$\Delta_{t_0}$ has smaller support than $\Delta$, thus
by induction, and the definition of flush, we have 
$g(t_0)>f(t_0)$ and $g(1) > f(1)$. Thus $g(0) > f(0)$,
or equivalently, $K_X + \Delta$ is flush at $F$. \qed \enddemo

\noindent
{\bf For the rest of (8.3) we work on the germ of a log
terminal surface $p \in S$.} 

Recall that throughout (8.3), $D = \rup \Delta.$, the support of $\Delta$.

\proclaim{8.3.2 Lemma} 
\roster 
\item If $(S,D)$ has normal crossings (in particular, $p$ is smooth) 
then $K_S + \Delta$ is level, and flush if $\Delta$ is a pure boundary.
\item Assume $(S,D)$ does not have normal crossings (that is either $p$ is singular, or $D$ has 
worse than a simple node). Let $f: \map T.S.$ be a log resolution of $(S,\Delta)$
(for the definition see \S 2). Assume $e(F,K_S+ \Delta) \leq 1$ for
any $f$-exceptional divisor $F$. Then for any exceptional divisor $V$
there is an $f$-exceptional divisor $F$ such that
$e(V,K_S + \Delta) \leq e(F,K_S + \Delta)$. In particular 
$K_S + \Delta$ is flush (resp. level) iff it is flush (resp. level)
at all exceptional divisors of $f$.
\endroster
\endproclaim
\demo{Proof} For (1) it is convenient to prove a little more. We consider only the
level case (the flush case is nearly identical). We show that if $(S,A+B)$
has normal crossings, $a,b  \leq 1$, where we allow negative numbers, then 
$e(V, K_S + a A + b B) \leq \text{ min }(a,b)$ for any exceptional divisor
over $p = A \cap B$. 

$V$ can be obtained by a sequence of smooth blow ups.
We induct on the number of these blow ups.
Let $f: \map T.S.$ be the blow up at $p$, with exceptional divisor $E$.
We have
$$
K_T + a \t{A} + b \t{B} + (b + a -1) E = f^*(K_S + \Delta).
$$
Thus $e =e(E,K_S + a A + b B) = (b + a -1) \leq \text{ min}(a,b)$.
If $V$ lies on a further blow up, then $V$ is centered over some
point $x$ of $E$. Since $A$ and $B$ are separated by $f$, we may
assume (by switching notation if necessary) that $x \not \in B$. Then by induction
$$
\align
e(V,K_S + a A + b B) &= e(V, K_T + a \t{A} + b \t{B} + eE)  \\
                     &=e(V,K_T + a\t{A} + e E) \leq e \leq \text{ min}(a,b).
\endalign
$$

For (2). Write 
$$
K_T + \Gamma = f^*(K_S + \Delta).
$$
Here $\Gamma$ is in general only a subboundary (for the definition, see \S 2). Assume 
$V$ is centered over some point $x$ on an
irreducible component, $F$ of $E$. Let $l$ be the coefficient 
of $F$ in $\Gamma$. By the proof of (1),
$$
e(V,K_S + \Delta) = e(V, K_T + \Gamma) \leq l = e(F,K_S + \Delta) \qed
$$
\enddemo

\proclaim{8.3.3 Lemma} Assume $p$ is singular. If the coefficient of every exceptional 
divisor of the minimal resolution $\pi :\map \tilde S.S.$ in $K_S + D$ is strictly less 
than one, then $K_S+D$ is log terminal.
\endproclaim
\demo{Proof} Suppose $(S,D)$ is not log terminal. $D$ is not empty, since $S$ is log 
terminal by assumption. Let $E$ be the (reduced) exceptional locus of $\pi$. By (8.3.2.2), 
$f$ is not a log resolution for $(S,D)$. By the classification of lt singularities, $E$ 
has normal crossings, thus $D$ does not have normal crossings with $E$. As the computation 
of coefficients is purely numerical, we may replace $\t{D}$ by a disjoint union of analytic
discs, each meeting $E$ normally, and such that at least two discs meet the same irreducible 
component of $E$,  without changing the coefficients of any of the exceptional divisors. 
But then, from the classification of log terminal singularities in Appendix $L$, the 
pushforward of the new configuration (for which $f$ is a log resolution) is not log 
terminal, and so by (8.3.2.2), the coefficient of some exceptional divisor must be at least 
one, a contradiction.
\qed\enddemo
\remark{8.3.4 Remark} If one only assumes in (8.3.3) that the coefficients are at most one, 
one cannot in general conclude that $K + D$ is log canonical at singular points. For example 
take a cyclic singularity with exactly one exceptional curve, $E$, and assume $\t{D}$ is 
simply tangent to $E$; numerically this is indistinguishable from a log canonical node. 
However using the same proof, and the classification in Appendix $L$ one can easily classify 
the exceptions. We note in particular that if $K + D$ is not log canonical at $p$, then there
are at most two exceptional divisors in the minimal resolution over $p$. 
\endremark

\proclaim{8.3.5 Lemma} Let $f:\map T.S.$ extract the irreducible divisor $E$. Assume 
$e(E,K_S + \Delta) \geq 0$ and let $\Gamma$ be the log pullback of $\Delta$. Let 
$\pi: \map \t{S}. S.$ be the minimal desingularisation.
\roster
\item If $p$ is singular, and  $K_S+\Delta$ is flush at every $\pi$-exceptional divisor 
$F$ then $K + D$ is lt and $K_S + \Delta$ is flush. 
\item If $K_S + D$ is lt and $K_S+\Delta$ is level at every $\pi$-exceptional divisor $F$ 
then $K_S+s\Delta$ is flush for any $s>1$ such that $s\Delta$ is a boundary. 
\item Suppose $p$ is singular. If $f$ is a $K_T$-contraction and $K_T+\Gamma$ is flush
at every  exceptional divisor of the minimal desingularisation of $T$, then $K_S + \Delta$ 
is flush.
\item Suppose $p$ is singular, and $K_S + D$ is log canonical. Then for any exceptional 
divisor, $V$, there is some $\pi$-exceptional $F$ with 
$e(V,K_S + \Delta) \leq e(F,K_S + \Delta)$.
\endroster
\endproclaim
\demo{Proof} We first prove (1). Suppose the smallest coefficient of $\Delta$ is $\lambda$. 
Consider the function $f(t)=e(K_S+tD,F)$. Now $f(0)\geq 0$ and as the coefficient of 
$K_S+\lambda D$ is at most the coefficient of $K_S+\Delta$, 
$$
f(\lambda) \leq e(F,K_S + \Delta) < \lambda,
$$
where the final inequality is from the definition of flush. As $f$ is an affine function, 
$f(t)<t$, for all $t\geq \lambda $. In particular $e(F,K_S + D) < 1$, so $K_S+D$ is log 
terminal by (8.3.3). Thus, by the classification of log terminal singularities, Appendix $L$, 
$\pi$ is a log resolution of $(S,\Delta)$ and so $K_S + \Delta$ is flush by (8.3.2.2). 
Hence (1).

Now we prove (2).  The minimal desingularisation is a log resolution of $K_S + t\Delta$, 
thus by (8.3.2.2), we need only prove flushness at each $F$. Note $D$ is irreducible, say 
$\Delta = a D$. Define $f(t)$ as in (1). $f(0) \geq 0$, $f(a) \leq a$, and $f(1) < 1$. Thus 
$f(t) < t$ for $1 \geq t > a$. Hence (2).
 
Now we prove (3). By (1) we only need to check flushness at each exceptional divisor, $F$, 
of $\pi$. Since $E$ is not an exceptional divisor of $\pi$ (it is $K_T$-negative), $F$ is 
an exceptional divisor of the minimal desingularisation of $T$. Thus 
$$
m(\Delta) \geq m(\Gamma) > e(F,K_T + \Gamma) = e(F,K_S + \Delta) 
$$
and $K_T + \Delta$ is flush. Hence (3).

Finally for (4), since $D$ is log canonical, $\pi$ is a log resolution.
Thus the result follows from (8.3.2.2).
\qed \enddemo

Here are a few concrete examples: 

\example{8.3.6 Examples} Let $p \in C \subset S$ be a curve germ on $S$.
\roster
\item If $C$ has an ordinary cusp at $p$, and $p$ is a smooth point of 
the surface $S$, then $K_S+aC$ is flush iff $a<4/5$, and level iff
$a\leq 4/5$.
\item Suppose $K_S+C$ is lt at $p$, and $p$ has type $(3)$ (in the sense of (8.0.5)).
Then $K_S+aC$ is 
flush at $p$ iff $a>1/2$ and level iff $a\geq 1/2$. 
\item Suppose $K_S + C_1 + C_2$ is lc at $p$, for two irreducible curves $C_1,C_2$.
Assume $a,b > 0$. Then $K_S+aC_1+bC_2$ is flush at $p$ iff $p$ is a 
smooth point and $\max (a,b)<1$, and it is level but not flush at $p$ iff
either $p$ is smooth and 
$\max (a,b)=1$, or $p$ is Du Val and $a = b$, or $a=b=1$.
\endroster
\endexample
\demo{Proof} For (1). By (8.3.2.2) we need only consider coefficients on 
a log resolution of $(S,C)$. Such a resolution is obtained by blowing
up $3$ times over $p$, always on the strict transform of $C$. One computes
the coefficient of each exceptional divisor. The determining coefficient
is 
$$
e(\Sigma_3,K_S + a C) = 6a - 4
$$
where $\Sigma_3$ is the exceptional divisor of the final blow up (so the unique
$-1$-curve over $p$). 

For (2). We need only consider $e(E,K_S + a C)$, where $E$ is the $-3$-curve
over $p$. One computes $e(E,K_S + a C) = 1/3 + a/3$.

For (3). The flush case follows from (8.3.5.1) and (8.3.2.1). The smooth part of the
level case is easy, using (8.3.2.1). So suppose $p$ is singular. By (8.3.5.4)
we need only consider coefficients of exceptional divisors from the minimal desingularisation.
Let $E$ be such a divisor. Assume $b \geq a$.
Let $f(t,s) = e(E,K_S + t C_1 + sC_2)$. $f(0,0) \geq 0$, and zero iff $p$ is Du Val. 
$f(1,1) =1$. $f(a,b) \leq f(b,b)$ with strict inequality unless $a = b$. By affineness,
$f(b,b) \geq b$, with equality iff $p$ is Du Val, or $b =1$. The result follows. \qed \enddemo

\remark{Remark} (8.3.6.1) shows that when we scale in the hunt, 
flushness might be lost. It will turn out this can only happen at smooth points of $S$.
The other way that flushness can be lost in the
hunt, is if at some stage a component of $\Delta_i$ is contracted, (8.3.1.2). 
\endremark

\proclaim{8.3.7 Lemma} Suppose $p$ is smooth and the pair $(S, \Delta)$ is flush. 
Let $\Delta = \sum a_i D_i$. Then
\roster
\item If $M_i$ is the multiplicity of $D_i$ then $\sum a_i M_i - 1 < m$. In particular 
$m < 1/(M -1)$ where $M$ is the multiplicity of $D$.
\item If $D$ has a node of genus at least two at $p$ (see (11.1)), and the coefficient of the
two branches of $B$ at $p$ are $a\geq b$, then $2a+b < 2$. 
\item If $p$ is a cusp of $D$ (see (11.2)) and $a$ is the coefficient of the branch
of $D$ at $p$, then $a <4/5$. 
\item If $m \geq 4/5$, then $D$ has normal crossings. 
\endroster
\endproclaim
\demo{Proof} Let 
$\pi :\map R.S.$ be the blow up of $S$ at $p$. Then the coefficient of 
the exceptional divisor is $\sum a_i M_i -1 $, which is less than $m$ by assumption.
Hence (1). 

(2) and (3) may be proved similarly. (4) is immediate from (1-3), since (2-3)
are the two possibilities for $D$ of multiplicity two, see (II.8) of
\cite{BPV84}.
\qed\enddemo
 
\demo{Proof of (8.0.8)} Let $s$ be the index of $\beta$ and $r_2$ the index of the chain 
obtained by deleting the leftmost integer from $\alpha$. Now $s=jr-r_2$, and 
the coefficient of the exceptional divisor of $\alpha$ on the left is 
$(r-r_2-1)/r$. The coefficient of $E_1$ is $((j-2)r+r-r_2-1)/s$. \qed\enddemo

\demo{Proof of (8.0.7)} $K_S + D$ is log terminal by (8.3.5.1). Let $E$ be the divisor of 
maximal coefficient, and let $\pi: \map T.S.$ extract $E$. By (8.3.5.4) and (8.3.1.1), 
$K_T+\Gamma=\pi^*(K_S+\Delta)$ is level. From the classification of log terminal 
singularities, it is enough to show $G = \rup \Gamma.$ is log terminal at the point 
$p = E \cap \t{D}$. By (8.3.5.1) we need only show $K_T + \Gamma$ is flush at $p$.  
The coefficient of $\t{D}$ in $\Gamma$ is strictly minimal, since $(S,\Delta)$ is flush, 
thus $K_T + \Gamma$ is flush at $p$ by (8.3.6.3). Hence (1). 

For (2). Suppose $p$ has index $r$, and $\Delta=aD$. $K_S+D$ is log terminal, by (8.3.5.1) 
and so the different of $K_S+D$ is $(r-1)/r$, see (L.2). Now the coefficient of $K_S+aD$ 
is $(k/r)+a(r-1-k)/r$, which is less than $a$. Thus $a>e>x$, where $x$ is the 
solution to the recursive equation
$$
x=(k+x(r-1-k))/r,
$$
which is $k/(k+1)$. Hence (2). \qed\enddemo 

We note that the proof of (8.0.7) also gives

\proclaim{8.3.8 Lemma} (in the notation of (8.0.5))
$$
e(E_1,K + \lambda D) = (k/r)+\lambda(r-1-k)/r = \lambda(r-1)/r+ (1 - \lambda)(k/r).
$$
\endproclaim

\proclaim{8.3.9 Lemma} Let $x$ be a non-chain singularity with
branches $\beta_1,\beta_2,\beta_3$. Assume the central curve is
a $-l$-curve.

$e(x)$ is the coefficient of the central curve, and is a rational
number of the form $\frac{k}{k+1}$, for some positive integer $k$.

Moreover if $\beta_1=\beta_2 =(2)$, and $\alpha = (l,\beta_3)$
then $k$ is the spectral value of $\alpha$ (see figure 4)
\endproclaim
\smallskip
\hskip 70pt
\epsfysize=120pt
\epsffile{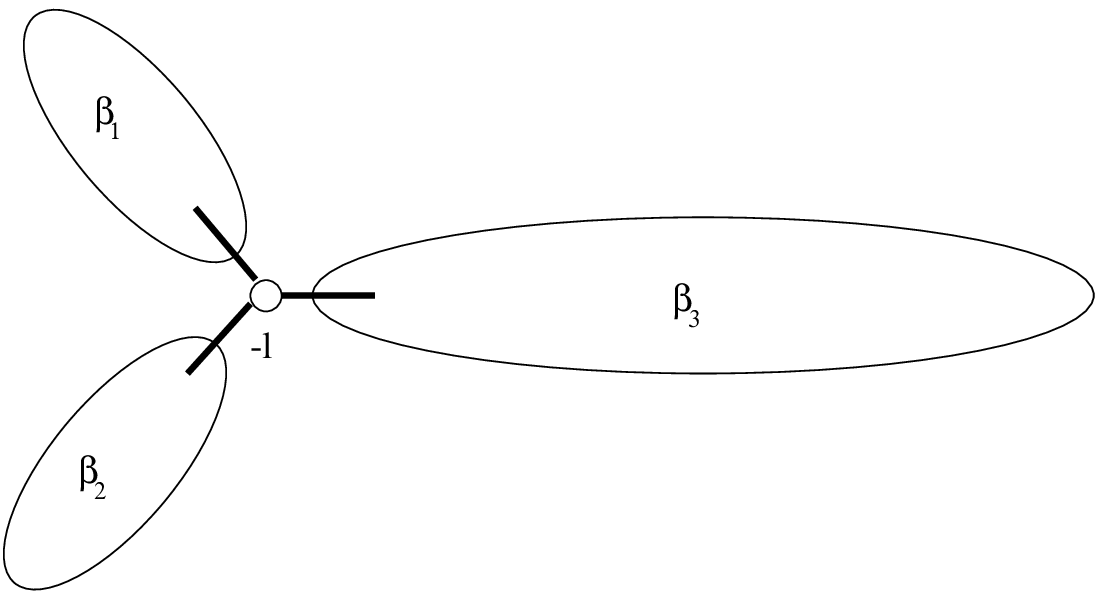}
\smallskip
\centerline{Figure 4}
\bigskip
Here ovals denote chains, as in Chapter 3 of \cite {Kollaretal}. 
\demo{Proof} We show first the central divisor has maximal coefficient. Let $f:\map T.S.$ 
extract a divisor $E$ of the minimal resolution of maximal coefficient, $e$. We assume $E$ 
is not central, thus it contains at most one non-chain singularity $y$. $K_T + E$ is not 
log terminal at $y$. By the choice of $E$, $K_T + e E$ is level at exceptional divisors
of the minimal resolution. Let $f(t) = e(F,K_T + t E) - t$ for $F$ an exceptional divisor 
of the minimal resolution over $y$. $f(0) \geq 0$, $f(e) =0$. Thus $f(1) \leq 0$. Now by 
(8.3.3), $f(1) =0$ (for some $F$), thus $f(0) =0$, and so $y$ is Du Val. Also by (8.3.4), 
$K_T + E$ is log canonical at $y$. Thus by Appendix $L$ two of the branches at $y$, are 
(2). Let the third branch be $(2,\dots,2)=A_t$ with the central curve meeting the left 
most $-2$-curve. Let the coefficients of the curves over $y$ be $e_1,\dots,e_t$. Let
$e = e_{t+1}$. Let $e_0$ be the coefficient of the central curve, $a$, $b$ the coefficients 
of the two $-2$-curves in the $(2)$ branches. One has (from the definition of the coefficient)
$$
\align
a = b &= e_0/2 \\ 
-2 e_0 + a + b + e_1 &= 0 \\
-2e_s + e_{s-1} + e_{s+1} &=0 \text{ for } t \geq s \geq 1 \\
\endalign
$$
Thus $e_0 =e_1 =\dots = e_t = e$.

Next we show $e(x)$ has the prescribed form. We use the notation and results of Chapter 3 of 
\cite{Kollaretal}. 

 Let $\beta_3$ have index $r$ and let $s$ be the index of the chain singularity
obtained by deleting the first integer on the left from $\beta_3$. 

 First suppose $\beta_1$ and $\beta_2$ both have index two. Thus
$$
\Delta=2.2.r.l-2.r-2.r-2.2.s=4(rl-r-s)
$$
and so the log discrepancy of $x$ is 
$$
\frac{2.2.r}{\Delta}(1/2+1/2+1/r-1)=1/(rl-r-s).
$$
But the index of $\alpha=rl-s$ and so the spectral value of $\alpha$ is 
exactly $rl-r-s-1$. 

 Thus we may suppose $\beta_1$ has index two, $\beta_2$ has index three,
and $3\leq r\leq 5$. Once again the log discrepancy is 
$$
\frac{2.3.r}{\Delta}(1/2+1/3+1/r-1)=\frac{6-r}{\Delta}.
$$
This deals with the case $r=5$. Otherwise we just need to prove that
$\Delta$ is divisible by three in the case $r=3$ and two in the case $r=4$. 
There are two cases, $\beta=(3)$ or $\beta=(2,2)$. 

 In the former,
$$
\Delta=2.3.r.l-5r-6s=6rl-5r-6s,
$$ 
and in the latter
$$
\Delta=2.3.r.l-4r-6s=6rl-4r-6s.
$$
The result now follows easily. \qed\enddemo

\subhead 8.4 Flush divisors in a hunt \endsubhead

\proclaim{8.4.1} Let $T$ be a surface. Let $G$ be a reduced 
curve on $T$, with $K_T + G$ lt. Let $\Gamma$ be a boundary with support $G$.
The following implications hold:
\roster
\item Assume
$K_T +\Gamma$ is level,
and $\Gamma$ has a unique component,
$E$, whose coefficient, $e$, is minimal.  Assume $e < 1$.
Let $\Gamma_{\epsilon} = \Gamma + \epsilon E$. Then 
$K_T + \Gamma_{\epsilon}$ is flush for all sufficiently small
$\epsilon > 0$. 
\item Assume any singularity of $G$ lies along a component of $G$ of minimal
coefficient. If $K_T + \Gamma$ is flush at all exceptional divisors of
the minimal desingularisation, then $K + t \Gamma$ is flush for any $t > 1$ such 
that $t \Gamma$ is a pure boundary.
\endroster
\endproclaim
\demo{Proof} 
For (1). Let $e = m(\Gamma)$ be the
coefficient of $E$. Let $e' = e + \epsilon$. Of course
$e' = m(\Gamma_{\epsilon})$. Let $V$ be any exceptional divisor. If $V$ does
not lie over $E$ then
$$
e(V,K_T + \Gamma_{\epsilon}) = e(V,K_T + \Gamma) \leq e < e'
$$
Since $e' = m(\Gamma_{\epsilon})$ it is thus enough to check flushness locally
at each point along $E$. At singular points this follows from (8.3.5.2),
while at smooth points (of $T$) we can apply (8.3.2.1). Thus (1).

For (2). Let $m = m(\Gamma)$. Note $tm = m(t\Gamma)$.
Let $V$ be an exceptional divisor, over the point
$x$ in $T$. Suppose first $x$ is a singular point. 
We can assume $V$ is an exceptional divisor of the minimal desingularisation
over $x$, by (8.3.5.4). Let
$$
g(t) = e(V,K + t \Gamma) - mt.
$$
We have $g(0) \geq 0$, $g(1) < 0$, so $g(t) < 0$ for $t \geq 1$.

Now suppose $x$ is a smooth point. If $G$ is smooth at $x$ then
$e(V,K + t \Gamma) < 0$. So assume $x$ is a node of $G$.
Then $m$ is the coefficient of one of the
components of $G$ through $x$, so we can check flushness locally around $x$,
and apply (8.3.2.1). \qed \enddemo

\definition{8.4.2 Definition} For a boundary $\Delta$ on a surface $S$, we let
$Z(\Delta)$ be the set of smooth points of $S$ where $\rup \Delta.$ (the support of
$\Delta$) is singular.
\enddefinition

\proclaim{8.4.3 Lemma} (Notation as in (8.2.4), (8.2.8)). Suppose $K_S + \Delta$ is 
anti-ample, and flush away from $Z(\Delta)$. Let $D \subset S$ be the support of $\Delta$.

Let $f:\map T.S.$ be the next step of the $K_S+ \Delta$ hunt. Assume $\Gamma'$ is a pure 
boundary. Then
\roster
\item $K_T+ \Gamma_{\epsilon}$ and $K_T + \Gamma'$ are flush away from 
$Z(\Delta)$.
\endroster
If $T$ is not a net, let $\overline{\Sigma} = f(\Sigma) \subset S$. The following hold:
\roster
\item[2] If $\overline{\Sigma}$ is a component of $D$, then 
$K_{S_1} + \Delta_1$ is not level at $q$, $q$ is a smooth point of $S_1$,
$D_1$ does not have normal crossings at $q$, and $K_{S_1} + \Delta_1$
is flush away from $Z(\Delta_1)$. Finally, if $K_S + \Delta$ is flush around
$\overline{\Sigma}$ then $\Sigma$ is the only exceptional divisor over
$q$ at which $K_{S_1} + \Delta_1$ fails to be level.
\item If $\overline{\Sigma}$ is not a component of $D$, and $\overline{\Sigma}$
is disjoint from $Z(\Delta)$ then $K_{S_1} + \Delta_1$ is flush away
from $\pi(Z(\Delta))$, and $K_T + \tilde{D} + E$ is log terminal in a neighbourhood
of $\Sigma$.
\endroster
\endproclaim
\demo{Proof} Note the final remark in (3) follows from previous statements of (3),
by (8.0.4) and the definition of $Z(\Delta)$. Note also that the assumptions and 
conclusions in (1) or (3) are such that we can shrink $S$, and assume $Z(\Delta)$ is 
empty, that is that $D \cap S^0$ is smooth. 

$K_{T_1} + \Gamma$ is level by (8.3.5.4) and (8.3.1.1). By the definition of flush, 
$e$, the coefficient of $E$ in $\Gamma$, is strictly smaller than the coefficient of 
any other component of $\Gamma$. So (1) follows from (8.4.1.1) followed by (8.4.1.2) 
(applied to $\Gamma_{\epsilon}$).

(3) is immediate from (8.3.1.2).

For (2). Suppose $\ol{\Sigma}$ is a component of $D$. Let $a$, $e'$ be the coefficients 
of $\Sigma$ and $E$ in $\Gamma'$. Then
$$
a = e(\Sigma,K_{S_1} + \Delta_1) > e' = e(\pi(E), K_{S_1} + \Delta_1).
$$
Thus $K_{S_1} + \Delta_1$ is not level at $q$. $K_T + \Gamma'$ is flush at any 
exceptional divisor of the minimal desingularisation of $T$, thus $q$ is not a
smooth point by (8.3.5.3). $D_1$ cannot have normal crossings at $q$, by (8.3.2.1). 
The final remark follows from (8.3.1.2). \qed\enddemo

\remark{8.4.4 Remarks} 
\roster
\item When we are hunting for a tiger, we can always assume $\Gamma'$ is a pure boundary 
by (8.2.5.8).
\item By definition of the hunt, $f(E)$ is a singular point, while $Z(\Delta)$ is contained
in the smooth locus, so we can think of $Z(\Delta)$ as a subset of $T$.
Note also that in the conclusion to (3), $q \not \in \pi(Z(\Delta))$.
\item Finally, in (3) the condition $\overline{\Sigma} \cap Z(\Delta) = \emptyset$
is satisfied whenever the multiplicity of $\Delta$ (that is the sum of the coefficients, 
weighted by the multiplicities of the components see (21.2)) is at least one at each 
point of $Z(\Delta)$. This will always be the case in any hunt we consider.
\endroster
\endremark

\proclaim{8.4.5 Lemma} In the hunt with $\Delta_0=\emptyset$, 
$K_{T_1} + e E_1$ is flush. If $T$ is not a net,  then $K_{S_1}+a_1A_1$ is flush. 
\endproclaim
\demo{Proof} Follow the proof of (1) and (3) of (8.4.3). \qed \enddemo

\proclaim{8.4.6 Corollary} In the hunt with $\Delta_0 = \emptyset$,
assume that 
$Z(\Delta_i) \cap \overline{\Sigma_{i+1}} = \emptyset$ for all $i > 0$ (see (8.4.4.3)).
Then $K_{S_i} + \Delta_i$ is flush away from $Z(\Delta_i)$ for all $i > 0$.
\endproclaim
\demo{Proof} Immediate from (8.4.3) and (8.4.5), by induction. \qed \enddemo

\proclaim{8.4.7 Proposition: General configuration of the first two steps of the
hunt with $\Delta_0 = \emptyset$} Notation as in (8.2.10).

For the first hunt step:

$K_{T_1} + E_1$ is log terminal, $K_{T_1}+\Gamma'$ is flush and one of the following holds:
\roster
\item $T_1$ is a net.
\endroster
otherwise $K_{S_1} + a_1 A_1$ is flush and one of the
following holds
\roster
\item[2] $g(A_1) > 1$ (arithmetic genus).
\item $g(A_1) =1$ and $A_1$ has an ordinary node at $q = q_1$.
\item $g(A_1) = 1$ and $A_1$ has an ordinary cusp at $q = q_1$.
\item $g(A_1) = 0$ and $K_{S_1} + A_1$ is log terminal.
\endroster

For the second hunt step one of the following holds:
\roster
\item[6] $T_2$ is a net
\item $A$ is contracted by $\pi_2$ (that is $A = \Sigma_2$), $K_{T_2} + \Gamma'$ is flush,
$K_{S_1} + A_1$ is log terminal,
$q_2$ is a smooth point of $S_2$, $B_2$ is singular at $q_2$, with a unibranch singularity,
and $K_{S_2} + \Delta_2$ is flush away from $q_2$, but is not level
at $q_2$. $\Sigma_2$ is the only exceptional divisor at which $K_{S_2} + \Delta_2$
fails to be flush.
\item $\Delta_2$ has two components.
\endroster

Suppose in (8) that $a_2+b_2 > 1$. Then:
\newline
(8.1): $\ol{\Sigma}_2 \cap \sg (A_1) = \emptyset$, $K_{T_2}+\Gamma'$ is flush away
$\sg(A_1)$. 
$K_{S_2}+\Delta_2$ is flush away from $\pi_2(\sg(A_1))$, and 
at least one of $-(K_{S_2}+A)$ or $-(K_{S_2}+B)$ is ample. Also, one of the 
following holds:
\roster
\item[9] $(S_2, A + B)$ is a fence. 
\item $(S_2, A + B)$ is a banana, $K_{S_2} + B$ is plt, and $x_1 \in A$.
\item $(S_2, A + B)$ is a tacnode, with tacnode at $q_2$. 
$K_{S_2} + B$ is plt. If $x_1 \in A$, $A \cap B = \{x_1,q_2\}$. 
If $x_1 \not \in A$ then $A \cap B = \{q_2\}$.
\endroster
\endproclaim
\remark{Remark} In (8.1), $\sg(A_1)$ (which is either empty, or $q_1$) is a smooth point 
of $S_1$ and so $f_1$ is an isomorphism around $\sg(A_1)$, thus we can think of 
$\sg(A_1)$ as a subset of $T_2$. 
\endremark

\demo{Proof of (8.0.14)} Everything through (5) is clear, using (8.4.3).

Consider the second hunt step, and suppose $T_2$ is not a net.

If $A$ is contracted, then obviously $A = \Sigma_2$, and we have case
(5) of the first hunt step. The rest follows from (8.4.3.2).

Now suppose $\Delta_2$ has two components, and $a_2+b_2 \geq 1$. 
In particular $a_2>1/2$, and so $\Sigma_2$ cannot meet a singular
point of $A$. All the flush remarks follow from (8.4.3.3). 
$A_2 + B_2$ cannot have a triple point at $q_2$, by (8.3.7). 
So $A \cap E_1 \cap \Sigma_2 = \emptyset$, and $A \cap \Sigma_2$, 
$E_1 \cap \Sigma_2$ are each at most one point (and the second cannot be empty 
as $B^2 > 0$). By (8.3.5.1), $K_{S_2} + A + B$ is lt at singular points of 
$S_2$. If $X$ is the smaller of $A$ and $B$, then $-(K + (a_2 + b_2) X)$
is ample. The rest follows from easy set theoretic considerations. \qed \enddemo

\heading \S 9 Bogomolov Bound \endheading

This section is logically independent from the rest of the paper and is
(we hope) of independent interest.

We refer to Chapter 10 of \cite{Kollaretal} for the definition and basic 
properties of $\Cal Q$-sheaves and $\hat \Omega ^1_S$.

The main result of this section is the following:

\proclaim{9.1 Theorem} Suppose $B \subset S$ is a reduced curve on 
a surface of Picard number one, and $K_S+B$ is lc. Then 
$\hat c_2(\hat \Omega ^1_S(\log B))\geq 0$, with strict inequality when $K_S+B$ is ample.
\endproclaim
\demo{Proof}  By (10.14) of \cite{Kollaretal} we can assume $K_S+B$ is 
anti-ample. We slightly modify an argument due to 
Kawamata \cite{Kawamata92}. We may assume $\hat \Omega ^1_S(\log B)$ is 
not semi-stable. Thus there is a saturated destabilising $\Cal Q$ subsheaf 
$\Cal L$, cf. the proof of (10.11) of \cite {Kollaretal}. Suppose 
$c_1(\Cal L)=t(K_S+B)$. The crucial observation 
is that $h^0(S, n\Cal L)<cn$ for some constant $c$, cf. the proof of (10.14) of
\cite{Kollaretal}. Thus $\Cal L$ is not ample, and $t\leq 0$. Now proceed, 
exactly as in \cite {Kawamata92}. \qed\enddemo

\proclaim{9.2 Corollary} Let $S$ be a log terminal surface of Picard number 
one and $B \subset S$ a reduced curve such that $K_S+B$ is lc. 
Let $U=S\setminus B$.
$$
\sum_{p\in U} \frac{r_p -1}{r_p} \leq  \ec{S} - \ec{B}
$$ 

If $S$ is rational, and $B$ has arithmetic genus $0$ and $s$ irreducible
components, then
$$
\sum_{p\in U} \frac{r_p -1}{r_p} \,
\left\{ \aligned \quad \leq 3 \qquad & \text {if $s=0$} \\
                 \quad \leq 1 \qquad & \text {if $s=1$ }\\
                 \quad =0 \qquad & \text {if $s=2$}.
\endaligned
\right.
$$
\endproclaim
\demo{Proof} The first inequality  follows from (10.7) and (10.8) of 
\cite {Kollaretal}. If $S$ is rational, then $\ec{S} =3$. 
$B$ is a connected tree of $\pr 1.s$. 
If $B$ is not empty, then $\ec{B} = 1 + s$. Hence the second set of
inequalities follows from the first inequality.
\qed\enddemo

\proclaim{9.3 Corollary} Fix $\epsilon > 0$, and a birational
equivalence class  $\Cal C$ of projective surfaces. The collection
$$
X=\{\, S\in \Cal C \, |\, \rho(S) =1, e(S) < 1 - \epsilon\,\}
$$
is bounded.
\endproclaim
\demo{Proof} Note in particular that $e(S) < 1$, so $S$ has quotient and thus rational 
singularities. For $S \in \Cal C$ with rational singularities, $\ec{S}$ depends only 
on $\rho(S)$. Thus (9.2) implies any surface $S\in X$ has at most $2 \ec{S}$ singularities. 
By the classification of quotient singularities, it is enough to bound $\rho(\t{S})$ for 
$S \in X$, and thus to bound $\{\, -K_{\t{S}}^2\, | \, S \in X\, \}$ from above, where 
$f:\map \t{S}.S.$ is the minimal desingularisation.

Write $f^*K_S=K_{\t{S}}+\sum e_i E_i$, and let $r_i=-E_i^2$. We have
$$
\align
-K_{\t{S}}^2 &= -K_S^2 + \sum e_i K_{\t{S}} \cdot E_i \\
             &\leq \sum e_i (r_i -2) 
\endalign
$$
Since $e_i < 1$, it is enough to bound the number of $E_i$ with $r_i \geq 3$.
This is bounded in terms of $\epsilon$ by (3.1.12) of \cite{Kollaretal} and 
simple coefficient calculations. \qed \enddemo

\subhead An interesting example in non-zero characteristic \endsubhead

Note in particular, the Bogomolov bound implies a rank one log Del
Pezzo can have at most five singular points (there is no rank one log del Pezzo 
with exactly six $A_1$ singularities). Here we show there is no bound on the number of 
singularities in characteristic two. 

 Take a strange conic in $\pr 2.$, and blow up the point where all the 
tangents meet. Now we have a curve $C$ in $\Hz 1.$, which is purely inseparable
over the base of degree two. 

 Pick any fibre and perform the following operation. Blow up along $C$ twice
and contract the two $-2$-curves. If we perform this operation $k$ times,
$C$ will now have self-intersection $4-2k$, and will lie in the smooth 
locus. Contracting $C$ we get a log del Pezzo of rank one, with $2k+1$ 
singular points. As a special case ($k=3$) we get a Gorenstein log del Pezzo
with seven nodes. 


\heading \S 10 Riemann Roch and Surfaces with Small Coefficient \endheading

We are seeking to partially classify those rank one log del Pezzo 
surfaces with coefficient less 
than $2/3$ which do not have tigers. Our method is simply to use 
Riemann Roch and the Bogomolov bound. In fact for our classification will require
much less than the absence of a tiger, see (10.8). 

It will be convenient (at various points in the paper) 
to have the following list of singularities of small 
coefficient:

\proclaim{10.1 Proposition} Let $x$ be an lt singularity, with
$e=e(x)$. If $0<e < 3/5$, the possibilities for $x$ are 
as follows:
\item{(1)} $e < 1/2$: $(3,A_j)$. $e= \frac{j+1}{2j+3}$.
\smallskip
\item{(2)} $e =1/2$:
\smallskip
\item\item{(a)} $(4)$
\smallskip
\item\item{(b)} $(3,A_j,3)$
\smallskip
\item\item{(c)} $(2,3,2)$
\smallskip
\item\item{(d)} $x$ is a non-chain singularity, with center $(2)$
and branches $(2)$, $(2)$, and $(A_j,3)$, with the central curve 
and the unique $-3$-curve meeting opposite ends of the $A_j$ chain.
This is the only non-chain singularity with $e < 2/3$.
\smallskip
\item{(3)} $1/2 < e < 3/5$:
\smallskip
\item\item{(a)} $(2,3,A_j)$ with $2 \leq j \leq 4$. $e = \frac{2j+2}{3j+5}$.
\smallskip
\item\item{(b)} $(4,2)$, $e=4/7$.
\endproclaim
\demo{Proof} (L.1) reduces the result to straightforward checking.
\qed\enddemo

\subhead 10.2 The following notation is fixed throughout the section \endsubhead

{\it $S$ is a rank one log del Pezzo, and
$\map \t{S}.S.$ is the minimal desingularisation. $\rho$ is the
relative Picard number of $\map \t{S}.S.$. 

For a singular point $p \in S$, $\rho(p)$ indicates the contribution of $p$
to $\rho$, that is the number of exceptional divisors lying over $p$.

Let $r:\map \overline S.S.$ be the minimal resolution of $S$ over the 
non Du Val locus. Thus $r$ resolves precisely the non Du Val singularities.
In particular, $\ol{S}$ 
has Gorenstein singularities. 

Let $D \subset \overline{S}$ be the (reduced) sum of the $r$-exceptional 
divisors, and $n$ the number of non Du Val singularities. Let 
$\Gamma=K_{\overline S}\cdot\sum(1-e_i)E_i$, where the sum is over the 
$r$-exceptional divisors, and 
$e_i = e(E_i,K_S)$. 

Since $H^0(2K_{\overline S} + D) = 0$ we have
by Riemann Roch
$$
\align
 h^0(-K_{\overline S} - D) \geq \chi (-K_{\overline S}-D)&=1-n+K_S^2+\Gamma \\
                          &=1-n+K_{\t S}^2+K_{\overline S}\cdot D\tag{10.3}
\endalign
$$

For a (non Du Val) singular point $t$ 
we let $\Gamma(t)$, $\delta(t)$, and $w(t)$ be the 
contribution of $t$ to $\Gamma$, $-K_{\t{S}}^2$, and $K_{\ol{S}} \cdot D$ respectively. 
Thus
$$
\align
\Gamma(t) &= \sum (1 - e_j) K_{\t{S}} \cdot E_j = \sum (1 - e_j)(-2 - E_j^2) \\
\delta(t) &= \sum e_j (-2 - E_j^2) \\
w(t)      &= \sum K_{\t{S}} \cdot E_j = \sum (-2 - E_j^2) \\
\endalign
$$
where in each case the sum is over those $j$ with $E_j$ lying over
$t$.
We also refer to $w(t)$ as the {\bf weight} 
of $t$. Let $w$ be the sum of the weights of the singular points, and 
$\beta$ be the weighted average of the coefficients, 
that is
$$
w \cdot \beta = \sum w(t) e(t).
$$
Observe
$$
\beta \leq (1/w) \sum w(t) e(S) = e(S).
$$
We let $\Delta(t)$ be the index of $t$ (that is the order of the local fundamental group).

In the hunt with $S = S_0$, we let $\alpha$ be sum of 
$\delta(t)$ for $t \not \in \sb$. Recall from (8.2.10) that $\sb \subset S_0$
is the image of $\Sigma_1 \subset T_1$.}

\proclaim{10.4  Lemma} (notation as in (10.2)) 
If $n\leq 1$, or $n=2$ and $e(S)\leq 1/2$, then $|-K_{\ol{S}} - D|$ is non-empty,
in particular any exceptional divisor over a non Du Val point is a tiger.
\endproclaim
\demo{Proof} By (10.1), $\Gamma(t)\geq 1/2$, if $e\leq 1/2$. Since $K_{S}^2 > 0$,  
$|-K_{\ol{S}} - D|$ is non-empty, by (10.3). \qed\enddemo

It is convenient to list all four-tuples of integers at least $3$, 
which can occur as indices on a rank one log del Pezzo (that is integers
satisfying the Bogomolov bound). We note in each case that by the 
Bogomolov bound there can be no other singularities on the surface.
\newline
\noindent
{\bf 10.5}
\roster
\item"(a)" $(3,3,3,m)$, $m \geq 3$, or
\item"(b)" $(3,3,4,m)$, $4 \leq m\leq 12$, or
\item"(c)" $(3,3,5,m)$, $5 \leq m \leq 7$, or
\item"(d)" $(3,3,6,6)$, or
\item"(e)" $(3,4,4,m)$ $4 \leq m \leq 6$, or
\item"(f)" $(4,4,4,4)$.
\endroster

\proclaim{10.6 Lemma} Notation as in (10.2). 
Let $S$ be a rank one log del Pezzo such that $|-K_{\ol{S}} - D|$ is empty.
Suppose $S$ has $m$ singularities of index at least three. 
  
If $m>3$, then the singularities of $S$ are, either
\roster
\item $(3)$, $(3)$, $(2,2,2)$ and $(3,2,2,2,2)$, or
\item $(3)$, $(3)$, $(3,2)$ and an $A_6$-point, or
\item  the indices are given by (10.5.a).
\endroster
\endproclaim
\demo{Proof} Suppose not. 

Note first that by the classification of quotient singularities
\cite{Brieskorn68} there is only one non cyclic singularity, $D_4$, of
index (that is order of local fundamental group) less than 16. 

Using (L.1) it is simple to classify singularities of small index.
For example the only non Du Val  singularities of index at most $7$ are
$$
(7),(6),(5),(4),(3),(4,2),(3,2),(3,2,2).
\tag{10.7}
$$
We will use this, as well as (10.1), repeatedly and without remark.

Clearly $m=4$, using the Bogomolov bound. We consider the
above list of possible indices. Let 
$z$ be the fourth singularity in the list. We first show that cases (d)-(f) 
above do not occur. Note $\rho \geq 10$, for otherwise
$K_{\t{S}}^2 \geq 0$ and $|-K_{\ol{S}} - D| \neq \emptyset$  by (10.3).

 In case (d), one of the index six points must be $A_5$ (or $\rho \leq 6$) 
the other type $(6)$ (otherwise $n \leq 2$ and $e_0 \leq 1/3$ violating 
(10.4)). But then both index three points are Du Val (or $\rho \leq 9$) and 
$n=1$, contradicting (10.4). 

In case (e), the points other than $z$ cannot all be Du Val (or $n \leq 1$ 
violating (10.4)). Thus $\rho(z) \geq 3$. Thus $z$ is Du Val. 
Since $n \geq 2$, $\rho \leq 9$ unless the singularities
are $(3)$,$(4)$, $A_3$ and $z = A_5$. This violates (10.4).

In case (f), $n \geq 3$ by (10.4), thus $\rho \leq 6$, a contradiction.

 In case (b), consider $w(z)$. If it is zero, then one checks $(n-1)-\Gamma \leq 0$, 
contradicting (10.3). So $w \geq 1$. $n > 1$ by (10.4). So $\rho - \rho(z)$ is at 
most $6$. Thus $\rho(z) \geq 4$. From $\Delta(z) \leq 12$ it follows that $z = (3,A_j)$. 
Thus by (10.4), $n \geq 3$. (1) is the only possibility with $\rho \geq 10$. 

In case (c), suppose first $z$ is non Du Val. Since $n \geq 2$, $\rho - \rho(z) \leq 7$, 
thus $\rho(z) \geq 3$, which implies $z = (3,2,2)$. It follows the singularities are $(3)$, 
$A_2$, $A_4$ and $(3,2,2)$. This violates (10.4). Thus $z$ is Du Val. Then one checks 
$(n-1)-\Gamma \leq 0$ unless we have (2). \qed\enddemo

\proclaim{10.8 Lemma} Notation as in (10.2). Assume $e(S) < 2/3$ and 
$|-K_{\ol{S}} - D| = \emptyset$. We have the following 
possibilities. Either
\roster
\item $n=4$, $w\geq 4$, three of the singularities are $(3)$, or
\item $K^2_{\t{S}}=-3$, $n=3$, $w=5$, $3/5 < \beta \leq e(S) $, or
\item $K^2_{\t{S}}=-2$, $n=3$, $w=4$, $1/2 < \beta \leq e(S) $, or
\item $K^2_{\t{S}}=-1$, $n=3$, $w=3$, $1/3< \beta \leq e(S) $, or
\item $K^2_{\t{S}}=-1$, $n=2$, $w=2$, $1/2 < \beta \leq e(S) $.
\endroster

Furthermore, if $e(S) < 1/2$ then $K^2_{\t{S}} =-1$, $e > 1/3$,  and 
either (4) holds, or
(1) holds, and the last singularity is $(3,A_6)$.
\endproclaim
\demo{Proof} 
By the Bogomolov bound, $n \leq 4$. By (10.3) and the 
definitions we have the following formulae:
$$
0 < - K_{\t{S}}^2 < 2w/3, \quad w \geq n
$$
$$
w \cdot \beta \geq \sum e_i K_{\ol{S}} \cdot E_i = K_S^2 - K_{\ol{S}}^2 > 
-K_{\ol{S}}^2
$$
$$
\Gamma > w/3, \quad \beta \leq e, \quad \Gamma < n -1, \quad -K^2_{\t{S}} \geq w - n + 1.
$$

If $n=4$ we have (1) by (10.6). 

Suppose $n=2$. Then $w/3 < \Gamma < 1$. Thus $w=2$ and $K^2_{\t{S}}=-1$. 

Suppose $n=3$. $w/3 < \Gamma < 2$ so $w < 6$. The corresponding values of 
$K_{\t{S}}^2$, and $\beta$  follow from the formulae above. 

Finally suppose $e(S) < 1/2$. Then by (10.1), $w = n \leq 4$, and 
$0< -K^2_{\t{S}}< w/2 \leq 2$. Thus $K^2_{\t{S}} = -1$, and $\beta \leq e < 1/2$. So if we 
do not have (4), then $n=w=4$. The last singularity has $\rho = 7$, and thus by (10.1) is
$(3,A_6)$. \qed\enddemo

 Here we house some simple results, which might otherwise be on the streets. 

\proclaim{10.9 Lemma} Notation as in (8.2.10).
$$
K_{S_0}^2 = \frac{(K_{S_0} \cdot \sb)^2}{\sb^2}
$$
\endproclaim
\demo{Proof} Immediate since $S_0$ has rank one. \qed \enddemo

\proclaim{10.10 Lemma} Let $(S, p)$ be a  lt germ, and $C\subset S$
an analytically irreducible curve germ at $p$. Suppose $p$ has
weight one, and $\t{C}$ meets only one exceptional divisor, $E$, of the 
minimal desingularisation over $p$, the unique $-3$-curve over $p$, and 
meets $E$ normally. Then the coefficient of $E$ in the pullback of $C$ is 
$e(p)$. \endproclaim
\demo{Proof} Follows since $\t{C}$ and $-K_{\t{S}}$ have the same 
intersection with each exceptional divisor of the minimal desingularisation 
over $p$. \qed \enddemo

\proclaim{10.11 Lemma} Let $x \in S$ be a cyclic non Du Val singularity. Then there is 
an exceptional divisor $E \subset \t{S}$ of the minimal desingularisation over $x$ of 
maximal coefficient (for $K_S$) such that $E^2 \leq -3$. \endproclaim
\demo{Proof} An easy coefficient calculation. \qed \enddemo


\heading \S 11 A partial Classification of $K_T$-contractions \endheading

The goal of this section is to give a partial classification for the
contractions $\pi_i$ in the hunt. We first consider birational contractions,
(11.0-4), and then $\pr 1.$-fibrations in (11.5).

\subhead 11.0 Birational Contractions \endsubhead

Our main objective is a local analytic classification of the possibilities
for 
$$
\pi_i:\map (T_i,\Sigma_i + D_i).(S_{i+1},D_{i+1}).
$$
for the first two hunt steps, $i=1$, $2$. By
local analytic, we mean locally analytically about $q_{i+1} \in S_{i+1}$. 

As such a classification will be useful at other points in the paper, 
we will work in the abstract (outside the context of the hunt).
\newline
{\it 11.0.1 We fix some notation. Let $\pi :\map T.S.$ be a proper birational contraction of a 
$K_T$-extremal ray. In particular $\rho(T/S) =1$.
Let $\Sigma$ denote the exceptional divisor of $\pi$, which is a $-1$-curve.
Let $q$ be the image of $\Sigma$. We will only be concerned with the local analytic description 
of $T$ about $\Sigma$, thus $p \in S$ is only an analytic surface germ. Let $W \subset T$ be a 
curve with smooth components crossing normally. We assume $W$ has at
most two irreducible components $X,Y$ ($Y$ is allowed to be empty). We assume $K_T+cX+dY$ is 
$\pi$-trivial and flush, with $0 < c,d < 1$. Thus by (8.0.4), $K_T + X + Y$ is log terminal.  In
particular $X$ and $Y$ cannot meet at any singular point of $T$.

We assume $\pi|_W$ is finite (that is $\Sigma$ is not a component of $W$). Thus
by (8.3.1.2), $K_S + cX + dY$ is flush. We assume $X$ and $Y$ have analytically
irreducible images at $q$, and thus each (is either empty or) meets $\Sigma$ at exactly one point. 

Let $D \subset S$ be the image (with reduced structure) of $W$.

At a chain singularity $(a_1,\dots,a_{\rho})$ a prime (resp. underline) indicates normal 
contact with $\Sigma$ (resp. $W$). Note this is the same notation we use in the hunt, see 
(8.2.10).

We let $h:\map \t{T}.\t{S}.$ be the induced map between the minimal desingularisations, 
a composition of smooth blow ups. }

We will also use the symbols $X$, $Y$ to indicate the images of these curves on $S$. This
is in keeping with our convention of occasionally using the same symbol 
to denote a divisor, and its strict transform under a birational map.

The starting point for our classification is the following simple observation:

\proclaim{11.0.2 Lemma} The centre of each blow up of $h$ lies on the exceptional divisor of 
the previous blow up. 
\endproclaim
\demo{Proof} $\t{\Sigma}$ is the only $-1$-curve lying over $q$. \qed \enddemo

We classify configurations by considering possibilities for the sequence of blow ups $h$. 
We think of $(T,W)$ as being obtained from $(S,D)$ by the sequence of blow ups $h$. 

We note next that (under some mild conditions) the coefficients $c$, $d$ do not decrease
with additional blow ups:

\proclaim{11.0.3 Lemma} If $(T',X',Y',\Sigma')$ is obtained from $(T,X,Y,\Sigma)$ by further
blow ups, that is if $h':\map \t{T'}.\t{S}.$ factors through $h:\map \t{T}.\t{S}.$
then 
$$
(c' X + d' Y) \cdot \Sigma \geq (c X + d Y) \cdot \Sigma
$$
In particular if $c = d$ and $c' = d'$, or $d'\leq d$, then $c' \geq c$. 
\endproclaim
\demo{Proof} Let $f$ be the induced map $\map \tilde T'.\tilde T.$. Let $p$, $p'$ be the 
minimal desingularisation maps. Note that 
$$
p_*f_* {p'}^*(K_{T'} + c' X' + d' Y') = K_T + c' X + d' Y + \alpha \Sigma
$$
for some $\alpha \geq 0$, and that the above divisor is pulled back from $S$, and
thus is $\Sigma$ trivial. Hence
$$
\align
(c'X + d' Y) \cdot \Sigma &= - \alpha \Sigma^2 - K_T \cdot \Sigma \\
                         &=  -\alpha \Sigma^2 + (cX + dY) \cdot \Sigma.	
\endalign
$$
As $\Sigma$ has negative self-intersection, the result follows. \qed\enddemo

When $D$ is smooth at $q$, then $K_S+D$ is log terminal at $q$, by (8.0.4). We will 
consider this case in (11.4).

Suppose $D$ is singular at $q$. By (8.0.4), $S$ is smooth. Let $g$ be the local 
arithmetic genus of $D$. We first classify those $D$ with a double point, (11.1-2).
Higher multiplicity is considered in (11.3). 

Double points for curves on a smooth surface germ are classified in
II.8 of \cite{BPV84}. There are two possibilities, and we consider each in turn.

\subhead 11.1 $D$ has two smooth branches at $q$ \endsubhead

In this case, $D$ has what is called an $A_{2g}$-singularity.
$X \subset S$ and $Y \subset S$ are both smooth, and meet to order $g$ at $q$. We call 
this a {\bf node of order $g$.}

Consider $h:\map \t{T}.S.$. If $g \geq 2$, we must first blow up repeatedly at 
$X \cap Y$ 
until $X$ and $Y$ cross normally, for otherwise $K + W$  will not be lt. We call 
this {\bf configuration $0$.} 
Here $X$ and $Y$ meet normally at a smooth point of $T$ along
$\Sigma$. $\Sigma$ contains a unique singular point, $(2',\dots,2)=A_{g-2}$
(and $g \geq 2$). Clearly $c+d=1$. 

\remark{Remark} Note that in this case $W$ is singular, so  
configuration $0$ cannot occur for $T = T_1$ of
the hunt, with $W = E_1$, since $E_1$ is smooth.
\endremark

If there is a further blow up it must be at $X \cap Y$, or otherwise
$X$ and $Y$ will meet at a singular point of $T$. We call the
resulting configuration,  
{\bf configuration I.}  Here $X$ and $Y$ each meet $\Sigma$ at a smooth point, 
there is one singular point $(2',2,\dots, 2)=A_{g-1}$ and once again
$c+d=1$. 

Consider the next blow up. It must be at the point where $\Sigma$ meets either
$X$ or $Y$, for if we blow up at any other point, then $X$ and $Y$ will meet 
at a singular point of $T$. Switching $X$ and $Y$ if necessary, we may assume
the next blow up is at the point of $X$. 

We call the resulting configuration, {\bf configuration $II$.} 
Here there is a unique
singularity, $A_g$ along $\Sigma$, $X$ meets $\Sigma$ at a smooth
point. On $\t{T}$, $Y$ and $\Sigma$ are disjoint, but meet the 
same $-2$-curve of the $A_g$ chain, an end of the chain. The $A_g$ point
has type $(\ul{2'},\dots,2)$. One computes 
$gd + (g+1)c= g+1$.

Now at every further step there are two possible point we could blow up,
either the point $x \in \Sigma \subset \t{T}$ lying over 
$\Sigma \cap X \in T$, or the point $y \in \Sigma \subset \t{T}$
lying over $\Sigma \cap Y \in T$ (if we blow up at any other point, then
$X$ and $Y$ will meet at a singular point of $T$). We call these possible
blow up points the points {\bf nearest} $X$ and $Y$ respectively.
We indicate the sequence of blow ups by
a string of $x$s and $y$s, for example $(II,x,x,y)$ means the configuration
obtained from $II$ by blowing up twice at $x$ and once at $y$. 

Suppose $g \geq 2$. Note that if any of the blow ups is at the point
nearest $X$, then $K_T + Y$ is not log terminal. Thus for $g \geq 2$ the
only possibilities are $(II,y,\dots,y)$, which we indicate by $(II,y^s)$, where
$s$ is the number of $ys$ in the string.

Figure 5 gives a picture of some possible configurations:
\smallskip
\epsfysize=370pt
\epsffile{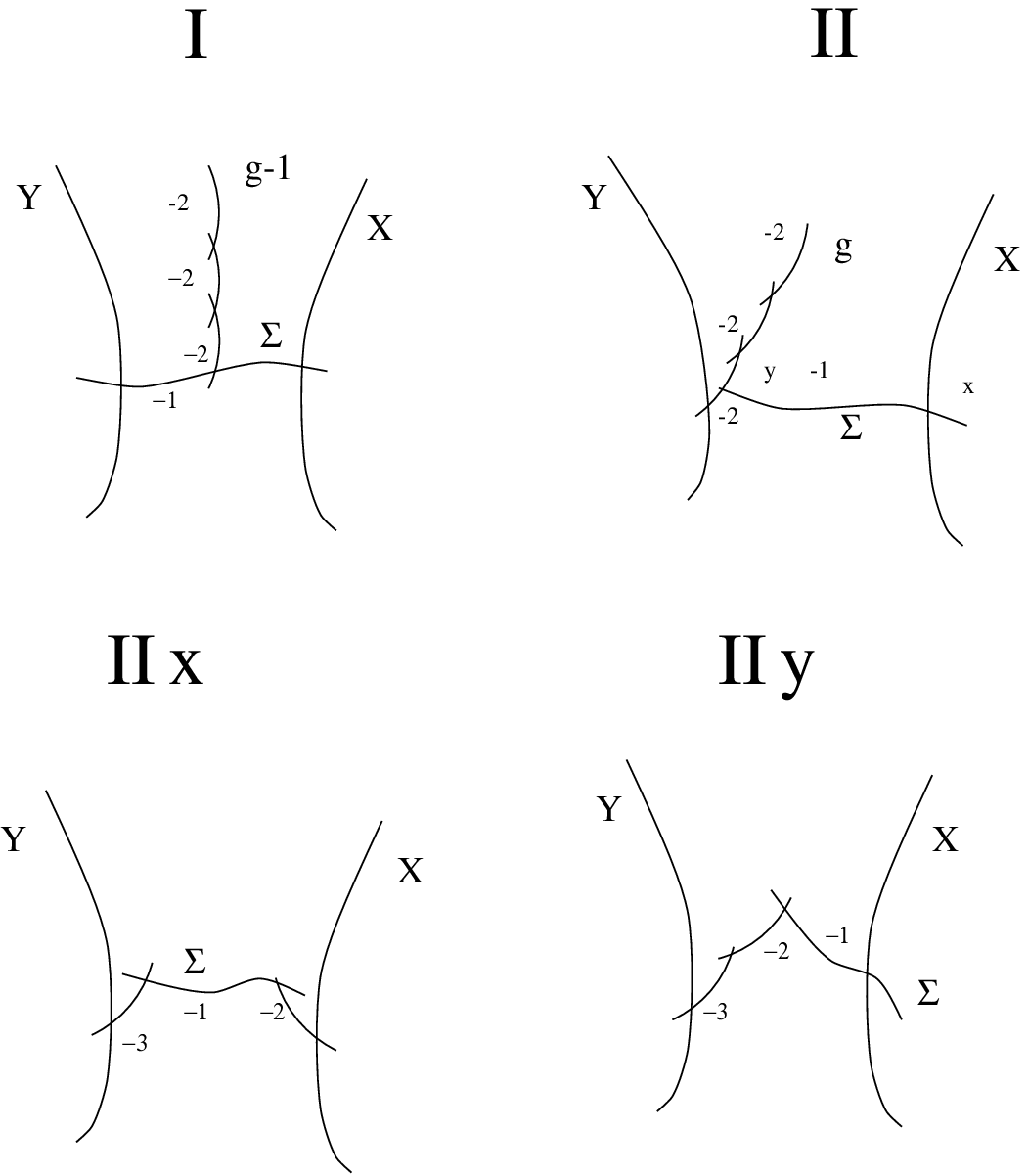}
\smallskip
\centerline{Figure 5}
\bigskip

\remark{11.1.0 Remark} The above classification of configurations (ignoring the values $c$ 
and $d$) requires only the assumption that $K_T+ X + Y$ is lt (that is it holds without assuming
either $K + c X + dY$ is flush, or $\pi$-trivial), as this is all we have used to this point.
Using the additional assumptions we can further restrict the possibilities:
\endremark

\proclaim{11.1.1 Lemma} Suppose $D$ has a node of order $g$. 

One of the following holds:
\roster
\item[1] $T$ has type I or 0, and $c+d=1$, or
\item[2] $g \geq 2$, $T$ has type II, and $gd+(g+1)c=(g+1)$, or 
\endroster
$g=1$, $K_T + \Sigma + X + Y$ is log canonical and either
\roster
\item[3]  $T$ has type $(II,x^{r-1})$, with $r \geq 1$ (for $r=1$ this means configuration 
$II$) there is a unique singularity, an $A_r$ point, on $\Sigma$, $X$ meets $\Sigma$
at a smooth point, $Y$ meets $\Sigma$ at the $A_r$ point, which has type $(\ul{2},\dots,2')$, 
and $c + \frac{d}{r+1} =1$, or
\item[4] $\Sigma$ meets $X$ and $Y$ each at a singular point of $T$, and those are the only 
two singularities along $\Sigma$.
\endroster
\newline
Moreover
\roster
\item[5]
$$
g(c+d-1)<\min (c,d) < \frac{g}{2g-1}.
$$
\endroster
\endproclaim
\demo{Proof} Let $E$ be the $-1$-curve of configuration I. One computes 
$e(E,K+ cX + dY)= g(c+d-1)$. (5) then follows from the flush condition. 

The equalities in the lemma are an easy computation.

Suppose $g \geq 2$. In configurations $(II,y^s)$, $Y$ has spectral value at least two, 
and so $c, d \geq 2/3$, by (8.0.7.2), which contradicts (5).

The rest is easy. \qed \enddemo
 
\subhead 11.2 $D$ has a unibranch double point at $q$ \endsubhead

We describe this as a {\bf $g^{th}$ order cusp.} In this case $Y$ is empty and 
$W=X$. We again consider $h$. The first series of blow ups must remove the 
singularity of $D$ (since $K_T + X$ is log terminal).
We call this {\bf configuration I.} 
Here $\Sigma$ has one singularity, 
$(2',\dots,2)=A_{g-1}$, and $X$ is simply tangent to $\Sigma$ at a smooth 
point. 

The next blow up must be at $X \cap \Sigma$. Call this {\bf configuration II.}
Here 
there is an $A_g$ singularity of type $(\underline 2', 2, \dots 2)$. 

The next blow up is again along $X$. Call this {\bf configuration III.} 
Then $\Sigma$ has two 
singularities $(A_{g-1},3')$ and an $A_1$ point, 
and $X$ meets $\Sigma$ at a smooth point. Let $u$, $v$, 
$w$ be the points where $\Sigma$ meets the $-3$-curve, $X$, and the $-2$-curve.
The next blow up is at one of these points, and we call the resulting 
configurations $u$, $v$ and $w$. 

In $u$ there are two singular point along $\Sigma$, 
$(A_{g-1},4')$ and $(\ul{2'},2)$.

In $v$ there is a single singular point $(A_{g-1},3,2',2)$ along
$\Sigma$, and $\Sigma \cap X$ is a smooth point of $T$.

In $w$ there are two singular point along $\Sigma$,
$(A_{g-1},3,\ul{2'})$ and $(3')$.

Beyond these there are two possible points to blow up. We indicate by a string of 
$ns$ and $fs$ whether we blow up at the point nearest, or farthest 
from $X$ (that is at a point of $\t{T}$ lying over $X \cap \Sigma \subset T$,
or a point of $\Sigma$ not lying over $X \subset T$).

Figure 6 gives a picture of some possible  contractions:
\smallskip
\epsfysize=370pt
\epsffile{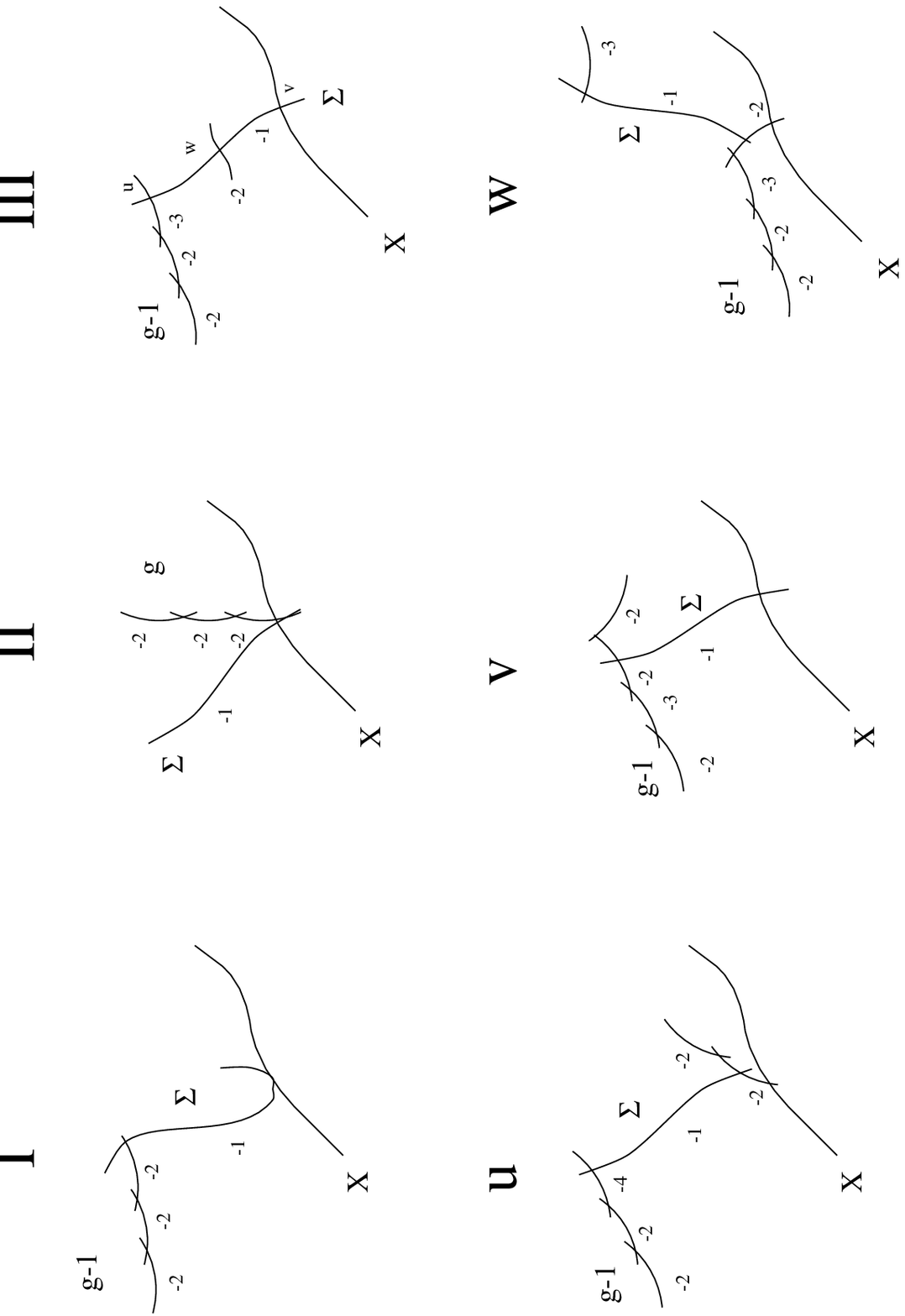}
\smallskip
\centerline{Figure 6}
\bigskip

\remark{11.2.0 Remark}The above classification requires only
the assumption that $K_T+ X$ is lt (that is  it holds without assuming either
$K + c X$ is flush, or $\pi$-trivial), as this is all we have
used to this point. Using the additional assumptions we can further restrict
the possibilities:
\endremark

\proclaim{11.2.1 Lemma} Suppose $D$ has a cusp of order $g$. Then, either
\roster
\item $T$ has type I, and $c=1/2$, or
\item $T$ has type II, and $c=(g+1)/(2g+1)$, or
\item $T$ has type III, and $c=(g+1)/(2g+1)$, or 
\item $g=1$, $T$ has type $u$ and $c=3/4$ or $g=2$ and $c=9/14$, or
\item $g=1$, $T$ has type $v$ and $c=5/7$ or $g=2$ and $c=7/11$, or 
\item $g=1$, $T$ has type $w$ and $c=7/9$, or
\item $g=1$, $T$ has type $(u;n)$ and $c=11/14$, or
\item $g=1$, $T$ has type $(v;f)$ and $c=10/13$, or
\item $g=1$, $T$ has type $(v;f^2)$ and $c=15/19$, or
\item $g=1$, $T$ has type $(v;n)$ and $c=3/4$, or
\item $g=1$, $T$ has type $(v;n^2)$ and $c=7/9$.
\endroster
Moreover 
\roster
\item[12] $c \leq \frac{g}{2g-1}$
\item If $g=1$, $c \leq 4/5$.
\endroster
\endproclaim
\demo{Proof} The stated values of $c$ are elementary calculations. Thus if 
we have type I, II or III, then we have the equalities of (1)-(3).
Let $E$ be the $-1$-curve of configuration I. One computes
$e(E,K+cX) = g(2c-1)$, so (12) follows from the flush 
condition. $(13)$ follows by computing the coefficient of
the $-1$-curve of configuration $III$. 

 Suppose we have type $u$. If $g\geq 3$, then we have a singularity at least $(4,2,2)$ 
which has coefficient $3/5$, which contradicts (12). Further blow ups are of form 
$(u;n^r)$ (or $K + X$ is not lt). If $g\geq 2$, then we have a singularity at least
$(2,4,2)$ which has coefficient $2/3$. If $g=1$ and $r \geq 2$ we get a singularity of 
spectral value four, impossible, by (13), and (8.0.7.2).

Suppose we have type $w$. If $g\geq 2$, then $E$ has a singularity of spectral
value at least $2$, impossible by (12). Further
blow ups are of the form $(w;n^r)$. Then $X$ has 
a singularity of spectral value at least four, violating (13).

Suppose we have type $v$. If $g \geq 3$ then we have a singularity at least
$(2,2,3,2,2)$ which has coefficient $3/5$, contradicting (12).
For $(v;n)$ or $(v;f)$ or beyond, if $g \geq 2$ then there is a point of 
coefficient at least $2/3$, contradicting (12).
So $g=1$. $c(v;f) = 10/13$. $c(v;f^2)=15/19$. In $(v;f,n)$, $(v;f^2,n)$,
$(v;f^3)$, $(v;n,f)$, $(v;n^2,f)$, or $v(n^3)$, $f$ has coefficient at 
least $4/5$, so these and further blow ups are not allowed. There are no 
remaining cases to consider. \qed\enddemo

\subhead 11.3 $D$ has multiplicity three \endsubhead
Our only application will be to the first hunt step when $A$ has multiplicity three.
(We do not need to consider higher multiplicities in the hunt, by (8.0.10)). 
A priori $A$ might have three analytic branches at $q$. Hence 
{\it We broaden the set up of (11.0.1) and allow $W$ to have
a third component, $Z$, and consider the log divisor $K_T+ cX + dY + f Z$, which we
assume is $\pi$-trivial and flush as before. } In fact we will show shortly that 
$Z$ is empty. 

Let $e$ be the smallest non-zero coefficient in $cX + dY + fZ$. $e < 1/2$ by 
(8.3.7). Thus by (8.0.7.2) and (8.0.8), $W$ is contained in the Du Val locus of $T$.

\proclaim{11.3.0 Lemma} Hunt notation as in (8.2.10). Suppose $S_0$ does not have a 
tiger, and $A \subset S_1$ has multiplicity at least three. Then $1/3 < a < 1/2$. 
Furthermore, if $E_1 \subset T_1$ contains
an $A_k$ point, then $x_0 = (3,A_k)$, $E_1$ is the $-3$-curve, and 
$e_0 = (k+ 1)/(2k + 3)$. \endproclaim
\demo{Proof} $S_0$ is not Du Val, for example by (8.0.9). $e_0 < a < 1/2$ by (8.3.7). Now the 
result follows from the list of singularities of small coefficient, (10.1). \qed\enddemo

In view of (11.3.0) we make the additional assumption:

{\it To simplify the classification we will assume that $e > 1/3$ and furthermore
if $W$ contains an $A_k$ point, then $e>(k+1)/(2k+3)$.} 

We again study the composition of blow ups $h:\map \t{T}.\t{S}.$. We 
write $\Sigma^i$ to indicate the $-1$-curve of the $i^{th}$ blow up.

\proclaim{11.3.1 Lemma} $Z$ is empty, that is $D$ has at most two branches.
\endproclaim 
\demo{Proof} As we have noted, $W$ is contained in the Du Val locus. Suppose $D$ has
three branches. $\Sigma$ can contain at most one Du Val point (or it will not be 
contractible) thus $\Sigma$ meets $W$ at two smooth points, and an $A_k$ point. But then 
$$
1< 2e + e/(k+1) \leq \Sigma \cdot (cX + dY + fZ) = - \Sigma \cdot K_T \leq 1
$$
a contradiction. \qed \enddemo

{\it In view of (11.3.1), we return to the notation of (11.0.1).}

\proclaim{11.3.2 Lemma} If $D$ has two branches at $q$, then one branch has a double 
point, a simple cusp, and the other is smooth. If $X$ is the branch with the cusp, then
$\Sigma$ meets $X$ normally at one smooth point, $\Sigma$ contains two singularities 
$(2)$ and $(3)$, and $Y$ meets the $-2$-curve, and is disjoint from $\Sigma$ on $\t{T}$. 
In terms of the classification of (11.2), this is configuration III for the genus one 
cusp $X$.
\endproclaim
\demo{Proof} Suppose $D$ has two branches. Since $D$ has multiplicity three, one has a 
double point, say $X$, necessarily a $g^{th}$ order cusp, and the other, $Y$, is smooth. 
The configuration with respect to the cusp $X$ belongs to the classification at the start 
of (11.2), by remark (11.2.0).

If $X$ and $Y$ both meet $\Sigma$ at smooth points, then from $e > 1/3$ it follows that each
meets $\Sigma$ normally. But then their images at $q$ are analytically isomorphic, a 
contradiction. So at least one contains an $A_k$ point. If each meets $\Sigma$ on $\t{T}$ 
then 
$$
(cX + d Y) \cdot \Sigma \geq e + e + e/(k+1) > 1
$$
a contradiction. 

Now consider the configuration for the cusp $X$.

If the configuration is I, $Y$ contains the $A_k$ point and we have 
$$
1 \geq 2c + d/(k+1) \geq e(2k+3/k+1) > 1
$$
a contradiction. If the configuration is II, then (since $X$ and $Y$ are disjoint on $T$), 
$Y$ meets $\Sigma$ normally at a smooth point. We have
$$
1 = d + c + c/(k+1) >1
$$
a contradiction. 

Since $K_T + X + Y$ is plt, $X + Y$ only contains Du Val points, and $Y$ is smooth at 
$q$ it follows that the configuration is III, $Y$ meets the $-2$-curve and is disjoint 
from $\Sigma$ on $\t{T}$. Thus $c + \frac d2 + \frac{g}{2g + 1} =1$, and so from the flush
condition, $g=1$. This is the configuration described. \qed \enddemo

\proclaim{11.3.3 Lemma} Suppose $Y$ is empty. Then either
\roster
\item $\Sigma$ has singularities $(3,2'),(3)$ and meets $X$ normally at a smooth point, or
\item $\Sigma$ has singularities $(3),(2)$ and on the minimal desingularisation
$X$ meets $\Sigma$ normally at the intersection of $\Sigma$ and the  $-2$-curve. 
\endroster
\endproclaim
\demo{Proof} Consider the multiplicity, $m$, of the strict transform of $D$ on the first 
blow up of $q$. 
 
Suppose first $m$ is three. Then by the flush condition, by considering the coefficient of 
$\Sigma$, $e <  2/5$. Thus, by our simplifying assumptions, $X$ lies in the smooth locus. 
Also by (10.1) the singularities of $T$ are either Du Val, or type $(3)$. Further, as 
$e > 1/3$, there can be at most one type $(3)$ singularity, and $\Sigma \cdot X \leq 2$. 
If $\Sigma$ is tangent to $X$, then $T$ is Du Val and thus by (3.3), $T$ has one Du Val 
singularity along $\Sigma$. But then $D$ has a double point, as this is configuration $I$ 
for the cusp, contradiction. So $X$ meets $\Sigma$ normally. $T$ cannot be Du Val
or $D$ is smooth. So there is a unique non Du Val singularity along $\Sigma$, a $(3)$ point. 
By (11.3.4) the singularities along $\Sigma$ are $(3)$ and $(2)$. This is configuration III 
for a genus one cusp, so $D$ has a double point, a contradiction.

If $m$ is one, then $X$ is triply tangent to $\Sigma^1$ and we must blow up at least twice 
more along $D$. At this point $T$ has a singularity of type $(\ul{3'})$, and $e =1/2$, a 
contradiction. By (11.1), further blow ups will only increase $e$.

 Thus $m=2$. As $D\cdot \Sigma^1=3$, $X$ must have an ordinary cusp on $S^1$. The next blow 
up must be along $X$. This configuration is not allowed, as $c=2/5$, and $k=1$. The two 
subsequent blow ups must be along $X$ and give (2) and (1).  On any further blow up there 
is a singularity of coefficient at least $1/2$, which is not allowed as $e < 1/2$.\qed\enddemo

\proclaim{11.3.4 Lemma} If $\Sigma$ contains a unique non Du Val singularity, a $(3)$ point, 
and $\pi(\Sigma)$ is smooth, then there is exactly one further singularity, an $A_1$ point, 
along $\Sigma$, that is $T$ is the surface of configuration III for a genus one cusp. 
\endproclaim
\demo{Proof} We consider the sequence of blow ups $h:\map \t{T}.S.$. Suppose the first 
$r\geq 2$ blow ups are over smooth points of $T$. In the resulting configuration the 
singularity is $(2',\dots,2) = A_{r-1}$. The next blow up (by the choice of $r$) is at the 
intersection of $\Sigma$ and the marked $-2$-curve. The resulting configuration has two 
singularities, $(3',A_{r-2})$ and $(2')$. On any further blow up there is a non Du Val
singularity, not of type $(3)$. Thus we have this configuration with $r =2$. \qed \enddemo

\subhead 11.4 $K_S + D$ log terminal at $q$  \endsubhead

Here $Y = \emptyset$, and $D = X$. We make an additional simplifying assumption:
\noindent
{\it We assume $\t{D}^2 > \t{W}^2$, or equivalently, that the first blow up
in $h:\map \t{T}.\t{S}.$ is at the point of $\t{D} \subset \t{S}$ lying
over $q$.}   
Note this holds for $W = E_i$ in
a step of the hunt, since (in the hunt notation) $\t{D_i}$ is $K$-negative,
while $\tilde{E}_i$ is $K$ non-negative.

There are infinitely many possible configurations, we now indicate how they can be generated. 
\noindent
{\it We will not assume in the analysis below that $K_T + cX$ is either
flush or $\pi$-trivial. We will only use our assumption on
$\t{D}^2$, and that $K_T + W$ and $K_S + D$ are log terminal}

Consider the possibilities for $h:\map \t{T}.\t{S}.$.
By assumption, the first blow up of $h:\map\t{T}.\t{S}.$ is along $X=D$.
Subsequent blow ups are of three
sorts. Either we blow up at the point of $\Sigma \subset \t{T}$ lying over
$\Sigma \cap X$, which we call the {\bf nearest} point,  or at a point,
lying over a singular point of $T$, different from the near point, which we
call the {\bf farthest} point (there will be at most one such point), or
at a point on $\Sigma$, not along $X$, 
lying over a smooth point of $T$. We call this
last case an {\bf interior blow up}. We indicate a series of 
non-interior blow ups by a sequence, $(n;f;n^2)$ for example, where $n$ 
indicates a blow up at the nearest point to $X$, and $f$ indicates a 
blow up at the farthest point. We count the first blow up of $h$ as
a blow up at the nearest point. 
We will indicate an interior blow up by $i$. After
an interior blow up, there is only one singular point of $T$ on $\Sigma$. A 
further interior blow up is impossible, for then $K_T+W$ is not lt, and so the 
possible sequence of blow ups is unique. We write $n^k$ if there are $k$ 
further blow ups. Note $K_T+X +\Sigma$ is lc iff there is no interior blow up,
and that $X \cap \Sigma$ is a smooth point of $T$ iff all the blow ups
are at the near point, that is iff the configuration is $(n^r)$ for some 
$r \geq 0$.  

For example if the singularity at $q$ is $(\ul{2},3)$, and we make the blow up
sequence $(n^2;f;i;n)$, the singularities on $T$ along $\Sigma$ are 
$(\ul{2},3',3,3,3)$ and $(2')$.

 Here is a picture:
\smallskip
\epsfysize=370pt
\epsffile{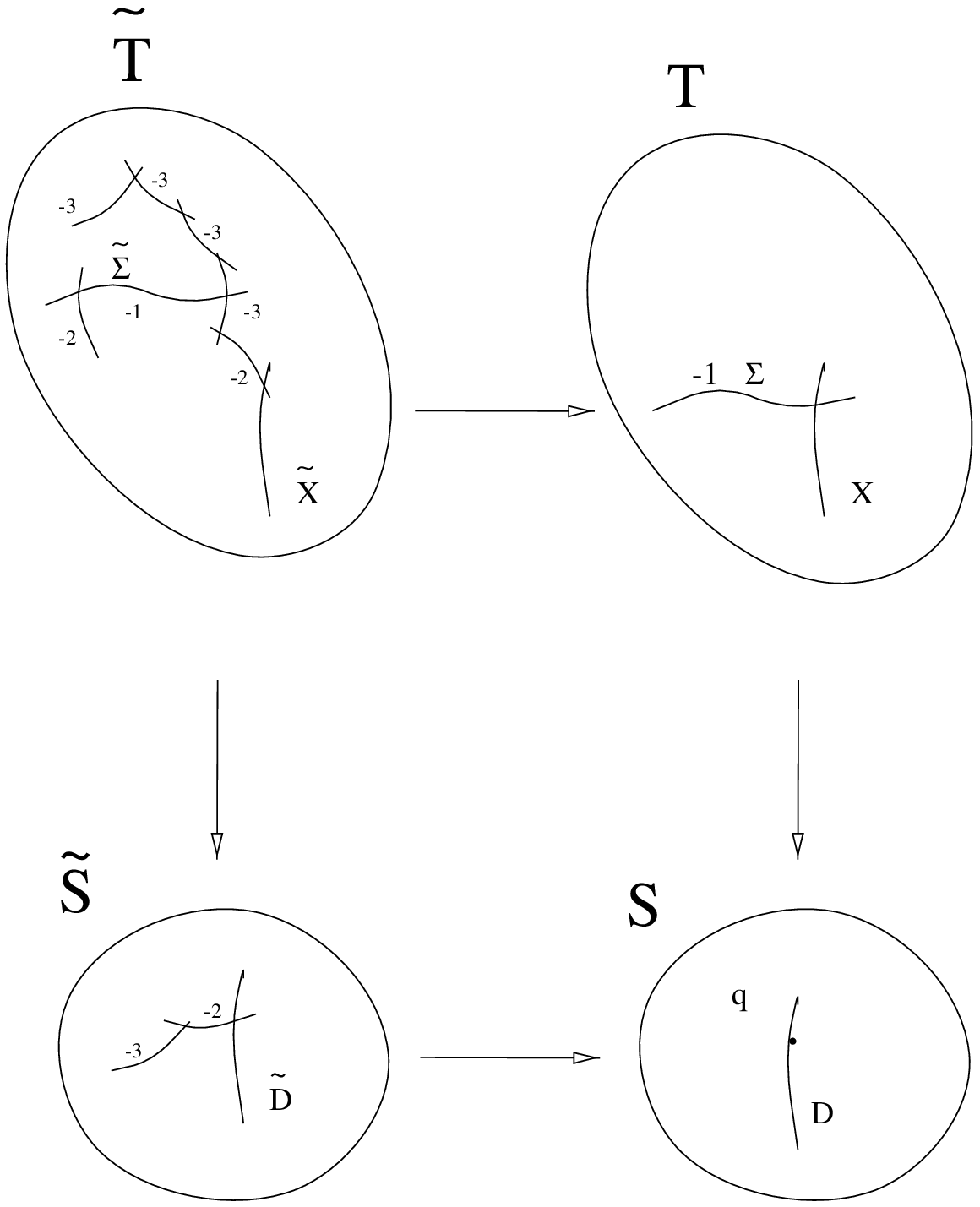}
\smallskip
\centerline{Figure 7}
\bigskip

\subhead 11.5 $\pr 1.$-fibrations \endsubhead
\newline
{\it 11.5.1 Notation: Let $(p,C)$ be a smooth curve germ, and let
$\pi:\map T.C.$ be a (proper) $\pr 1.$-fibration of relative 
Picard number one, with $T$ a normal surface.  
Let $\t{\pi}$ be the composition
$$
\t{\pi}: \t{T} \rightarrow T @>{\pi}>> C
$$
(where, as throughout the paper, $\map \t{T}.T.$ is the minimal desingularisation).
Let ${\Cal F} \subset \t{T}$, ${\Cal G} \subset T$ be the scheme-theoretic fibres
(over $p$), with reductions $F,G$. }

\proclaim{11.5.2 Lemma} The following are equivalent:
\roster
\item $\pi$ is smooth.
\item $T$ is smooth.
\item ${\Cal G}$ is generically reduced.
\endroster
\endproclaim
\demo{Proof} The equivalence of (1-2) is well known, and obviously (1) implies (3). Assume
(3). ${\Cal G} \subset T$ is Cartier, and hence Cohen-Macaulay. Thus it has no embedded
points, and so is integral. It follows that 
$$
h^1(\ring G.) = 1 - \chi(\ring G.) = 1 - \chi(\ring {\pr 1.}.) = 0
$$
thus $G = \pr 1.$. Since ${\Cal G}$ is Cartier, (2) follows. \qed \enddemo

{\it From now on we assume $T$ is singular.} In view of (11.5.2) we will sometimes
abuse notation and call $G$ a {\it multiple} or {\it singular} fibre of $T$ (even though
$G$ is by assumption reduced, and is isomorphic to $\pr 1.$ ).

\proclaim{11.5.3 Lemma} $\t{G} \subset F$ is the unique $-1$-curve in the fibre. 
\endproclaim
\demo{Proof} This holds since the exceptional locus of $\map \t{T}.T.$ is $K_{\t{T}}$ 
nef. \qed\enddemo

\noindent
{\it 11.5.4 Further Notation:
Let 
$$
\CD
\t{T} @>{h}>> W \\
@V{\t{p}}VV  @V{f}VV  \\
C @=     C
\endCD
$$
be a relative minimal model over $p$, so $f$  is a smooth $\pr 1.$-bundle, and $h$ is a 
composition of smooth blow ups. Let $H$ be the fibre of $f$ over $p$. 
In the case 
where the dual graph of the inverse image of $F$ is a chain, we will use 
the following notation to describe the scheme-theoretic fibre ${\Cal F}$: The sequence 
$k(-a)+l(-b)+\dots +m(-c)$, is taken to mean that the first rational curve 
in the chain has self-intersection $-a$ and has multiplicity $k$ inside 
the fibre, the second has self-intersection $-b$ and multiplicity $l$ 
and so on. }

Consider the possibilities for $h$. The first blow up is at some point along
the fibre,  after which the fibre is $(-1)+(-1)$. One of the components is
(the strict transform of) $H$, the other, $V$, is $h$-exceptional. 

There must be a further blow up. Necessarily, by (11.5.3), 
at the intersection point of the two $-1$-curves, after which the fibre is 
$(-2)+2(-1)+(-2)$. All subsequent blow ups are at some point along the 
unique $-1$-curve, by (11.5.3). At each stage one either blows up
at a point where the $-1$-curve meets another fibral component, or at a point of
the $-1$-curve which is not on any other component. We call the latter an
{\bf interior blow up}. Note in particular that whatever the
sequence of blow ups there are exactly two reduced components in the fibre,
the strict transforms of $V$ and $H$. We collect these observations as:

\proclaim{11.5.5 Lemma} Notation as in (11.5.1), (11.5.4). ${\Cal F}$ has exactly two reduced
components. The collection of multiple components is contractible. Let $h':\map T.W'.$ 
contract these components. $W'$ is smooth, and $f':\map W'.C.$ has fibre $(-1)+(-1)$. 
$h':\map T. W'.$ blows down (at each stage) the unique $-1$-curve. $f'$ has two relative 
minimal models $g_1$, $g_2$ given by contracting one, or the other $-1$-curve. $\t{\pi}$ has 
exactly two relative minimal models $g_1 \circ h'$ and $g_2 \circ h'$, and the relative 
minimal model is determined by choosing which reduced component of $F$ not to contract. This 
component is the strict transform of $H$, the other component is the strict transform of $V$.
 $T$ has rational singularities, and at most two singular points along $G$. 
\endproclaim
\demo{Proof} The last two statements are the only new observations, and follow easily from the 
description of the morphism $h'$. \qed\enddemo

Possible $\pi$ are classified by the choices of the blow ups in $h$. Using the classification of 
lt singularities, one can further classify possibilities when $T$ is log terminal. One easily 
checks in particular the following:
\proclaim{11.5.6  Lemma} Notation as in (11.5.1). $K_T + G$ is log terminal iff
there are no interior blow ups.
\endproclaim

Note after the first 
interior blow up, $F$ consists of the $-1$-curve, and a chain of $K_T$-positive
curves. We call this chain (and its strict transform after further blow ups), 
the {\bf principal} chain of the fibre. Its two ends are the two reduced 
components of the fibre. 

\subhead 11.5.7 Fibrations with a section \endsubhead
\newline
{\it Additional Notation: Let $E \subset T$ be a reduced curve germ, finite over $C$, 
such that $K_T+E$ is lt, and $\t{E}\cdot {\Cal F}= 1$.} 

Necessarily $E$ is a section of $\pi$, $\t{E}$ meets a unique components of ${\Cal F}$, a 
reduced component, and $\t{E}$ has normal contact. By (11.5.5) there is a unique choice of 
$h$ which is an isomorphism in a neighbourhood of $\t{E}$. There can be at most one interior 
blow up (otherwise $E$ contains a non-chain singularity). Let $x=E\cap G \in T$. If $K_T+G$ is 
lt then $F$ (and hence $h$) is uniquely determined by the marked singularity $(x,E)$, and 
otherwise it is determined by the marked singularity $(x,E+G)$. We indicate the singularities
as in (11.0.1), with an underline (resp. prime) indicating normal contact with $\t{E}$
(resp. with  $\t{G}$).

\remark{Remark} Of course one can classify global $\pr 1.$-fibrations $\pi:\map T.C.$, 
for a smooth proper curve $C$, by using the above local classification, so long as one
has a classification of smooth $\pr 1.$-bundles $\map W.C.$. For example one
can carry this out for $C = \pr 1.$. \endremark

\remark{11.5.8 Example} Suppose $C = \pr 1.$, $E \subset T$ is a $-k$-curve, and a section, with
$K_T + E$ log terminal, and $E$ contains a single singularity $x=(\ul{2},2,2)$. Let $G$ be the 
(reduction of the) fibre through $x$. 

Then the above analysis shows: $\pi$ is smooth away from $G$. There is a unique relative
minimal model $h:\map \t{T}.W.$ which is an isomorphism in a neighbourhood
of $\t{E} \subset \t{T}$. $W = \Hz k.$ and (the
image of) $E \subset W$ is the negative section. If $K_T + G$ is log terminal, then $F$ is 
$(-\ul{2})+2(-2)+3(-2)+4(-1)+(-4)$.  If $K_T+G$ is not lt, then $x$ is the 
principal chain, the only singularity along $G$, and the marked fibre is
$(-\ul 2)+2(-2')+(-2)$. \endremark

 Here is a picture:
\smallskip
\epsfysize=350pt
\epsffile{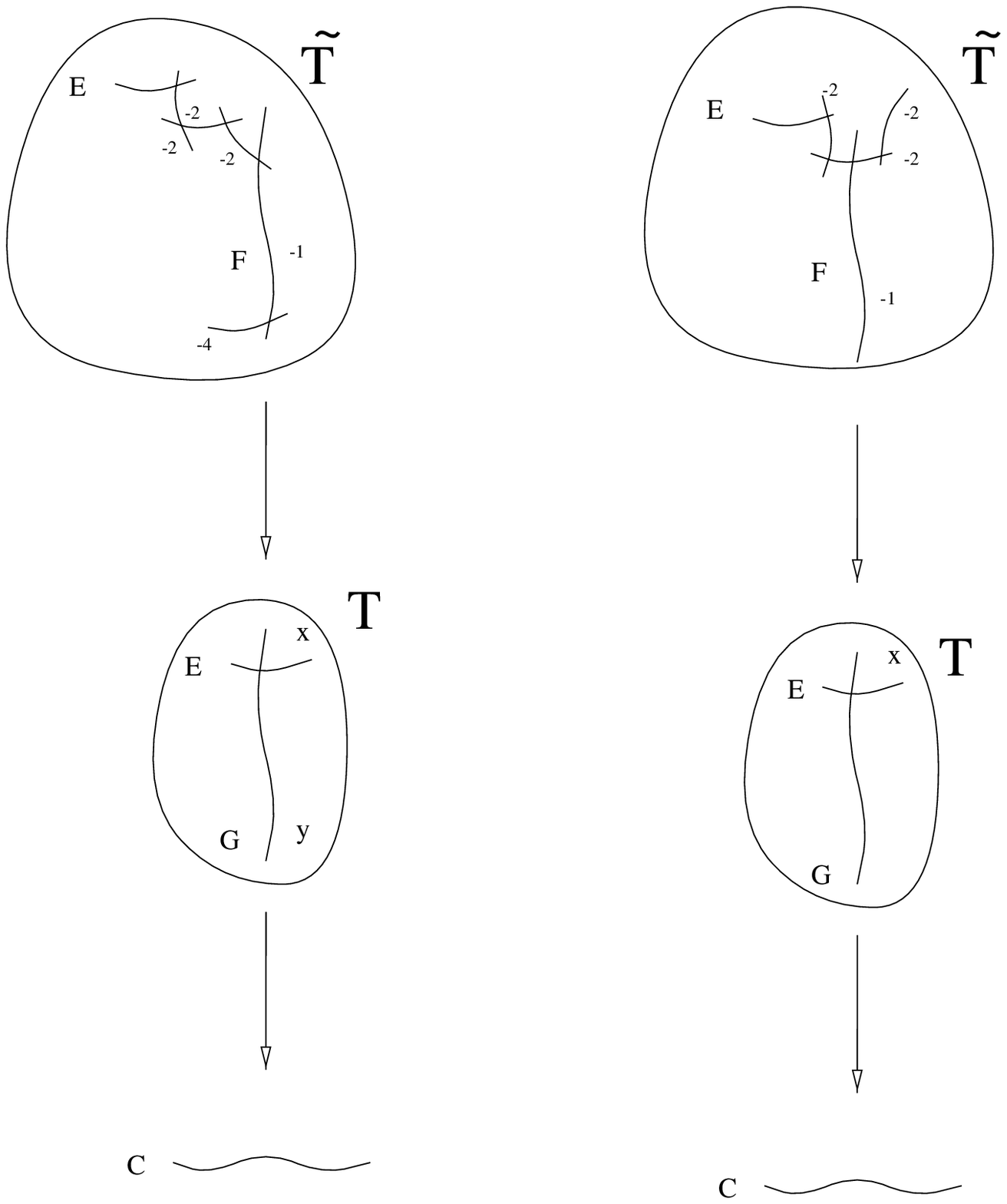}
\smallskip
\centerline{Figure 8}
\bigskip

\proclaim{11.5.9 Lemma} Notation as in (11.5.1). Assume $T$ is log terminal, 
$G$ is a multiple fibre of $\pi$, of
multiplicity $m$, and $G$ contains a cyclic singularity, either Du Val or 
almost Du Val. If $e(T) < 2/3$, then $G$ is one of the following: \newline
If $K_T+G$ is not lt at any singular point:
\roster 
\item $(2, 2', 2)$, $m=2$ or 
\item $(3,2',2,2)$, $m=3$.
\endroster
If $K_T+G$ is log terminal at one singular point, but not lt:
\roster 
\item[3] $(2';z)$, $m=4$. $z$ is a non-chain singularity, with centre $-2$ 
and branches $(2)$, $(2)$, and $(2,\dots,3')$ (or $(3')$).
\item $(2,3',2;2')$ $m=4$. 
\endroster
If $K_T+G$ is log terminal:-
\roster
\item[5] $(A_k; (k+1)')$, $k\leq 4$, $m=k+1$. The fibre is
$$
-(k+1) + [k+1](-1) + k(-2) + [k-1](-2) + \dots + (-2).
$$
\item $(2,3';2',3)$, $m=5$. The fibre is
$(-2) + 2(-3) + 5(-1) + 3(-2) + (-3)$.
\item $(3,2,2';4',2)$ $m=7$. The fibre is
$(-3) + 3(-2) + 5(-2) + 7(-1) + 2(-4) + (-2)$.
\item $(4,2';3',2,2)$, $m=7$. The fibre is
$(-4) + 4(-2) + 3(-1) + 7(-3) + 2(-2) + (-2)$.
\endroster
\endproclaim
\demo{Proof}
One checks that each of (1-8) is an allowed configuration. 
 Suppose $K_T+ G$ is lt. One checks that for each configuration
(5-8), and either of the two possible non-interior blow ups, the resulting 
configuration is either one of (5-8), or disallowed by the hypotheses. Since 
the fibre $(-2)+2(-1)+(-2)$ is (5), $k=1$, we have considered all 
possibilities.

Now suppose $K + G$ is not lt. Thus there is some interior
blow up. Consider the sequence of blow ups. We begin with a sequence
of non-interior blow ups, giving one of the fibres (5-8). Call this
configuration $J$. Then we make a string $i^r$ of $r \geq 1$
interior blow ups.

Suppose there are two singular points along $G$. Then we make in
addition at least one non-interior blow up.
$e(T) \geq 2/3$ unless $J$ is given by (5) with $k=1$. $r \geq 2$
gives (3), and $r=1$ gives (4). Any further blow up has $e(T) \geq 2/3$.

Suppose there is only one singular point along $G$, cyclic
and Du Val or almost Du Val by assumption. Then $r=1$ (otherwise
the singularity cannot be cyclic). $J$ given by (5) with
$k=1,2$ are the only possibilities in which the singularity
is Du Val or almost Du Val. This gives (1) and (2). No further
blow up is possible (the next blow up has to be non-interior),
as the singularity will not be Du Val or almost Du Val. \qed\enddemo

\proclaim{11.5.10 Lemma} Suppose $D\subset \Hz n.$ is an irreducible curve of degree 
$d\geq 2$, arithmetic genus $g$. If $g\leq 1$, then either 
\roster
\item $g=0$, $n=0$ and $D^2=2d$, or
\item $g=0$, $d=2$, $n=1$ and $D^2=4$, or
\item $g=1$, $d=3$ and $n=1$, or
\item $g=1$, $d=2$ and $n\leq 2$.
\endroster
\endproclaim
\demo{Proof} $(K_S+D)\cdot D\leq 0$. But $(K_S+D)\cdot F\geq 0$, where
$F$ is a smooth fibre, and $(K_S+D)\cdot E_{\infty}> 0$, unless $n\leq 2$,
where $E_{\infty}$ generates the other extremal ray. Write $D$ numerically
as $dE + r F$. One computes $r$ using $(K + D) \cdot D = 2g -2$, and
then checks the result. \qed\enddemo

\proclaim{11.5.11 Lemma} Suppose $\pi :\map T.{\pr 1.}.$ is a log Fano fibration of relative 
Picard number one, general fibre $\pr 1.$, $C$ is a section, $D$ is rational, an irreducible 
double section, and $C$ and $D$ meet at only smooth points of $T$. Let $F_1$ and $F_2$ be the 
fibres, where $\pi |_D$ has a ramification point. Let $W$ be the unique minimal model 
$h:\map \tilde T.{W=\Hz n.}.$, of (11.5.5), (11.5.7), such that $h$ is an isomorphism in a 
neighbourhood of $\t{C}$. If $K_T+D$ is lt, then 
\roster 
\item Any multiple fibre is either $F_1$ or $F_2$ (though either $F_i$ might be
smooth).
\item Neither $C$ nor $D$ contains more than  two singularities of $T$.
\item $h|_D$ is an isomorphism.
\endroster
\endproclaim
\demo{Proof} Let $G$ be a multiple fibre. Suppose $G$ is neither $F_1$ nor $F_2$. Then each 
point of $G\cap C$ and $G\cap D$ is a singular point, contradicting (11.5.5). Hence (1). 

 (2) follows easily, using (11.5.5).

 Suppose (3) does not hold. Then one of the ramification points of $D$ is smooth on $T$, and 
has a cusp on $W$. Considering the first two blow ups of $h$, we see that $D$ does not meet 
the exceptional locus of $\map \tilde T.T.$ transversally, a contradiction (since we assume 
$K_T + D$ is lt). \qed\enddemo

\proclaim{11.5.12 Lemma} Notation as in (11.5.1), (11.5.4). Suppose $\pi: \map T. {C}.$ is a 
$\pr 1.$-fibration, and $G \subset T$ is an irreducible fibre contained in the Du Val locus 
of $T$. One of the following holds:
\roster
\item $T$ is smooth along $G$.
\item There are exactly two singularities, $A_1$ points, along $G$. $K_T + G$
is lt. The fibre is $(-2) + 2(-1) + (-2)$ (notation as in (11.5.4)).
\item There is a unique singularity, an $A_3$ point along $G$, 
$(2,2',2)$. 
\item There is a unique singularity, a $D_n$ point, along $F$, its branches
are $(2)$, $(2)$ and \newline 
$(\ul{2},\dots,2')$, where the central curve meets the
underlined curve, and $G$ meets the primed curve.
\endroster
\endproclaim
\demo{Proof} The only possibilities for $h$ which yield Du Val singularities
are $(-2) + 2(-1) + (-2)$, or this configuration followed by a sequence of 
interior
blow ups. \qed \enddemo

We will use the following ad-hoc result in \S 14.

\proclaim{11.5.13 Lemma} Notation as in (11.5.1). If $G$ is a multiple fibre of multiplicity 
$m=3$, and the coefficient $e(T) < 2/3$ then $G$ is one of the fibres in the classification 
(11.5.9).
\endproclaim
\demo{Proof} Consider the possibilities for $h$ as below (11.5.4). It is easy to see 
that the only way to get a fibre of multiplicity three is to begin with 
$(-3)+3(-1)+2(-2)+(-2)$ and make a series of $r \geq 0$ interior blow ups. 
After two interior blow ups there is a non-chain point of coefficient $2/3$. 
For $r\leq 1$, the fibres are in (11.5.9). \qed \enddemo

\heading \S 12 The linear system $|K_S+A|$ \endheading
\subhead The following notation is fixed throughout the section \endsubhead
{\it Let $S$ be a rank one 
log del Pezzo surface. Let  $A \subset S$ a reduced and irreducible divisor, such that 
$-(K_S+(1/2)A)$ is ample and $K_S+A$ is lt at singular points of $S$. Note in 
particular $A$ has planar singularities, and is thus Gorenstein. Assume the 
arithmetic genus $g$ of $A$ is at least one. $r: \map \t{S}.S.$ denotes the minimal 
desingularisation of $S$, and tildes denote strict transforms. Let $D = F + M$ of (12.0) 
below. Let $\ol{D}$ be the reduction of $D$. Let $\map Y.S.$ extract exactly the 
exceptional divisors (of the minimal desingularisation) adjacent to $A$. }

Note the assumptions imply $r|_{\t{A}}$ is an isomorphism. 

\proclaim{12.0 Definition-Lemma} 
$h^1(\omega_{\t{S}})=h^0(\omega_{\t{S}})=0$ and 
$H^0(\omega_{\t{S}}(\t{A})) \simeq H^0(\omega_{A})$.
There are
integral, effective, Weil divisors $F,M \subset S$, $E \subset \t{S}$ with
the following properties:
\roster
\item $\t F+\t M+E \in | K_{\t S}+\t A |$, $F + M \in |K_S + A|$.
\item The support of $E$ is $\map \t{S}.S.$ exceptional.
\item $|\t{F}|$ has no fixed components. No component of $\t{M}$ moves.
\endroster
\endproclaim
\demo{Proof} 
The first equalities hold since $\t{S}$ is rational.
The second follow from the first and 
the exact sequence
$$
\ses {\omega_{\t{S}}}. {\omega_{\t{S}}(\t{A})}. {\omega_{A}}. .
$$
By assumption, $|\omega_A|$ is non-empty. Thus $|K_{\t{S}} + \t{A}|$ 
is non-empty. We can write a general member as 
$\t{F} + \t{M} + E$, where
$\t{F}$ is the moving part, $\t{M} + E$ is the fixed part, and $E$
is the $\map \t{S}.S.$ exceptional part. Now let
$F,M$ be the images of $\t{F},\t{M}$. Note
$\t{F}$ and $\t{M}$ are indeed the strict transforms of $F$ and $M$, by construction.
\qed \enddemo

\proclaim{12.1 Lemma} Notation as in (12.0).
\roster
\item $K_S+D$ is negative (that is anti-ample) and $\ol{D}$ has arithmetic genus zero.
\item If $D$ is reducible, then the union of the pairwise intersection of 
components of $D$ consists of a single point $b$. 
\item $F \neq \emptyset$ iff $g \geq 2$. $| \t{F}|$ is basepoint free.
\item $\ol{D} \cap \sg (S)$ contains any non Du Val singularity
of $S$, and any singularity of $S$ which lies on $A$.
\item If $\ol{D}$ is reducible, then $K_S + \ol{D}$ is plt away from $b$.
\item $\t{A}$ is disjoint from $\t{M}$ and $E$.
\item The support of $\t{M}$ is a disjoint union of $-1$-curves. 
$\t{M} \subset Y$ is contractible by a series of $K$-negative
contractions, and the number of irreducible components 
of $M$ is at most the number of singularities along $A$.
\endroster
\endproclaim
\demo{Proof} 
$K_S+D=2(K_S+ 1/2A)$ is negative. (1), (2) and (5) follow by adjunction.

$h^0(\omega_A) \geq 1$ iff $g(A) \geq 2$. (3) follows, since a moving
smooth rational curve on a surface with $H^1(\ring S.)=0$ is basepoint 
free.
 
(6) follows from (12.7).

 From $K_{\tilde S}+\t{A}-\t{M}-\t{F}=E$, we see that $E$ is nef,
and hence by (2.19) of \cite{Kollaretal}, empty, over any 
point not contained in $M+F$. (4) follows.

For (7) we use the following:

\proclaim{Claim} Let $\gamma: \map Y. Z.$ be proper birational
map, which is an isomorphism in a neighbourhood of $\t{A} \subset Y$.
Let $N$ be an effective divisor on $Y$, disjoint from $\t{A}$,
such that $\gamma$ and $\map Y.S.$ are finite on $N$. Then
$$
0 > \gamma(N) \cdot (K_Z + \gamma(\t{A})) = K_Z \cdot \gamma(N).
$$
\endproclaim
\demo{Proof} Since $-K_S$ is ample, there is an effective
${\Bbb Q}$-divisor $J$ supported on the exceptional locus of $\map Y.S.$
such that $-(K_Y+J)$ is ample. We have
$$
0 > (K_Z + \gamma(J)) \cdot \gamma(N) \geq K_Z \cdot \gamma(N)
= (K_Z + \gamma(\t{A})) \cdot \gamma(N). \qed
$$
\enddemo
Now let $\gamma: \map Y. Z.$ contract some components of
$\t{M}$ (we allow $\gamma$ to be the identity). Let $N$ be the
reduction of some component of $M$ which is not contracted.
By (6) and the claim
$$
0 > \gamma(N) \cdot (K_Z +\gamma(\t{A}))= \gamma(N) \cdot K_Z =
\gamma(N) \cdot \gamma(\t{M}+\t{F}+E)\geq \gamma(N)^2 
$$ 
Thus $\gamma(N)$ can be contracted, in fact by a $K_Z$-negative contraction.
Thus, $M$ can be contracted by a series of $K$-negative contractions.
(7) follows. \qed\enddemo

\proclaim{12.2 Lemma} Suppose $S=S_1$ of the hunt. If $S_0$ does not have a 
tiger, then every component of $M$ contains at least two singularities of $S$. 
\endproclaim
\demo{Proof} Suppose $N$ is a component of $M$ passing through a unique
singular point $z$ (it must contain a singular point, or it
would be contractible on $S$). 
As $\t{M}$ is disjoint from $\t{A}$, by (12.1.6), $z$ is a point of $A$. 
$M$ cannot meet $\Sigma_1$ on $T_1$, for otherwise $z = \pi(\Sigma_1)$,
a singular point of $A$. Thus $M$ meets $E_1$ at a singular
point of $T_1$ and $x_0$ is the only singular point of $S_0$ 
along $M$. 

Let $D$ be the unique exceptional curve of $\t S_0$ which meets $M$
(there is a unique $D$ by (4.12.2) and (12.1.1)), and let  
$g: \map T'.S_0.$ be the extraction of $D$. Then $M$ is a $-1$-curve in the 
smooth locus of $T'$, and we may contract $M$ on $T'$ to obtain $S'$. Let 
$D'$ be the image of $D$ in $S'$. 

By the classification of lt points, if $K_{S'}+D'$ is lc then it is anti-nef 
and we have a tiger. So we can assume that $x_0$ is a non-chain singularity
and one of the singular points along $E_1$ other than $z$ has index at least 
three. Thus $z$ has index at most $5$. As $M^2>0$, it follows that $z = A_4$ 
(just consider the possibilities, see for example the proof of
(10.6)), and $D$ is the primed curve $(2,2',2,2)$. 
But then $(K_S+A)\cdot M \leq -1 + 3/5 < 0$, and
$-(K_S+A)$ is ample, a contradiction.  \qed \enddemo

\proclaim{12.3} Notation as in (12.0). Suppose $g=2$. Then 
\roster
\item $\t{F}$ is the fibre of a ruling on $\t{S}$, $\t{A}$ is a double
section.
\item $\t{F}$ is disjoint from $\t{M}$ and $E$ on $\t{S}$.
\item If $p \in F$ is a singular point, then either $F$
meets the same curve over $p$ as $A$, or $F$ meets a curve of 
self-intersection at most $-3$.
\endroster
If $p \in F \setminus F \cap M$ is a singular point, then $E$ is empty over $p$ and 
\roster
\item[4] If $p \in A$ then $p$ is Du Val and $F$ meets the same
exceptional curve at $p$ as $A$.
\item If $p \not \in A$ then $p$ has weight one, and $F$ meets 
the unique $-3$-curve. 
\endroster
\endproclaim
\demo{Proof} By (12.1.3) and (12.1.6) we have

$$
\align
0 \leq \t{F}^2 \leq (\tilde F+\tilde M+E)\cdot \tilde F&=
(K_{\tilde S}+\tilde A)\cdot \tilde F \leq 
-2+\tilde A\cdot \tilde F \\
&=-2 + \t{A} \cdot (K_{\t{S}} + \t{A}) = 2(g-2).
\endalign
$$
(1-2) follow. 

For (4) and (5), (2) and the equality $K_{\t {S}}+\t{A}-\t{F}=E$ over $p$ 
imply $E$ is nef over $p$, thus empty, by (2.19) of 
\cite{Kollaretal}. (4) and (5) follow. 
(3) is similar.\qed\enddemo

\proclaim{12.4 Lemma-Definition} Notation as in (12.0). Suppose $g=1$. By (12.1.7)
$\t{M} \subset Y$ is contractible, let
$\pi:\map Y.W.$ be the contraction.
Let $G = \sum G_i$ be the exceptional locus of $f:\map Y.S.$. By the
definition of $Y$,
$G$ is the union of the exceptional divisors (of  $\t{S}$) adjacent to $A$.
Let $D_i = \pi(G_i)$. Let $\tilde{A} \subset Y$ be the strict transform of $A$.
We have 
\roster
\item The rank of $W$ is one more than 
the difference between the number of singular points
on $A$ and the number of irreducible components of $M$.
\item $\pi$ is an isomorphism in a neighbourhood of $\t{A}$.
\item $A \in |-K_W|$, $A \subset W^0$, and $W$ is Gorenstein.
\item $D_i\cdot A =1$. 
\item If $A$ contains only one singular point, then $W$ has rank one and $\tilde A^2\geq 1$. 
\endroster
\endproclaim
\demo{Proof} The relative Picard number of $f:\map Y.S.$ is the number of singularities
along $A$, and to get from $Y$ to $W$, we contract $M$. Hence (1). 

(2) is immediate from (12.1.6).

By (12.1.3) and (12.0.1), $K_Y + \t{A} = \t{M}$. (3) follows from this, and (2). 

(4) follows from (3) since $G_i \cdot \t{A} =1$, as $K_S + A$ is lt at singular points

$\tilde A^2=K_W^2$, by (2) and (3). If $W$ has rank one, $K^2_W>0$. Hence (5).  \qed\enddemo

\remark{Remark 12.4.6} Suppose $S = S_1$, $A = A_1$ of the hunt, and $g(A) =1$. The transformation
$S_1 \dasharrow W$ of (12.4) is not in general the next hunt step. However, since $W$ is 
Gorenstein, it is often a useful alternative, especially if one can argue that
$W$ has rank one. \endremark

\proclaim{12.5 Lemma} Notation as in (12.0). If $A$ has genus one and $K_S+2/3A$ is negative 
then no irreducible component of $M$ passes through more than two singular points of $S$. 
\endproclaim
\demo{Proof} Suppose not. Let $N$ be the reduction of an irreducible component of $M$ containing 
three singular points, say of index $p$, $q$ and $r$. As $K_S+(2/3)A$ is negative, $K_S+2N$ is 
negative. But by adjunction
$$
\align
(K_S+2N)\cdot N& \geq (p-1)/p+(q-1)/q+(r-1)/r-2+N^2 \\
               &\geq (p-1)/p+(q-1)/q+(r-1)/r-2+(-1+1/p+1/q+1/r)=0,
\endalign
$$
a contradiction. \qed\enddemo

\subhead 12.6 Canonical map of an integral Gorenstein curve \endsubhead

Here we prove the following, as we were unable to find it in
the literature:
 
\proclaim{12.7 Lemma} Let $A$ be a proper, integral, Gorenstein
curve. If $A \neq \pr 1.$ then $|\omega_A|$ is basepoint free.
\endproclaim

\proclaim{12.8 Lemma} Let $A$ be a proper, integral, Gorenstein
curve. There is a differential form
$F \in H^0(\omega_A)$ which does not vanish at any singular point
of $A$.
\endproclaim
\demo{Proof} Let $\pi:\map \t{A}.A.$ be its normalisation. Define the
Cartier divisor $D$ by $\omega_{\t{A}}(D) = \pi^*(\omega_A)$ (thus 
$D$ is cut out by the conductor ideal on $\t{A}$). We have a commutative
diagram with exact rows
$$
\CD
0 @>>> \pi_*(\omega_{\t{A}}) @>>> \pi_*(\omega_{\t{A}}(D)) @>>>
\pi_*(\omega_D) @>>> 0 \\
@. @VVV                         @|                        @V{g}VV   @. \\
0 @>>> \omega_A @>>> \omega_A \otimes \ring \t{A}. @>>> \omega_A \otimes
(\ring \t{A}.)/ (\ring A.) @>>> 0 .
\endCD
$$
This induces a short  exact sequence
$$
0 @>>> H^0(\omega_{\t{A}}) @>>> H^0(\omega_A)  @>>> \text{ker}(g) @>>> 0.
$$

Fix a singular point $q \in A$ and a local generator $F_q\in \omega_{A,q}$.
Fix a neighbourhood $q \in U$ missing the other singularities, and
let $V = \pi^{-1}(U)$.
The polar divisor of $F_q$, thought of as a rational form on $V$,
is exactly $D|V$, and the image of $F_q$ in 
$\omega_D|V$ is an element $G_q \in \text{ker}(g)$. Any lifting
of $G_q$ to $H^0(\omega_A)$ has polar divisor $D|V$ on $V$, and
so is a generator for $\omega_{A,q}$. \qed \enddemo

\demo{Proof of (12.7)} By (12.8) we need only show
there are no smooth basepoints. By Duality
$p$ is not a basepoint for $|\omega_A|$ iff $\Hom(m_p,\ring A.)$
is one dimensional.  If $p$ is a smooth point
$\Hom(m_p,\ring A.) = H^0(\ring .(p))$, so if the dimension is
at least two, then $|p|$ is basepoint free, and defines a birational
map to $\pr 1.$, and so an isomorphism. \qed\enddemo


\heading \S 13 Classification of Fences and Bananas \endheading

\noindent
{\bf 13.0 Notation, fixed throughout \S 13:}
{\it $S$ is a rank one log del Pezzo surface.
$X,Y \subset S$ are irreducible rational curves such that $K_S +X + Y$
is log terminal at singular points of $S$.}

For the definitions of Fence and Banana see (8.0.11-12). 

\proclaim{13.1 Lemma} Notation as in (13.0). If $X$ and $Y$ are
smooth then each of $X$, $Y$ contains at most one
singular point of $S$. If furthermore there is a singularity off $X \cup Y$,
then $S = \Hbz n.$, and $X,Y \in |\sigma_n|$ (so in particular, $X + Y \subset S^0$).
\endproclaim
\demo{Proof} We argue first that each of $X,Y$ contains at most one singular point. 

Let $B$ be either $X$ or $Y$. By (9.2), $S \setminus B$ contains either a unique 
singularity, or two $A_1$ points. Since there is no $S(4A_1)$ in the list (3.1), there is at
least one component which contains a unique singular point. 

Thus we can (switching if necessary) assume $X$ contains a unique singular point, and
(arguing by contradiction), $Y$ contains exactly two singular points, each $A_1$ points.
Note $K_S+X$ and $K_S+Y$ are both anti-ample, by adjunction.

Consider the hunt with respect to $K_S+ X$, which is clearly flush. Obviously $x$, for the first
hunt step, is the unique singularity along $X$. Note that after this we will 
not use flushness.

Note $X \subset T^0$.

 Suppose $T$ is a net. If $X$ is fibral, then it is a smooth fibre and by negativity $Y$
is a section, which meets $X$ once. But then $T$ has two fibres which both contain one $A_1$ singularity,
impossible. Thus, since $K_S+ X$ is negative, $X$ is a section, in the smooth locus of $T$. Thus $T$ is 
smooth, a contradiction (there are two singularities along $Y$).

 Hence we have a contraction $\pi:\map T.S_1.$. Since $X \subset T^0$, and 
$K_S+X$ is negative, $\Sigma\cap X=\emptyset$, and $X \subset S_1^0$. Note
$K_{S_1} + Y$ is negative, so $\Sigma \cap Y$ cannot contain both of the
$A_1$ points (otherwise $Y \subset S_1$ has a singularity at $q$). Thus by (13.7),
$S_1 = \Hbz 2.$,  $C + X \subset S_1^0$, $X$ and $C$ are linearly equivalent and 
$X^2= C^2 = 2$, where $C = \pi(E)$.  A contradiction, since $X \cdot C =1$.

So each of $X$, $Y$ contains at most one singular point of $S$. Now suppose
there is a singularity $z$ off $X + Y$, and a singularity $w$ on $X + Y$.
Then by (9.2), both of $z$ and $w$ are $A_1$ points. Thus $S$
has at least two singular points, and all of its singularities are $A_1$ points.
This contradicts the list (3.1). Now apply (13.7). \qed\enddemo

\proclaim{13.2 Proposition (The Smooth Banana)} If $K_S+X+Y$ is a banana, with
$X$ and $Y$ smooth, and $S$ is singular at some point of $X \cup Y$, 
then each of 
$X$ and $Y$ contains a unique singular point. Moreover one component is 
a $2$-curve and contains an $(\ul{r+1},2)$ point, the other is a $1$-curve and
contains an $A_r$ singularity. After extracting the $-(r+1)$-curve, $E$, there is a 
$K$-negative contraction down to $\Hbz 2.$, and the image of $X+Y+E$ is two 
sections (in $|\sigma_2|$) and a fibre. After extracting the adjacent $-2$-curve, $E'$, 
of the $A_r$ point, there is a $K$-negative contraction down to $\pr 2.$, and the image of the 
(component of $X + Y$ which is a) $2$-curve is a conic, the image of the $1$-curve
is a secant line to the conic, 
and 
the image of $E'$ is a tangent line to the conic.

 Moreover if $K_S+aX+bY$ is flush and anti-ample and does not have 
a tiger, then the next hunt step for $K_S + aX + bY$ extracts $E$, and
$S_1 = \Hbz 2.$. 
\endproclaim
\demo{Proof} Suppose  $Y$ is contained in
the smooth locus. Then $(K_S+X+Y)\cdot Y=0$, and $K_S+X+Y$ 
is numerically trivial. But then $(K_S + X + Y) \cdot X =0$ which
implies $X \subset S^0$, a contradiction. Of course the same
applies to $X$, thus either component contains a singular point. By (13.1) each contains 
a unique singular point, and $S$ is smooth off $X + Y$.

Switching components if necessary, we can assume $X$ contains the singular point of highest 
index. Note $-(K_S+X)$ and $-(K_S+Y)$ are both ample, by adjunction.

Consider the hunt with respect to $K_S+ X$. Obviously $x$, for the first hunt step, is the unique 
singularity along $X$. We note that after scaling some coefficient could be greater than one 
(since we are not assuming there is no tiger). This will not effect our argument, which will 
not use flushness at all.

$T$ cannot be a net, as in the proof of (13.1).

 Hence we have a contraction $\pi:\map T.S_1.$. Since $X \subset T^0$, and 
$K_S+X$ is negative, $\Sigma\cap X=\emptyset$, and $X \subset S_1^0$. 
Then by adjunction $K_{S_1}+X+Y$ is numerically trivial. 
Thus $X+Y \subset S_1^0$, and $S_1$ is Du Val. By (13.7) $S_1$ is either $\pr 2.$ or $\Hbz 2.$.

Let $C \subset S_1$ be the image of $E$. 

If $S_1=\Hbz 2.$, $C$ is a fibre, thus $C$ contains a singular point, and so 
(since $E$ contains at most one singular point), $\Sigma$ meets $E$ at a 
smooth point, and $Y$ at a singular point. $Y+C$ has an ordinary  node at 
$q$. The given description now follows from (11.1.0).

If $S_1=\pr 2.$, then $X$ and $C$ are lines and $Y$ a conic tangent to $C$ at 
$q$. Now consider the composition of blow ups $h:\map \t{T}.S_1.$, as in \S 11.
The first two blow ups must be at $Y$ (to separate $Y$ and $C$). In this 
configuration, $C$ is a $-1$-curve. Since $\t{C}^2 \leq -2$, the next blow up 
is also at $C$. Now there is only one possible choice for any further blow up 
such that $K_T + Y$ is plt. This yields the given description.

  Now suppose $K_S+aX+bY$ is flush and anti-ample, we do not have a tiger, and the next hunt 
step yields  $\pr 2.$. We may assume (by switching components if necessary) that $x\in X$. 
As we get a tacnode, $e < 2/3$ by (11.1.1). Then by (8.0.7.2), the spectral value of 
$(\ul{r+1},2)$ is at most one, thus $r=1$, $x$ is an $A_1$ point and the singularity on $Y$ 
is an $A_2$ point. By adjunction $b>3/4$ (or $K_S + X + b Y$ is anti-nef) and so, by 
flushness, $a/2 = e> (3/4)(2/3)=1/2$. Thus $a > 1$, a contradiction. \qed\enddemo

\proclaim{13.3  Proposition (Fence with a smooth side)} Let $S$ be a log del Pezzo
surface of rank one and suppose $X+Y$ is a fence, with $Y$ smooth.

If $X$ is smooth, then $-(K_S+X+Y)$ is ample, and $S$ is smooth away from $X+Y$. 

If $X$ is singular, then $Y$ is a $-1$-curve and contains at least two 
singular points. If furthermore $Y$ contains at least $3$ singular points then either
\roster
\item $K_S+X$ is ample and $K_S+ 1/2(X+Y)$ is nef, or
\item $S=S(A_1 + A_2 + A_5)$, $X \subset S^0$ is a genus one curve with a 
simple node, and $K_S+Y$ is numerically trivial.
\endroster
\endproclaim
\demo{Proof} If $X$ is smooth, then (13.1) applies. In particular $-(K_S+X+Y)$ is ample 
and $S$ is smooth away from $X+Y$.  

Now suppose $X$ is singular. If there is at most one singular point along 
$Y$, then $K_S+X+Y$ is negative by adjunction, a contradiction. If $\t{Y}$ 
moves then by (5.11), $\kappa(K_{\t{S}}+\t{X})=-\infty$. However since 
$h^1(K_{\t{S}})=h^0(K_{\t{S}}) =0$, 
$h^0(K_{\t{S}} + \t{X})=h^0(K_{\t{X}}) > 0$. Thus $Y$ is a $-1$-curve. 

Suppose $K_S+X$ is ample and $K_S + 1/2(X + Y)$ is negative. By (12.1.4) $K_S+X=D$ for some 
effective $D$ containing all the singularities of $S$ along $X$, and having a non-empty
moving part if $g(X) \geq 2$. It follows that the support of $D$ is not contained in $Y$. 
$K_S + Y + D$ is negative, so by adjunction, $Y$ can contain at most two singular points.

Finally, if $K_S+X$ is numerically trivial, then by adjunction (since $X$ is singular),
$X \subset S^0$, $X$ has genus $1$, and $S$ is Gorenstein. If $X$ has a cusp, then $S^0$ is
simply connected by (7.4.1), and so from the list (3.2), $Y$ contains exactly two singular 
points. Otherwise $K_S + X$ is log canonical, and by (9.2), 
$\sum_{p \in S} \frac{r_p -1}{r_p} \leq 2$. 

Suppose $Y$ contains at least three singular points. It follows
(since from the list, (3.1), $S$ cannot contain $4$ index two points)
that $Y$ contains exactly three singularities, and $S$ is
smooth off $Y$. Note 
$$
K_S^2 = \frac{1}{Y^2} = \frac{1}{-1 + \sum \frac{r_p-1}{r_p}}.
$$
It follows (since $K_S^2$ is an integer) 
that $\sum \frac{r_p-1}{r_p} =2$, $K_S^2 =1$,  and $\sum r_p = 11$.
Thus (from the list) $S=S(A_1+A_2+A_5)$. 
\qed \enddemo

\subhead Fences with one smooth side, and one side of genus $1$ \endsubhead

\proclaim{Proposition 13.4} Let $K_S + X + Y$ be a fence, with $Y$ smooth and $g(X) =1$. 
Suppose $X$ has an ordinary cusp, and $K_S + 1/2 X$ is negative. Further assume either that $X$ 
contains at most one singular point, or $K_S + 2/3 X$ is negative.

Then $X \subset S^0$, $S=S(A_1+A_2)$, and the configuration is  obtained from $\pr 2.$ by taking 
a flex line to a cubic, blowing up to separate the line and the cubic, then blowing 
down all the $-2$-curves. In particular $X^2=6$, $Y$ contains two 
singularities, an $A_1$ and an $A_2$, and there are no other singularities 
on $S$. Moreover, if $\Delta$ is any effective divisor supported on
$X + Y$, and $K_S + \Delta$ is anti-nef, then $K_S + \Delta$ has a tiger.
\endproclaim
\demo{Proof} We begin by showing that $X \subset S^0$. Suppose not. 
We consider the transformation
$S \dasharrow W$ of (12.4) (this is where we use the assumption that $K_S + 1/2 X$
is negative).
\proclaim{Claim} $W$ has rank one.
\endproclaim
\demo{Proof of Claim} Note that by (9.2), either $X$ contains a single singular point of 
$S$, or exactly two $A_1$ points and $S$ is smooth off $X + Y$. Now suppose the claim fails. 
Then by (12.4.1), $M$ is irreducible, and $X$ contains two $A_1$ points. Since $Y$ is a 
$-1$-curve, $S$ is not Du Val along $Y$. Thus by (12.1.4), $M$ contains at least three 
singular points of $S$. This contradicts (12.5). \qed \enddemo

Let $G$ be an exceptional divisor adjacent to $X$. $W$ is a rank one Gorenstein
del Pezzo, $X \subset W^0$, and $K_W \cdot G = K_W \cdot Y = -1$. $W^0$
is simply connected by (7.4.1). Thus $W = S(E_8)$ by (3.9), and 
$X$, $G$, $Y$ are all members of $|-K_W|$. But any two members of $|-K_{S(E_8)}|$ intersect
in a single point, the unique basepoint of the linear system, by (3.6). But clearly
$G \cap X$ and $X \cap Y$ are different points of $W$. Thus we have a contradiction.

So we may assume $X \subset S^0$. $S^0$ is simply connected by (7.4.1). By (13.3) there are 
two chain singularities on $Y$, so from the list $S=S(A_1 + A_2)$. Blowing up the $A_1$ 
point and blowing down $Y$ gives the described configuration in $\pr 2.$. 

For the final statement, let $\Delta = aX + b Y$. If $a\geq 5/6$, then there is a tiger over 
the singular point of $X$ (the $-1$-curve of (11.2) configuration III), else $Y$ itself is a 
tiger, as $K_X+(5/6)X+Y$ is numerically trivial. \qed\enddemo

\proclaim{Proposition 13.5} Let $K_S + X + Y$ be a fence, with $Y$ smooth, $g(X) =1$ and $X\subset S^0$. 
Assume $X$ has an ordinary node. Then one of the following holds:
\roster
\item $S=S(A_1 + A_2)$ with the construction as in (13.4). This
is the only possibility with $S^0$ simply connected.
\item $S=S(A_1 + A_5)$. The configuration is obtained from 
a nodal cubic $X$ and a smooth conic $C$ meeting $X$ to
order $6$ at a point $p$. Blow up $6$ times at $p$, along $C$. $Y$ is
the final exceptional divisor over $p$. 
\item $S=S(3A_2)$. Begin with a nodal cubic $X$ and two
flex lines $L$ and $M$, meeting $A$ at $l$ and $m$. Blow up
$3$ each times at $l$ along $L$, and $m$ along $M$. $Y$ is the final 
exceptional divisor over $l$. The singularities along $Y$ are $A_2$
and $(\ul{2}^L,2^M)$. There is additionally one $A_2$
point away from  $X + Y$.
\item $S=S(A_2 + A_5)$. 
\item $S=S(A_1 + 2A_3)$, $2A_3 \in Y$.
\item $S=S(2A_1 + A_3)$, $A_1 + A_3 \in Y$.
\item $S = S(A_1 + A_2 + A_5)$, $A_1 + A_2 + A_5 \in Y$.
\endroster
\endproclaim
\demo{Proof} Assume we do not have (7). Then by (13.3) there are exactly two
singularities along $Y$, say $A_s$ and $A_t$. $S$ is Du Val. We have
$1< K_S^2 = \frac{1}{Y^2}= \frac{(t+1)(s+1)}{st-1}$. Now one goes
through the list (3.1)  to check these are the given possibilities for
$s$, $t$. It remains now to show that the configurations are as
described. For (1-2) one blows up the $A_1$ point and
contracts $Y$ just as in the proof of (13.4). 

For (3): Suppose first $S$ is any $S(3A_2)$. If $f: \map T.S.$
extracts a $-2$-curve, $E$, then by (3.3), (3.4), 
and the list, the $K_T$-negative ray gives a birational contraction
$\pi:\map T.S(A_1 + A_2).$. Let $\Sigma$ be the $-1$-curve
contracted by $\pi$. It follows that 
$K_T + \Sigma$ is log terminal, $\Sigma$ contains two
$A_2$ points, and meets $E$.

Now apply the above with $E$ one of the $-2$-curves in the $A_2$ point off
of $Y$. Since $K^2 = 3$, $K + Y + \Sigma$ is anti-ample. It follows
that $Y \cap \Sigma$ is a single point, one of the $A_2$ points
along $Y$. Suppose $Y$ and $\Sigma$ meet the same $-2$-curve. Then
when we extract this curve, both $Y$ and $\Sigma$ are contractible,
a contradiction. Thus $K_S + Y + \Sigma$ is log canonical.
Resolve the $A_2$ point $Y \cap W$, and contract $Y$ and $W$
(which are disjoint on the resolution). This gives the described
configuration. \qed \enddemo

\remark{Remark} One can also give explicit constructions in 
(4-6), along the lines of (2-3). 
We have not done so as we will not need it. \endremark

It will be convenient to have a version of the final statement of 
(13.4) for a slightly more general situation:

\proclaim{13.6 } Let $X$, $Y \subset S$ be two integral rational curves such that
$X \subset S^0$ has arithmetic genus one, and an ordinary cusp, and 
$Y$ meets $X$ at only
one point, and normally (but unlike in (13.0), we do not assume
$K_S + Y$ is lt at singular points). Then $K_S + \Delta$ has a tiger
for any effective $\Delta$ with support contained in $X + Y$ such that
$K_S + \Delta$ is anti-nef.
\endproclaim
\demo{Proof} $X \in |-K_S|$, so $S$ is Gorenstein. $S^0$ is simply connected by (7.4.1). 
$K_S \cdot Y = -1= - X \cdot Y$.

Observe that $Y$ is the unique $-1$-curve of (3.8). 
Indeed if $K_S^2 \geq 2$ then $K_S + Y$ is anti-ample, so $Y$ is smooth and 
thus a $-1$-curve. Otherwise, by the list, $S = S(E_8)$, $Y \in |-K_S|$, and
since $X$ has a cusp, $Y$ is the unique $-1$-curve on $S$, by (3.6). 
We note that in this case, $Y$ itself is not smooth.

Now we consider the possibilities in (3.8). In (3.8.1-2), $K_S + Y$ is log
terminal, so $X + Y$ is a fence and (13.4) applies.

In the remaining cases, (3.8.2-5), $Y$ meets a unique curve, $V$ of the minimal 
desingularisation. Let $f:\map T. S.$ extract $V$. $Y \subset T^0$ is a 
$-1$-curve. Let $\pi:\map T.S_1.$ contract $Y$.
Let $Y_1 = \pi(V)$. Observe $(S_1,Y_1,X)$ again satisfies the assumptions of (13.6).
$K_{S_1}^2 = K_S^2 -1$. Scaling as in (8.2.4), we can find $\Delta_1$ supported in 
$X + Y_1 \subset S_1$, so that a tiger for $K_{S_1} + \Delta_1$ implies a tiger for 
$K_S + \Delta$. Thus we may induct on $K_S^2$. \qed \enddemo

\proclaim{13.7 Lemma} Suppose $C \subset S^0$ is a smooth rational
curve. Then either $S = \pr 2.$, or $S = \Hbz n.$ and $C \in |\sigma_n|$.
\endproclaim
\demo{Proof} If $S$ is smooth, $S = \pr 2.$. Otherwise let
$f :\map T.S.$ extract any exceptional divisor of the minimal 
resolution. Since $K_S + X$ is negative, and $X \subset T^0$, $T$ must
be a smooth net, and $X$ a section. \qed \enddemo


\heading \S 14  $T_1$ a net \endheading
Here we prove
\proclaim{14.1 Proposition} If $T_1$ is a net, and $\alg (S_0^0)=0$, then $S_0$
has a tiger.
\endproclaim

{\bf We assume throughout \S 14 that $S_0$ does not have a tiger, and
$\alg (S_0^0)=0$.}
We will derive a contradiction.

We follow the notation of (8.2.10). Since there are no further hunt steps
we drop the subscripts, writing for example $T,S,E,e$ for $T_1,S_0,E_1,e_0$.

We have the $\pr 1.$-fibration $\pi:\map T.{\pr 1.}.$. 
Recall in particular that $a$ is the coefficient of $E$ in $\Gamma'$, that is
$K_T + a E$ is $\pi$-trivial.

$K_T + a E$ is flush, klt, and anti-nef by (8.4.5).

\proclaim{14.2 Lemma} $\alg (T^0) = 0$. $E$ is not fibral. Let $d$ be its degree. $d \geq 3$.
$e< a =2/d$. Any singularity along $E$ has spectral value at most one.
\endproclaim
\demo{Proof} The first claim follows from (7.3). $E^2 < 0$, so $E$ cannot be fibral. 

Let $F$ be a smooth fibre. Note $F$ and $E$ generate the two edges of the cone of curve 
$\ol{NE}_1(T)$. 
$$
(K_T + E) \cdot E < (K_T + e E) \cdot E = f^*(K_S) \cdot E = 0
$$
since $E^2 < 0$. 
$$
(K_T + E) \cdot F = -2 + d
$$
where $d$ is the degree of $E$. $-(K_T+E)$ is not nef (otherwise
$E$ is a tiger), so $d \geq 3$.

$$
0 = (K_T + a E)\cdot F = -2 + a \cdot d 
$$
Thus $a = 2/d \leq 2/3$. Now by (8.0.7.2) any singularity along $E$ has spectral 
value at most $1$. \qed \enddemo

We get some control from the triviality of the algebraic fundamental group:

\proclaim{14.3 Lemma} If $f:\map T.{\pr 1.}.$ is a $\pr 1.$-fibration with 
irreducible fibres and with three or more multiple fibres, or two multiple 
fibres of non-relatively prime multiplicity, then $\alg (T^0)\neq 0$.
\endproclaim
\demo{Proof} If two fibres have non-relatively prime multiplicity, then there 
is a non-trivial torsion Weil divisor, and hence a finite \'etale cover of
$T^0$. If there are at least three multiple fibres, then there is a branched cover 
$\map C.{\pr 1.}.$, with ramification matching the multiplicities of the fibres (see (6.5) 
of \cite{Kollar91}). The normalisation of the pullback of $T$ gives a cover, \'etale in 
codimension one (see (III.9.1) of \cite{BPV84}).  \qed \enddemo

So by (14.3), $\pi$ has at most two multiple fibres, and if there are 
two multiple fibres, their multiplicities are relatively prime.

\proclaim{14.4 Lemma} $e \geq 1/2$.
\endproclaim
\demo{Proof} Assume $e < 1/2$. By the final remark of (10.8), $S$ has at least three
non Du Val singularities, and $1/3 < e < 1/2$. Thus  by (10.1), $x = (3,A_r)$ for $r > 1$ 
and $E$ contains a unique singularity, $y$, an $A_r$ point. 
$e = \frac{r +1}{2r + 3} \geq 2/5$. Thus by (14.2), $d \leq 4$. Let $G_1$ be the multiple 
fibre containing $y$. By (11.5.5), there is exactly one other multiple fibre, $G_2$, and 
it contains a non Du Val point. 

The possible non Du Val fibres, by (11.5.9), are (11.5.9.2), (11.5.9.6), and (11.5.9.5) with 
$k=2$, which have respective multiplicities three, five, and three. 

(11.5.9.6) is not possible, for $E$ must meet such a fibre at a singular point (otherwise 
$d$ is at least $5$), a contradiction, as this fibre contains no Du Val singularity (and 
the only singularity on $E$ is Du Val).

Thus any non Du Val fibre has multiplicity three. In particular this holds for $G_2$. Thus 
by (14.3), $G_1$ is Du Val.  $G_2$ then must contain two non Du Val points, contradicting the
descriptions (11.5.9.2) and (11.5.9.5). \qed \enddemo

By (14.2) and (14.4), $d =3$.

We let $h:\map \t{T}.W.$ be a relative minimal model, as in (11.5.4). By (11.5.5) there are 
two choices for $h$ at each multiple fibre.

\proclaim{14.5 Lemma} Let $G$ be a multiple fibre of $\pi$. 
We may choose $h$ so that 
$h(E)$ is smooth at the image fibre, and $\t{E}^2$ goes up by at most $4$. 
Further if $\t{E}^2$ changes by $4$, then $G$ is uniquely determined. 
\endproclaim
\demo{Proof} We use the classification (11.5.9). We check first that it applies. By (14.2), the 
singularities along $E$ have spectral value at most $1$, so if $G$ meets $E$ in a singular point, 
then by (8.0.8.b), (11.5.9) applies. Otherwise $G$ meets $E$ in a single smooth point, and 
so the multiplicity of $G$ is three, and we can apply (11.5.13).

We will indicate the type of the fibre (using the notation of (11.5.4)),  the change in $\t{E}^2$ 
under $h$, and leave the choice of $h$ to the reader. We let $c$ be the change in $\t{E}^2$ (for 
optimal choice of $h$). If $G$ meets $E$ in two singular points, then the fibre is 
$$
(-\ul{3})+(-1)+(-\ul{2})+(-2), \quad c=2.
$$
If $G$ meets $E$ in one singular point, and one smooth point, then the fibre 
is 
$$
(-\ul{2}) + (-\ul{1}) + (-2), \quad c=2,
$$
or given by (11.5.9.1) and $c=3$. If $G$ meets $E$ in one point, a singular point, and $\t{E}$ 
and $\t{G}$ meet, then the fibre is again $(-2) + (-1) + (-2)$, $c=2$. If $G$ and $E$ meet only 
once, at a smooth point, then $m=3$, and the fibre is given by (11.5.9.2) or (11.5.9.5), with 
$c= 4$, or $3$ respectively. Finally if $\t{E}$ and $\t{G}$ are disjoint, and $G$ meets $E$
in a single singular point, then the fibre is given by
(11.5.9.5) with $k=3$, or (11.5.9.6), and in either case $c=3$.\qed\enddemo

\demo{Proof of (14.1)} As remarked below (14.3), $\pi$ has at most two multiple fibres, and 
cannot have two of the same multiplicity. Thus by (14.5) we may choose $h$ so that $h(E)$ 
is a smooth triple section of self-intersection at most $5$ (corresponding to the  case when 
$\t{E}^2 = -2$, there are two multiple fibres, and, under $h$, $\t{E}^2$ goes up by three at 
one multiple fibre, and four at the other). But by (11.5.1)0, $h(E)^2 = 2d =6$. Thus we have
a contradiction. \qed\enddemo


\heading \S 15 $g(A) > 1$. \endheading
Here we prove
\proclaim{15.0 Proposition} If $S=S_0$ does not have a tiger and 
$g(A_1) \geq 2$, then $S$ is the surface in (15.2), and $S^0$ is 
uniruled.
\endproclaim
{\bf Throughout \S 15 we assume $S=S_0$ has no tiger.}

$1/3 \leq e_0 < a$, by (10.8), and so  by (8.3.7) $A$ either has multiplicity three, or a 
double point. By (8.4.7), $K_{T_1} + a_1 E_1$ is flush.  The possible configurations 
for $\map T_1.S_1.$ (in a neighbourhood of $\Sigma_1$) are classified in (11.1-3). 

We will repeatedly use the hunt notation (8.2.10), the notation of (11.0.1) to
describe the contraction $\pi:\map T_1.S_1.$, and the notation of (10.2). 

We begin by ruling out multiplicity three:

\proclaim{15.1 Lemma} $A_1$ cannot have multiplicity three.
\endproclaim
\demo{Proof} Suppose $A_1$ has multiplicity three. The possible local pictures for $\pi$ are 
given in (11.3), with $W = E_1$. By (8.0.9.1), $e_0<1/2$, so $x_0= (3,A_s)$ for some $s$, by 
(10.1). By (10.8), $K_{\t{S}_0}^2=-1$.

Consider what happens when $\pi$ is given by (11.3.2) or (11.3.3.2) (numerically these cases 
are indistinguishable). In the notation of (11.0.1), $X$, $Y$ are the two analytic discs 
where $E$ meets $\Sigma$. From the description of the configuration, $x_0=(3,2)$. By (10.8) 
there is a non Du Val point away from $\sb$. In particular $\alpha \geq 1/3$ ($\alpha$ of 
(10.2)). By (10.9)
$$
K_{S_0}^2 = -1 + 1/3 + 2/5 + \alpha = 
\frac{(K_{S_0} \cdot \Sigma)^2}{\Sigma^2} = 
\frac{(-1 + 1/3 + 3/5)^2}{(-1 + 1/3 + 7/5)}.
$$
and so $\alpha = 3/11 < 1/3$, a contradiction.

 Now consider (11.3.3.1). Since there is at least a $(2,3)$ singularity, $x_0=(3,A_s)$ with 
$s \geq 1$. By (10.6) there are no singularities away from $\sb$. As $\rho=10$, $s=6$. But 
then the singularities along $\sb$ are exactly $(3')$, $(2',3)$ and $(3',A_6)$ and so 
$K_{S_0} \cdot \sb = 0$, a contradiction. \qed\enddemo

 We now turn to double points of genus $g \geq 2$.

We consider in turn each possible configuration on $T_1$ from the classification (11.1-2).   
We note that configuration $0$ for a node cannot occur (since $E$ is smooth). In the
case of the node we apply the classification with $X$ and $Y$ (in the notation of (11.0.1)) 
the two branches of $A$ at $q$ (that is the two branches of $E$ along $\Sigma$), and $c =d =a$. 
For the cusp, we take $X = A$, $c = a$. We note that the node and cusp configurations
labeled by the same roman numeral are numerically indistinguishable.

\subhead   Configuration $I$ \endsubhead
In this subsection we also allow $g =1$, since most of the analysis in this case is the same 
(and this will save a little work in \S 16-17).

By (11.1-2.1), $a = 1/2$. Thus $1/3 \leq  e_0 < 1/2$. By (10.8), $K^2_{\t{S_0}}=-1$, 
$e_0 > 1/3$. By (10.1), $x_0 = (3,A_r)$ with $r \geq 1$, and 
$e_0 =  \frac{r +1}{2r + 3} \geq 2/5$. By the description of Configuration I, 
$A_1 \subset S_1$ contains a unique singular point, an $A_r$ point. By (10.9) we have 
$$
-1+e_0+\alpha = \frac{(K_{S_0} \cdot \ol{\Sigma})^2}{\ol{\Sigma}^2}=
\frac{(-1 + 2e_0)^2}{4e_0+\frac{g-1}g -1}
$$
so that $\alpha = \frac{e_0+g-1}{4e_0g -1}$. By (10.8), $n \geq 3$, so there
are at least two non Du Val singularities off $\sb$. Thus $\alpha \geq 2/3$.
Since $e_0 = \frac{r+1}{2r+3} \geq 2/5$, 
$\alpha \leq \frac{5g-3}{8g-5} \leq 2/3$. Thus we must have equality in each
of the inequalities, $r=g=1$ and the non Du Val singularities of $S_0$ are exactly 
$(2,3),(3),(3)$. By (10.6.2) there is exactly one other singular point, an $A_6$ point. Now 
consider the transformation $S_1 \dasharrow W$ of (12.4). $W$ has rank one by (12.4.1). The
$-1$-curve $M$ (of (12.0)) is irreducible, $K_{S_1}+M$ is negative, and $M$ contains the singular 
points $(3),(3),(2)$ of $S_1$, by (12.1.3-4). By adjunction $M$ cannot contain any other singular 
points, and $K_{S_1}+ \ol{M}$ is log terminal, where $\ol{M}$ is the reduction of $M$. Thus 
$W = S(A_2 + A_6)$. But there is no such surface in the list (3.1). \qed

\subhead Configuration $II$ \endsubhead

By (11.1.1) and  (11.2.1) $e_0 < a = \frac{g+1}{2g+1} \leq 3/5$, and so by (10.1) either 
$x_0=(2,3,2,2)$ and $g=2$, or $x_0=(A_g,3)$.  

Suppose $x_0=(3,A_g)$. $e_0 < 1/2$, so the final remark of (10.8) holds.
Let $V$ be the $-2$-curve which
meets $\Sigma$. One checks $\ol{\Sigma}$ has coefficient $\frac{4g}{2g+3}$ on 
$V$ and $\frac{2g+1}{2g+3}$ on $E$. Also $e_0 = \frac{g+1}{2g + 3}$. By
(10.9) we have 
$$
-1 + \frac{g+1}{2g+3} + \alpha = 
\frac{(-1 + \frac{2g+1}{2g+3})^2}{\frac{6g+1}{2g+3} -1}
$$
so that $\alpha=\frac {g}{2g-1}$. Since $n \geq 3$, $\alpha \geq 2/3$, so it follows that 
$g=2$ and $\alpha=2/3$. Thus there are exactly two other non-Du Val points, each $(3)$. 
Since $\rho = 10$ there is another singular point, $z$, the only other singularity, by 
(9.2). Thus $\rho(z) = 5$. But then $z$ has index at least $6$, and we have singularities 
of index at least $(3,3,6,7)$, which violates (9.2).

Now suppose $x_0 =(2,2,3,2)$. By (10.11), $E$ is the $-3$-curve. 
Then, by the description of configuration II,
$A$ contains a unique singular point, $w$, of $S_1$, an $A_1$ point. 
$M$ (of \S 12) is irreducible or empty by (12.1.7).
$K_{S_1}+M+F$ is negative, and so by (12.1.4) and adjunction, $n < 4$.
$e_0 = 6/11$, so we have one of (10.8.3-5).
In particular, $\kts = -1$ or $-2$. We have
$$
\kts +6/11+\alpha =\frac{(-1+6/11+4/11)^2}{-1+24/11}
$$
thus $\alpha =6/13 < 1/2$ or $\alpha =1+6/13$, depending on whether
$\kts$ is $-1$ or $-2$. 

Suppose we have the second possibility. Then by (10.8), $n=3$, $w = 4$, thus 
there are two more singular points, one, $y$ of weight one, and the other, 
$z$, of weight two. The only weight two singularities of coefficient
at most $e_0 = 6/11$ are given by (10.1.2). In particular, $\delta(z) = 1$. Hence
$\delta(y) = e(y) = 6/13$. By (10.1), $y=(3,A_5)$. Since $\rho = 11$,
$\rho(z) =1$, $z = (4)$, and there are no other singular points.
By (12.3.5) $M \subset S_1$ is nonempty, and must pass through $(2)$ and 
$(4)$. By adjunction $M$ is smooth. But then $M$ contracts on $S_1$,
a contradiction. 

Thus we have the first possibility. Necessarily there is exactly one more non Du Val 
singularity, $y$, and $w(y) =1$. So 
$$
\alpha = \delta(y) = e(y) = 6/13
$$
and thus $y = (3,A_5)$ by (10.1). Since $\rho  = 10$, there are no other singular points.
We now show there is exactly one possibility for $S = S_0$. 

{\bf Notation:} {\it For a conic $X \subset \pr 2.$ and points $w,z,p \in X$, we let
$L_{wz}$ and $M_p$ indicate the secant through $w,z$, and the tangent line through
$p$. }

\proclaim{15.2 Definition-Lemma} Let $A$ and $B$ be two smooth conics meeting to 
order four in $\pr 2.$ at a point $c$. Pick a point $b\neq c$ of $B$. Pick 
a point $a$ of $A\cap M_b$. Let $b'$ be the other point of $L_{ac}\cap B$. 
Let $a'$ be a point of $A\cap L_{bb'}$. 

 Let $f:\map \tilde S.{\pr 2.}.$ blow up once at $a$, $a'$, $b'$, twice along 
$M_b$ at $b$, and five times along $A$ at $c$.

 We have 
\roster
\item $\tilde S$ is the minimal resolution of $S$.
\item All configurations $(A, B, b)$ are projectively isomorphic. In 
particular $S$ is unique.
\item $L_{aa'}$ meets $B\setminus \{b,b',c\}$ in two distinct points. 
\item $L_{bb'}$ is not tangent to $A$. Equivalently $q$ is not a unibranch
singularity of $A_1$. 
\item $S^0$ is uniruled.
\endroster
\endproclaim
\demo{Proof} 
Let $\map Z.S_1.$ extract the divisors adjacent to $F$ (of \S 12). 
Then $F$ is the general fibre of a $\pr 1.$-fibration on $Z$, by (12.3). 

We argue first that $y$, $w \in F$. Otherwise $Z$ has Picard number
two, and $\map Z. S_1.$ extracts a unique divisor $E$, a section
of the fibration. $E$ cannot be over $w = A_1$, for otherwise 
$E \subset Z^0$, and so (since $E$ is a section) $E$ is smooth, a 
contradiction, as $y \in Z$. Then by (12.3.3), $E$ is the $-3$-curve
of $y$, so $Z$ is Du Val, and has an $A_5$ singularity, contradicting (11.5.4).

Now we claim $M$ is empty: As we have observed, it has at most one irreducible
component. Since $y$, $w \in F$, and $K_{S_1} + F + M$ is negative, $M$ can only contain
one singular point by adjunction, and so by (12.1.6), $w \in M$. But then $M$ contracts
on $S_1$, a contradiction. So $M$ is empty. Then by (12.3), $E$ is empty
as well, and $K + \t{A} = \t{F}$. 

$\map Z. S_1$ extracts (by (12.3.3)) the divisors $M_b$, the $-2$-curve of $w$, and $B$, the 
$-3$-curve of $y$. $Z$ has Picard number three. There is a unique reducible fibre $F_1+F_2$, 
and $A \subset Z^0$. $Z$ is Du Val, with unique singular point $v=A_5 \in B$. $F_1$ and $F_2$ 
are contractible $-1$-curves, and so the singularities along either are described
by (3.3). $K_Z + A = F$, so $A$ meets each $F_i$ normally in a single (smooth) point.
$v$ cannot lie on an irreducible fibre, by (3.4). Now $F_i$ is not in the smooth locus, or 
after contracting it, $v$ lies on an irreducible fibre. It follows that $v = F_1 \cap F_2$ 
and $K_Z + F_1 + F_2$ is lc at $v$.

Let $F_1$ be the component meeting $M_b$  (there can be only one as $M_b$ is a section 
contained in $Z^0$). If we contract $F_2$ we obtain $\Hz 2.$, $M_b$ becomes $\sigma_{\infty}$ 
and $B$, as it remains disjoint from $M_b$, becomes $\sigma_2$. It follows that $K_Z+B+F_1$ 
is lc at $v$, and $F_2$ and $B$ meet the same curve over $v$.

Let $t \in \Hz 2.$ be the image of $F_2$. Then $A$ and $B$ meet only 
at $t$, with fifth order contact.

Consider extracting $V$, the $-1$-curve of the second blow up at $q$ along 
$A$, and blowing down the fibre (through $q$). The resulting surface is 
$\pr 1. \times \pr 1.$, where $M_b$ becomes a zero curve, and
$A$ a section (in the direction for which $M_b$ is a fibre). Let $L$ be the 
unique curve disjoint from $M_b$ passing through $t$.  $L$ is a $2$-curve 
on $\Hz 2.$ through $t$ and $q$, and disjoint from $M_b$. $A$  meets $L$ to order 
four at $q$. $A$ and $L$ are transverse at $t$ (as is clear on $\pr 1. \times \pr 1.$)
thus $B$ meets $L$ transversally at $t$ as well. 

Let $L_{bb'}$ be the marked curve $(2,2^{L_{bb'}},3^A,2)$ over $x_0$ and 
$L_{ac}$ the marked curve $(3^B,A_4,2^{L_{ac}})$ at $y$. Let $G$ be the fibre 
of the ruling (on $Z$) through $q$. 

 Then $f$ is given by contracting $\Sigma$, $G$, $F_1$, $F_2$, $L$ and all exceptional 
curves over $S$, other than $A$, $B$, $L_{bb'}$, $L_{ac}$, $M_b$. Hence (1).

Pairs $(C,p)$ of a point on a smooth conic are projectively equivalent, 
so we may assume $A$ is the conic $yz = x^2$ and $c$ is the origin
$(0,0,1)$. Then $B$ is in the pencil of conics with $4^{th}$ order contact
at $c$, $yz - x^2 + ty^2 =0$, for $t \neq 0,\infty$. The automorphisms
$\map (x,y,z).(\lambda x , \lambda^2 y, z).$ (for $\lambda \neq 0$)
act transitively on the pencil elements with $t \neq 0,\infty$. Thus we may
assume $B: yz -x^2 + y^2=0$. The automorphisms
$$
\phi_d:\map (x,y,z).
(x + \frac{d}{2} \cdot y,y,d \cdot x + \frac{d^2}{4} \cdot y + z).
$$
for $d \in \af 1.$ preserve $A$ and $B$. Acting on a general point of
$B$ gives a non-constant map $\map {\af 1.}.B. \setminus \{c\} = \af 1.$,
which is necessarily a surjection. Hence (2). 

By (2) we can assume $b$ is the point $(0,-1,1)$. Now (3) and (4)
are straightforward calculations, which we leave to the reader.

By (3) and (4), $\ol{\Sigma}$ and $L_{aa'}$ are each $-1$-curves on $S$,
meeting $\sg( S)$ twice. $\ol{\Sigma}$ has branches $(2,2',3,2)$
and $(2,2,3',2)$, both at $x$, while $L_{aa'}$ has two branches,
each $(A_5,3')$ at $y$. Every branches is smooth by (4.12.2).
The indices of $x$ and $y$ are $11$ and $13$. Thus $S^0$ is uniruled
by (4.10.3) and (6.6). \qed \enddemo

\remark{15.2.6 Remark} It is tempting to think that (15.2.2) follows
immediately from the fact that $A+ B$ is canonically embedded. However, the scheme structure on
$A + B$ is not so obvious. If one constructs a scheme $Z(\phi)$ by gluing together two copies 
of $\pr 1.$ along an automorphism $\phi$ of the length $4$ curvilinear scheme supported at 
some point, then the isomorphism class of $Z(\phi)$ depends in general on $\phi$. For example 
if one takes the identity for $\phi$, then the canonical map for the resulting Gorenstein 
curve (of arithmetic genus $3$) is the two to one map to $\pr 1.$ which is the identity on 
each component. (15.2.2) can be viewed as saying that any two $Z(\phi)$ for which the 
canonical map is an embedding, are isomorphic.
\endremark

Note by (11.1-2) that for further configurations, $q \in A$ is a unibranch singularity. 
The possibilities are classified in (11.2.1.3-5). We consider each in turn.

\subhead Configuration $III$ \endsubhead 

Let $y$ be the $(A_{g-1},3)$ point along $\Sigma$. We have by (10.1.1) and (11.2.1.3) 
$$
e(y) =\frac g{2g+1} \leq e_0 < a = \frac{g+1}{2g+1}\leq 3/5.
$$ 
Thus 
$$
\frac{2g}{2g+1}  \leq e_0+e_y <1. 
$$
One easily rules out (10.8.1), (10.8.2) and (10.8.5). So we may assume there is exactly one 
other non Du Val singularity, $z$, of $S_0$, of weight one or two, and $w\leq 4$. 

 Suppose $z$ has weight two. Then we have (10.8.3). As $\beta >1/2$, $1/2 < e(z) < 3/5$. 
Thus by (10.1), $z=(4,2)$. $e(z) = 4/7 < a = \frac{g+1}{2g + 1}$, thus $g=2$. 
By (9.2) there are no other singularities. Since $\rho = 11$, $\rho(x) = 6$. $x$ has weight 
one, and  $1/2 < e(z) < 3/5$, so by (10.1), $x = (2,3,A_4)$ which has index 17. But now we 
have singularities with indices $(2,5,7,17)$, which contradicts (9.2). 

Thus $z$ has weight one. Suppose first $x$ has weight two. From $1/2 < e(z) < 3/5$,
it follows from (10.1) that $x = (4,2)$. Exactly as in the previous paragraph, $g =2$
there are no other singularities, and $\rho(z) = 6$. But then $z$ has index at least
$13$, and we have indices of at least $(2,5,7,13)$, which contradicts (9.2).

Thus $x$ has weight one as well, and we have (10.8.4). $\kts =-1$, and $\rho = 10$. 
By (10.9) we have
$$
-1 + e_0 + e_y + \alpha 
= \frac{(-1+e_y+e_0)^2}{-1+e_y+e_0+1/2}
$$
so that $\alpha = \frac{1-e_0-e_y}{2e_0+2e_y-1}$. As $\alpha \geq 1/3$, it 
follows that $e_0+e_y \leq 4/5$. Thus $\alpha =1/3$, $e_0=e_y=2/5$, and thus 
the non Du Val singularities are exactly 
$$
x,y,z = (2,3),(2,3),(3).
$$
By (9.2) there is at most one other singularity, and any additional singularity can have 
index at most $3$. But then $\rho \leq 7$, a contradiction. \qed

\subhead Configurations $u$ or $v$ \endsubhead

The final possibilities are the genus two unibranch singularities given by (11.2.1.4) and 
(11.2.1.5). In the first case, configuration $u$, there are two singularities along 
$\Sigma \subset T$, $(2,4')$ and $(\ul{2}',2)$. In the second case, configuration $v$, there 
is one singularity along $\Sigma$, $(2,3,2',2)$, and $E$ meets $\Sigma$ normally at a smooth 
point.

\proclaim{15.3 Lemma} If we have configuration $u$ or $v$ then $S_1$ has at least
two singular points along $F$ ($F$ of (12.3)).
\endproclaim
\demo{Proof} $\t{F} \subset \t{S}$ is a fibre, so obviously $F \not \subset S_1^0$.
Suppose there is a unique singular point, $b$ of $S_1$ along $F$.

Let $f:\map Z.{\pr 1.}.$ be 
the $\pr 1.$-fibration of (12.3) obtained  by extracting the unique divisor $E$, adjacent to 
$F$. Note $E$ is a section, and by (12.3.1) $A$ is a double section. 
As $\rho(S_0)\geq 10$, $\rho(Z)\geq 4$, and so there is at least one multiple 
fibre $F'$.  

Let $\Gamma$ be the fibre passing through the singularity $q$ of $A$,
$\Gamma_q$ the local analytic branch of $\Gamma$ at $q$, and $\Gamma_e$ the local 
branch where $\Gamma$ meets $E$. 

 Then $\Gamma$ is a smooth fibre and $q$ is a ramification point of $f|_A$. 

Since $A \cdot \Gamma =2$, $\Gamma_q$ and $A$ are separated by the first blow up of 
$\map \t{T}.S_1.$ at $q$. $x_0 \not \in \Gamma_q$ and by (11.2.1) the singularity of 
$\Gamma_q$ is respectively $(2',4)$ for $u$, and  $(2',3,2,2)$ for $v$. Thus 
$K_{S_0}+\Gamma_q$ is lt.  $\Gamma_e$ and $\Gamma_q$ pass through different points of
$S_0$, and are the only branches of $\Gamma$ through singular points of $S_0$.
Thus $K_{S_0} + \Gamma_e$ is not lc, for otherwise $\Gamma$ is a tiger.

Suppose $b \not\in A$. Thus $K_{S_1} + \Gamma_e$ is not lc at $b$. A multiple
fibre meets $A$ in only one point, for otherwise it will contain
three singularities (the two points where it meets $A$ and the point 
where it meets $E$), contradicting (11.5.5). As $f|_A$ ramifies at only two 
points, $G$ is the unique multiple fibre, and there is a unique singular 
point $b'$ on $E$. Thus $b$ is not a chain singularity (otherwise
$K_{S_1} + \Gamma_e$ is lt). By (12.3.3) $E$ is 
not a $-2$-curve, and as $e_0<2/3$, $b$ is given by (10.1.2.d), and
$G$ meets the unique $-3$-curve. But then $K_{S_1} + \Gamma_e$ is
lc (see Appendix $L$), a contradiction.

Suppose $b \in A$. $x_0$ is a chain singularity by (10.1) since $1/2 < e_0 < 2/3$. 
Thus $\Gamma_e$ does not meet an end of the chain at $x_0$. As $a < 2/3$, by
(8.0.7.2), any singularity along $A$ has spectral value at most $1$ (with respect to $A$), and so 
by (8.0.8) and (12.3.3), $K_{S_1} + F$ is plt at $b$. Thus $E$ is an end of the chain 
at $b$, and not at $x_0$. Thus $A \cap E$ is a smooth point of $Z$, there is a unique
singular point $b'$ on $E$, and $G$ is the unique multiple fibre. As $b'$ is a chain 
singularity, Du Val or almost Du Val, $G$ is described by (11.5.9). If $m\neq 2$, since 
$A$ is a double section, by (11.5.5), $G$ meets $A$ at a unique point, $b'' \neq b$, so 
$G$ contains two chain singularities, each either Du Val or almost Du Val. This rules out 
cases (2-4), (7-8), and (5) for $k \geq 3$. $A$ meets a curve over $b''$ of multiplicity 
two in the fibre, thus in (6) the marked singularity $(A,b'')$ is $(3',2)$, which has 
spectral value $2$, a contradiction. In the remaining cases, which are (1), and (5) with 
$k \leq 2$, $\rho (Z) \leq 3$, a contradiction. \qed\enddemo

\proclaim{15.4 Lemma} If we have type $u$ or $v$ then $S_0$ has a tiger.
\endproclaim
\demo{Proof} By (15.3) we may suppose $F$ passes through at least two singular points.

We show first that $M$ (of \S 12)  is empty. Suppose not. Then by (12.1.2), $M \cap F$ is a single
point $b$. So there is a singularity along $F \setminus M$. Let $r$ be the index of one such 
singularity.

Let $N$ be an irreducible component of $M$. Note that as $K_S+(7/11)A$ is negative, so is 
$K_S+(7/4)(M+F)$. Let $V$ be the curve which $F$ meets at $b$, and let $\eta = e(V,K_{S_1})$. 
Note $e(V, K_{S_1}+N+F)\geq 1$ and so the coefficient of $V$ in $N+F$ is at least $1-\eta$. 

We have
$$
0 > (K + (7/4)(N + F))\cdot F \geq -2 + \frac{r-1}{r} +1 +3/4(N + F) \cdot F
\geq  -1/r + 3/4(1 - \eta + 1/r).
$$

It follows that $2/3 > a > \eta > \frac{3r-1}{3r} \geq 5/6$, a 
contradiction. 

Hence $M$ is empty. The singularities along $A$, as
well as any non Du Val singularities of $S_1$ away from $A$, are 
described by (12.1.4) and (12.3.3-4). Let
$r_1,\dots,r_t$ be the indices of the singular points of
$A \cap F$ (when there are no singular points, take $t=1$, $r_1 = 1$).
Since $K_{S_1} + A$ is lt at singular points (by (8.4.5) and (8.0.4)),
and $A$ is in the Du Val locus
$$
A \cdot F = 2 + \sum \frac{r_i -1}{r_i}.
$$

Note that $K \cdot F = -2 + \alpha$ ($\alpha$ as in (10.2)).

By (11.2.1.4-5)
$$
-2 + \alpha + 7/11( 2 + \sum \frac{r_i -1}{r_i}) = (K_S+7/11A)\cdot F<0. \tag{15.4.1}
$$
Suppose first $t \geq 2$. Then (15.4.1) implies $t =2$, $r_1 = r_2 =2$ and $\alpha = 0$.
Thus $S_1$ is Du Val and (by the configuration on $T$) $n=2$. Hence by (10.8), $w=2$. Since 
$x_0$ has weight at least one, the configuration cannot be $u$ (which has a $-4$-curve away
from $x_0$), thus the configuration is $v$, $E$ is a $-3$-curve, and $x_0$ is a non-chain 
singularity with branches $(2,2)$, $(2)$ and $(2)$. But then $e_0 \geq 3/4$, a contradiction.

Next suppose there are at least two singular points along $F \setminus F \cap A$. Then 
(15.4.1) and (10.1) imply that there are exactly two singularities on $F$, each  of type
$(3)$. Thus $n=4$, and so by (10.8), three of the singularities on $S_0$ are $(3)$, clearly 
impossible.

Thus we can assume there is exactly one singularity, along $A$ (necessarily in $A \cap F$) 
and exactly one non Du Val singularity, $b$ (necessarily in $F \setminus F \cap A$) of $S_1$.
(15.4.1) implies $r_1=2$, and $b=(3)$ or $(3,2)$. $x_0$ cannot have weight one (or there is
another singularity of higher coefficient). If the configuration is $u$, then 
$e_0 \geq e(2,2,4,2) > 2/3$, a contradiction. Thus the configuration is $v$, and 
$x_0 = (2,4)$. Thus $n=3$ and $w=4$, and so by (10.8), $\rho(S_0) = 11$. By (9.2) there 
is at most one other singular point. Since $\rho(S_1) = \rho(S_0) -5 = 6$, there must be 
exactly one more singular point, Du Val, of index at least $4$. But now we have indices on 
$S_0$ of at least $(4,3,7,11)$, which violates (9.2). \qed\enddemo



\heading \S 16 $g(A)=1$: a cusp. \endheading

This section is devoted to proving:

\proclaim{16.1 Proposition} If $A_1$ has a cusp of genus one, and $\alg (S^0)=0$, then $S$ has a 
tiger.
\endproclaim

\noindent
{\bf Throughout \S 16 we assume $S$ has no tiger, and $\alg (S^0) =0$.}
Our goal will be to obtain a contradiction.

$K_{T_1} + a \Sigma$ and $K_{S_1} + a A$ are flush by (8.4.7). The possible configurations 
for $\map T_1.S_1.$ are given by (11.2). As Case I has been ruled out in (\S 15), we may 
assume by (11.2.1) that $a\geq 2/3$. Note by (7.3.1), since $q=q_1=\pi_1(\Sigma_1)$ 
is a smooth point of $S_1$, $\pi_1^{alg}(S_1^0)$ is trivial. $A \not \subset S_1^0$, for 
otherwise $K_{S_1} + A$ is numerically trivial, and $A$ is a tiger.

By (11.2.1), $2/3\leq a_1 < 4/5$. By (8.3.8), $e_1 \geq 1/3$, so $a_1+e_1\geq 1$. Since 
$A \neq \Sigma_2$ ($\Sigma_2$ is smooth) it follows from (8.4.7) that:
\remark{16.2} The possibilities for the next hunt step are:
\roster
\item $A_2+B_2$ is a tacnode
\item $A_2+B_2$ is a banana
\item $A_2+B_2$ is a fence
\item $T_2$ is a net
\endroster
\endremark

Recall that $E_1$ is a $-k$-curve, and $E_2$ is a $-j$-curve.

For the next lemma (and its proof) we use the notation of \S 12. 

\proclaim{16.3 Lemma} The following hold:
\roster
\item $S_2$ is not a fence,
\item $W$ (of (12.4)) has rank at least two,
\item $S_1$ is singular along $A$ in at least two points,
\item $S_1$ is Du Val outside $A$,
\item $S_1$ is not Gorenstein, and 
\item $S_2$ is not a tacnode.
\endroster 
\endproclaim
\demo{Proof} Suppose $S_2$ is a fence. Then by (8.4.7.8.1), $-(K_{S_2}+B)$ is ample, 
so $B$ is smooth by adjunction. But then by (13.4), there is a tiger, a contradiction.

 Suppose $W$ has rank one.  

As $A$ is simply connected, $W$ is also, by (7.4.1). 

Claim: There is only one singularity along $A$: Otherwise by (12.4) and (3.9), $W = S(E_8)$,
(for some $i$) $D_i \subset W^0$, $D_i \in |-K_W|$. But as $A \in |-K_W|$ has a cusp, this 
contradicts (3.6).  

Now $M$ is irreducible by (12.4.1).

Now define $\Gamma$ by $K_Y + \Gamma = f^*(K_S + a A)$, and define $\Gamma' = \lambda \Gamma$ 
such that $K_Y + \Gamma'$ is $\pi$ trivial (as in the definition of the hunt), let
$\Delta' = \pi(\Gamma')$. $K_W + \Delta'$ has a tiger by (13.6). But then $K_S + a A$ has a 
tiger, a contradiction.

 (3) follows from (2) and (12.4.1). 

 If there is a non Du Val singularity outside $A$, then $M$ passes through that point by (12.1.4).
Say there are $l$ singular points along $A$. Then by (2) and (12.4.1), $M$ has at most $l-1$ 
components, and contains at least $l+1$ points. By (12.1.2) some component of $M$ contains at
least three singular points, contradicting (12.5). Hence (4).

 Suppose $S_1$ is Du Val. By (3) and the classification of simply connected 
Gorenstein surfaces $S_1=S(A_1+A_2)$. Thus 
$$
\tilde{A}_1^2= K^2_W>K^2_Y=K^2_{S_1}=6,
$$
which contradicts (16.4) below. Hence (5). 

 Thus $S_1$ has a non Du Val point along $A_1$. (6) follows exactly as in
the second paragraph of the proof of (18.5). \qed\enddemo

\proclaim{16.4 Lemma} If $\tilde A_1^2\geq 7$, then $A$ has three 
singularities along it. 
\endproclaim
\demo{Proof} By (11.2.1) the self-intersection of $E$ can only go up 
by nine (corresponding to $(v,n^2)$). Thus $k=-2$ and $x_0$ is not a 
chain singularity, by (10.11). \qed\enddemo

\proclaim{16.5 Lemma} (16.2.2) does not occur. 
\endproclaim
\demo{Proof} Suppose not.
 
 $B \not \subset S_2^0$ (or $B$ is a tiger), so there is a singularity $y \in B$. 
$x_1 \in A$ (the only way to get a Banana), so there is a unique singularity along $E_2$, 
by (8.0.7.1), and $\Sigma$ meets $E$ at a smooth point. Since $a_1 + e_1 \geq 1$, 
$\Sigma$ must meet $A$ at a singular point $z$. The configuration along $\Sigma_2$ is 
given by (11.1.1.3). In particular $z$ is a Du Val point.

 As there is only one singularity along $B$, $B$ is not a $-1$-curve. On the other
hand $-K_S\cdot B\leq A\cdot B=2$. It follows that $B$ is a $0$-curve. 

Suppose $S_2$ is singular along $A_2$. Then $y$ is not a Du Val point, $x_0$ is not a 
chain singularity (as there are at least $3$ singularities along $A_1$, and hence along 
$E_1$). $x_1$ is at least $(2,3)$. Since $e_0 < 4/5$, one checks using (8.3.9) that the 
only possibility is that $k=j=2$, and the branches of $x_0$ are $(2)$, $(2)$ and 
$x_1 = (A_r,3)$, $r \geq 1$. Then $z = A_1$, $A_2$ contains an $A_1$ point, 
$y = (A_{r-1},3)$, and by (9.2) there are no other singularities. Let $\pi:\map T.S_2.$ 
extract the $-2$-curve, $E_3$, adjacent to $B$ (this is in fact the next hunt step, but 
we will not use this). $T$ is a net, $B$ is a smooth fibre, and $E_3$ is a section, 
containing a unique singularity, $(A_{r-2},3)$ (by which we mean $(3)$ if $r \leq 2$). 
Thus the fibre through the $(A_{r-2},3)$ point is the unique multiple fibre, and so this 
fibre contains the  $A_1$ point as well. But there is no such net by (11.5.9). 

 Thus $A_2$ lies in the smooth locus and $S_2$ is Gorenstein. As $y$ is a chain 
singularity and $S_2$ is simply connected, the possibilities for $S_2$ are $S(A_1)$, 
$S(A_1+A_2)$, or $S(A_4)$. Moreover $x_1=(j,y)$, $z=A_{j-1}$ (since $B$ is a $0$-curve) 
and $\t{A}_1^2=K_{S_2}^2-1$. Now $j\geq 3$, by (16.3.5). Since the spectral value is at 
most three by (8.0.7.2), $y=(2)$ or $(2,2)$ and $j=3$, by (8.0.8). If $y=(2,2)$, then 
$e_1>3/4$ and $a/3+b>1$, so $(K_{T_2} + \Gamma')\cdot \Sigma_2 > 0$, contradiction.

 Thus $y=(2)$ and by (16.4), and the simply connected list, $S_2=S(A_1+A_2)$. But then if 
$f:\map T.S_2.$ is the blow up at $y$, $T$ is a Du Val net, with exactly one singularity, 
an $A_2$ point, contradicting (3.4). \qed \enddemo

For the next Lemma, we allow the possibility that $A$ has a node (as the 
analysis will be the same in that case).

\proclaim{16.6 Lemma} (16.2.4) does not occur. 
\endproclaim
\demo{Proof} Suppose $T_2$ is a net. As $a_1 \geq 2/3$, necessarily $E_2$ is a section and
$A$ is a double section. By (11.5.11.2), $E$ contains at most two singularities
and so $x_1$ is a chain singularity. In particular $j \geq 3$.
But by (11.5.10), $E_2$ is a $-2$-curve, a contradiction. \qed \enddemo

As we have eliminated every case of (16.2), this completes the proof of (16.1). 


\heading \S 17 $g(A_1)=1$: a node \endheading

Our goal is to prove:
\proclaim{17.1 Proposition} If $\alg (S^0)=0$ and $g(A_1)=1$, and $A_1$ has a node, then 
either $S_0$ has a tiger, or $(S_2,A_2 + B_2)$ is a fence, with $B_2$ smooth. 
In any case $S_0$ is log uniruled. \endproclaim

For the cases when $S_0$ does not have a tiger, we will essentially classify the possible fences 
$(S_2,A_2 + B_2)$, see (17.3), as well as the original $S_0$, see (17.5-14); we will stop the
analysis once we have enough information to apply the criteria of \S 6, but further 
analysis of the same sort would yield an explicit classification. The fences
$(S_2,A_2 + B_2)$ are in the classification of (13.5), with one exception, (17.3.2), which is
partially classified in the proof of (17.3). 
 
{\bf We assume throughout \S 17 that $S_0$ does not have a tiger, and $\alg (S^0) = 0$. }

$K_{T_1}+a_1E_1$ is flush by (8.4.7). Configurations on $T_1$ are classified in (11.1). 
We showed in \S 15 that there is a tiger if the configuration is I. If $A_1 \subset S_1^0$ 
then $K_{S_1}+A_1$ is trivial and we have a tiger, a contradiction. 
Thus $A_1$ contains a singular point, and the configuration on 
$T_1$ is II or higher. Note also that $S_1^0$ is simply connected by (7.3.1). 

{\bf 17.1.0 Notation:}
{\it  Let $C$ and $D$ be the two branches of $A_1$ at the 
node, and $c$, $d$ the points of $T_1$ where the branches meet $E_1$. We may 
assume (by switching $C$ and $D$ if necessary) 
that the first two blow ups of $h:\map \t{T}_1.S_1.$ are along $C$. 
Let $r+1$ ($r \geq 1$) be the number of initial blow ups along $C$. Note
$d$ is necessarily singular. }

\proclaim{17.2 Lemma} 
\roster
\item $K_{T_1}+\Sigma_1 + E_1$ is log canonical.
\item $\sb$ has two smooth branches through $x_0$ and meets no other 
singularities.
\item If $c$ is smooth, then $d$ is an $A_r$ point, $r \geq 1$ and 
$a_1=\frac{r+1}{r+2}$.
\item $T_1$ is singular at some point of $E_1 \setminus E_1 \cap \Sigma_1$.
\item $S$ has exactly two non Du Val points and $e_0>1/2$. 
\item $a_1 \geq 2/3$, and $a\geq 4/5$, unless we have (3) with $r  \leq 2$.   
\item $A$ contains exactly one singularity.
\endroster
\endproclaim
\demo{Proof} (1)-(3) are immediate from (11.1.1) (since we have
ruled out configurations $0$ and I).

(4) holds, for otherwise $A \subset S_1^0$. 

For (5), one shows there is at most one non Du Val point away from $A$ as in
the proof of (16.3.4). So by (2), there are at most two non Du Val points
on $S_0$. Now (5) follows from (10.8). 

For (6), suppose we are not in the situation of (3). $a \geq 4/5$ by (3) and (11.0) unless 
$r=1$ or $r=2$. By (11.0) it's enough to consider the configurations after one more blow up. 
Suppose $r =1$. After the next blow up $c =(2)$ and $d = (3)$. One computes $a = 4/5$. 
Suppose $r =2$. After the next blow up $c =(2)$ and $d = (\ul{2},3')$. One computes $a=6/7$. 
Hence (6).

Finally suppose (7) fails. Then, since $\pi$ removes one of the singularities along $E$, 
$A$ contains two singularities, $x_0$ is a non-chain singularity, and we have (3). We 
will use (8.3.9) to compute $e_0$. If $E_1$ is not a $-2$-curve, then 
$e_0 \geq \frac{r+1}{r+2}$, a contradiction. Thus we may assume $k=2$. As $e_0> 1/2$, 
(by (5)), $e_0 \geq 2/3$, by the classification of non-chain singularities, and $r\geq 2$.

One checks that $e_0 \geq \frac{r+1}{r+2}$ unless $r=2$ and the other two 
singular points are $(2)$ and $(3)$. Thus $\t{A}_1^2 =4$. As $S_1$ has a
non Du Val point outside $A_1$, $M$ is reducible by (12.5). Thus $W$ of (12.4) has 
rank one, and $K_W^2 = \t{A}^2 =4$. $\pi$ (of (12.4)) has two 
exceptional divisors. $G_1$ and $G_2$ of (12.4) are the
$-2$ and $-3$-curves adjacent to $A$, and $G_i \subset Y^0$. By (3.9)
$W = S(2A_1 + A_3)$, and each $D_i = \pi(G_i)$ contains two
singular points, thus each $G_i$ must meet both $\pi$ exceptional divisors. But
then $D_1 \cap D_2$ contains two points, contradicting (3.9). \qed \enddemo

By (17.2), there are two non Du Val points, $x_0$ and $y$, of $S_0$ and
$S_1$ has exactly one singularity, $z$, along $A$.
Let $s$ be the index of $z$. We will first show that $S_2$ is a fence (17.3), 
and then partially classify each possible $S_0$, and show in each case,
using \S 6, that $S_0^0$ is uniruled. (With a little more analysis one 
can completely classify the possibilities, but we stop once we have  
enough information to apply the criteria of \S 6).

Recall from (8.2.10) that $E_2$ is a $-j$-curve, and $E_1$ is a $-k$-curve. 

\proclaim{17.3 Proposition (Outcome of the hunt)} $(S_2, A_2+B_2)$ is 
a Fence. $g(A_2) =1$, $B_2$ is smooth. 
\roster
\item If $x_1\in A_1$ then $\Sigma_2$ meets $E_2$ at a smooth point,
and contains a unique singular point $(A_t,3,A_{j-2})$ (for some $t$).
$K_T + \Sigma_2$ is log terminal, and $\Sigma$ meets the end of the
$A_{j-2}$ chain. $(S_2, A_2+B_2)$ is given by (13.5). $q_2$ is an $A_{t+1}$ point.
\item If $x_1 \not \in A_1$, and $A_2 \not \subset S_2^0$, then 
$S_0^0$ is uniruled. 
\item If $x_1 \not \in A_1$ and $A_2 \subset S_2^0$, then $(S_2,A_2+B_2)$ is 
given by (13.5.1), and $a<6/7$.
\endroster 
\endproclaim
\demo{Proof} As $a+e_1\geq 2/3+e(y) \geq 1$, we have the same possibilities for the 
outcome of the hunt as in (16.2). 

By (16.6) $T_2$ is not a net. 

 Suppose $x_1\in A_1$. Then $S_2=W$ of (12.4), and $M = \Sigma$. As $M$ is disjoint from $A$ 
(on $T = Y$ of (12.4)) we must go to a fence, with $B$ smooth by (11.1), since
$a_2 + b_2 > 1$. By (13.3) (with $B$ playing the r\^ole of $Y$ of (13.3)) 
there are at least two singularities along $B$, and $B$ is a $-1$-curve.
Thus $\Sigma_2$ meets $E$ at a smooth point, and contracts to a cyclic
Du Val singularity. Now (1) follows easily by (11.4). 

Thus we may suppose $x_1\notin A_1$. Suppose $A_2\cup B_2$ has a node of 
genus $g\geq 2$. Since $a_2 + b_2 > 1$, the configuration on $T_2$
is given by (11.1.1.2).

 If $z\notin \Sigma$, then $K_{S_2}+A_2$ is not trivial and $b<1/2$ by 
(11.1.1.2). Thus $x_1=(3,A_g)$, by (10.1). But then $B$ is contained in the smooth
locus and as $B$ is not a tiger, $g\geq 3$. (11.1.1.2) implies $b\leq 4/9$. 
On the other hand, $e_1\geq e(3,A_3)=4/9$, a contradiction. 

If $z\in \Sigma_2$, then $A_2$ is contained in the smooth locus, and $S_2$ is simply 
connected Gorenstein. As $z=A_g$, $e_1\geq ag/(g+1)$ and by (11.1.1.2) $2ag<(g+1)$, 
that is $2/3 \leq a<1/2+1/(2g)$. Thus $a=2/3$, $g=2$. By (17.2.6), $r=1$, and $x_0$ is a chain 
singularity. In particular $k \geq 3$ by (10.11). Thus 
$$
S_2^2 = A_2^2 = r + 4 + g -k \leq 4.
$$
But $B$ contains
a chain singularity (it is not a tiger), which contradicts the list (3.2). 

Thus we must go to a fence. $B$ is smooth, and by (13.3) contains at
least two singularities, and is a $-1$-curve.

If $A_2\subset S_2^0$, then $S_2$ is Gorenstein, and $S_2^0$ is simply connected, 
by (7.3.1), since $q_2$ is smooth. Thus $(S_2,A_2 + B_2)$ is given by (13.5.1). 
Blowing up the $A_1$ point and contracting $B$ gives a flex line to a nodal cubic, 
and the inequality $b/2 + 3 a < 3$. Since $a > b$, (3) follows.

Now suppose $A_2\not \subset S_2^0$. Thus $\Sigma_2$ meets $A$ at a smooth point. 
As $a+e_1\geq 1$, by (11.1) $\Sigma_2$ meets $E_2$ at an $A_s$ point, for some $s \geq 1$. 
$E_2$ is a $-2$-curve, so $x_1$ is a non-chain point. Since $A_1$ contains only one singular 
point, $\tilde A_1^2\geq 1$, by (12.4.6), and so (from the configuration on $T_2$), 
$\tilde A_2^2 = \t{A}_1^2 + 1 + s$. 

Suppose $\t{A}_2^2 \geq 4$ (for example if $s \geq 2$, or $\tilde A_1^2 \geq 2$).

We  apply (12.4) to $(S_2,A_2)$. Since $A_2$ contains a unique singularity, $W$ has rank 
one, by (12.4.1). $K_W^2 = \t{A}_2^2 \geq 4$, so $W$ is given by (3.9.2). $G$ (of (12.4)) 
is irreducible, the curve over $z$ adjacent to $A_2$. $M$ is irreducible. $D =\pi(G)$
and $B \subset W$ are each smooth $-1$-curves, and thus their configuration is given by 
(3.9.2). In particular they meet exactly once, at opposite ends of the $A_3$ point, which is 
necessarily $\pi(M)$, and each contains in addition an $A_1$ point. Since $G$ contains at 
most one singularity, $M$ meets $G$ at a smooth point. Now consider the sequence of blow ups
$h:\map \t{Y}.\t{W}.$. The first blow up is along $D$ (since $D$ is a $-1$-curve, and $G$ is 
$K_Y$ non-negative). All further blow ups must also be along $G$, for otherwise $M$ meets 
$G$ at a singular point. Say $G$ is a $-v$-curve. Then there is a unique singularity 
$(\ul{2},2,3,A_{v-2})$ on $M$, at $M \cap B$, $B$ meets the underlined curve, and 
$K + B + M$ is log canonical (so $M$ meets the end of the $A_{v-2}$ chain). Thus the branches
of $x_1$ are $A_s$, $(2)$ and $(2,2,3,A_{v-2})$. Since this last branch has index at least 
$7$, $s \leq 1$. Note also that $z = (\ul{v},2)$.

Hence $s=1$. Then $a_2 + b_2/2 =1$. By (10.1), $b_2 > e_1 \geq 1/2$,
thus $a \leq 3/4$. So by (17.2.6), $a=2/3$, $r =1$. Then $b < 2/3$ so $x_1$ is given
by (10.1.2.d). As $e(2,4,2)=2/3$, $k=3$. Thus $\t{A}_1^2 = 4 + r -k = 2$. Thus
the singularities are as in the previous paragraph. Since $e_1 = 1/2$, $k=2$, by 
(8.0.7.2). Thus $x_0=(2,3^A,2^G,2)$. The index of $x_0$ is eleven. $\sb$ has two smooth 
branches, $(2',3,2,2)$ and $(2,3',2,2)$. One computes $-K_{S_0}\cdot \sb =2/11$. Hence 
$S_0^0$ is uniruled by (6.5). \qed\enddemo

\proclaim{17.4 Lemma} Either $x_0$ is a chain singularity, or $S_0^0$ is uniruled. 
\endproclaim
\demo{Proof} Suppose $x_0$ is a non-chain singularity, and $S_0^0$ is not uniruled. We will 
use (8.3.9) to compute $e_0$. 
 
By (17.2) $c$ and $d$ are singular points. 

 Suppose $c$ has type $(2)$. By (11.1), the two points of $E_1 \cap \Sigma_1$ 
are $c=(\ul{2}')$ and $d=(\ul{2},\dots,3')$ (where there are $r-1$, 
$-2$-curves). In this case one checks that $a$ rescales to $\frac{2r+2}{2r+3}$. 
On the other hand, if $E_1$ is a $-3$-curve, then $e_0\geq \frac{2r+2}{2r+3}$. 
Thus $E_1$ is a $-2$-curve. As $e_0> 1/2$, it follows that the third branch is 
not an $A_1$-singularity. Thus (by the classification of non-chain singularities)
either $r=1$ and the last point has index at most $5$, or $r=2$, and the
last point has index three. One checks the only case is $r=1$, $e_0=2/3$, 
${\tilde A}^2=3$, $a=4/5$ and $A$ contains an $A_2$-singularity. 

 Suppose first $x_1 \not \in A_1$. Let $w = \Sigma_2 \cap E_2$. We have
(17.3.3). By (11.1), if $w$ is smooth, then since $B$ is a $-1$-curve,
$j =2$ and $x_0$ is Du Val, a contradiction. Thus $w$ is singular, and
$x_1$ is a non-chain point with branches $(2)$, $(2,2)$ and $w$. 
$z = \Sigma_2 \cap A$. Since $\t{A}_2^2 =6$, the first three
blow ups of $h:\map \t{T_2}.\t{S_2}.$ are along $A$. Thus
by (11.1), $w = (4',\ul{2})$, $w = (4',\ul{2})$, which has 
index greater than five, a contradiction.

We conclude $x_1=z$, and $(S_2,A_2 + B_2)$ is given by (17.3.1). $B_2$ 
contains an $A_1$ point, and $A_2^2=3$. It follows from (13.5) that
$S_2$ is given by (13.5.2). Let $R'$ be given by (17.5.2). Then
$R'$ has a node at $x_0$. One computes each branch has index nine and
$-K_{S_0} \cdot R' = 1/3 > 2/9$. Hence $S_0^0$ is uniruled by (6.5). 

 Now suppose $c$ does not have type $(2)$. Then both $c$ and $d$ have
index at least $3$, and the indices of singularities along $E_1$
are $2,3,n$ with $n \leq 5$. Consider the possibilities for $\map \tilde T_1.\tilde S_1.$, 
classified in (11.1). Index consideration shows $r=1$, and we are exactly one blow up beyond 
$c=(2)$ and $d=(3)$. There are two possibilities, $c=(\ul{2},2')$, $d=(4)$, $a=6/7$ or 
$c=(3)$, $d=(\ul{3},2')$, $a=7/8$. But if $x_0$ has centre a $-2$-curve and branches $(2)$, 
$(2,2)$ and $(4)$, then $e_0=6/7$ and so the former case is impossible. $e_0\geq 16/17$, in 
the latter case, again  impossible. \qed\enddemo

 Thus from now on we will assume $x_0$ is a chain singularity. In particular
$k\geq 3$, and there are exactly two singularities along $E_1$, $z$ and $A_r$. 

\heading 17.5 Proving $S_0^0$ is uniruled  \endheading

We now classify possibilities for $S_0$ leading to $(S_2,A_2+B_2)$ given in (17.3.1) and 
(17.3.3). In each case we show that $S_0$ is uniruled using the methods of \S 6. This 
involves finding certain rational curves on $S_0$. One candidate, by (17.2.2) is $\sb$. The 
others we will use are given by the following lemmas. We encourage the reader to first go 
through example (6.10), as the analysis to follow is very similar, but less detailed.

The notation for the rational curves below is fixed for the remainder of \S 17.

\proclaim{17.5.1 Lemma} Let $A \subset \pr 2.$ be a nodal cubic, with node at $q$ and a flex 
point at $p$. Let $r\neq p$ be the unique point of $A^0$, whose tangent line meets $p$. 
\roster
\item There is a nodal cubic $R$ through $p$ and $q$, smooth
at $p$ and $q$, with seventh order contact with $A$ at $p$.
\item There is a nodal cubic $\ol R$ through $p$, $q$, and $r$, smooth at 
these points, with forth order contact with $A$ at $p$, and third order contact
with $A$ at $r$. 
\endroster
\endproclaim
\demo{Proof} We prove the first case, and the leave the second to the reader, which is 
analogous. Let $S=S(A_1+E_7)$. Blow up once at $q$ and seven times 
at $p$ along $A$, to obtain $\t{S}$. It is enough to show there is a nodal 
elliptic curve $R \subset S^0$. For this we study the elliptic pencil 
$|-K_S|$. \cite{HW81} implies there is a smooth member $C \in |-K_S|$, 
with $C \subset S^0$, and if $V$ is the $-1$-curve over $q$ then 
$V\in |-K_S|$. 

Thus we have an elliptic fibration $\pi :\map T.{\pr 1.}.$ after blowing up 
the basepoint of $|-K_S|$. Then by (16.4) of \cite{BPV84} we have 
$$
e(\t{T}) = e(\pr 2.) + 9 = e(I_2) + e(\t{E}_7) + n
$$
where $e$ is the Euler characteristic, $n$ is the contribution to the sum from
the other singular fibres and $\t{E}_7$, $I_2$ are as on page 150 of 
\cite{BPV84} (The fibre through the $A_1$ point contains $V$ and so we get 
the graph $I_2$). Thus $n=1$ and one fibre of $\pi$ is a nodal cubic 
$R$.\qed\enddemo
\remark{17.5.1.1 Remark} By a dimension count one expects in
(17.5.1.1) a one dimensional family of elliptic curves passing through
$p$ and $q$ meeting $A$ to order $7$ at $p$. However finding a 
{\it rational} elliptic curve in the family, smooth at $p$ and $q$,
is more subtle. For example the proof of (17.5.1) shows that there is
no such curve if $R$ has an ordinary cusp. \endremark

\proclaim{17.5.2 Lemma} Let $A$, $C \subset \pr 2.$ be a nodal
cubic, with node $q$, and a smooth conic, meeting $A$ to order six at a point
$p\neq q$. Then there is a nodal cubic $R'$ which meets $A$ at $p$ to order 
seven. \endproclaim
\demo{Proof} If we blow up 3 times at $p$ along $A$ and then blow down the
tangent line and the two $-2$-curves, then $C$ becomes the flex line to $p$,
and the result follows from (17.5.1.2). \qed \enddemo

\proclaim{17.5.3 Lemma} Let $A$, $C \subset \pr 2.$ be a nodal cubic, with 
node $q$, and a smooth conic, meeting to order six at a point
$p\neq q$. Then there is a rational quartic $R''$, with an oscnode at $p$
(that is two smooth branches meeting to order $3$), such that one branch of 
$R''$ at $p$ meets $A$ to order $5$, and the other meets $A$ to order $3$,
and such that $R''$ meets one branch of $A$ at $q$ to order $3$. 
\endproclaim
\demo{Proof} After a birational transformation (cf. the proof of (17.5.2))
$p$ becomes a flex point and the problem is to find a smooth conic $Q$ 
tangent to $A$ at $p$ with the described contact at $q$. As all pairs $(A,C)$ 
are projectively isomorphic, it is enough to construct $Q$ for some pair. 

Begin with a smooth conic $E \subset \pr 2.$, two tangents $C$ and $ G$ at 
$c$ and $g$, and a secant line $A$ through $c$, $g \not \in A$. Let 
$w:\map T.{\pr 2.}.$ be obtained by blowing up $3$ times at $c$ along $C$, 
and $3$ times at $E$ along $g$. Let $Q$ be the $-1$-curve over $g$, and $V$ 
the $-1$-curve over $c$. Let $h:\map T.{\pr 2.}.$ be given by contracting all 
the $w$-exceptional curves except $Q$ and $V$, together with $E$ and $G$ 
(which are both $-1$-curves on $T$). $Q$ is the required conic. \qed\enddemo

\proclaim{17.5.4 Definition} Suppose $(S_2,A_2+B_2)$ is given by (13.5.1). Let $Z,W \subset S$ be 
the strict transforms of the tangent lines to the two branches $C$, $D$, of (17.1.0), of 
$A \subset \pr 2.$ at the node $q$. Note $W$ has a unibranch singularity at $x_0$ of form 
$(A_{r-1},2',k',\dots)$ (on $\t{T}_1$, $\t{W}$ meets the intersection of $\t{E}_1$ with the first 
$-2$-curve of the $A_r$ point). Note if $r \geq 2$, $Z$ is smooth at $x_0$, and has form 
$(A_{r-2},2',2,k,\dots)$, while if $r =1$, $Z$ meets $\Sigma_1$ at a smooth point, and does not 
pass through $x_0$. 
\endproclaim

Note $-k + r + 4 = \t{A}_1^2$. We have

$$
0 > (K + a A_1) \cdot A_1 = \frac{s-1}{s} - \frac{A_1^2}{r+2} \tag{17.6 }
$$

Suppose $z$ has type $\Gamma_n$. Then $A_1^2 = \t{A}_1^2 + \frac{s-1-n}{s}$
and one finds from (17.6) 
$$
\t{A}_1^2 > (s-1) (k-3) + n \tag{17.6.1}
$$

In 17.7-12 we classify possibilities where $x_1 \in A_1$. (17.3.1) 
describes the configurations on $S_2$ and $T_2$. We consider each possibility 
for $S_2$ in turn. In any case $s \geq 3$. Note $\t{A}_1^2 = K_{S_2}^2$. 
\subhead 17.7 $S_2=S(A_1 + A_2)$ \endsubhead
Then $\t{A}_1^2 =6$, $r = k + 2$. Note $E_2$ has either an $A_1$ point or 
$A_2$ point along it. 
\subhead 17.7.1 If $E_2$ has an $A_2$ point \endsubhead
Thus $z=(j,2,2)$, $y=(A_{j-2},3)$, $x_0=(2,2,j,k,A_{k+2})$.
$s=3j-2$. $n = 3j-6$. By (17.6.1), $j < 4$.
\subhead 17.7.1.1 If $j=3$ \endsubhead
(17.6.1) gives $k=3$. $\Delta(x_0) =73$.
$-K_{S_0} \cdot \sb = 3/73$. Log uniruled
by (6.5). \qed
\subhead 17.7.1.2 If $j=2$ \endsubhead
(17.6.1) gives $k=3$, $4$, $\Delta(x_0) =34,67$.
One checks the smooth branch of $\sb$ has index $17$ 
in the first case. $-K_{S_0} \cdot \sb = 6/34,3/67$.
Log uniruled by (6.5). \qed

\subhead 17.7.2 If $E_2$ has an $A_1$ point \endsubhead
Thus $z=(j,2)$, $y=(A_{j-2},3,2)$, $x_0=(2,j,k,A_{k+2})$. $s= 2j-1$,
$n=2j-4$. (17.6.1) gives $j < 5$.
\subhead 17.7.2.1 If $j=4$ \endsubhead
(17.6.1) gives $k=3$. $\Delta(x_0) =79$. 
$-K_{S_0} \cdot \sb = 2/79$. \qed
\subhead 17.7.2.2 If $j=3$ \endsubhead
(17.6.1) gives $k=3$. $\Delta(x_0) = 53$. 
$-K_{S_0} \cdot \sb = 4/53$. \qed
\subhead 17.7.2.3 If $j=2$ \endsubhead
(17.6.1) gives $k=3,4,5$. $\Delta(x_0)=27$, $52$, $83$. Note the
two branches of $Z$ at $x_0$ are
$(2',j,k,A_{k+2})$ and $(2,j,k,2,2',A_k)$.
One checks the smooth branch of $Z$ has index
$27$, $52$ in the first two cases.
$-K_{S_0} \cdot Z = 9/27$, $8/52$, $5/83$. \qed

\subhead 17.8 If $S_2 = S(A_1 + A_5)$ \endsubhead
Then $\t{A}_1^2 =3$, $r = k -1$. 
Note $E_2$ has either an $A_1$ point or $A_5$ point along
it. 
\subhead 17.8.1 If $E_2$ has an $A_5$ point \endsubhead
Thus $z=(j,A_5)$, $y=(A_{j-2},3)$, $x_0=(A_5,j,k,A_{k-1})$.
$s=6j-5$, $n=6j-12$. (17.6.1) gives $j=2$, $k=3$. 
$\Delta(x_0)=31$, $-K_{S_0} \cdot \sb = 3/31$, so
log uniruled by (6.5). \qed

\subhead 17.8.2 If $E_2$ has an $A_1$ point \endsubhead
$z=(j,2)$, $y=(A_{j-2},3,A_4)$, $x_0 =(2,j,k,A_{k-1})$.
$s=2j -1$, $n=2j-4$. (17.6.1) gives $j=3,2$.

\subhead 17.8.2.1 If $j=3$ \endsubhead
(17.6.1) gives $k=3$. $\Delta(x_0) =29$, $\sb$ has two smooth branches at $x_0$. 
$R''$, (of 17.5.3) has two smooth branches at $x_0$.  Thus log uniruled by (6.5), using
$R''$ and $\sb$.
\subhead 17.8.2.2 If $j = 2$ \endsubhead
(17.6.1) gives $k=3,4$. For $k=3$, $\Delta(y) =11$.
$R''$ has two smooth branches at $y$, $(A_4,3')$ 
and $(A_2,2',2.3)$. $-K_{S_0} \cdot R'' = 3/11$. log uniruled
by (6.5). For $k=4$, $\Delta(x_0)=31$. $R'$ has one smooth, and
one singular branch at $x_0$. The singular branch is
$(2,2,4',2',A_2)$, which one checks has index $31$. 
Thus (with $R'$ and $\sb$) log uniruled by (6.5). \qed

\subhead 17.9 If $S_2=S(3A_2)$ \endsubhead
$\t{A}_1^2 =3$. $r = k-1$. 
$y=(2,3,A_{j-2})$, $x_0 =(A_2,j,k,A_{k-1})$. 
$s = 3j-2$, $A_1^2 = 3 + \frac{3}{3j-2}$. Then (17.6)  gives
$j=2$, $k=3$. $\Delta(x_0)=19$, and $-K_{S_0} \cdot \sb = 3/19$,
log uniruled by (6.5). \qed 

\subhead 17.10 If $S_2 = S(A_2 + A_5)$ \endsubhead
$\t{A}_1^2 =2$. Thus by (17.6.1), $k=3$, $r = 1$,
$j=2$, and $z$ is Du Val .

Note $E_2$ contains either an $A_2$ or an $A_5$ point.
\subhead 17.10.1 If $E_2$ contains an $A_2$ point \endsubhead
Thus $y=(A_4,3)$ and $x_0=(2,3,A_3)$. $\Delta(x_0) =14$. One
checks the smooth branch of $\sb$ has index 14, and 
$-K_{S_0} \cdot \sb =2/14$, thus log uniruled by (6.5). \qed

\subhead 17.10.2 If $E_2$ contains an $A_5$ point \endsubhead
$y=(2,3)$, $x_0=(2,3,A_6)$. $\Delta(x_0) = 23$. 
$-K_{S_0} \cdot \sb =2/23$, thus log uniruled by (6.5). \qed

\subhead 17.11 If $S_2=S(A_1 + 2A_3)$ \endsubhead
By (13.5.5) $E_2$ contains an $A_3$ point. As in (17.10) $\t{A}_1^2 =2$ implies $k=3$, $r=1$, 
$j=2$. We have $y=(A_2,3)$ and $x_0=(2,3,A_4)$. $\Delta(x_0) =17$, $-K_{S_0} \cdot \sb = 2/17$, 
so log uniruled by (6.5). \qed

\subhead 17.12 If $S_2 = S(2A_1 + A_3)$ \endsubhead
Note $E_2$ contains either an $A_3$ or an $A_1$ point. We have $\t{A}_1^2 = 4$. $r=k \leq 4$ by 
(17.6.1). Suppose $k=4$. (17.6.1) gives $n=0$, and $5>s$, thus $j=2$, and $E_2$ contains
an $A_1$ point, $y=(A_2,3)$, $x_0=(A_4,4,A_2)$, $\Delta(x_0) =38$, one checks the smooth branch 
of $\sb$ had index $38$ and $-K_{S_0} \cdot \sb = 2/38$, so log uniruled by (6.5). 

Now assume $k=r=3$. 

\subhead 17.12.1 If $E_2$ contains an $A_3$ point  \endsubhead
By (17.6.1), $n<4$ implies $j=2$. Thus $x_0=(A_3,3,A_4)$,
$\Delta(x_0) =29$, $-K_{S_0} \cdot \sb = 4/29$, log uniruled 
by (6.5). \qed
\subhead 17.12.2 If $E_2$ contains an $A_1$ point \endsubhead
$z=(\ul{j},2)$, which has type $2j-4$. Thus (17.6.1) gives
$j=2,3$. For $j=2$, $x_0=(A_3,3,A_2)$, $\Delta(x_0) =19$,
$-K_{S_0} \cdot \sb = 4/19$. For $j=3$, $x_0=(A_3,3^A,3,2)$,
$\Delta(x_0)=37$, $-K_{S_0} \cdot \sb = 2/37$. \qed

In 17.13-14 we classify possibilities where  $x_1 \not \in A_1$.
Thus we have the setup of (17.3.3). We divide the analysis into
two cases depending on whether or not $\Sigma_2 \cap E_2$ is 
a smooth point. Note $y = x_1$. 

\subhead 17.13   If $\Sigma_2$ meets $E_2$ at a smooth point \endsubhead
So $z$ is a Du Val point, $A_t$, $j=t+2$, $\t{A}_1^2=5$,
$r=k+1$,$x_0=(A_t,k,A_{k+1})$, $x_1=(2,t+2,2,2)$. $s = t+1$,
so (17.6.1) gives

$$5> t(k-3) \tag{17.13.1}$$

In particular $k < 8$.

\subhead 17.13.2 If $k=3$ \endsubhead
One checks $t=1$ (or $e_1 > e_0$). $x_0=(2,3,A_4)$, $\Delta(x_0) =17$.
$-K_{S_0} \cdot \sb = 5/17$. $S_0$ is log uniruled by
(6.5). \qed
\subhead 17.13.3 If $k=4$ \endsubhead
16.16.1 give $t=1,2,3,4$. $\Delta(x_0) = 32,45,58,71$.
One checks in each case that the singular branch of $R$, which
is $(A_{t-1},2',4^3,A_5)$ has index $\Delta(x_0)$, and that
the smooth branch, which is $(A_t,4,2',A_4)$, has index $\Delta(x_0)$
except in the case $t=2$ when the index is $9$. One computes
$-K_{S_0} \cdot R = \frac{6(5-t)}{13t + 19}$. Thus $S_0$ is
log uniruled by (6.5) in every case. \qed
\subhead 17.13.4 If $k=5$ \endsubhead
16.16.1 gives $t=1,2$. $\Delta = 51,73$.
One checks in the first case that the smooth
branch $\sb$ has index $51$, and $-K_{S_0} \cdot \sb = 3/51$.
In the second case $R$ has two branches at $x_0$, one
smooth and one singular. One checks the singular branch has
index $73$. The singular branch is $(A_6,5^3,2',2)$.
 The smooth branch is $(A_5,2',5,2,2)$. One
computes $-K_{S_0} \cdot R = 7/73$. So log uniruled
by (6.5). \qed
\subhead 17.13.6 If $k=6$ \endsubhead
$t=1$. $\Delta(x_0) =74$. One checks both branches
of $R$ have index $74$. $-K_{S_0} \cdot R =16/74$. \qed
\subhead 17.13.7 If $k=7$ \endsubhead
$t=1$. $\Delta(x_0) = 101$. One checks the
singular branch of $R$ has index $101$.
$-K_{S_0} \cdot R = 9/101$. \qed

\subhead 17.14 If $\Sigma_2$ meets $E_2$ at a singular point \endsubhead
Let $w = E_2 \cap \Sigma_2$. $x_1$ is necessarily a 
non-chain point, with branches $w,(2),(2,2)$. From an easy classification 
(since $x_1$ is not Du Val), $e_1 \geq 2/3$. Consider 
$h:\map \t{T}_2.\t{S}_2.$. 
\subhead 17.14.1 If the second blow up of $h$ is along $A$ \endsubhead
Note  $j=2$. The only possibility for $w$ (with $e_1 < 6/7$, and $x_1$ non 
Du Val) is $(3)$. Then $z=(2)$,$\t{A}_1^2 =4$. $R$ (on $\t{S}_0$) is a 
$3$-curve. Its singular branch at $x_0$ is $(2',k^2,A_k)$ and its smooth
branch is $(2,k,2',A_{k-1})$. Also after scaling $b = 2 -3/2 a$. Then
from the configuration in $\pr 2.$ it follows that $a < 8/9$.
Thus $r=k \leq 6$. If $k=3$ then $e_0 < e_1=2/3$. Thus
$k=4$, $5$, $6$. $\Delta(x_0) = 27$, $44$, $65$. For $k=4$ one
checks the smooth branch of $\sb$ has index $27$, and
$-K_{S_0} \cdot \sb =3/27$. For $k=5$ one checks
both branches of $R$ have index $44$ and
$-K_{S_0} \cdot R = 12/44$. For $k=6$ the singular branch
of $R$ has index $13$, the smooth branch has index $65$,
while $-K_{S_0} \cdot R = 7/65 > 1/13 + 1/65$. So (6.5) applies
in each case. \qed

\subhead 17.14.2 If the second blow up of $h$ is along $E_2$ \endsubhead
Note $j \geq 3$. If $j \geq 4$ then $e_1 \geq 6/7$. Also if $j \geq 3$ and
$w$ has index at least $3$, $e_1 \geq 6/7$. We conclude
$w=(2)$, $j=3$, $z=(3)$, $\t{A}_1^2=5$, $r = k+1$, $e_1 =3/4$,
$x_0=(3,k,A_{k+1})$.
By (17.6.1) $k \leq 4$. If $k=3$ then
$e_0 = 21/28 < 3/4$, a contradiction. So $k=4$. $\Delta(x_0) =51$. 
$R$ has branches $(3',4^3,A_5)$ and $(3,4,2',A_4)$. One checks each
has index $51$. $-K_{S_0} \cdot R = 4/17$. \qed

\heading \S 18 $A_1$ smooth \endheading
In this and the following section we prove

\proclaim{18.1 Proposition} If $A_1$ is smooth, then either $S$ has a tiger, or
$(S_2,A_2 + B_2)$ is a smooth banana, classified in (13.2). In either case $S_0^0$ is uniruled.
\endproclaim

{\bf We assume throughout \S 18--19 that $S$ has no tiger.}

In this section we will consider all possibilities for the next hunt step, and rule out all 
but the smooth banana. The smooth banana will be considered in \S 19

Recall from (8.2.10) that $E_1$ is a $-k$-curve, and $E_2$ is a $-j$-curve. We first put 
together some elementary observations.

$K + a A$ is flush and $K_{S_1} + A$ is lt by (8.4.7).

\proclaim {18.2 Lemma} The following hold:
\roster
\item $A$ is ample and $\t{A}^2 \geq -1$.
\item $A$ contains exactly three singularities.
\item There is a singularity of index at least four along $A$. 
\item $a >2/3$, $e_1 > 1/2$.
\endroster
\endproclaim
\demo{Proof} (1) follows, since $K_{S_1}\cdot A<0$. 

Next we consider (2). As $A$ is not a tiger, $A$ must contain at least three singularities,
by adjunction. Suppose it contains at least $4$. Since $E$ can contain at most three, 
and the only additional singularity can be at $q_1$, we conclude $x_0$ is a non-chain 
singularity, $\Sigma \cap E$ is a smooth point, and there are exactly four
singularities along $A$. By adjunction one has index at least three.

$e_0\geq 1/2$ by (10.1). $e_1 \geq 2/3 a$ (since there is a singularity
along $A$ of index at least three), by (8.3.8). If $e_0 \geq 2/3$, then 
$a_2 + 3b_2 >2$. If $e_0 < 2/3$, then $x_0$ is given by (10.1.d), and
$A$ contains a non-Du Val singularity, so $e_0 \geq e_1 \geq 1/2$ by (8.0.7). Thus
again $a_2 + 3b_2 > 2$. In particular $a_2 + b_2 > 1$.

Consider the next step of the hunt; the possibilities are given by (8.4.7).

Suppose $T_2$ is a net. $A$ is not fibral by (11.5.5). $A + E$ has
degree at least three, for otherwise $-(K_T + A + E)$ is nef and we
have a tiger. Since $a_2 + 3b_2 > 2$, it follows that one of $A$ or
$E$ is a section, and the other is a double section. This contradicts
(11.5.11.2). Thus $\pi_2$ is birational. 

The possibilities for $(S_2,A_2+ B_2)$ are given by (8.4.7). $A_2$ must contain at
least two singularities. Thus $B$ is not smooth, by (13.1). Hence we have a fence. 
$A_2$ contains three singularities, thus by (13.3), $S_2 = S(A_1 + A_2 + A_5)$, 
$B \subset S_2^0$, $B$ has arithmetic genus one, and a simple node, and $B^2 =1$.
By (11.1.1), $T_2$ is smooth along $\Sigma_2 \setminus \Sigma_2 \cap E_2$, and
$\Sigma \cap B$ is an $A_r$ point, and 
$$
1 = B_2^2 = r +1 - j.
$$ 
Thus $r = j$. $x_1 = (\ul{j},A_j)$. $S_1$ is smooth away from $A_1$, and its 
singularities are exactly $x_1,A_1,A_2,A_5$. $S_1$ cannot be Gorenstein 
(since $\t{S_1}$ has Picard number at least nine). One of the singular points is $q_1$. 
This cannot be $x_1$(or $x_0$ has branches of index $(2,3,6)$). 
Thus $x_1$ has index at most five, either $x_1=(3)$, $x_1=(4)$ or $x_1=(3,2)$. This 
contradicts our expression $x_1=(j,A_j)$. Hence (2).

(3) follows immediately from (2) and adjunction.

Now we prove (4). The second inequality follows from the first, and (3) by (8.3.8).

Suppose first $x_0$ is a chain singularity. Then $\Sigma$
meets $A$ at a smooth point and $q_1$ is a singular point. The configuration
on $T$ is of form $(n^s)$ of (11.4). There is a unique singularity,
$y$, along $\Sigma$, and $K_T + \Sigma$ is lt. $k \geq 3$, so the
self-intersection of $E$ goes up by at least two, and thus $y$ is of the form 
$(2',s,\dots)$ (for some $s\geq 3$) where $x_1$ has type $(\ul{s-1},\dots)$. 
Let $g$ be the coefficient of the $-s$-curve. Then the
marked $-2$-curve has coefficient $g/2$, and 
$a+g/2=1$. $K_T + a E$ is flush by (8.4.7), thus $a > g$, so $a>2/3$. 

Now assume $x_0$ is not a chain singularity, and
$e_0 < 2/3$. Then by (10.1.1.d), two of the singular points along $E$ are index 
two, $E$ is a $-2$-curve, and $e=1/2$. We may assume $\Sigma \cap E$ is one of 
the $(2)$ points, for otherwise the indices along $A$ are $(2,2,m)$ for 
some $m$, and we have a tiger by adjunction. The configuration on $T$ is 
$(n^s,f)$ of (11.4), for some $s \geq 1$ (in no other configuration does
$\Sigma$ meet $E$ at an $A_1$ point). There is one singularity
$y$ along $\Sigma \setminus A \cap \Sigma$.
Let $g$ be the coefficient of the curve over $y$ which meets 
$\Sigma$. We have $a/2+g=1$. Arguing as before, $a > 2/3$. 
Hence (4). \qed\enddemo

Define $t>0$ so that $K_{S_1} + tA$ is numerically trivial.
$K_{S_1} + t A$ is flush by (8.3.2.1) and (8.3.5.2).

\proclaim{18.3 Lemma} Suppose $A$ is a $(l-2)$-curve, and $A$ has only Du Val
singularities along $A$. If $\epsilon=(K_S+A)\cdot A$, then
\roster
\item $t=l/(l+\epsilon)$, $\epsilon=l(1-t)/t = A^2 -l$
\item $S_1$ is not Gorenstein.
\endroster
\endproclaim
\demo{Proof} We have 
$$
\align
0=(K_S+tA)\cdot A&=K_S\cdot A+tA^2 \\
                 &=-l+t(l+\epsilon)=-l+lt+t\epsilon
\endalign
$$
(1) follows. 

Suppose $S_1$ is Gorenstein. By (18.2.2), $A$ contains three Du Val chain singularities, say of 
index $s+1$, $v+1$, $r+1$. By adjunction 
$$
\frac{s}{s+1} + \frac{v}{v+1} + \frac{r}{r+1} > 2.
$$
One checks from the list (3.1) that no rank one Gorenstein log del Pezzo contains three such 
singularities. \qed\enddemo

\remark{18.4 Remark} Since $a_1 + e_1 >1$, by (8.4.7) one of the following occurs.
\roster
\item"(a)" $(S_2, A_2 + B_2)$ is a tacnode.
\item"(b)" $S_2$ is a fence.
\item"(c)" $S_2$ is a banana.
\item"(d)" $A$ is contracted, that is $A =\Sigma_2$.
\item"(e)" $T_2$ is a net.
\endroster

Furthermore $K_{T_2} + \Gamma'$ is flush, and $K_{S_2} + a_2 A_2 + b_2B_2$
is flush in cases (a-c). 
\endremark

\proclaim{18.5 Lemma} (18.4.a) never holds. 
\endproclaim
\demo{Proof} Suppose (18.4.a) holds. By (11.1.1) the configuration on $T_2$ is type $II$,
and $2a_2+3b_2\leq 3$.

 Suppose $S_1$ has a non Du Val Point along $A$. Then by (8.3.8)
$$
b_2 > e_1\geq (1/r)+(r-2/r)(2/3)
$$ 
and $r>2$, where $r$ is the index of this point. 
Thus 
$$
2a_2+3b_2> 2(2/3)+(2r-1)/r=2+4/3-1/r\geq 3,
$$
a contradiction.

So $S_1$ is Du Val along $A$. Let $w = e(E_2,K_{S_1} + t A)$. Then $t> w > e_1$, and 
$K_{T_2} + t A + w E$ is $\pi_2$-trivial, so $2t + 3w \leq 3$ as in (11.1.1.2). Thus 
$t < 3/4$ by (18.2.4). By (18.3.1), $\epsilon > l/3 \geq 1/3$. Applying the Bogomolov bound, 
we deduce that either $S_1$ is smooth away from $A$, or it has exactly one singular point 
away from $A$, a point of index $2$ (the contribution to the sum in (9.2) for the 
singularities away from $A$ is less than $2/3$). But then $S_1$ is Gorenstein, contradicting 
(18.3.2). \qed\enddemo

\proclaim{18.6 Lemma} (18.4.b) does not hold.
\endproclaim
\demo{Proof} Suppose (18.4.b) holds. By (13.3), $B_2$ is singular, and $A_2$ is a $-1$-curve.
Since $B_2$ is singular, $x_1 \in A_1$. By (8.0.7.1) and (8.0.4), there is at most one 
singularity of $S_2$ along $B_2$. By (13.4) and (18.2.4), if $B$ has genus one, then $B$ has 
a node.

 Suppose $A_2$ lies in the Du Val locus. Then $(K + B) \cdot A_2 =0$, so $B\subset S_2^0$, 
$B$ has a genus one, and $S_2$ is Gorenstein. By (11.1.1) we have a contraction of type 
(II, $x^{r-1})$. As $S_1$ is not Gorenstein, by (18.3.2), $x_1$ has type $(j,A_r)$, 
$j\geq 3$, which has spectral value at least $r+1$. By (8.0.7) this contradicts (11.1.1.3). 

Thus we may assume $A_2$ has at least one non Du Val point along it. Thus $e_1 > 1/2$ by 
(8.0.7). If $B$ has genus one, then it has a node and does not lie in the smooth locus (or 
$S_2$ is Du Val). Thus $\Sigma_2$ meets $E$ at smooth points, and we have a node of type I. 
This contradicts (11.1.1.1).

Thus the genus of $B$ is $g \geq 2$. By (11.1.1.5) and (11.1.1.12), and (8.0.7) the spectral 
value of $\Delta_1$ is at most one. Thus by (8.0.8), $A_1$ has Du Val or almost Du Val
points along it. We consider the possible types of contraction, see (11.1) and (11.2). 

 By (18.2.4), $e_1 > 1/2$, so we cannot have a contraction of type I. Suppose we have one of 
type II. Then 
$$
e_1< b_2 \leq (g+1)/(2g+1)\leq 2/3
$$ 
and so the point we blew up must be an $A_{g+1}$ (otherwise $x_1$ has spectral value at least
two). But $A$ has a singularity with spectral value one. Since we still decided to blow up 
the $A_{g+1}$ point, we have that $a$ is at least the solution to the following equation in 
$x$
$$
1/3+x/3=x(g+1)/(g+2).
$$

Thus $a \geq (g+2)/(2g+1)$ and $e_1 \geq (g+1)/(2g+1)$, a contradiction.

It follows by (11.1.1) that $B_2$ has a unibranch singularity. Now type III is 
impossible, as $S_1$ would violate the Bogomolov bound (9.2).

 This only leaves types $u$ and $v$. 

In type $u$, $x_1 = (\ul{j},2,2)$, and there is in addition a point $(4,2)$ (index $7$) away 
from $A_1$. By spectral value, (8.0.8), $j=2$.

In type $v$ there is a $(2,3,2,2)$ point (index $11$) away from $A_1$.

Now consider the possible indices of 
the singular points along $A_1$. Let $p$ be the largest index of a non Du Val 
point along $A_2$ (we have already noted there is at least one non Du Val
point along $A_2$), let $q$ be 
the index of the other point on $A_2$.  By (8.3.8) the index of 
$x_1$ is at least $p$ (or the hunt would not choose $x_1$).

 Suppose $q$ is at least three. Then we have four singularities of index at 
least three on $S_1$ and so we may use the list (10.5).
Since we have a point off $A_1$ of index at least $7$, it
follows from the list that $p$ and $q$ both have index exactly three. But then
since $A_2$ is a $-1$-curve, $A^2_2\leq 0$, a contradiction.

 Thus $q$ has index two. As $A_2$ is not contractible, $p$, and hence
$x_1$, have index at least five. Thus we have configuration $u$ (in
$v$ the index of $x_1$ is four) and $S_1$ has singularities of index
at least $2,5,5,11$, which violates (9.2). \qed\enddemo

\proclaim{18.7 Lemma} $A$ is not contracted. 
\endproclaim
\demo{Proof} Suppose $A$ is contracted, in other words $A=\Sigma_2$. 

 By (8.3.1.2) and (8.3.5.1), $q_2$ is a smooth point of $S_2$,
$B_2$ has a unibranch singularity and $\Sigma_2$ is the only divisor whose 
coefficient in $K_{S_2}+bB$ is greater than $b$. Now $\Sigma_2$ cannot 
be in the Du Val locus, or $B_2$ would be smooth. Thus $e_1 > 1/2$. 

 Suppose $B_2$ has a triple point. Note that the first exceptional divisor
lying over $q_2$ has coefficient $3b-1=(2b-1)+b>b$. Thus $\Sigma_2$ must be
this divisor, and in particular $\Sigma_2$ lies in the smooth locus, 
a contradiction. 

 Thus $B_2$ has multiplicity two at $q_2$. Since $T_2$ has exactly two 
singularities along $\Sigma_2$, and meets $E_2$ at a smooth point, we must 
have type III of (11.2). 

 Suppose $g\geq 2$. Then $b>2/3$, since $A_1$ has spectral value at least
two. But if we blow up $q_2$ twice then the resulting exceptional divisor, 
which is not $\Sigma_2$, has coefficient $2(2b-1)=b+(3b-2)>b$, a 
contradiction. 

 Thus $g=1$. Exactly as in the proof of (16.3.1), we see that $S$ has a 
tiger, a contradiction. \qed\enddemo

\proclaim{18.8 Lemma} $\pi _2$ is not a net.
\endproclaim
\demo{Proof} Suppose $\pi _2$ is a net. Let $d$ be the degree of $E_2$ over $\pr 1.$. 

Suppose $A$ is a fibre of $\pi_2$. $d >2$ (or $K_{T_2}+A+E_2$ is anti-nef) and so $e_1<2/3$. 
Thus $A$ and $E$ have singularities of spectral value at most one. $A$ contains two singular 
points of $T_2$. Since $e_1 \geq 1/2$, $d=3$. It follows that $A$ is given by (11.5.9.5) 
with $k=2$. Thus by (18.2.3) $x_1$ has index at least $4$. Since the spectral value is at most 
$1$, $E$ must be a $-2$-curve, and contain a unique singularity $z$, either Du Val or almost 
Du Val. Let $h:\map \t{T}_2.W.$ be a $K_T$ relative minimal model, an isomorphism at the generic 
point of the $-3$-curve meeting $A$. $h$ is determined by selecting which reduced component of 
each multiple fibre becomes the fibre of $W$, see (11.5.5).

$A$ and $G$, the fibre through $z$, are the only multiple fibres: $\map E.{\pr 1.}.$ is 
necessarily ramified over the image of any multiple fibre, in particular it is ramified at $z$, 
and totally ramified at $E \cap A$. If there were another multiple fibre, then it would meet
$E$ at a smooth point, and so $E$ would be totally ramified there as well, a contradiction. 

We consider in turn the two possibilities, either $E$ meets $G$ only
at $z$, or at $z$ and also at some smooth point $y$.

In the first case $E$ meets a curve over $z$ of multiplicity $3$ in the fibre. From (11.5.9) 
the possibilities for $G$ are (5), with $k=3$, or (6). Choose $h$ by selecting in the first 
case the $-4$-curve, and in the second the $-3$-curve over $z$. Then $h(E) \subset W$ is a 
smooth rational triple section, and $h(E)^2 =4$, contradicting (11.5.10).
 
In the second case, necessarily $m(G) = 2$. The only possibility is (11.5.9.1) (in (5) with 
$k=1$, $x_1$ has index $3$). Choose $h$ by selecting the $-2$-curve meeting $E$ at $q$. 
$h(E) \subset W$ is a smooth triple section of self-intersection $4$, contradicting (11.5.10).

Thus  $A$ dominates $\pr 1.$. Let $\delta$ be the degree of $A$ over $\pr 1.$. $d+\delta>2$ 
(otherwise $K_{T_2} + A + E$ is anti-nef and we have a tiger). By (18.2.4), one of $E$ or 
$A$ is a section and the other is a double section. Thus $x_1 \in A$ by (11.5.11.2). The 
spectral value of $x_1$ is either one or zero. 

Suppose $E$ is a section. Each singular point of $A$ lies on a different multiple fibre, by 
(11.5.5), since each must meet $E$ at a singular point. Thus there are at least two multiple 
fibres, and so at least two singularities along $E$, a contradiction. Thus $A$ is a section, 
and $E$ is a double section. There are exactly two multiple fibres by (11.5.11.1).

Suppose $E \subset T^0$. Then $E$ is a $-3$-curve (by spectral value $j \leq 3$, and the hunt 
would not choose an $A_1$ point as there is a point of index $4$ along $A$, by (18.2.3)). Both 
multiple fibres are of multiplicity two, each given by (11.5.9.1) (if either is (11.5.9.5), with 
$k=1$, then the indices of $S_1$ along $A$ are at most (2,3,4) and $K_{S_1} + A$ is negative, a 
contradiction). Then $h(E) \subset W$, of (11.5.4), is a smooth double section of 
self-intersection $3$, contradicting (11.5.10). Thus there is a unique singular point, $z$, along 
$E$, and $E$ is a $-2$-curve.

Let $F_1$, $F_2$ be the two multiple fibres of (11.5.11). Each meets $A$ at a singular point.
One of the $F_i$, say $F_2$ meets $E$ at $z$. $E$ meets a curve of multiplicity two at $z$. By 
spectral value (with respect to both $A$ and $E$) the only possibility, from (11.5.9) is 
$$
(-\ul{3}) + 3(-1) + 2(-\ul{2}) + (-2).
$$
Thus $x_1 = A_3$, and one other singularity of $A$ is $(3)$. Thus the last 
singularity has index at least $3$. $F_1$ has multiplicity two. Since it
meets $A$ at a point of index at least $3$, the only possibility from 
(11.5.9) is (1), and so the third singularity along $A$ is $(\ul{2},2,2)$. 
Now if we choose $h:\map \t{T}_2.W.$ an isomorphism along $A$, $h(E)$ is a 
smooth double section of self-intersection $3$, violating (11.5.10).
\qed\enddemo

\heading \S 19 The smooth banana \endheading
Here we prove:
\proclaim{19.1 Proposition} If $S_2$ is a smooth banana then $S_0$ is log uniruled.
\endproclaim

{\bf We assume throughout \S 19 that $S_0$ does not have a tiger. }

 Smooth bananas are classified in (13.2), where we find $S_3 = \Hbz 2.$, 
$x_2 = (\ul{r+1},2)$. We will explicitly classify the possibilities. 
Note $\Sigma_1$ meets $E_1$ at a smooth point, and so by (11.1.1) meets $A_1$ 
at a Du Val point, say $A_s$.

$A+B+C \subset {\Hbz 2.}$ 
is a configuration of two sections and a fibre, thus since 
$K_{\Hbz 2.}+aA+bB+cC$ is negative
$$
2a_2 + 2b_2 + c < 4. \tag *
$$
In particular $e_2 < 4/5$, so by (8.0.7.2), $x_2$ 
has spectral value at most three. Thus $r \leq 2$.

\proclaim{19.2 Lemma} There are two possibilities:
\newline
{\bf Case 1:} $x_2 \in B$, the singularities of $S_1$ along $A_1$ are
$x_1= (\ul{s-1},r+1,2)$, $A_s$, $A_r$, $s \geq 3$, and $A_1$ is a $0$-curve. 
$(s,r)$ is either $(3,2)$ or $(4,1)$.
\newline
{\bf  Case 2:} $x_2 \in A$, the singularities along $A_1$ are 
$x_1 =(\ul{s},A_r)$, $A_s$, $(\ul{r+1},2)$, $A_1$ is a $1$-curve ($s \geq 2$).
$(s,r)$ is either $(3,2)$, $(4,1)$ or $(4,3)$.
\endproclaim
\demo{Proof} The descriptions of the singularities follow from (11.1.1) 
and our 
description of $(S_3, A+B+C)$. Thus we need only determine $s$ and $r$.
We compute the spectral value of $x_1$ and use $(*)$:

In Case 1: Suppose first $r=1$. If $s=3$ then $K_{S_1}+A_1$ is numerically 
trivial and we have a tiger. If $s \geq 5$ then $x_1$ has spectral
value at least $6$, and so $a_2 > b_2 \geq 6/7$ and $c \geq 2/3 \cdot 6/7$,
which
contradicts $(*)$. If $s=3$ then $K_{S_1} + A$ is numerically trivial,
and $A$ is a tiger.

Now suppose $r$ is two, we show $s$ is at 
most $3$. If $s \geq 4$ then $x_1$ has spectral value at least $7$, 
$a_2 > b_2 \geq 7/8$, and $c \geq 4/5 \cdot 7/8$ 
which again violates $(*)$. 

In Case 2: Suppose first $r=2$. $s$ cannot be two, (or the hunt would choose 
$x_1 = (3,2)$). If $s\geq 4$ then $x_1$ has spectral value of at least
$6$, contradicting $(*)$ as above. Hence $s =3$. 

If $r=1$ then $s$ can be at most $4$, for otherwise the spectral
value of $x_1$ is at least $6$ which gives a contradiction as before.
If $s=2$ then $K_{S_1} + A_1$ is numerically trivial. \qed\enddemo

We will now classify the possible surfaces $S_0$ and in each
case prove log uniruledness using \S 6.

\subhead 19.3  If $x_0$ is a chain singularity \endsubhead

\proclaim{19.3.1 Lemma} If $x_0$ is a chain singularity, then there is one singularity of 
$T_1$ on $\Sigma_1$, $\Sigma_1$ meets $E_1$ at a smooth point, and $K_{T_1} + \Sigma_1$ is 
log terminal.
\endproclaim
\demo{Proof} Since there are only two singular points along $E_1$, $\Sigma_1$ must meet 
$E_1$ at a smooth point. The rest now follows from the fact that $A_1$ is smooth, see (11.4).
\qed \enddemo

By (19.3.1) there are two singular points $x=x_0$ and $y$ on $S = S_0$. $y$ is the unique 
singular point on $\Sigma=\Sigma_1$. We say that we are `adding' the singularity $y$. Recall 
that $E=E_1$  is a $-k$-curve.

Let $v$ be the coefficient of the curve meeting $\Sigma$ at $y$. Let $e=e_0$. Then 
$K_S \cdot \sb =-1+e+v$. Thus $S$ is a log del Pezzo iff $e+v< 1$.

 We note in all cases there is a map $j:\map \t{S}.{\pr 2.}.$, a composition of smooth blow 
ups. To describe germs of curves on $S$ we will follow the notation of example (6.9).

\subhead 19.3.2  Case 1 \endsubhead
Fix a configuration in $\pr 2.$ of a conic $B$, a secant
line $A$, a tangent line $C$. Let $C \cap B = c$, $A \cap C =t$,
and $A \cap B = \{a,b\}$. Let $L_{wz}$ indicate the line in $\pr 2.$ from $w$ 
to $z$, and $M_p$ the tangent line to the conic at $p$.

\subhead  19.3.2.1 Case 1: $(s,r) = (3,2)$ \endsubhead

\subhead  19.3.2.1.1 Add $(2,3,2)$ \endsubhead
$y=(A_{k-1},3',3,2)$. $e_0 = \frac{12k -24}{12k-17}$
$v = \frac{7}{8k+5}$. $e(3') = \frac{7k}{8k+5}$. For 
$S$ to be del Pezzo, and this to be a step of the hunt we have
$e+v < 1$ and $e \geq e(3')$. The possibilities are 
$k=4,5$. $j$ is given by blowing up 3 times at $c$ along
$C$, and once more over $c$, at the point where the $-1$-curve meets
an exceptional $-2$-curve, 
4 times at $b$ along $B$, and $k$ times at $a$ along
$A$. $x=(2,2^C,k^A,2,2,2)$. $y=(2,3,3^B,A_{k-1})$. Here
for example $k^A$ indicates that the strict transform of the line
$A$ on $\t{S}$ is a $-k$-curve.
$\Delta(x) =31,43$. $\Delta(y)=37,45$. Let $W=L_{ac}$ (on
$S$). $W$ is a $-1$-curve and has two smooth branches at $y$,
they are $(A_{k-2},2',3,3,2)$ and $(A_{k-1},3,3,2')$.
For $k=4$, $-K_S \cdot W = 3/37$, thus $K_S$ is log uniruled
by (6.5). Let $Z=M_b$. $Z$ is a $-1$-curve, and has
two smooth branches at $x$. They are
$(2,2',k,2,2,2)$ and $(2,2,k,2,2',2)$. When $k=5$ one computes
$-K_S \cdot W = 1/45$, $-K_S \cdot Z = 1/43$. $S$ is log
uniruled by (6.5).

\subhead  19.3.2.1.2 Add $A_3$ \endsubhead
$y=(A_{k-1},3,2,2)$. $e_0 = \frac{24k-42}{24k-31}$
$v = \frac{3}{4k+3}$. $k=3,4$. $j:$ Blow up 4 times at $C,c$, 
$4$ times at $B,b$ then $k$ times (above $b$) along $A$.
Let $Z=L_{ac}$. $-K_S \cdot Z = 9/41$. Log
uniruled by (6.5). For $k=4$ let $W=M_a$. $W$ has two branches
at $x$. One checks that each branch
has index 65 =$\Delta(x)$. $-K_S \cdot W = 5/65$,
log uniruled by (6.5).

\subhead  19.3.2.1.3 Add $A_2$ \endsubhead
$y=(A_{k-1},3,2)$. $e = \frac{32k-56}{32k-44}$
$v = \frac{2}{3k+2}$. $e + v <1$ gives $k <4$, $k=2$ not
a hunt step, so $k=3$.  $j:$ Blow up 4 times at $B,b$,
4 times at $C,c$ and then 3 times above $c$ along $A$.
$\Delta(x) =52$. Let $Z=L_{bc}$. $Z$ has two branches 
at $x$, one lt. One checks the non lt branch has index 52.
$-K_S \cdot Z = 4/52$. 

\subhead  19.3.2.2 Case 1 $(s,r)=(4,1)$ \endsubhead

\subhead  19.3.2.2.1 Add $(3,2,2)$ \endsubhead
$y=(A_{k-1},4',2,2)$. $e = \frac{10k-21}{10k-13}$
$v = \frac{6}{7k+3}$. $e + v <1$ gives $k < 9$. 
$e \geq e(4') = \frac{6k}{7k+3}$ gives $k \geq 5$.
$k=5,6,7,8$. Corresponding indices are
$\Delta(x) = 37,47,57,67$. $\Delta(y) = 38,45,52,59$.
$j:$ Blow up 5 times at $B,b$, 3 times at $C,c$, $k$ times
at $A,a$. Let $Z= M_{b}$. $-K_S \cdot Z = \frac{k-9}{10k-13}$.
Let $W=L_{ac}$. $-K_S \cdot W= \frac{k-9}{7k+3}$. $Z$ 
has two branches at $x$, one lt, one smooth; 
$W$ two at $y$, one lt, one smooth. Apply (6.5) to $Z$
for $k \neq 7,8$, and (6.5) to $Z+W$ for $k=8$. 
For $k=7$ let $H$ be the conic of (6.9.1), having
$3^{rd}$ order contact with $B$ at $c$ and tangent to
$A$ at $a$.  $H$ has
two branches at $y$, one singular, one smooth,
and one checks that each has 
local index $52$. $-K_S \cdot H = 4/52$, so we can apply
(6.5). \qed

\subhead  19.3.2.2.2 Add $A_1$ point \endsubhead

$y=(A_{k-1},3)$. $e = \frac{35k-55}{35k-43}$
$v = \frac{1}{2k+1}$. $e + v <1$ gives $k \leq 4$.
$k=3,4$. $\Delta(x) = 62,97$. $\Delta(y)=7,9$.
$j:$ 5 times at $B,b$, 3 times at $C,c$, $k$ times
at $A,t$. For $k=3$ take $H$ smooth conic
tangent to $A$ at $a$ and with $3^{rd}$ order contact
with $B$ at $c$. $H$ has two branches at $x$, 
(each singular) each
of index 62. $-K_S \cdot H = 21/62$. 
For $k=4$ let $Z = L_{bc}$. Two branches at $x$,
one lt, one smooth. $-K_S \cdot Z = 3/97$. \qed

\subhead  19.3.2.2.3 Add $A_4$ point \endsubhead
$y=(A_{k-1},3,2,2,2)$, $e = \frac{14k-22}{14k-13}$,	
$v = \frac{4}{5k+4}$. $e + v <1$ gives $k \leq 7$,
$k=3,4,5,6,7$. $\Delta(x) =29,43,57,71,85$. 
$\Delta(y)=19,24,29,34,39$. $j$: 3 times at $C,c$, 5 times
at $B,b$ then $k$ times over $b$ along $A$. $Z=L_{ac}$,
$W=M_a$. Each has two branches at $x$, one of which is
singular.
$-K_S \cdot Z = \frac{2(8-k)}{14k-13}$.
When $\Delta(x)$ is prime, the singular branch of
$Z$ is non-Cartier by (4.12.2), and we can apply (6.5). One checks 
that the non lt branch of $Z$ has index 57 when $k=5$ and
the non lt branch of $W$ has index 85 when $k=7$. Since
$-K$ is not a generator we can apply (6.5). \qed

\subhead \S 19.3.3 Case 2 \endsubhead
Here fix the configuration as in Case 1, but with
$A$ the conic and $B$ the secant line. 

\subhead 19.3.3.1 $(r,s)=(2,3)$ \endsubhead

\subhead 19.3.3.1.1 Add $(3,2,2)$ point \endsubhead
$y=(A_{k},4,2,2)$, $e = \frac{21k-32}{21k-23}$,	
$v = \frac{6}{7k+10}$. $e + v <1$ gives $k \leq 3$. 
But for $k < 4$ this is not a hunt step (the $-4$-curve has
higher coefficient). So this case does not occur.

\subhead 19.3.3.1.2 Add $(3,2)$ point \endsubhead
$y=(A_{k},4,2)$, $e = \frac{28k-44}{28k-33}$,	
$v = \frac{4}{5k+7}$. $e + v <1$ gives $k \leq 3$.
$k=2$ not a hunt step so $k=3$. $\Delta(x)=51$,
$\Delta(y)=22$.$j:$ 4 Times at $B,b$, $4$ times
at $c$, to get exceptional locus over $c$,
$(-2,-3,-1,-2)$, with $C$ a $-2$-curve, then blow up
4 times above $c$ along $A$. $x=(2,2^C,3^B,3^A,2,2,2)$.
Let $W=M_b$. $W$ has two branches at  $x$,
$(2,2',3,3,2,2,2)$ and $(2,3,3,2,2',2)$, each of
index $51$. $-K_S \cdot W =5/51$. \qed

\subhead 19.3.3.1.3 Add $A_3$ point \endsubhead
$y=(A_{k},3,2,2)$, $e = \frac{35k-41}{35k-29}$,	
$v = \frac{4}{5k+7}$. $e + v <1$ gives $k \leq 2$.
$k=2$ not a hunt step. So this case does not occur. \qed

\subhead 19.3.3.2 Case 2, $(r,s)=(1,4)$ \endsubhead

\subhead 19.3.3.2.1 Add $(4,2)$ point \endsubhead
$y=(A_{k},5,2)$, $e = \frac{15k-30}{15k-22}$,	
$v = \frac{6}{7k+9}$. $e + v <1$ gives $k \leq 5$.
$e \geq e(5) =\frac{6k+6}{7k+9}$ gives $k=5$. 
$\Delta(x) =53$. $j:$ 6 times at $A,a$,3 times at
$C,c$, 5 times at $B,b$. $Z=L_{bc}$. $Z$ has
two branches at $x$, one lt, one smooth.
$-K_S \cdot Z =2/53$. \qed

\subhead 19.3.3.2.2 Add $A_2$ point \endsubhead
$y=(A_{k},3,2)$, $e = \frac{35k-50}{35k-28}$,	
$v = \frac{2}{3k+5}$. $e + v <1$ gives
$k \leq 3$. $k=2$ not hunt step so $k=3$.
$\Delta(x)=67$. $j$: 5 times at $B,b$, 3 times
at $C,c$ then 4 times at $A,c$. $x=(2^C,4^B,3^A,A_4)$.
$W=M_b$. Two branches at $x$, one lt, one smooth.
$-K_S \cdot W = 8/67$. \qed

\subhead 19.3.3.2.3 Add $A_4$ point \endsubhead
$y=(A_{k},3,A_3)$, $e = \frac{21k-30}{21k-21}$,	
$v = \frac{4}{5k+9}$. $e + v <1$ gives $k \leq 4$.
$k=3$ not a hunt step so $k=4$. $\Delta(x) = 64$.
$j$: 3 times at $C,c$, 5 times at $B,b$ then 5 at
$A,b$. $W=L_{ac}$. $W$ has two branches at $x$,
on lt, one singular. The singular branch has index 64.
$-K_S \cdot W = 4/64$. \qed

\subhead 19.3.3.3 Case 2 $(r,s)=(1,3)$ \endsubhead

\subhead 19.3.3.3.1 Add $(3,2)$ point \endsubhead
$y=(A_{k},4,2)$, $e = \frac{12k-24}{12k-17}$,	
$v = \frac{4}{5k+7}$. $e + v <1$ gives $k \leq 8$.
$e \geq e(4)=\frac{4k+4}{5k+7}$ gives $k=4,5,6,7,8$.
$\Delta(x)=31,43,55,67,79$. $j$: 3 times at $C,c$,
$k+1$ times at $A,a$, 4 times at $B,b$. $W=L_{bc}$.
$W$ has two branches at $x$, on lt, one smooth.
$-K_S \cdot W = \frac{9-k}{12k-17}$. One checks when
$k=6$, the smooth branch of $W$ has index $55$. Thus
(6.5) applies when $k \neq 8$. For $k=8$ let
$Z=M_{a}$. $Z$ has two branches at $y$, one lt, other
smooth. $\Delta(y)=47$ (so prime). $-K_S \cdot W =1/47$.
Now apply (6.5). \qed 

\subhead 19.3.3.3.2 Add $A_3$ point \endsubhead
$y=(A_{k},3,2,2)$, $e = \frac{15k-24}{15k-16}$,	
$v = \frac{3}{4k+7}$. $e + v <1$ gives
$k=3,4,5,6,7$. $\Delta(x)=29,44,59,74,89$.
$j$: 3 times at $C,c$, 4 times at $B,b$ then
$k+1$ at $A,b$. $Z=L_{ac}$. $Z$ has two
branches at $x$, one lt, one singular. The
singular is not Cartier by (4.12).
$-K_S \cdot Z = \frac{2(8-k)}{15k-16}$. One checks
for $k=4$, $6$ that the singular branch has index
$44,74$. \qed

\subhead 19.3.3.3.3 Add $A_2$ point \endsubhead
$y=(A_{k},3,2)$, $e = \frac{21k-32}{21k-23}$,	
$v = \frac{2}{3k+5}$. $e + v <1$ gives
$k=3,4,5,6$. $\Delta(x) =37,57,77,97$. 
$j$: 4 at $B,b$, 3 at $C,c$ then $k+1$ at $A,c$.
$Z=M_b$. $Z$ has two branches at $x$, on lt,
one singular. The singular branch is $(2^C,3^B,k',2',2,2)$.
By (4.12) the singular branch is not Cartier.
For $k=4$ on checks the index of the singular
branch is $37,57$. $-K_S \cdot Z =\frac{3(7-k)}{21k-23}$.
Thus (6.5) applies for $k \neq 5$. For $k=5$ let
$W=M_a$. $W$ has one lt branch, one singular branch,
at $x$. One checks the singular branch has index $77$.
Since $K_S$ is not a generator, (6.5) applies. \qed

\subhead \S 19.4 If $x_0$ is a non-chain singularity \endsubhead

Thus there are three singularities along $E$, and along
$A_1$. If $\Sigma$ meets $E$ at a smooth point, then $\pi(E)$ 
is smooth. It follows (since $A_1$ is smooth), that $\Sigma$
is in the Du Val locus, thus $a=1$, a contradiction.

So we can assume $E$ meets $\Sigma$ at a singular point $z$.
Now the possibilities for $\map {\t{T}}.\t{S}_1.$ are 
given in (11.4). Since by (4.12.3) there must be a non-Du Val 
singularity along $\Sigma \setminus \Sigma \cap E$, it follows
(in the language of (11.4)) that there are no interior
blow ups, and $K_T + E + \Sigma$ is lt. Let
$w$ be the other singularity along $\Sigma$.

\subhead  19.4.1.1 Case 1, $(s-1,r+1) =( 2,3)$ \endsubhead
By the classification of lt singularities, 
necessarily the singularities
along $E$ are $A_3$, $A_2$ and $z=A_1$. Thus (since $x_0$ is 
not Du Val ), $k \geq 3$. By (11.4), $w=(3',A_{k-2},3,3,2)$. 
$e_0 =\frac{12k -24}{12k-23}$. One easily checks that 
$K_{S_0} \cdot \Sigma > 0$, a contradiction. \qed 

\subhead 19.4.1.2 Case 1, $(s-1,r+1)=(3,2)$ \endsubhead
The singularities along $A_1$ are $(3,2,2)$, $A_4$,
and $A_1$. 

Suppose first that $(3,2,2)$ is not added.
Then necessarily the singularities along $E$ are 
$(3,2,2)$,$A_1$, and $z=A_1$, and $w=(3',A_{k-2},3,2,2,2)$,
$k \geq 2$. One checks that if $k \geq 3$ then $K_{S_0}$ is
ample. Thus $k=2$: The singularities are a non-chain
$x$ with center $(-2)$ and branches $(3',2,2),(2),(2)$, and
$y=w=(3,3,2,2,2)$. The map $\map \t{S_0}.{\pr 2.}.$
is obtained from a configuration with conic $B$, secant
line $A$ and tangent line $C$. Blow up $3$ times at $c$ along
$C$, $5$ times at $b$ along $B$ and then $3$ more times over
$b$ to give the prescribed configuration. $A$ becomes
$E$, $C$ is a $-2$-curve, one of the branches, and
$B$ is the $-3$-curve of the $(3,2,2)$ branch. Let $W$ be the
tangent line at $a$. $W$ is a one curve on $S_0$ with
two branches at $x$. One meets $C$, and has local index
(one computes) $16$, the other is tangent to $B$ and
transversal to $A$ (at $a$). Its local index is $8$. 
One computes 
$$
-K_{S_0} \cdot W = 3 - 3/8- 2 \cdot 3/4 -3/4 =3/8.
$$ 
Thus $S_0$ is log uniruled by (6.5).

Now we assume $(3,2,2)$ is added. Then $z$ has index either
$2$ or $3$. Suppose first it is $2$. Then $k \geq 3$,
and $y=(3',A_{k-2},4,2,2)$. $e_0= \frac{5k-10}{5k-9}$.
One checks $K_{S_0}$ is ample.

Suppose $z=(2,2)$. Then $k \geq 3$, $w=(4',A_{k-2},4,2,2)$,
and $e_0=\frac{30k -60}{30k-59}$. One checks
$K_{S_0}$ is ample.

Finally suppose $z=(3)$. Then $k \geq 2$, $w=(2',3,A_{k-2},4,2,2)$,
and the coefficient of the $-3$-curve at $x_0$ is
$\frac{21k -33}{30k -49}$. The coefficient of $2'$ is
$\frac{7k+2}{21k+2}$. One checks $K_{S_0}$ is nef.
\qed

\subhead 19.4.2.1 Case 2: $(r,s) =(2,3)$ \endsubhead
Singularities along $A_1$ are $(3,2,2),A_3,(3,2)$.
This case is ruled out by index consideration. \qed

\subhead 19.4.2.2 Case 2: $(r,s) =(1,4)$ \endsubhead
Singularities along $A_1$ are $(4',2),A_4,(2,2)$. 
Necessarily $(4,2)$ is added, $z=(2)$, $k \geq 3$, and
$y=(3',A_{k-1},5,2)$. $e_0 = \frac{30k -60}{30k -59}$.
One checks $K_{S_0}$ is ample. \qed

\subhead 19.4.2.3 Case 2: $(r,s) = (1,3)$ \endsubhead
 Singularities along $A_1$ are $(3',2),A_3,(2,2)$.
In any case $z=(2)$.

Suppose first that $(3,2)$ is added. Then $k \geq 3$
and $y=(3',A_{k-1},4,2)$. $e_0=\frac{12k -24}{12k-23}$.
One checks $K_{S_0}$ is ample. \qed 


\heading \S 20 Proof  of  (1.1) and corollaries \endheading

We begin with the formal proof of (1.1).

\proclaim{20.1 Theorem} Let $S$ be a rank one log del Pezzo surface. Assume
$S^0$ has trivial algebraic fundamental group. If $S$ does not have a tiger
then $S^0$ is uniruled. Moreover, exactly one of the following holds for
the hunt with $(S_0,\Delta_0) = (S,\emptyset)$:
\roster
\item $A_1 \subset S_1$ is smooth, and $(S_2,A_2 + B_2)$ is a smooth 
Banana and included in the classification (13.2).
\item $g(A_1) = 2$ and $S$ is the surface given in (15.2).
\item $g(A_1) =1$ and $A_1$ has a simple node. $(S_2,A_2 + B_2)$ is a fence with
$B_2$ smooth. Either $(S_2, A_2 + B_2)$ is included in the classification (13.5), or we have (17.3.2).
\endroster
\endproclaim
\demo{Proof} The possibilities for the first hunt two hunt steps are described in (8.4.7).
By (14.1), $T_1$ cannot be a net. Now apply (15.1), (16.1), (17.1), (17.1.1) and (18.1). \qed\enddemo

\demo{20.2 Proof of (1.1)} Notation as above (1.1). As above (1.3), the $(K_X + D)$-MMP reduces (1.1) 
to (1.3). By (6.2) and (7.2), we may assume $D$ is empty, $S$ has no tiger, and $S$ has trivial 
algebraic fundamental group. Now apply (20.1). \qed\enddemo

Here we give proofs of the corollaries in the introduction, in those cases where the result does 
not follow immediately from (1.1) and (1.3).

\demo{Proof of (1.4)} (3) and (4) are instances of (6.1).  Case (2) reduces to case (1)
after passing to a log terminal model, as in the proof of (6.3). Thus its enough
to consider (1). Now pass to a log resolution, replace $B$ by the total transform,
and apply (1.1). \qed \enddemo

\demo{Proof of (1.6)} Let $B = \rdown \Delta .$. 
If $S$ has Picard number one the result follows from 
(7.5) and (7.10). Otherwise, running $(K_S+\Delta)$-MMP we may assume $S$ has Picard
number two, and each edge gives a $\pr 1.$-fibration, and $B$ has degree at
most one for either fibration. As in the proof of (7.10), we can deform the
union of the two general fibres $F_1 + F_2$ to get a dominating family of
integral rational curves $D \subset S^0$ which deform with a general
point fixed. $D \cdot B  \leq 2$, so $D \cap U$ gives a connecting family of
images of $\aos$. \qed \enddemo

\demo{Proof of (1.7)} This follows from (3.3-3.4) of \cite{Batryev89}. 
Batryev works with terminal three folds, but his arguments apply
without change to log surfaces. The
description of the edges follows from his definition of total
pullback, and our (1.3). \qed \enddemo

\demo{Proof of (1.8), due to Mori} By (2.1) and (3.1) of \cite{MR82} it
is enough to show that if 
$$
\map T_S(- \log(D)).L.
$$
is a generically surjective map to a line bundle $L$, 
then $c_1(L) \cdot H \geq 0$, for any ample $H$. By (1.7) it is enough to show
$c_1(L) \cdot [C] \geq  0$, for $C$ as in (1.7). This follows from (5.5.1).
\qed \enddemo

\demo{Proof of (1.8.1)} Immediate from (1.8), and (10.8), (10.12) of 
\cite{Kollaretal}. \qed \enddemo

\demo{Proof of (1.12)} We can assume $C$ is not rational.

Suppose first $\text{dim}(X) =2$. Since $K_X \cdot C < 0$,
we have $C^2>0$ by adjunction. Thus $K_X + D$
is not effective, and so by (1.1) $X$ is dominated by rational
curves meeting $D$ at most once. The dominating family
meets $C$ by the Hodge Index Theorem, and has no basepoints
along $C$ by (5.5).

Now suppose $\text{dim}(X) =3$.
Let $W \in |m(K_X + D)|$ for some $m > 0$. Let $\t{C}$ be the
normalisation of $C$. As in the proof of (2.13) of
\cite{Kollar91}, there is a proper smooth map $f:\map T.B.$
from a smooth surface to a smooth curve, with $\t{C}$ a fibre of
$f$, and a map $g:\map T.X.$ with two dimensional image,
such that $g|_{\t{C}}$ is the normalisation map (in short $C$ sweeps out
a surface). Let $S = \overline{g(T)}$. 
Since $S$ is covered by $W$ negative curves, $S$
is a component of $W$, and $S \cdot C < 0$. Thus $K_S \cdot C < 0$.

Let $p:\map S'.S.$ be the normalisation. $g$ lifts to $g':\map T.S'.$.
Let $C' = g'(\t{C})$, 
$\pi:\map S''.S'.$ be the minimal desingularisation, 
$C'' \subset S''$ be the strict transform of $C'$, $V = D|S$,
$V' = p^*(V)$, and $V'' = \pi^*(V')$.

There is an effective integral Weil divisor
$\Gamma$ such that $K_{S'} + \Gamma = p^*(K_S)$.
Recall that the intersection product with $\Bbb Q$-coefficients
is defined on any normal surface.
We have 
$$
\align
(K_{S''} + V'') \cdot C'' &\leq (K_{S'}+ V') \cdot C' \\
                          &\leq (K_{S'} + \Gamma + V') \cdot C' =  (K_S + V) \cdot C < 0
\endalign
$$
where the third inequality holds since $C'$ sweeps out $S'$.

Thus by the
surface case there is a rational curve through
a general point of $C''$, meeting $V''$ at most once. Of course 
its image passes through a general point of $C$ and
meets $D$ at most once.
\qed \enddemo

\heading \S 21 $S$ with $\alg(S^0) = \{1\}$ but no tiger \endheading

Our first goal is to give an example of $S$ without tiger. We again use the hunt,
but of course, now we cannot assume $S=S_0$ does not have a tiger. None the less, 
following (8.2), (8.2.5), and (8.2.8), we construct $(S_i,\Delta)$
starting with $(S_0,\Delta_0)= (S,\emptyset)$, by exactly the same procedure used in 
hunting for a tiger. We continue to use the notation (8.2.10). 

By (8.2.5.8), if $K_{S_{i+1} + \Delta_{i+1}}$ has a tiger, so does $K_S$. To get an 
implication in the other direction, we need to look at a smaller boundary.

\subhead 21.1 The boundaries $\gamma_i$ \endsubhead
We introduce some more notation. We recall our convention of using the same symbol for a 
curve, and its strict transform under some birational map. Inductively define 
\newline
$\gamma_i = \sum_{j\leq i} m_j E_j$ and $m_{i+1} = e(E_{i+1},K_{S_i}+\gamma_i)$.
We let $\gamma_0 = \emptyset$. Thus $\gamma_1=e_0A_1$, $m_2=e(E_2,K_{S_1}+e_0A_1)$, and
$\gamma_2 = e_0 A + m_2 B$. 

For the next lemma, assume $T_j$ is not a net for $j \leq i$.
\proclaim{21.1.1 Lemma} If $K_S$ has a tiger, then so does $K_{S_i}+\gamma_i$.
\endproclaim
\demo{Proof} We prove inductively that if 
$K_{S_i} + \gamma_i$ has a tiger, then so does
$K_{S_{i+1}} + \gamma_{i+1}$. 

Suppose $K_{S_i}+\gamma_i + \alpha$ is numerically trivial, where $\alpha$ is
effective. Define effective $\beta$ by
$$
\align
f_i^*(K_{S_i} + \gamma_i + \alpha) &= K_{T_{i+1}} + {\gamma_i}
+ m_{i+1} E_{i+1} + \beta \\
&= \pi_{i+1}^*(K_{S_{i+1}} + \gamma_{i+1} + \pi_{i+1}(\beta)).
\endalign
$$
Then $K_{S_i} + \gamma_i + \alpha$ is klt iff
$K_{S_{i+1}} + \gamma_{i+1} + \pi_{i+1}(\beta)$ is klt. \qed \enddemo

\remark{21.1.2} We note that (21.1.1), and its proof, hold for
any sequence of relative Picard number one extractions and contractions
$f_i:\map T_{i+1}.S_i.$ and $\pi_{i+1}:\map T_{i+1}.S_{i+1}.$,
with $f_i$ $K$-positive. Indeed, all we use in the proof is that 
$\beta$ is effective.
\endremark

\subhead 21.2 A condition to ensure there is no tiger \endsubhead

Let $\Gamma$ be an effective $\Bbb Q$-Weil divisor. For a point $p$ on a 
smooth surface, let $m_p(\Gamma)$ be the coefficient of the exceptional 
divisor of the blow up of $p$ in the inverse image of $\Gamma$ (when $\Gamma$ 
is reduced this is the usual multiplicity).

\proclaim{21.2.1 Lemma} Let $X$ and $Y$ be divisors meeting normally at a 
smooth point $p$ of a surface $S$. If 
$$
\max(a,b)+m_p(G)<1
$$ 
then $K_S+aX+bY+G$ is klt at $p$. 
\endproclaim
\demo{Proof} We may assume $a \geq b$. Blow up $S$ at $p$. The exceptional 
divisor $E$ has coefficient $a+b+m_p-1<b$, and the strict transform of $X$ 
and $Y$ union $E$ has normal crossings. Since $m_p$ does not go up under
blow ups, the result follows by induction. \qed\enddemo

\proclaim{21.2.2 Proposition} Let $B$, $A$, $C$ be a conic, a secant line
and a tangent line in $\pr 2.$, with $A \cap B \cap C = \emptyset$. 
$K_{\pr 2.}+(60/67)A+(381/469)B+(30/67)C$ has no tiger.
\endproclaim
\demo{Proof} Let $\alpha$ be an effective $\Bbb Q$-Weil divisor
such that
$$
K + \Delta = K_{\pr 2.}+(60/67)A+(381/469)B+(30/67) + \alpha
$$
is numerically trivial.

Let $p$ be a point of $\pr 2.$, and let $L$ be a line.
Note that 
$m_p(\alpha)\leq \alpha\cdot L=3-60/67-2.381/469-30/67=15/469$. Let
$\Sigma^i$ be the $-1$-curve of the $i^{th}$ blow up of 
$\pr 2.$ at $z=B \cap C$ along $B$. Then
$$
\align
e(\Sigma^1,K + \Delta) &\leq 137/469. \\
e(\Sigma^2,K + \Delta) &\leq 274/469.
\endalign
$$
Now blow up twice at $z$ along $B$. Applying (21.2.1), to the
total transform, which has normal crossings, we conclude
that $K + \Delta$ is klt. \qed\enddemo

\remark{Remark} The coefficients in (21.2.2) may seem a bit bizarre. We chose them to make
the analysis in the next section as simple as possible. \endremark

\subhead 21.3 An example of a smooth Banana with no tiger \endsubhead

Throughout (21.3) let $S=S_0$. 
\proclaim{21.3.1 Lemma} If $S_2$ is a smooth banana, then $\pi_1(S^0)$ is
Abelian. \endproclaim
\demo{Proof} By (13.2), after a blow up and blow down of $S_2$ we obtain
$\Hbz 2.$. We will use the hunt notation (although it is possible
this last transformation is not part of the hunt). We have by (13.2)
a rational map $h:S \dasharrow {\Hbz 2.}$, and the exceptional
locus of $h^{-1}$ is contained in $A+B+C\subset \Hbz 2.$, for
$A$, $B$ two sections and $C$ a fibre, such that $A+B+C$
has normal crossings. Blowing up a point of $p\in A\cap B$ and
blowing down the fibre through $p$, gives a rational map 
$g:S \dasharrow {\pr 2.}$ such that the exceptional locus of $g^{-1}$ is contained 
in $W \subset \pr 2.$, a union of four lines, with normal crossings. Then 
$U = \pr 2. \setminus W$ is an open subset of $S^0$, thus 
$\pi_1(U) \twoheadrightarrow \pi_1(S^0)$. $\pi_1(U)$ is Abelian by Zariski's 
conjecture.\qed\enddemo

\proclaim{21.3.2 Lemma} If $S_2$ is a Banana of type given in (19.3.2.2.1), 
then $S^0$ is simply connected. \endproclaim
\demo{Proof} By (21.3.1), $\pi_1(S^0) =H_1(S^0,{\Bbb Z})$. If $D$ is the 
exceptional locus of $\map \t{S}.S.$, then according to (2.1) of \cite{MZ88}
there is an exact sequence
$$
\ses r:\Pic(\t{S}).H^2(D,{\Bbb Z}).H_1(S^0,{\Bbb Z}).
$$
Let $M$ be the cokernel. 

If we let ${\Cal D} \subset \Pic(\t{S})$ be the subgroup generated
by $D$, then (since $x$ and $y$ are chain singularities), 
we have
$$
\ses \Cal D.H^2(D,{\Bbb Z}).\frac{{\Bbb Z}}{\Delta(x) \cdot \ZZ} \oplus 
\frac{{\Bbb Z}}{\Delta(y) \cdot \ZZ}.
$$
with either end component mapping to the generator in 
$\frac{{\Bbb Z}}{\Delta(p) \cdot \ZZ}$ (for $p=x$, $y$). 

Let $E_c$, $E_b$  be the $-1$-curve over $c$, $b$ (in the notation of 
(19.3.2)). One easily checks that the image of $E_b$ in $M$ is $(1,3)$ and the
image of $E_c$ in $M$ is $(1,2)$. Thus $r$ is surjective. \qed\enddemo

\proclaim{21.3.3 Example} If $S$ is the surface from (19.3.2.2.1) with 
$k=8$ then $S^0$ is simply connected, and $S$ has no tiger.
\endproclaim
\demo{Proof} Let $f_1$ and $f_2$ be (as usual) the first two hunt steps. $(S_2,A + B)$ is a smooth 
banana. Let $f_3$ blow up the $A_1$ point along $A \subset S_2$ (this is not a hunt step), so 
$(S_3,B+A+ C)$ is as in (21.2.2). $m_1=e_0= \frac{60}{67}$. $x_1$ has the form (with respect to $A_1$) 
$(\ul{3},2,2)$. Thus by (8.3.8), $m_2 =3/7+3e_0/7=381/469$, and $m_3 = (1/2)(60/67)=30/67$.
By (21.2.2), $K_{S_3}+\gamma_3$ has no tiger. We conclude $S$ has no tiger by (21.1.1). $S^0$ is 
simply connected by (21.3.2). \qed  \enddemo

\subhead 21.4 Counter-Example to a conjecture of Miyanishi \endsubhead

\definition{21.4.1 Definition} We say that a surface $V$ is
{\bf affine ruled} if there is an open subset of $V$ which is isomorphic
to a product $U \times \af 1.$.
\enddefinition

According to \cite{Zhang94a}, Miyanishi has conjectured that given any
rank one log del Pezzo $S$ there is a finite \'etale cover of $S^0$ which is 
affine ruled. We show here that any tiger-free $S$ with $S^0$ simply connected
(for example the surface of (21.3.3) ) is a counter-example.

\proclaim{21.4.2 Lemma} Let $S$ be a rank one log del Pezzo, and 
$V \subset S^0$ an open subset dominated by images of $\af 1.$. If
the complement of $V$ contains a divisor, then $S$ has a tiger.
\endproclaim
\demo{Proof} Assume that there is a dominant map $f:\map {U \times \af 1.}.V.$,
and that $D \subset V^c \subset S$ is an irreducible curve. 
Let $C$ be the closure of the image of a general $\af 1.$. 
If $C$ meets $D$ at a general point of $D$, then 
$(K_S + D)\cdot C < 0$ by (5.11), and thus $D$ is a tiger. 

 Otherwise every $C$ passes through a common point $p = C \cap D$. 
Let \map $S'.S.$ extract the divisor $E$ over $p$, so that $C'$ cuts out $E$.
Then $C'$ and $D'$ are disjoint on $T$ (here the prime indicates strict
transform). Define $\lambda$ so that 
$K_T+E+\lambda D' = f^*(K_S+\lambda D)$. Note $\lambda > 0$ since $S$ is lt. 
Now $C'$ meets $E$ in a single smooth point, and so by (5.11), 
$(K_T + E) \cdot C' = (K_S + \lambda D) \cdot C < 0$. Thus 
$K_S + \lambda D$ is anti-ample. It is not klt by construction, and so
$\lambda D$ is a special tiger for $K_S$. \qed\enddemo

\proclaim{21.4.3 Corollary} An affine ruled rank one log del Pezzo surface
has a tiger. \endproclaim

\proclaim{21.4.4} Let $S$ be a rank one log del Pezzo. If $S$ does not have a 
tiger and $\alg (S^0)=\{1\}$, then no open subset of $S$ has a finite \'etale cover which is 
affine ruled.
\endproclaim
\demo{Proof} Let $S$ be a simply connected rank one log del Pezzo, with no tiger. Suppose 
$V \subset S$ has a finite \'etale cover which is affine ruled. Dropping points from $V$, 
we may assume $V \subset S^0$. By (21.4.2) $V^c$ has codimension at least two, and thus $V$ 
is simply connected (see the proof of (7.3)). Thus $S$ is affine ruled, contradicting 
(21.4.3). \qed \enddemo

\remark{Remark} Definition (21.4.1) is due to Miyanishi. Koll\'ar has suggested that in view of 
the definition of ruled, (see for example (IV.1.1) of \cite{Kollar96}) it is perhaps more natural to 
assume only the existence of a birational map $f:\map {U\times \af 1.}.S^0.$. We do not know if 
(21.4.3) holds under this weaker assumption. 
\endremark


\heading \S 22 Tigers, complements and toric surfaces\endheading

 The following section prepares the way for the classification that appears
in \S 23. First we introduce the notion of complements and give the 
connection between tigers and complements. 

 We first recall the definition of a complement, cf. Chapter 19 of 
\cite{Kollaretal}, for more details.  

\definition{22.1 Definition} Let $X$ be a normal variety and $\Delta=D+B$ a 
boundary, where $D$ is reduced and $B$ is a pure boundary. We say $K_X + \Delta$ 
is 
$n$-complemented, for a positive integer $n$, 
if there is a divisor $D'$, such that 
\roster
\item $nD'\in |-nK_X|$,
\item $nD'\geq -nD-\rdown (n+1)B.$ and 
\item $K_X+ D'$ is lc.
\endroster
\enddefinition

On a first encounter, (22.1) may seem a bit bizarre. For our classification in
\S 23, we need only a simple form. For the readers convenience, we restate
the definition in this case:

\definition{22.1.4 Definition ($1$-complement of a reduced divisor)} $(X,D)$ as in (22.1). We say 
$K_X + D$ is $1$-complemented if there is a reduced divisor $Y$ such that $K_X + D + Y$ is 
Cartier, trivial, and log canonical. In the notation of (22.1), $D' = D + Y$.\enddefinition

The following result, and its proof, are due to Shokurov. We learned of them from Alessio Corti:

\proclaim{22.2 Lemma} Let $S$ be a Fano surface, and $\Delta$ an effective ${\Bbb Q}$-divisor, 
such that $K_S+\Delta$ is log canonical. If $K_S + \Delta$ has a tiger, then $K_S+\Delta$ has 
a $1$, $2$, $3$, $4$, or $6$-complement.
\endproclaim
\demo{Proof} After scaling and applying (19.2) of \cite{Kollaretal}
we may assume $K_S + \Delta$ is anti-nef and maximally log canonical.
There is an extraction $\pi :\map T.S.$ such that 
$K_T+\Gamma=\pi ^*(K_S+\Delta)$ is lc, and $D = \rdown \Gamma.$ is non-empty, 
and contained in the smooth locus. By (19.2) of \cite {Kollaretal} it is 
enough to show $K_T + \Gamma = K_T + D + B$ has a $1$, $2$, $3$, $4$ or 
$6$-complement.

Look at the natural restriction exact sequence
$$
\multline
0\longrightarrow \ring T.(-nK_T-(n+1)D-\rdown (n+1)B.) \longrightarrow 
\ring T.(-nK_T-nD-\rdown (n+1)B.) \longrightarrow \\
\longrightarrow\ring D.(-nK_D-\rdown (n+1)B|_D.)\longrightarrow 0.
\endmultline
$$

Now $K_D+B|_D$ is $n$-complemented, for some $n \in \{1,2,3,4,6\}$ by
(19.4) of \cite {Kollaretal}. We argue next that 
$$
H^1(T,\ring T.(-nK_T-(n+1)D-\rdown (n+1)B.))
=H^1(T,K_T+\rup -(n+1)(K_T+\Gamma).)=0.
$$

If $K_S + \Delta$ is anti-ample, this holds by Kawamata-Viehweg vanishing 
\cite {KMM87}, so we may assume $K_S + \Delta$ is numerically trivial.
Let $G = K_T+\rup -(n+1)(K_T+\Gamma).$. $R^1\pi_*(G) =0$ by Kawamata-Viehweg 
vanishing. Hence we are reduced to showing $H^1(S,\pi_*(G))=0$.

Now $\pi_*(G) = \ring S.(K_S+\rup -(n+1)(K_S + \Delta).)$. 
Since the Picard number of $S$ is one, $\rup -(n+1)(K_S + \Delta).$
is ample, and Kawamata-Viehweg vanishing applies, unless
$-(n+1)(K_S + \Delta)$ is integral, in which case (22.2.1) applies.

 Thus we may find a divisor $D'$ on $T$, such that $nD'\in |-nK_T|$
and $K_T + D'$ is lc in a neighbourhood of $D$. Now the divisor 
$K_T+(1-\epsilon)D'+\epsilon (D+B)$ is not klt in a neighbourhood of $D$. By
the connectedness Theorem, (17.4) of \cite {Kollaretal}, it follows that 
this divisor is klt away from $D$. Thus $K_T+D'$ is lc and $K_T+\Gamma$ 
has a complement. \qed\enddemo

\proclaim{22.2.1 Lemma} Let $X$ be a $\Bbb Q$-factorial 
log Fano variety and $H$ a numerically trivial integral Weil
divisor on $X$. Then $H^i(X,\ring X.(K_X + H)) =0$ for all
$i > 0$. \endproclaim
\demo{Proof} $H$ induces a cyclic Galois cover, \'etale in
codimension one, $g:\map Y.X.$ such that $K_Y = g^*(K_X + H)$. In particular, 
$Y$ is log Fano, and so $H^i(\ring Y.(K_Y)) =0$. $\ring X.(K_X + H)$ is a 
direct summand of $g_*(\ring Y.(K_Y))$. \qed \enddemo

 The following Lemma, refines (22.2), in a special case. All we need for this 
paper is the case (22.3.1), which will be used in (23.7).

\proclaim{22.3 Lemma} Let $S$ be a log del Pezzo of rank one. Suppose $K_S+A$
is plt and $-(K_S+A)$ is ample, where $A$ is integral and non-empty. $A$ contains
at most three singularities of $S$. Let 
$p\geq q\geq r$ be the indices of these singularities (we allow index $1$).

 Then $K_S+A$ has an $n$-complement $X$. Moreover 
\roster
\item $n=1$, and $K_S+C$ is plt for any component $C$ of $X$, in the case
$r=1$ (that is when there are at most two singularities along $A$).
\endroster
 
 Otherwise $A$ has three singularities, $X$ has one other component $B$, 
$K_S+A+B$ is lc in a neighbourhood of $A$, and $A$ meets $B$ at one point. 
Moreover
\roster
\item[2] $n=2$ and $X=A+B$, when $q=r=2$.
\endroster
Otherwise $X=A+(B/n)$. Moreover
\roster
\item[3] $n=3$, when $p=q=3$. 
\item $n=4$, when $p=4$ and $q=3$,
\item $n=6$, when $p=5$ and $q=3$. 
\endroster
\endproclaim
\demo{Proof} $A$ can contain at most three singularities by adjunction. By the proof of (22.2), 
we know that if $K_A + \Diff_A(\emptyset)$ (see (L.2)) has an  $n$-complement, then so does 
$K_S+A$. Moreover since we are lifting an element of $|-nK_A|$, we have a local picture around 
$A$. Now $A$ is a copy of $\pr 1.$, and the description of $X$ follows from the division into 
cases in the  proof of (19.4) of \cite{Kollaretal}. Note in particular that if $R=X-A$, then 
$R\cdot A=1/p$ when $n=2$ and $R\cdot A=1/sn$ when $n\geq 3$, where $s$ is the index of the point 
where $A$ meets $B$. On the other  hand $B\cdot A\geq 1/s$, with equality in the case that 
$K_S+A+B$ is lc in a neighbourhood of $A$ and $R=(k/n)B$, for some $k\geq 1$. The result now 
follows, except in the case $n=2$. Here we make the ad hoc argument that we have a classification 
cf. (23.6) and that $B$ is the image of the obvious fibre. \qed\enddemo


\subhead Toric Pairs \endsubhead

 We now recall the definition of a toric variety. As we shall see in \S 23, 
they are fundamental to the classification of all log del Pezzos. 
Also we present a result about toric varieties which is rather surprising, 
in light of the boundedness results of Alexeev \cite{Alexeev94}.

\definition{22.4 Definition} A pair $(S,X)$ is said to be toric, 
if $S\setminus X$ is a copy of a torus, and the natural action of the torus 
on itself, extends to an action on $S$. 
\enddefinition

Toric pairs are very easily understood. Either one may use the approach of fans, in which case a 
pair $(S,X)$ is essentially specified by two pairs of relatively prime integers (see for example 
\cite{Fulton93}) or one may untwist a toric pair, using a series of transformations of 
type (a) and thus create $(S,X)$ inductively, see  (23.9) and the proof of (23.12). The results 
below illustrate this idea. 

Several people have asked the following question: Does the set
$$
\{\, K^2_S\, |\, S \text { is a log del Pezzo surface of rank one}\, \}
$$
satisfy ACC for bounded sequences ?

  We prove
\proclaim{22.5 Proposition} 
\roster
\item The set
$$
\{\, K^2_S\, |\, S \text { is a toric surface of rank one, with two 
singularities}\, \}
$$
is dense in the set of real numbers at least four, and 
\item The set
$$
\{\, K^2_S\, |\, S \text { is a toric surface of rank one}\, \}
$$
is dense in the set of positive real numbers. 
\endroster
\endproclaim

\proclaim{22.6 Lemma} Let $S$ be a toric surface of Picard number one, with 
three singular points of index $p$, $q$ and $r$. 

 Then 
$$
K^2_S=\frac {(p+q+r)^2}{(pqr)}.
$$
\endproclaim
\demo{Proof} Now there are three toric divisors, $D_1$, $D_2$ and $D_3$,
which form a triangle, such that each vertex is one of the singular points.
Moreover $K_S+D_1+D_2+D_3$ is lc and linearly equivalent to zero. 
Thus 
$$
K^2_S=(D_1+D_2+D_3)^2.
$$

 But as $K_S+D_1+D_2+D_3$ is lt, 
$$
\{\,D_i\cdot D_j\,|\, 1\leq i<j \leq 3\, \}=\{\,p,q,r\,\}.
$$

 The result is now an easy computation. \qed\enddemo

\proclaim{22.7 Lemma} 
\roster
\item The set
$$
\{\, \frac {(x+1)^2}x\, |\, \text {$x$ is a positive real number}\,\}
$$ 
is dense in the set of real numbers at least four.
\item Fix positive integers $N$ and $r$. The set
$$
I_{N,r}=\{\, \frac {p'}q |\, \text {$c$ and $q>N$ are coprime 
integers, $p'=r(c+1)-q>0$}\,\}
$$
is dense in the set of positive real numbers. 
\endroster
\endproclaim
\demo{Proof} (1) is elementary calculus and (2) is clear. \qed\enddemo

\demo{Proof of Proposition} Pick a positive integer $r$ and coprime integers
$c$ and $q$. The three vectors $\{\, (1,0), (1,r), (-c,q)\,\}$ 
specify a fan in the lattice $\Bbb Z^2$, and give rise to a toric surface 
$S$ of Picard number one, with indices $p=rc-q$, $q$ and $r$. 

 By (22.6) 
$$
\align
K^2_S&=\frac {(p+q+r)^2}{(pqr)} \\
     &=(1/r)\frac {((p/q)+1+(r/q))^2}{(p/q)}\\
     &=(1/r)\frac {((p'/q)+1)^2}{(p'/q)}(1+\epsilon)
\endalign
$$
where $|\epsilon|=o(1)$ (where we think of fixing $r$ and let $q$ go to 
infinity) and $p'=p+r$. The result now follows easily from (1) 
and (2) of (22.7). \qed\enddemo
 
\heading
\S 23  Classification of all but a bounded family of rank one log del Pezzo surfaces
\endheading

Here we will give a classification of all but a bounded family of 
rank one log del Pezzos. The proofs use the methods of the hunt, introduced in \S 8, but
as we expect there will be readers who are interested
in the statement of the classification, independent of the rest of the paper, we will state the
results in a self contained way. This will involve some repetition of definitions from
\S 8. The analysis of this section is independent from most of the paper. It uses results
and notations from \S 2-3 and \S 8-13 only.

Our classification is of the following sort. Starting with a rank one log del Pezzo
we will make an essentially canonical birational transformation (the first hunt
step) from $S$ to either a $\pr 1.$-fibration or to a second rank one log del Pezzo. We 
will explicitly classify the resulting surfaces. We will also indicate how the inverse
transformations can be classified, the procedure is easy and elementary, but we do not
write down an explicit list as it would be notationally too involved.

\subhead 23.0 The First Hunt Step \endsubhead
Given a rank one log del Pezzo $S$, let $E$ be an exceptional
divisor of the minimal desingularisation with maximal coefficient (that is minimal discrepancy).
Let $f:\map T.S.$ be the extraction (of relative Picard number one) of $E$. We assume
$E$ is chosen so that $K_T + E$ is log terminal (such a choice is always possible, see
(8.3)) $T$ has a unique $K_T$-negative contraction $\pi$. $\pi$ is either a 
$\pr 1.$-fibration, or birational. In the second case we let $S_1$ be the image, and let 
$A_1 \subset S_1$ be the image (with reduced structure) of $E$. $S_1$ is again a rank one 
log del Pezzo. For details see (8.3).

\proclaim{23.1 Lemma} (Notation as in (23.0))
For all but a bounded family of rank one log del Pezzos, $-(K_T + E)$ is effective,
and either $\pi$ is a $\pr 1.$-fibration or $K_{S_1} + A_1$ is log canonical and anti-nef.
\endproclaim

\proclaim{23.2 Theorem} (Notation as in (23.0)) Let $S$ be a rank one log del Pezzo surface,
such that  $-(K_T + E)$ is effective. One of the following holds:
\roster
\item $\pi$ is a  $\pr 1.$-fibration, and $-(K_T+E)$ is nef, $E$ is a multi-section of degree 
at most two and $-(K_T+E)$ is ample if $E$ is a section. 
\item $K_{S_1}+A_1$ is lt and anti-ample and $A_1$ contains at most two 
singular points. Possibilities for $(S_1,A_1)$ are classified in (23.12).
\item $K_{S_1} + A_1$ is lt and anti-ample and $A_1$ contains exactly three
singularities. Pairs $(S_1,A_1)$ are classified in (23.5.1). They are in one
to one correspondence with non cyclic quotient singularities (of dimension two)
\endroster

Otherwise $S_1$ is Du Val, and one of the following holds:
\roster
\item[4] $A_1 \subset S_1^0$ is a rational curve of arithmetic genus one. Possibilities for
$S$ such that $A_1$ has a cusp are bounded. 
\item $K_{S_1}+A_1$ is lt and numerically trivial, and $A_1$ contains exactly three 
singular points. Pairs $(S_1,A_1)$ fall into $5$ families,
listed in (23.5.2-4). 
\endroster
\endproclaim

\proclaim{23.3 Lemma} Let $(S_1,A_1)$ be a pair of a rank one log del Pezzo with
a reduced rational curve such that $K_{S_1} + A_1$ is anti-nef. If $\pi: \map T.S_1.$ is a 
contraction of relative Picard number one, such that $E$, the strict transform of $A_1$, 
contracts $f:\map T.S.$ to a log terminal surface, then $S$ is a rank one log del Pezzo.
\endproclaim

We now make remarks regarding the above:

\remark{23.4 Remarks} 
\roster
\item Rank one Du Val del Pezzo surfaces are bounded, and there is a short list
of possibilities. See \S 3, or
\cite{MZ88}.
\item It is a simple matter to classify possible $\pr 1.$-fibrations in (23.2.1), along 
the  lines of (11.5). One can classify pairs in (23.2.4) by considering the 
anti-canonical series $|-K_{S_1}|$, see for example (3.6).
\item Given a pair $(S_1,A_1)$ contractions $\pi:\map T.S_1.$
such that $K_T+E$ is lt are classified in (11.4). Using the classification 
of log terminal singularities, one can then further classify possibilities 
such that $E$ contracts to a quotient singularity. In this sense (23.1-3) 
give a classification of all but a bounded family of log del Pezzo 
surfaces. 
\item Note the assumption in (23.2), is exactly that $E$ is a tiger.
We believe it would be possible, using the methods of this 
paper, to completely classify those log del Pezzos $S$, such that $E$ is not a 
tiger. Indeed, in the body of the paper we have already constructed an explicit collection
$\goth F$ which includes all $S$ with no tiger, 
under the additional assumption that $\alg(S^0) = \{1\}$.
\endroster
\endremark

\subhead Proofs of (23.1-3) \endsubhead

\demo{Proof of (23.1)} By (9.3), it is enough to find $e'< 1$ such that  if $e=e(S)>e'$, 
then $E_1$ is a tiger. Obviously we may assume $a$ (of the hunt) is less than one.  

If $T$ is a net, and $e>2/3$, then $E$ has degree at most $2$, thus $K_T + E$ is 
anti-nef. Now suppose $\pi$ is birational. $K_{S_1} + a A_1$ is flush by (8.4.5). 
If $e > 4/5$ then $K_{S_1} + A_1$ is lc by (8.0.9.3) and, by (5.3) of \cite{Kollar92b}, 
if $e>41/42$, then $K_{S_1} + A_1$ is anti-nef. \qed \enddemo

\demo{Proof of (23.3)} Suppose $S$ is not a log del Pezzo. Then $K_S$ is nef, and pulling 
up and pushing down $K_{S_1}+eA_1$ is nef. Since $S$ is lt, $e<1$, and we get
a contradiction. Hence (2). \qed\enddemo 

\demo{Proof of (23.2)} Clearly (1) holds if $T$ is a net. So assume
$\pi$ is birational. We will show that $K_{S_1} + A_1$ is anti-nef,
log terminal at singular points, and, except for a bounded collection
of $S$, log canonical. The division into cases then follows easily by adjunction.

By assumption $-(K_{T_1} + E)$ is effective, thus
$K_{S_1} + A_1$ is anti-nef. Suppose $a$ (from the hunt) is at least
one. Then $K_{S_1} + E_1$ is $\pi$ non-positive, thus, since $K_{T_1} + E_1$ is
log terminal, so is $K_{S_1} + A_1$. 

So we may assume $a < 1$. $K_{S_1} + a A_1$ is flush by (8.4.5), thus $K_{S_1} + A_1$
is log terminal at singular points, by (8.0.4). If $e>4/5$ then $K_{S_1}+A_1$ is lc by 
(8.0.9.3). But if $e<4/5$ then $S$ is bounded by (9.3). \qed \enddemo

\proclaim{23.5 Proposition} Let $(S_1,A_1)$ be a pair of a rank
one log del Pezzo surface, $S_1$ and a reduced irreducible curve
$A_1 \subset S_1$ such that $K_{S_1} + A_1$ is anti-nef and 
log terminal. Assume there are at least three singularities 
of $S_1$ along $A_1$.

If $-(K_{S_1}+A_1)$ is ample then 
\item{(1)} $S_1 \setminus A_1$ has exactly one singular point,
a non-cyclic singularity, $z$.  
If $\map Z.S_1.$ 
extracts the central exceptional divisor of $z$, then $Z$ is a $\pr 1.$-fibration 
and $E$ and $A$ are sections. $(S_1,A_1)$ is uniquely determined by
$z$, and all non-cyclic singularities, $z$,  occur in this way
for some pair $(S_1,A_1)$.
\newline
\smallskip
If $K_{S_1}+A_1$ is numerically trivial then $S_1$ is Du Val, and $(S_1,A_1)$ 
is in one of five families: 
\item{(2)} $S = S(A_1 + 2A_3)$. $(S_1,A_1)$ is given by (19.2), Case 1, 
with $s=3$ and $r=1$. 
\item{(3)} $S = S(3A_2)$. $(S_1,A_1)$ is given by (19.2), Case 2, 
 with $s=2$ and $r=1$. 
\item{(4)} $K_{S_1}^2 =1$, $A_1$ is a $-1$-curve and $S_1$ is
one of $S(A_1 + A_2 + A_5)$, $S(2A_1 + A_3)$, $S(4A_1)$. 
The pairs are obtained from (13.5.1), (13.5.3) and (13.5.6) respectively,
by the inverse of a transformation of type (11.1.1).
\endproclaim
\remark{Remarks} 
\roster
\item For a construction, in (1), of $(S_1,A_1)$ from $z$,
see the proof of (23.5.1.1)
\item (23.5) is a classification of abstract pairs
$(S_1,A_1)$, but our principal interest comes from (23.2), where the pair $(S_1,A_1)$ is
the first hunt step, and this is why we use this subscripted notation. In fact,
it is easy to check using (23.3) and (11.4) that each of the pairs $(S_1,A_1)$ actually 
occurs as the first hunt step for some $S$. 
\endroster
\endremark

\demo{Proof} 
By (23.5.1.1) below, $z$ in case (1) determines $(S_1,A_1)$, and 
$K_{S_1} + A_1$ is negative.
In (2-4) one can compute directly that $K_{S_1} + A_1$ is numerically
trivial. We now show (without distinguishing between negative
and numerically trivial) that any pair $(S_1,A_1)$ falls into one of 
the cases (1-4). 

Let $r$ be the maximal index of any singularity along $A_1$.

Choose a singular point $y \in S_1$ as follows. If
there is a point $y \not \in A_1$ with $e(y) \geq \frac{r-1}{r}$, choose
such a $y$ of maximal coefficient (for $K_{S_1}$).
Otherwise take $y$ a point of $A_1$ of index $r$, which has
maximal spectral value (among points of index $r$). One
checks using (8.3.8), that $y$ is (one of the possible
choices for) $x$ for the next hunt step with respect to 
$K_{S_1} + a A$ for any $1 - \delta < \geq a < 1 $, for
fixed $\delta > 0$ sufficiently small. (Note that some care is needed to 
define $x$, to ensure that it is independent of $a$.)

Now we run the hunt with respect to $K_{S_1} + a A$ for $a$ close to one, and
$x = y$ as above. Note $K_{S_1} + aA$ is flush.

Here the situation is slightly different than in (8.4.7). We cannot assume 
that there is no tiger, indeed $A$ is obviously a tiger. None the
less, the division into cases (8.4.7) will be the same, using 
essentially the same argument. This will require a bit of preliminary analysis.

To avoid having any coefficients larger than one, we alter slightly the scaling 
convention of (8.2.4). Define $b \geq e_1$  such that $K_T + aA_1+b E$ is $R$-trivial 
($R$ is the $K_T$-negative extremal ray). Let $\Gamma' = a A_1 + e E$. 

\proclaim{Claim} $b < a$, unless $T$ is a net, $K + A_1$ is negative and 
$x \not \in A_1$, in which case $b = 2 - a$.
\endproclaim
\demo{Proof} 
If $K_{S_1}+ A_1$ is numerically trivial, then $b = e_1 < a$, so we can assume
$K + A_1$ is negative. 

Now suppose $T$ is a net, and $x \not \in A_1$. Since (for $a=1$) $e_1 \geq 1/2$, $A$ 
and $E$ are sections, so clearly $b = 2-a$. 

We will show in the remaining cases that $b < a$. It is enough to check for $a =1$. 

If $A \cap E \neq \emptyset$ then since $A_1 \subset T$ contains at least two 
singularities, adjunction gives
$$
0 > (K_T + A + b B) \cdot A \geq -1 + b.
$$
So we can assume $x \not \in A_1$. Then $T$ is not a net.
$K_{S_2} + A + b B$ is negative, so again $b < 1$ by adjunction. \qed \enddemo

For the inequalities and equalities for coefficients below, as well as
the analysis in the exceptional case of the claim, we take $a =1$.

\proclaim{Claim} (8.4.7.6-11) hold 
\endproclaim
\demo{Proof} $K_T + e_1 E + a A_1$ is level by (8.3.5.4). Since $b \leq a$, we can check 
flushness of $K_T + \Gamma'$ locally at each point of $E$, where it follows by (8.3.5.2). 
(8.4.7.6-11) now follow exactly as in the proof of (8.4.3).  \qed \enddemo 

Note that $e_1 \geq 1/2$, and $e_1 \geq 2/3$ unless $A_1$ contains $4$ singularities of 
index $2$.

Consider the next hunt step. Let $j =-\t{E}^2$ (as in the usual hunt).

As $2a+b>2$, we cannot have a tacnode, or a triple point, by (8.3.7).

  Suppose $x\notin A$. If we go to a log del Pezzo, then $B$ is lt and smooth,
since $A+B$ cannot have a triple point. Hence we must go to a fence. But as
$A \subset S_2$ contains at least two singular points, $K_{S_2} + A + B$
is nef by adjunction, which contradicts (13.3). If $T$ is a net, then $A$ and $E$ are sections, 
for otherwise, since $2a + b > 2$, one is a section, and the other is 
a double section, which contradicts (11.5.11.2). There are at least
$3$ multiple fibres, so $x$ is a non-chain point. Now (1) follows from (23.4.1.1).

 Suppose $x\in A$. We consider the various possibilities, (8.4.7.6-11).

If we have a Banana, then from the classification (13.2), (2-3) follow.

 Now suppose $S_2$ is a fence. Since $A_1$ contains at least two singular points, by (13.3), 
$B$ is singular, $x \in A_1$ and $A_1$ is a $-1$-curve. Suppose $A_2$ does not lie in the Du 
Val locus. By considering the possible indices for singularities along $A_1$, and using the 
fact that $x$ has maximal index, it follows that $A_2$ contains either $(3)$ or $(4)$. One 
checks in all cases that $A_2^2 \leq 0$, a contradiction. Thus $B \subset S_2^0$, $g(B) =1$, 
and $S_2$ is Du Val. Now the possibilities are classified in (13.4) and (13.5).

 If $B$ has a cusp, then by (13.4) the singularities are $A_1+A_2$. Consider 
$h:\map \t{T}.{S}_1.$. Since $\t{E}_2^2 \leq -2$, the first five blow ups must be along $B$. 
This gives configuration $(v;n)$ of (11.2.1), and $x=(2)$. But then $r =2$, and $A_2$
contains three $A_1$ points, contradicting (13.3). By (11.2.1) the configuration on $T$ is 
$(v;n^2)$, $x=(3)$ and $r =3$. But then there is a point outside $A$ of coefficient $2/3$, a 
contradiction, since by our conventions the hunt would choose that point.

Suppose $B$ has a  node. Then $E$ contains an $A_r$ point, and if $E$ is a 
$-j$-curve, $-j+r+4=K_{S_2}^2$. (4) follows by comparing the possible indices,
with the possibilities of (13.5). 

 Suppose $x \in A$, and $T$ is a net. Then $A$ must be a fibre (otherwise by (11.5.11.2), as 
above, both $A$ and $E$ are sections, and $A$ contains at least two singular points, while 
$E$ contains at most one). By (11.5.5), $A_1$ contains exactly three singular points, so 
$e_1 \geq 2/3$. 

Suppose $E$ has degree at least three. Then $E$ has degree exactly three, $A$ is a fibre
of multiplicity three, $e_1 = 2/3$ and $r =3$. $A \subset T$ is the fibre (11.5.9.5), with 
$k=2$. In particular $A \subset T$ contains singularities $A_2$ and $(3)$. It follows 
from our convention for $x$, that $x = (3)$, and $E \subset T^0$. If we extract the 
$-3$-curve adjacent to $A$ and contract $A$, $E^2$ becomes zero, thus we obtain $\Hz 0.$. 
This contradicts (11.5.10), since $E$ remains smooth and has degree three. 

It follows $E$ is a double section and the singularities of $T$ along $A$ are each $A_1$. 
Consider $h:\map \t{T}.W.$, for $W$ a smooth relative minimal model (see (11.5.4)). Note 
$K_T+E+\lambda F$ is trivial for some $\lambda > 0$. Thus $(K_W+E)\cdot E<0$ and $E$ remains 
a smooth double section. By (11.5.10), $E^{2}=4$ on $W$, and $W=\Hz 1.$ or $\Hz 0.$. 
Let $p$, $q \in W$ be the ramification points of $\map E.{\pr 1.}.$, with $p$ the point on 
the fibre corresponding to $A$. Clearly over $p$, $h$ is given by blowing up twice along $E$.
Thus (since $\t{E}^2 \leq -2$ on $T$)  $q$ must correspond to a singular fibre, $G$, of $T$. 
There can be no other singular fibres. The first $j+2$ blow ups of $h$ over $q$ must be along
$E$. At this point $G$ is in the Du Val locus, and has a $D_{j}$ singularity. The next blow 
up (if it exists) is the non-interior blow up away from $E$ (see (11.5.4)), and subsequent 
blow ups are non-interior, determined by $x$. $S_1$ has a non-chain singularity $z$, with 
centre $-2$, and two branches $(2)$. The last branch is determined by $x$. For example if 
$x =(\ul{2},3)$, then the last branch is $(3',2)$, with the marked curve meeting the central 
curve. If $\map Z.S_1.$ extracts the central divisor at $z$, then $Z$ is a net, $G$ is a 
fibre, $A$ and the central divisor are sections, and we have case (1).

The remaining possibility is that $A=\Sigma$. Thus $A$ is a $-1$-curve, and
we have (8.4.7.7). $A_1 \subset T$ cannot contain two $A_1$ points
(or it would not contract), so $b > 2/3$. It 
follows as in the proof of (18.7) that $B$ acquires an ordinary
cusp, and the singularities along $\Sigma$ are $(2)$ and $(3)$.

If $B \subset S_2^0$, then $S_2$ is Du Val and $S_2^0$ is simply
connected. $x=(j)$, $j \geq 3$. Then $S_2^2 = 6-j\leq 3$. It follows from
the simply connected list that
$S_2=S(E_{j+3})$.  Let $Q$ be an end, opposite the central divisor,
in an $A_{j-1}$ chain of the $E_{j+3}$ point. By
(3.8), there is a unique $-1$-curve, $G$, meeting $Q$ 
(and no other curves over the singular point). Since
$B \in |-K_{S_2}|$, $G$ meets $B$ normally. If
$\map Z.S_1.$ extracts the central divisor at $z$, then
$G$ is a fibre, $A$ and the central divisor are sections, and we have
case (1). In this case $G$ is the fibre (11.5.9) with $k = j-1$.

Otherwise $B$ contains exactly one singular point $y$. Note by adjunction, $b = 5/6$. We 
consider the next hunt step with respect to $K_{S_2} + b B$. 

We argue that $(S_3,B+ C)$ is given by (13.4). If $x_2 \in B$, then $T = Y$ of (19.4), 
$M = \Sigma$, $B \subset S_3^0$ and the result is clear. So we can assume $x_2 \not \in B$. 
Since $2b + c > 2$ we cannot have a net, or, by (8.3.7), a tacnode. $b + c > 1$, so we 
cannot have a triple point. Thus $S_3$ is a fence, and the result follows from (13.4).

Note $\t{B}_2^2 = 6 - j \leq 4$, while $B_3^2 = K_{S_3}^2 = 6$. Thus the first $j$ blow ups 
of $h:\map \t{T}_3. S_3.$ are along $B$. In particular $x_2 \not \in B$, and $x_2$ is a 
non-chain singularity with center $-2$ and two branches $(2)$ and $(2,2)$. Let 
$w = E_3 \cap \Sigma_3$. Note (see (11.1)) that $w$ is uniquely determined by $j$, and the 
marked point $y \in B$, and thus by $x \in A$. For example if $x = A_t$, then $w = (t + 1)$.
$x$ has index at most $6$, and there are at least $2$-curves over $x$ (since $B$ contains a 
singular point), and if $x$ is Du Val it has index at least $4$ (or the hunt chooses the 
$(3)$ point). Thus the possibilities for $x$ are $(\ul{3},2),(\ul{2},3)$ and $A_t$, with 
$3 \leq t \leq 5$. One checks for each possibility that $e(x_2) \geq \frac{r-1}{r}$, so with 
our convention, the hunt would choose $x_2$ instead of $x$, a contradiction.
\qed\enddemo

\proclaim{23.5.1.1 Lemma} Let $(S_1,A_1)$ be as in (23.5).
Suppose there is a non-chain
singularity $z \not \in A_1$ such that if $f:\map Z.S_1.$ extracts
the central curve, $E$, then $Z$ is a net and $E$ and $A_1$ are sections.
Then:
\roster 
\item $(S_1,A_1)$ is uniquely determined by $z$. 
\item $K_{S_1} + A_1$ is negative.
\item There are exactly three multiple fibres, each log terminal.
\item Every non-chain singularity $z$ occurs for some pair $(S_1,A_1)$.
\endroster
\endproclaim
\demo{Proof} Suppose $f:\map Z.S_1.$ is as in the statement. Clearly there
are three multiple fibres. Let $G$ be a multiple fibre. By (11.5.5), 
$G$ contains exactly two singular points, $G \cap E$ and $G \cap A_1$.
The reduced components of ${\Cal F}$ (in the language of (11.5.1) and (11.5.5)), each lie
over one of the singular points. It follows there are no interior blow ups
so $G$ is log terminal. 

Let $h:\map \t{Z}.W.$ be a relative minimal model which is an isomorphism in a neighbourhood
of $\t{E} \subset \t{Z}$. By (11.5.7), $h$ is uniquely determined by $z$. If $E$ is a $-j$-curve, 
then $W = \Hz j.$, and $E \subset W$ generates one of the extremal rays, and $A_1 \subset E$ is 
an $j$-curve, disjoint from $E$. The first blow up of $h$ (for each multiple fibre) is along $A$, 
and then there is a sequence of interior blow ups uniquely determined by the marked singularity 
of $E$ at $E \cap G$. The singularities $A \cap G$ and $E \cap G$ have the same index, thus 
$K_{S_1} + A_1$ is negative (since $E$ is contractible). Any two $j$-curve in $\Hz j.$ disjoint
from $E$ are permuted by an automorphism which is trivial on $E$. Thus $(S_1,A_1)$ is uniquely 
determined by $z$. \qed \enddemo

\remark{23.6 Remark} In the case of (23.5.1) 
where $A$ contains two $A_1$ singularities, let $B \subset S_1$ 
be the image of the fibre
which passes through the third singularity along $A$. Then 
$K_{S_1}+ A + B$ is numerically trivial, log canonical, and $2(K_{S_1}+A+B)$
is linearly equivalent to zero.  Thus
$B$ is a two complement for $K_{S_1} + A$, see (22.3.2).
\endremark

\subhead 23.7 The case (23.2.1) \endsubhead
  Suppose $A_1$ contains at most two singularities. By (22.3.1) $K_{S_1}+A_1$ 
has a $1$-complement, $X$. That is, there is a reduced curve $X$, containing $A$ as an
irreducible component, with the following properties: \newline

(*) $K_{S_1} + X$ is trivial, Cartier, log canonical,
and $K_{S_1} + C$ is log terminal and anti-ample for every irreducible component 
$C$ of $X$. 

We will forget about $A_1$, drop the subscript, and classify pairs $(S,X)$, satisfying
(*).

\remark{23.8 Remarks} 
The conditions in (*), together with adjunction, imply the following:
\roster
\item The only possible singularities of
$S$ along $X$ are at the intersection points of two irreducible components.
\item $X$ is necessarily reducible, and has either two irreducible components, meeting
in two points, or has three irreducible components forming a triangle.
\item $S \setminus X$ is Du Val.
\endroster
\endremark

We first introduce two useful birational transformations. We start with a 
pair $(S,X)$ and obtain a new pair $(S',X')$. $S'$ is obtained by one blow up and one
blow down, and $X'$ consists of the pushforward of the total transform of $X$. 

We remind the reader of our sporadically employed 
convention of using the same notation to indicate a curve,
and its strict transform under a birational transformation.

\proclaim{23.9 Definition-Lemma} 
\item{(a)} Suppose a component, $C$, of $X$ is a $-1$-curve. Let $\map T.S.$ extract a 
divisor, $E$, adjacent to $C$. $C \subset T$ is contractible. Let
$\map T.S'.$ be the contraction of $C$. Let $q \in S'$ be the image of $C$.
\item{(b)} Suppose there a $-1$-curve $\Sigma$ with the following properties:
\roster
\item $\Sigma\cap X$ is a 
single point, $p$, a singular point of $S$, and $\Sigma$ meets a 
unique curve of the minimal desingularisation over $p$.
\item The strict transforms of $\Sigma$ and $X$ on the minimal desingularisation are
disjoint.
\item $\Sigma$ contains at most
one other singular point of $S$, and $K_S + \Sigma$ is log terminal away from $p$.
\endroster
Let $\map T.S.$ extract the divisor, $E$, over $p$ which is adjacent to $\Sigma$. 
$\Sigma \subset T$ is contractible. Let $\map T.S'.$ be the contraction of $\Sigma$.
Let $q \in S'$ be the image of $\Sigma$.

 In both cases the new pair $(S',X')$ also satisfies (*). In case (a), $X'$ 
has the same number of components as $X$, but in case (b)  the number of 
components goes up by one. 
\endproclaim

The inverse transformations $(S',X') \dasharrow (S,X)$ are easy to describe:

\proclaim{23.10 Definition-Lemma } Let $(S',X')$ satisfy (*). Let 
$q \in X'$ be a singular point (resp. a smooth point) of $X'$, lying on an 
irreducible component $E$ of $X$ (there are two choices of $E$ when $q$ is a singular
point, and only one choice when $q$ is a smooth point). Let 
$\map \t{S}.\t{S'}.$ be a composition of blow ups, with the center of each
blow up at the unique point of (the strict transform of) $E$, lying over $q$. 
Blow up at least enough times so
that the self-intersection of $E \subset \t{S}$ is at most $-2$. Let
$C \subset \t{S}$ (resp. $\Sigma \subset \t{S}$) 
be the unique $-1$-curve over $p$, that is the exceptional divisor of
the last blow up. Let $\map \t{S}.S.$
contract all the $K$ non-negative curves, one of which, by construction, is $E$. 
Let $p \in S$ be the image of $E \subset \t{S}$.
$(S,X)$ satisfies (*), and $(S',X')$ is obtained from
$(S,X)$ by making a transformation of type (a)  (resp (b)) extracting $E$, and 
using $p \in C \subset S$
(resp. $p \in \Sigma \subset S$).
\endproclaim
\demo{Proofs of (23.9-10)} These are easy, see for example (11.4). \qed \enddemo

We will refer to the transformations of (23.9-10) as being of type
(a),(b),(a)$^{-1}$ and (b)$^{-1}$.

\proclaim{23.11 Lemma} Let $\map \ol{S}.S.$ be the minimal desingularisation in a
neighbourhood of $X$ (that is resolving only the singularities along $X$). 
If $\ol{S}$ has a $K_{\ol{S}}$-negative birational contraction, 
then $(S,X)$ has a transformation of type (a) or (b).
\endproclaim
\demo{Proof} Note by (23.8.3), $\ol{S}$ is Du Val. Suppose $\ol{S}$ has a 
$K_{\ol{S}}$-negative birational contraction. If the associated $-1$-curve is a component, 
$C$ of $X$, then clearly we have a transformation of type (a). Otherwise call the $-1$-curve
$\Sigma$.  By (3.3), $K_{\ol S} + \Sigma$ is log terminal, and $\Sigma$ contains at most one 
singularity. Let $X' \subset \ol{S}$ be the total transform of $X$. 
$K_{\ol{S}}+X'=\pi ^*(K_S+X)$, thus $K_{\ol S}+X'$ is numerically trivial. Thus 
$X'\cdot D=-K_{\tilde S}\cdot \tilde \Sigma=1$. Thus $\Sigma$ meets a unique irreducible 
component of $X'$. This must be an exceptional component, or $\Sigma \subset S$ is 
contractible. \qed \enddemo

\proclaim{23.12 Proposition} Let $(S,X)$ satisfy (*). One of the following holds:
\item{(1)} The pair $(S,X)$ is toric (see (22.4)), or obtained from a toric pair
by a single transformation of type $(b)^{-1}$.
\item{(2)} $S=\Hbz n.$, $X$ consists of a fibre union a section which passes
through the singular point. (As a special case we get $\pr 2.$ and $X$ is a 
conic union a line.)
\item{(3)} $S=\Hbz 2.$ and $X$ consists of two sections in the smooth locus.
\item{(4)} $(S,X)$ is obtained from the pair in (2) by a sequence of transformations
of type $(a)^{-1}$. Moreover, exactly one component, $C$, of $X$ is a $-1$-curve and $C$ has 
a 1-complement $Y$ (see \S 22) such that $(S,Y)$ is toric. 
\endproclaim
\remark{23.13 Remarks} 
\roster
\item In case (4), the proof gives a bit more information. 
For the transformations, the section is never contracted. $C$ is the other component. 
\item Note in particular that the surface $S$ is toric, except possibly in the second case 
of (1).
\endroster
\endremark
\demo{Proof} Note there can only be a finite sequence of transformations (of any type),
since each improves the singularities. Note also, in any sequence of transformations,
there can be at most one of type (b), since otherwise $X'$ would have at least four
components, violating (23.8.1). Finally note that if $(S',X')$ has a transformation of
type (b), then so does $(S,X) = (a)^{-1}(S',X')$ (an $(a)^{-1}$ transformation will not
destroy a curve of type $\Sigma$ in (23.9)). 

Now start with a pair $(S,X)$. Perform a
transformation of type (b) if one exists. The only further possible transformations are
of type (a).

Now assume we have $(S,X)$ and there is no transformation of type (a) or (b).
Then by (23.7), either $S$ is smooth along $X$, or $\ol{S}$, the minimal desingularisation 
of $S$ in a neighbourhood of $X$, is a net, of Picard number two. Since $K_{\ol{S}} + X'$
is trivial, and $X' \subset \ol{S}^0$, it follows that $\ol S$ is smooth. Now if $X$
has two components it is easy to see we have (2) or (3),
and if $X$ has 
three components, it is easy to check that $S=\Hbz n.$ and $X$ is 
two fibres union a section in the smooth locus (as a special case we get 
$\pr 2.$ and three lines), which are all toric. 

Note $(a)^{-1}$ transformations preserve toric pairs. So if we end with three components,
we began in the situation of (1).

So we can assume we end with (2) or (3), and we obtain the original pair by a sequence of
$(a)^{-1}$ transformations.

In (3), after one transformation of type $(a)^{-1}$, there is a curve of type $\Sigma$ (the
strict transform of the fibre through $q$), that is the new pair admits a transformation
of type (b). We ruled out this possibility above.

So we can assume our original pair is obtained from $(a)^{-1}$ transformations applied to (2).
We need only check the extra conditions in (4). 
Let $X_1$ be the section of self-intersection 
$n+2$, and $X_2$ the fibre. 

Let $y \in \Hz n.$ be the intersection point
of the strict transforms of $X_1$ and $X_2$. Let $x \in \Hz n.$ be the point where
$X_1$ meets the negative section. Let $G$ be the fibre of $\Hbz n.$ through $x$. There is
a unique section 
$\Sigma_n \subset \Hz n.$ which is disjoint from the negative section, and has
${(n+1)}^{\text {st}}$ order contact with $X_1$ at $y$. Scheme-theoretically,
$X_1 \cap \Sigma_{n} = (n+1) y$. As remarked above, the pair
$(X_2 + G + \Sigma_{n},\Hbz n.)$ is toric.

We will argue that after any blow up
$\t{X}_1^2$ (the self-intersection of the strict transform on the minimal
desingularisation) remains positive. In particular, $X_1$ will never be chosen
as $E$ (in the language of (23.10)). Thus each of the transformations is 
toric with respect to $(X_2 + G + \Sigma_n)$ (and its transforms). (4) then follows immediately.

Note that unless $q$ (in the language of (23.10)) corresponds to either $x$ or $y$, there is a 
neighbourhood of $\t{G} + \t{\Sigma}_n + \t{X_1}$ (on the minimal desingularisation) which is not 
effected. Further there is a neighbourhood of $\t{G}$ which is effected only if $q = x$, and any
choice of $q$ has the same effect on numbers $\t{X}_1^2$, $\t{\Sigma}_n^2$, and 
$X_1 \cdot \Sigma_n$ (the order of contact at $y$). 

Now if we choose $q = x$, then after the transformation
$G$ is a $-1$-curve of type $\Sigma$,  a possibility we have ruled out.

If after a blow up (in the construction of $\map \t{S}.\t{S}'.$ in (23.10))
$\t{X}_1$ becomes non-positive, then $\t{\Sigma}_n$ is separated from $\t{X}_1$, and
on the new surface, $\Sigma_n$ is a $-1$-curve of type $\Sigma$, a contradiction
as before. \qed \enddemo

\heading Appendix $L$: Log terminal surface singularities and adjunction \endheading

We note for the readers convenience (or edification) 
that a normal surface is log terminal iff it has quotient
singularities, and that for $C$ an analytically irreducible 
germ of a curve at $p$,
$p\in C \subset S$, if $S$ is given locally analytically
about $p$ as the quotient $q:\map V.S.$ with $V$
smooth, then $K_S + C$ is log terminal (resp. log canonical)
at $p$ iff
$q^{*}(C) \subset V$ is smooth (resp. has normal crossings). For
a proof see (20.3) of \cite{Kollaretal}.

Let $p \in C \subset S$ be a point of $S$ on the (non-empty)
reduced curve $C$. We give the criteria for log terminal and
log canonical at $p$ in this case. 
For proofs see chapter 3 of \cite{Kollaretal}.

Suppose first $p$ is a singular point of $S$:

$K_S+C$ is lt iff $p$ is a cyclic singularity,
and the resolution graph is a chain, with $\t{C}$ meeting
one end of the chain, and normally. In particular $C$ is smooth
at $p$. We will refer to the unique exceptional divisor which $\t{C}$ meets
as the divisor {\bf adjacent to} $C$ over $p$. 

If $C$ is singular at $p$, then $K_S + C$ is log canonical (lc) at
$p$ iff $C$ has two analytic branches at $p$, each of which if 
smooth, $p$ is
a chain singularity, the two branches of $\t{C}$ are disjoint, and
meet opposite ends of the chain, and normally. 

If $C$ is smooth at $p$  then $K_S + C$ is lc but not lt, iff $p$
is a non-chain singularity, with two branches $(2)$, $(2)$, and $\t{C}$
meets normally the end of third branch, at the opposite end from 
the central divisor. 

Note in particular that the criteria above are local analytic. 
When $p$ is smooth this is not the case: 

$K + C$ is lt at $p$ iff either $C$ is smooth at $p$, or $C$ has
two irreducible components (in a Zariski neighbourhood), 
each smooth, meeting normally at $p$.

$K + C$ is lc at $p$ but not lt, iff $C$ is irreducible (in a Zariski
neighbourhood of $p$), and has a simple node at $p$. 

Thus for example if $C$ is the union of three general lines in
$\pr 2.$, $K + C$ is lt, but if $C$ is a nodal cubic, it is lc but
not lt.

Because of this peculiarity of the definition, we will sometimes use
the variant purely log terminal (plt). For the general definition see
chapter 2 of \cite{Kollaretal}. For our purposes it is enough
that $K + C$ is plt at a singular point
$p$ iff it is lt at $p$, and plt at a smooth point $p$ iff $C$ is
smooth at $p$. In short, $K + C$ is plt iff it is lt and $C$ is smooth.

We often use the following standard result implicitly in this paper:

\proclaim{L.1 Lemma} Let $(S,x)$, $(T,y)$ be germs of log terminal surface singularities, and 
let $\Gamma',\Gamma$ be the dual graphs of their minimal  resolutions. Assume 
$\Gamma$ is obtained from $\Gamma'$ either by increasing the weight of one of the vertices,
or by adding a single vertex and edge. The following hold:
\roster
\item The index of $y$ is strictly bigger than the index of $x$.
\item The coefficient of each exceptional divisor of the minimal resolution of $y$
is strictly bigger than the coefficient of the corresponding divisor over $x$, 
unless $y$ is a Du Val singularity (in which case all the coefficients are zero).
\endroster
\endproclaim
\demo{Proof} The cyclic case of (1) 
follows from (3.1.8) of \cite{Kollaretal}. The non-cyclic case (and indeed the
cyclic case as well) can be checked from the explicit lists in \cite{Brieskorn68}. 

For (2) we will consider the case when the weight goes up. The other case
is handled analogously, or alternatively, follows from (3.1.3) of \cite{Kollaretal}.

Of course we may assume the weight goes up by exactly one.
Since the singularities are determined by the resolution graph, we
can assume $(T,y)$ is obtained from $(S,x)$ by  
a smooth blow up $h:\map \t{T}.\t{S}.$ at a point lying on a unique
exceptional divisor, the divisor corresponding to the vertex whose weight increases.
Let $\Sigma$ be the $-1$-curve extracted by $h$. Note $e=e(\Sigma,K_S) < 0$ since
$y$ is log terminal. Let $E$ be an exceptional divisor (of the minimal
desingularisation) over $x$. Let $f:\map T'.T.$,
$g:\map S'.S.$ extract $E$. We have an induced map $h':\map T'.S'.$. 
$$
(g \circ h')^*(K_S) = {h'}^*(K_{S'} + e(E,K_S) E) = K_{T'} + e(E,K_S) E + (e-1) \Sigma. 
$$
Thus, computing intersections on $T'$, 
$$
-e(E,K_S) E^2 = (K_{T'} + (e-1) \Sigma) \cdot E < K_T' \cdot E = 
(f^*(K_T) -e(E,K_T) E) \cdot E = -e(E,K_T) E^2. \qed
$$
\enddemo

\subhead L.2 Adjunction and Inversion of Adjunction \endsubhead

In chapter 3 of \cite{Shokurov93}, Shokurov introduces the
following:
\definition{Definition-Lemma (The Different)}
Let $S$ be a normal variety and $C \subset S$ a reduced effective
divisor. Let $\ol{C}$ be the normalisation, and $i:\map \ol{C}.S.$
the induced map. Let $B$ be any effective $\Bbb Q$-divisor on $S$. 
Assume $K_S + C + B$ is $\Bbb Q$-Cartier:

There is a canonically defined effective $\Bbb Q$-Weil divisor 
$\Diff_C(B)$ on $C$ such that
$$
i^*(K_S + C + B) = K_C + \Diff_C(B)
$$
\enddefinition

(The assumption that $K_S + C + B$ is $\Bbb Q$-Cartier is not necessary, but
without it the definition of $i^*$ is more technical. The above is sufficient
for our needs.)

When $(S,C)$ is log terminal the different has a simple interpretation
via Koll\'ar's Bug-Eyed cover. See (4.15).

Shokurov obtained the following quite surprising result:

\proclaim{L.2.1 Theorem (Inversion of Adjunction)} Assume $\di(S) \leq 3$. 

$K_S + C + B$
is log canonical (resp. plt and $B$ is a pure boundary) 
in a neighbourhood of $C$ iff
$K_C + \Diff_C(B)$ is log canonical (resp. klt).
\endproclaim
For the proof see (3.4) of \cite{Shokurov93}, 

By an elegant extension of an idea  of Shokurov, Koll\'ar proved the plt case of (L.2.1)
in all dimensions, see (17.6) of \cite{Kollaretal}.

\proclaim{L.2.2 Lemma} Let $C \subset S$ be a reduced irreducible curve
on a $\Bbb Q$-factorial projective surface. If $(K_S + C) \cdot C < 0$,
then $C = \pr 1.$. 
\endproclaim
\demo{Proof} By the subadjunction formula, (5.1.9) of \cite{KMM87}, for
sufficiently divisible $m > 0$, 
$\omega_{C}^{[m]}$ is a subsheaf of the line bundle $\ring C.(m(K_S + C))$.
As the latter has negative degree, it follows that
$H^1(\ring C.)=H^0(\omega_C) =0$. \qed \enddemo

\remark{Remark} We expect there is also a local version of (L.2.2), namely that
$\Diff_C(0)$ is at least as big as the conductor, which is either empty,
or has degree at least two. \endremark

The next result is our principal application of (L.2.1):

\proclaim{L.2.3 Corollary} Let $S$ be a $\Bbb Q$-factorial projective 
surface, and $C,D \subset S$ 
two curves, with $C$ integral.
If $(K_S + C + D) \cdot C \leq 0$, and $C \cap D$ contains at
least two points, then $K_S + C + D$ is log canonical, $C \cap D$
is exactly two points, and $(K_S + C + D)\cdot C = 0$.
\endproclaim
\demo{Proof} $C = \pr 1.$ by (L.2.2). By (L.2.1), since $K + C + tD$
is not lc at $D \cap C$, for any $t > 1$, the coefficient 
of $Q$ in $\Diff_C(D)$ is at least one, for any $Q \in D \cap C$. It follows
that $D \cap C$ is exactly two points $\{Q_1,Q_2\}$  and
$\Diff_C(D) = Q_1 + Q_2$. In particular $K_C + \Diff_C(D)$ is lc, and thus
$K_S + C + D$ is lc in a neighbourhood of $C$ by (L.2.1). \qed \enddemo

{\bf Warning:} In this paper 
we frequently use (L.2.1-3) without reference, or under the catch all {\it by adjunction}.

\heading Appendix N: Normalisation of an algebraic space \endheading

Here we give a construction of the normalisation of an algebraic
space, as we have not been able to find it in the literature.

Note that the property of being normal is preserved by \'etale
covers, and so makes sense for algebraic spaces.

Recall for a reduced scheme $X$ there is a unique finite morphism
$\map \t{X}.X.$ (the normalisation) which is universal for dominant maps 
from normal schemes. 

\proclaim{N.1 Lemma} Let $f:\map Y.X.$ be a dominant \'etale map of 
schemes, with $X$ reduced. Then there is a natural fibre diagram
$$
\CD
\t{Y} @>>> Y \\
@VVV       @VVV \\
\t{X} @>>> X .
\endCD
$$
\endproclaim
\demo{Proof} Immediate from the universal property of the 
normalisation.\qed\enddemo

Now let $j:\map R.Y \times Y.$ be an \'etale equivalence
relation. By the universal property there is an induced map
$\t{j}:\map \t{R}.\t{Y} \times \t{Y}.$.

\proclaim{N.2 Lemma} (assumptions as above)
$\t{j}$ is an \'etale equivalence relation. \endproclaim
\demo{Proof} 
By the definition of an equivalence relation, 
we need to show that $j$ is a monomorphism, and 
$\Hom(T,\t{R}) \subset \Hom(T,\t{Y}) \times \Hom(T,\t{Y})$
is an equivalence relation for all $T$. 

It follows from (N.1) that $\t{j}$ is a monomorphism, and the
projections $\map \t{R}.\t{Y}.$ are \'etale. Also, by
the universal property, for any normal scheme $T$,
$\Hom_d(T,\t{R})= \Hom_d(T,R)$ and 
$\Hom_d(T,\t{Y} \times \t{Y}) = \Hom_d(T,Y \times Y)$, where
$\Hom_d$ indicates the subset of dominant maps. It follows that 
$\Hom_d(T,\t{R}) \subset \Hom_d(T,\t{Y}) \times \Hom_d(T,\t{Y})$ is
an equivalence relation for any normal $T$. 

For the transitivity axiom. Let
$T = \t{R} \fb {(p_2,p_1)}.\t{R}$, where $p_1,p_2:\map \t{R}.\t{Y}.$
are the two projections. It is enough to show that
$(p_1\circ pr_1,p_2 \circ pr_2) \in \Hom(T,\t{R})$. But
since $T$ is normal, and $p_i$ are dominant, this follows from the 
remarks above. The other axioms are analogously obtained. \qed\enddemo

\proclaim{N.3 Definition-Lemma: Normalisation} 
Let $W$ be a reduced algebraic space. Then there is a unique finite dominant 
map $\map \t{W}.W.$ from a normal algebraic space, universal
for dominant maps from normal algebraic spaces. \endproclaim
\demo{Proof} Let $f:\map Y.W.$ be an \'etale cover by an
affine scheme, let $R = Y \fb W. Y$, so that $W = Y/R$.  Let
$\t{W} = \t{Y}/\t{R}$, which is an algebraic space by (N.2), and
clearly normal. By construction there 
is an induced map $\map \t{W}.W.$. (N.1) and the
construction implies $\t{Y} = \t{W} \fb W. Y$. Now the
universal property of $\map \t{Y}.Y.$
implies the universality of $\map \t{W}.W.$.
\qed \enddemo

\remark{N.4 Remark} Note that Lemma (N.1) holds for algebraic spaces,
just copy the proof of (N.1) and use (N.3). 
\endremark

\heading Index \endheading
\settabs\+indent &Coefficientindexblahblah \quad &\cr
\+ &$\alpha$, $\beta$, & {\bf \S 10}\cr
\+ &almost Du Val, & {\bf \S 2}\cr
\+ &adjacent to, & {\bf \S 2} \& {\bf L}\cr
\+ &banana, & {\bf (8.0.11)}\cr
\+ &branches, & {\bf \S 2}\cr
\+ &bug-eyed cover, & {\bf (4.0)}\cr
\+ &central divisor, & {\bf \S 2}\cr
\+ &chain, & {\bf \S 2}\cr
\+ &coefficient, & {\bf (1.16)}\cr
\+ &complement, & {\bf \S 22}\cr
\+ &dominated by, & {\bf (1.1)}\cr
\+ &effective, & {\bf \S 2}\cr
\+ &fence, & {\bf (8.0.12)}\cr
\+ &flush, & {\bf (8.0.2)}\cr
\+ &furthest, f& {\bf (11.2), (11.4)}\cr
\+ &$\Hbz n.$, & {\bf \S 2}\cr
\+ &hunt, & {\bf (1.15)}\cr
\+ &interior, i & {\bf \S 11}\cr
\+ &level, & {\bf (8.0.2)}\cr
\+ &list, simply connected list, & {\bf \S 3}\cr
\+ &log uniruled, & {\bf \S 2}\cr
\+ &near, n & {\bf \S 11}\cr
\+ &net, & {\bf (1.16)}\cr
\+ &$\rho$, & {\bf \S 2}\cr
\+ &spectral value, & {\bf (8.0.6)}\cr
\+ &tacnode, & {\bf (8.0.13)}\cr
\+ &tiger, & {\bf (1.13)}\cr
\+ &$\tilde S$, & {\bf \S 2}\cr
\+ &type $-(\list a.n.)$, & {\bf (8.0.5)}\cr
\+ &toric, & {\bf \S 22}\cr
\+ &$X$, $Y$, & {\bf \S 11}\& {\bf \S 12} \& {\bf \S 13}\cr
\+ &weight, w& {\bf \S 10}\cr

\heading References \endheading

\bibliography{big}

\end